\documentclass{memo-l}
\usepackage[english]{babel}
\usepackage[utf8]{inputenc}
\usepackage[T1]{fontenc}
\usepackage{faktor}

\usepackage{epsfig}
\usepackage{graphicx}
\usepackage{color}
\usepackage{url}
\usepackage{hyperref}
\usepackage{amssymb}
\usepackage{esint}
\def\R{\mathbb R}
\def\Z{\mathbb Z}
\def\N{\mathbb N}

\def\C{\mathbb C}
\def\J{\mathbb J}
\def\n{\|{\hspace{-0.12em}}|}

\def\loc{{\text{\upshape loc}}}
\def\negquad{\!\!\!\!}
\def\negqquad{\!\!\!\!\!\!\!\!}
\def\eps{\varepsilon}
\newcommand{\cal}[1]{{\mathcal #1}}
\DeclareMathOperator\essinf{essinf}

\DeclareMathOperator\DivGamma{div_{_\Gamma}}
\DeclareMathOperator\DivGammaPrimo{div_{_{\Gamma'}}}
\DeclareMathOperator\DivGammaZero{div_{_{\Gamma_0}}}
\DeclareMathOperator\DivGammaUno{div_{_{\Gamma_1}}}
\DeclareMathOperator\Det{det}
\DeclareMathOperator\Tr{Tr}
\numberwithin{equation}{chapter}

\DeclareMathOperator{\Ker}{Ker}
\DeclareMathOperator{\Real}{Re}
\DeclareMathOperator{\Ima}{Im}

\theoremstyle{plain}
\newtheorem{thm}{Theorem}[section]
\newtheorem{lem}[thm]{Lemma}
\newtheorem{prop}[thm]{Proposition}
\newtheorem{cor}[thm]{Corollary}
\newtheorem{rem}[thm]{Remark}

\theoremstyle{definition}
\newtheorem{definition}[thm]{Definition}
\newtheorem{example}[thm]{Example}

\numberwithin{section}{chapter}
\numberwithin{equation}{chapter}
\makeindex

\begin{document}
\frontmatter

\title[The wave equation with acoustic...]
{The  wave equation with acoustic boundary conditions on non-locally reacting surfaces}
\author{Delio  Mugnolo}
\address{D.~Mugnolo, Lehrgebiet Analysis, Fakult\"at Mathematik und Informatik, Fern\-Universit\"at in Hagen, D-58084 Hagen, Germany}
\email{delio.mugnolo@fernuni-hagen.de}

\author{Enzo Vitillaro}
\address{E.~Vitillaro, Dipartimento di Matematica e Informatica, Universit\`a di Perugia\\
       Via Vanvitelli,1 06123 Perugia ITALY}
\email{enzo.vitillaro@unipg.it}

\date{\today}



\thanks{The authors would like to thank the anonymous reviewer for the helpful comments .The work of the first author was partially supported by the Deutsche Forschungsgemeinschaft (DFG), Grant 397230547. The work of the second author  was realized within the auspices of the INdAM -- GNAMPA Projects
{\em Equazioni alle derivate parziali: problemi e modelli} (Prot\_U-UFMBAZ-2020-000761), and it was also supported by {\em Progetto Equazione delle onde con condizioni acustiche,  finanziato  con  il Fondo  Ricerca  di Base, 2019, della Universit\`a degli Studi di Perugia} and by {\em Progetti Equazioni delle onde con condizioni iperboliche ed acustiche al bordo,  finanziati  con  i Fondi  Ricerca  di Base 2017 and 2018, della Universit\`a degli Studi di Perugia}.}

\begin{abstract} The aim of the paper is to study  the problem
$$
\begin{cases} u_{tt}-c^2\Delta u=0 \qquad &\text{in
$\R\times\Omega$,}\\
\mu v_{tt}- \DivGamma (\sigma \nabla_\Gamma v)+\delta v_t+\kappa v+\rho u_t =0\qquad
&\text{on
$\R\times \Gamma_1$,}\\
v_t =\partial_\nu u\qquad
&\text{on
$\R\times \Gamma_1$,}\\
\partial_\nu u=0 &\text{on $\R\times \Gamma_0$,}\\
u(0,x)=u_0(x),\quad u_t(0,x)=u_1(x) &
 \text{in $\Omega$,}\\
v(0,x)=v_0(x),\quad v_t(0,x)=v_1(x) &
 \text{on $\Gamma_1$,}
\end{cases}$$
where $\Omega$ is a open domain of $\R^N$ with uniformly $C^r$ boundary ($N\ge 2$, $r\ge 1$),
$\Gamma=\partial\Omega$, $(\Gamma_0,\Gamma_1)$ is a relatively open partition of $\Gamma$ with $\Gamma_0$ (but not $\Gamma_1$) possibly empty.
Here $\DivGamma$ and $\nabla_\Gamma$ denote the
Riemannian divergence and gradient operators on $\Gamma$, $\nu$ is the outward normal
to $\Omega$, the coefficients $\mu,\sigma,\delta, \kappa, \rho$ are suitably regular functions on $\Gamma_1$ with $\rho,\sigma$ and $\mu$ uniformly positive while $c$ is a positive constant.
This problem have been proposed long time ago by Beale and Rosencrans, when $N=3$, $\sigma=0$, $r=\infty$, $\rho$ is constant,  $\kappa,\delta\ge 0$,  to model acoustic wave propagation with locally reacting boundary.

In this paper we first study well-posedness in the natural energy space and give regularity results. Hence we give precise qualitative results for solutions when $\Omega$ is bounded and   $r=2$, $\rho$ is constant, $\kappa,\delta\ge 0$.
These results motivate a detailed discussion of the derivation of the problem in Theoretical Acoustics and the consequent proposal of adding to the model  the integral condition $\int_\Omega u_t=c^2\int_{\Gamma_1}v$.
\end{abstract}

\maketitle
\tableofcontents

\chapter{Introduction and main results} \label{intro}
\section{Presentation of the problem and literature overview}
We deal with the wave equation posed in a suitably regular open domain of
$\R^N$, supplied with an acoustic  boundary condition on a part of the boundary and an homogeneous Neumann boundary condition on the (possibly empty) remaining part of it. More precisely we consider the
initial- and boundary-value problem
\begin{equation}\label{1.1}
\begin{cases} u_{tt}-c^2\Delta u=0 \qquad &\text{in
$\R\times\Omega$,}\\
\mu v_{tt}- \DivGamma (\sigma \nabla_\Gamma v)+\delta v_t+\kappa v+\rho u_t =0\qquad
&\text{on
$\R\times \Gamma_1$,}\\
v_t =\partial_\nu u\qquad
&\text{on
$\R\times \Gamma_1$,}\\
\partial_\nu u=0 &\text{on $\R\times \Gamma_0$,}\\
u(0,x)=u_0(x),\quad u_t(0,x)=u_1(x) &
 \text{in $\Omega$,}\\
v(0,x)=v_0(x),\quad v_t(0,x)=v_1(x) &
 \text{on $\Gamma_1$,}
\end{cases}
\end{equation}
where $\Omega\subset\R^N$, $N\ge
2$, is an open domain with boundary $\Gamma=\partial\Omega$ uniformly of class $C^r$ in the sense of \cite{Stein1970}, where
the value of $r\in\N\cup\{\infty\}$ will be further specified when needed, so $r=1$ when nothing is said (in most of the paper we shall take $r=1$ or $r=2$).
We assume $\Gamma=\Gamma_0\cup\Gamma_1$, $\overline{\Gamma_0}\cap\overline{\Gamma_1}=\emptyset$,
so $\Gamma_0$ and $\Gamma_1$ are clopen in $\Gamma$, and $\Gamma_1\not=\emptyset$.
All these properties of $\Omega$, $\Gamma_0$ and $\Gamma_1$ will be more formally stated in \S~\ref{section2.2} and we shall refer to them  as assumption (A0).

Moreover $u=u(t,x)$, $v=v(t,y)$, $t\in\R$, $x\in\Omega$, $y\in\Gamma_1$,
$\Delta=\Delta_x$ denotes the Laplace operator with respect to the space
variable, while $\DivGamma$ and $\nabla_\Gamma$ respectively denote the
Riemannian divergence and gradient operators on $\Gamma$. By
$\nu$ we denote the outward normal to $\Omega$ and $c$ is a fixed positive constant.

Acoustic boundary conditions as those in problem \eqref{1.1} have been introduced by Beale and Rosencrans, for general domains, in \cite{beale, beale2,BealeRosencrans} to model acoustic wave propagation, motivated by \cite[pp.~259--264]{morseingard}. The physical derivation of \eqref{1.1} will be the  subject of \S~\ref{section7}, so we here  briefly recall the one presented in \cite{BealeRosencrans}.

In it $N=3$, $\Gamma_0=\emptyset$, and $\Omega$ is either a bounded or an external domain filled
with a fluid which is at rest but for acoustic wave motion. Since the fluid is assumed to be non-viscous, one denotes by $u$ the
velocity potential, so   $-\nabla u$ is the particle velocity, and $u$ satisfies the wave equation $u_{tt}-c^2\Delta
\phi=0$ in $\R\times\Omega$, where $c>0$ is the sound speed in the fluid.  In the Beale--Rosencrans model, one supposes that $\Gamma$ is not rigid but subject to small oscillations, and that each point of it reacts to the excess pressure of the acoustic wave like a
(possibly) resistive harmonic oscillator or spring, so there is no transverse tension between neighboring points of $\Gamma$. These surfaces are called
locally reacting in \cite[pp.~259--264]{morseingard}.

Now, $v$ denotes the normal displacement of the surface $\Gamma$  into the domain,
$\delta=\delta (y)\ge 0$ the surface resistivity, $\mu=\mu(y)>0$ the surface mass density, $\kappa=\kappa(y)>0$ the spring constant, $\rho=\rho_0>0$ the  constant unperturbed gas density. Since the excess pressure
 is given by $\rho_0u_t$,
 \begin{footnote}{We shall explain in detail this point in \S~\ref{section7}. }\end{footnote}
 one gets that $v$ satisfies the equation \eqref{1.1}$_2$ with $\sigma\equiv 0$. The boundary condition $v_t=\partial_\nu u$ on $\Gamma_1$ follows since fluid particles cannot penetrate $\Gamma$, so they have to move with it.

The same model (for $\delta=0$) has been proposed in \cite{WeitzKeller} when $\Omega$ is a strip in $\R^2$ or $\R^3$, and in \cite{Peters}
when $\Omega$ is the half-space in $\R^3$, to describe acoustic wave propagation in an ice-covered ocean. We refer the interested reader to~\cite{Belinsky} for a historical overview of these and related problems in mathematical physics, especially focusing on the results by the soviet school.

It is worth observing that, in all the above recalled models, only $\nabla u$ and $u_t$ have a physical meaning, while the velocity potential $u$ is merely a mathematical object.

After their introduction acoustic boundary conditions for locally reacting surfaces have been the subject of several papers. See for example \cite{Alcantara, CFL2004,CavFrota, FL2006, FG2000, GGG,  graber, Hao, JGPhD, JGSB2012, KJR2016, KT2008, LLX2018,  Maatoug2017,  mugnolo, Shomberg}.
When one dismisses the simplifying assumption that neighboring point do not interact (in the terminology of \cite[p.266]{morseingard}) such surfaces are called of {\em extended reaction}.
We shall call those which react like a membrane {\em non-locally reacting} (other types of reactions can be considered), and in this case one has to take $\sigma>0$ in \eqref{1.1}.

The simplest case in which $\sigma$  is constant and  the operator $\DivGamma(\sigma \nabla_\Gamma)$ reduces (up to $\sigma$) to the Laplace--Beltrami operator $\Delta_\Gamma$ was briefly considered in \cite[\S 6]{beale} and then studied in \cite{FMV2011,  FMV2014, VF2013bis,VF2016, VF2017} and in the recent paper \cite{becklin2019global}.
In all of them the authors assume that $\Gamma_0\not=\emptyset$ and that the homogeneous Neumann boundary condition on it is replaced by the (mathematically more attracting) homogeneous Dirichlet boundary condition. Such a  boundary condition, as it will be clear from \S~\ref{section7}, does not appear to have a meaning in the physical model.

Several papers in the literature also deal with the wave equation with porous acoustic boundary conditions, where $\mu\equiv 0$ in \eqref{1.1}, which has a different nature. We refer to
\cite{AN2015bis, AN2015, BB2017, BB2018, graber2010}. See also \cite{Shoubo, LP2018} for  related problems.

In the present paper we shall study the case of non-locally reacting boundary conditions without any simplifying Dirichlet boundary condition, which reflects the original physical model.

The coefficients $\mu,\sigma,\delta, \kappa, \rho$ are given real  functions on $\Gamma_1$
satisfying the following assumptions depending on the value of $r$:
\renewcommand{\labelenumi}{{(A\arabic{enumi})}}
\begin{enumerate}
\item \label{A1} $\mu,\sigma\in W^{r-1,\infty}(\Gamma_1)$ with $\mu_0:=\essinf_{\Gamma_1}\mu>0$ and $\sigma_0:=\essinf_{\Gamma_1}\sigma>0$;
\item \label{A2} $\rho\in W^{r,\infty}(\Gamma_1)$ with $\rho_0:=\essinf_{\Gamma_1}\rho>0$;
\item \label{A3} $\delta, \kappa \in W^{r-1,\infty}(\Gamma_1)$.
\end{enumerate}
We notice that in (A1--3) we denoted $W^{\infty,\infty}(\Gamma_1)=\bigcap\limits_{n=1}^\infty W^{n,\infty}(\Gamma_1)=C_b^\infty(\Gamma_1)$.
\begin{footnote}{Here and in the sequel the subscript ``$b$'' in spaces of type $C^r$ means that all derivatives up to order $r$ are (not necessarily uniformly when $r=\infty$) bounded. Clearly  Morrey's Theorem is used.}\end{footnote}

The meaning of the Sobolev spaces used above is the standard one when $\Gamma_1$ is compact, while it will be made precise in  \S~\ref{section3} in the non-compact case.

\section{Main results I: well-posedness}\label{mainresultsI}
Our first aim is to study well-posedness for problem \eqref{1.1} and regularity of its solutions, including in our treatment all types of domains already considered in the literature above when $\sigma\equiv 0$, i.e., bounded or external domains, that is domains with compact boundary, and
strips or half-spaces, that is domains with non-compact but trivial boundary.
To avoid considering two different cases we are forced to  consider domains with possibly non-compact and non-trivial boundary. This choice also allows to consider types of domains not yet treated in the literature but of interest in applications. An example of them are cylindrical domains, that is domains of the form $\Omega=\R\times \Omega_1$, where $\Omega_1\subset\R^2$ is an open domain with compact $C^r$ boundary. When $\Omega_1$ is bounded such a  domain models an infinite tube of constant section. Another example is given by
domains of type $\Omega=\{(x,y,z)\in\R^3: \sqrt{x^2+y^2}<{\cal R}(z)\}$, where ${\cal R}\in C^r_b(\R)$, ${\cal R}>0$,
modeling an infinite pipeline of non-constant section.
We shall try to optimize our regularity assumption, so we shall not take $r$ larger than needed.

 For this reason we are considering, as already stated, domains with uniformly $C^r$-boundary, for which a trace operator is available, as shown in \cite{LeoniSobolev2}. This choice causes a first difficulty which we are now going to describe.

 Sobolev spaces of real nonnegative order on a manifold $\Gamma$ and boundedness properties of  differential operators in  \eqref{1.1} are classical objects when $\Gamma$ is compact (see for example \cite{lionsmagenesIII} for the smooth case and \cite{grisvard} for $r<\infty$). They are well-known objects also
when $\Gamma$ is  possibly non-compact but smooth  and the order is an integer (see for example \cite{aubinmanifolds, hebey, jost}). As shown in \cite{hebey}  in the non-compact case many usual properties fail unless $\Gamma$ satisfies good geometrical properties, essentially related to curvature bounds.

When $\Gamma$ is possibly non-compact but smooth, and the order is a real number, a satisfactory theory is given in \cite{triebel} and related papers, provided $\Gamma$ has bounded geometry. Unfortunately this approach does not extend to $r=1,2$, since curvatures bounds are meaningless.
A different approach has been recently proposed by Amann in \cite{Amann2013}, where curvature bounds are replaced by an uniform regularity assumption, formulated in a form which is slightly different from the one we are using in this paper. Indeed in our case $\Gamma$ is not an abstract manifold.  This assumption  is equivalent to the bounded geometry assumption, as shown in \cite{Amann2013, Disconzi2016}. When $\Gamma$ is smooth our assumption essentially reduces to it.

The only reference found by the authors concerning these spaces when $r\in \N$  is \cite[Chapter 8]{LeoniSobolev2}. Our assumption originates from it.  In the quoted book the author only gives the definition of  Besov spaces and proves  Trace Theorem, without studying the  relationship with the Riemannian structure on $\Gamma$ and the boundedness of differential operators. Such a study is mandatory in this paper.

Fortunately the general approach in \cite{Amann2013} can be generalized to $r\in\N$. The first task in the paper will be to prove this assertion and giving boundedness properties of differential operators.   It will be accomplished in  \S~\ref{section3}.  All cases considered in the references quoted above will be included as particular cases. In the sequel of this introduction we shall freely use Sobolev spaces $H^s(\Gamma_1)$, referring to \S~\ref{section3} for more details and making the reader aware that all familiar properties of the compact case continue to hold.
We hope that \S~\ref{section3} could be of independent interest for scholars dealing with problems with relevant boundary terms.

To state our first result we introduce the phase space
\begin{equation}\label{1.2}
{\cal H}=H^1(\Omega)\times H^1(\Gamma_1)\times L^2(\Omega)\times L^2(\Gamma_1).
\end{equation}
We make the reader aware that all functions spaces considered in the paper will be complex, as they are commonly used in Acoustics. On the other hand the corresponding spaces of real-valued functions are trivially invariant under the semigroup governing problem \eqref{1.1}.
\begin{thm}[\bf Well-posedness]  \label{theorem1.1}
Under assumptions (A0--3), for any choice of data $U_0=(u_0,v_0,u_1,v_1)\in {\cal H}$, problem \eqref{1.1} has a unique weak solution
\begin{footnote}{I.e., a solution in a suitable distributional sense, see Definition~\ref{Definition4.2}  below.}\end{footnote}
\begin{equation}\label{1.2bis}
(u,v)\in C(\R;H^1(\Omega)\times H^1(\Gamma_1))\cap C^1(\R;L^2(\Omega)\times L^2(\Gamma_1))
\end{equation}
continuously depending on data.
Moreover, when we also have
\begin{footnote}{Here and in the sequel $\Delta$, $\partial_\nu$ and $\DivGamma$ here are taken in suitable distributional sense, see  \S\ref{subsection3.3} below}\end{footnote}
\begin{equation}\label{1.2ter}
\begin{aligned}
&u_1\in H^1(\Omega),\quad v_1\in H^1(\Gamma_1),\quad \Delta u_0\in L^2(\Omega), \quad\partial_\nu {u_0}_{|\Gamma_1}\in L^2(\Gamma_1),\\
&\partial_\nu {u_0}_{|\Gamma_0}=0,\quad \text{and}\quad \DivGamma(\sigma\nabla_\Gamma v_0)\in L^2(\Gamma_1),
\end{aligned}
\end{equation}
then
\begin{equation}\label{1.2quater}
(u,v)\in C^1(\R;H^1(\Omega)\times H^1(\Gamma_1))\cap C^2(\R;L^2(\Omega)\times L^2(\Gamma_1))
\end{equation}
and \eqref{1.1}$_1$--\eqref{1.1}$_4$ hold a.e., respectively in $\R\times\Omega$, on $\R\times\Gamma_1$ (twice) and on $\R\times\Gamma_0$.
Finally, when $\rho(y)\equiv \rho_0>0$, solutions satisfy for all
$s, t\in \R$
the energy identity
\begin{equation}\label{energyidentity}
\begin{split}
&\frac {\rho_0}2 \int_\Omega |\nabla u|^2\!+ \frac {\rho_0}{2c^2}\int_\Omega |u_t|^2
+\frac 12 \int_{\Gamma_1} \sigma |\nabla_\Gamma v|_\Gamma^2 \\
&\qquad +\frac 12
\int_{\Gamma_1} \mu |v_t|^2+\frac 12
\int_{\Gamma_1} \kappa|v|^2\Big|_s^t=-
\int_s^t\int_{\Gamma_1}\delta|v_t|^2.
\end{split}
\end{equation}
\end{thm}
The proof of Theorem~\ref{theorem1.1} is based, beside boundedness properties of differential operators mentioned above, on standard linear semigroup theory and on the characterization  of  generalized solutions (in the semigroup sense) as  weak solutions.

Our second main result concerns optimal regularity of solutions. As usual in hyperbolic problems higher regularity
requires corresponding regularity of $\Gamma$ and data, as well as compatibility conditions.
To state next result we set, for $1\le n\le r$,
\begin{equation}\label{1.5}
{\cal H}^n:=H^n(\Omega)\times H^n(\Gamma_1)\times H^{n-1}(\Omega)\times H^{n-1}(\Gamma_1),
\end{equation}
and  we introduce, when $r=\infty$, the spaces
\begin{equation}\label{Frechet}
\begin{split}
C^\infty_{L^2}(\overline{\Omega})&:=\{u\in C^\infty(\Omega): d^n u\in L^2_n(\Omega)\quad\forall n\in\N_0\},\\
C^\infty_{L^2}(\Gamma_1)&:=\{v\in C^\infty(\Gamma_1): D_\Gamma^n v\in L^2_n(\Gamma_1)\quad\forall n\in\N_0\},
\end{split}
\end{equation}
where $d$ and $D_\Gamma$ respectively denote the differential on $\Omega$ and the covariant derivative on $\Gamma$, while
$L^2_n(\Omega)$ and $L^2_n(\Gamma_1)$ stand for the spaces of $n$-times covariant tensor fields with square integrable norm.

\begin{rem}
Since they look quite uncommon we make some remarks on the first one, $L^2_n(\Omega)$. They also apply, {\em mutatis mutandis},  to the second one, $L^2_n(\Gamma_1)$.
By Morrey's Theorem all elements of $C^\infty_{L^2}(\overline{\Omega})$ continuously extend, with all their derivatives, to $\overline{\Omega}$, so motivating the notation we used. By the same reason $C^\infty_{L^2}(\overline{\Omega})=\bigcap_{n\in\N_0} H^n(\Omega)$,
so it is a Fréchet space with respect to the associated family of seminorms. Hence the notation $C^\infty(\R;C^\infty_{L^2}(\overline{\Omega}))$
is meaningful in the sense of the G\^{a}teaux derivative (see \cite[pp.~72--74]{Hamilton}), and trivially $C^\infty(\R;C^\infty_{L^2}(\overline{\Omega}))\subseteq C^\infty(\R\times\overline{\Omega})$.
Applying Morrey's Theorem again we have the continuous (and possibly strict) inclusions  $C^\infty_{L^2}(\overline{\Omega})\subseteq C^\infty_b(\overline{\Omega})=C^\infty_b(\Omega)\subseteq C^\infty(\overline{\Omega})$. Trivially $C^\infty_{L^2}(\overline{\Omega})= C^\infty_b(\overline{\Omega})$ when $\Omega$ has finite measure, and $C^\infty_{L^2}(\overline{\Omega})= C^\infty(\overline{\Omega})$ when $\Omega$ is bounded.
\end{rem}

\begin{thm}[\bf Optimal regularity]
\label{theorem1.2} Let (A0--3) hold, $r\ge 2$ and $n\in\N$ such that $2\le n\le r$.
Suppose that $(u_0,v_0,u_1,v_1)\in {\cal H}^n$ and
the following compatibility conditions, where $[\cdot]$ stands for the integer part, hold:
\begin{equation}\label{1.6}\left\{
\begin{aligned}
&\partial_\nu \Delta^iu_0=0,\qquad \text{on $\Gamma_0$ \qquad for $i=0,\ldots, [n/2]-1$,}\\
&\partial_\nu \Delta^iu_1=0,\qquad \text{on $\Gamma_0$ \qquad for $i=0,\ldots, [(n-1)/2]-1$, when $n\ge 3$,}\\
&\partial_\nu u_0=v_1,\qquad \text{on $\Gamma_1$,}\\
&\mu \partial_\nu u_1=\DivGamma (\sigma\nabla_\Gamma v_0)-\delta \partial_\nu u_0-\kappa v_0-\rho u_1,\qquad \text{on $\Gamma_1$, \quad when $n\ge 3$,}\\
&\begin{split}
c^2\mu \partial_\nu \Delta^i\!u_0\!=\!\DivGamma (\sigma\nabla_\Gamma \partial_\nu\Delta^{i-1}\!u_0)\!
-\!\delta \partial_\nu\Delta^{i-1}\!u_1\!-\!\kappa\partial_\nu\Delta^{i-1}\!u_0\!-\!c^2\rho\Delta^i \!u_0\\
\text{\quad on $\Gamma_1$  \quad for $i=1,\ldots, [n/2]-1$, \quad when $n\ge 4$,}&
\end{split}\\
&\begin{split}
c^2\mu \partial_\nu \Delta^i\!u_1\!=\!\DivGamma (\sigma\nabla_\Gamma \partial_\nu\Delta^{i-1}\!u_1)\!
-\!c^2\delta \partial_\nu\Delta^{i-1}\!u_0\!-\!\kappa\partial_\nu\Delta^{i-1}\!u_1\!-\!c^2\rho\Delta^i \!u_1\\
\text{on $\Gamma_1$  \quad for $i=1,\ldots, [(n-1)/2]-1$, \quad when $n\ge 4$.}&
\end{split}
\end{aligned}\right.
\end{equation}
Then the corresponding weak solution $(u,v)$ enjoys the further regularity
\begin{equation}\label{1.7}
(u,v)\in \bigcap\limits_{i=0}^n C^i(\R, H^{n-i}(\Omega)\times H^{n-i}(\Gamma_1)).
\end{equation}
Moreover, when $r=\infty$, if $u_0,u_1 \in C_{L^2}^\infty(\overline{\Omega})$,
$v_0,v_1\in C_{L^2}^\infty(\Gamma_1)$
and  \eqref{1.6} hold for all $i\in\N$, we have
\begin{equation}\label{1.8}
(u,v)\in  C^\infty(\R;C^\infty_{L^2}(\overline{\Omega})\times C^\infty_{L^2}(\Gamma_1)).
\end{equation}
\end{thm}

The conditions \eqref{1.6} are quite involved but unavoidable. Indeed, by \eqref{1.7} and the Trace Theorem (see \S~\ref{subsection3.3.2}) one gets $u_{|\Gamma}\in \cap_{i=0}^{n-2} C^i(\R, H^{n-2-i}(\Gamma))$. When $n=2$, using \eqref{1.7} and Lemma~\ref{Lemma4B} one easily gets that equations \eqref{1.1}$_1$--\eqref{1.1}$_4$ hold true in the spaces $C(\R,L^2(\Omega))$, $C(\R,L^2(\Gamma_1))$ (twice) and $C(\R,L^2(\Gamma_1))$, from which \eqref{1.6} follow since $v_{tt}(0)=\partial_\nu {u_1}_{|\Gamma_1}$. When $n\ge 3$ we then get \eqref{1.6} by deriving $(n-2)$-times with respect to time equations \eqref{1.1}$_2$ and \eqref{1.1}$_4$   and using \eqref{1.1}$_1$ and  \eqref{1.1}$_3$ to transform time derivatives in spatial operators.
Conditions \eqref{1.6} were never pointed out before in their general form, neither when $\sigma>0$ nor
$\sigma\equiv 0$.
In the standard physical applications, however, one models the evolution after an acoustic perturbation originating  in a compact subset of $\Omega$, with $\Gamma_1$ being at rest,  so one can suppose that $\nabla u_0$ and $u_1$ are compactly supported and $v_0=v_1=0$. In this case, \eqref{1.6} trivially holds.
Actually the regularity result in \cite[Theorem~2.2]{beale} (when $\sigma\equiv 0$) is stated  in this situation.

The proof of Theorem~\ref{theorem1.2} is based on boundedness properties of differential operators mentioned above and standard linear semigroup theory. The main difficulty is technical and it consists in extending standard regularity theory for elliptic problems on $\Omega$ and $\Gamma_1$ to the present case in which $\Gamma_1$ can be non-compact and non-smooth. This extension, of possible independent interest, is given in  Theorems~\ref{lemma4.5}--\ref{lemma4.5BISS} below.

\section{Main results II: solutions' behavior}\label{mainresultsII} The second aim of the paper is to study
the long-time behavior of solutions of \eqref{1.1} when $\Omega$ is bounded. We shall also assume that $\Gamma$ is of class $C^2$ and that the following
properties of $\rho,\kappa,\delta$, coming from the physical model, hold:
\begin{itemize}
\item[-] the unperturbed fluid density $\rho$ is constant; this property is mandatory in the physical derivation of problem \eqref{1.1};
\item[-] the surface resistivity $\delta$  and the spring constant $\kappa$ are non-negative.
\end{itemize}
More formally, in the sequel we shall assume, additionally to (A0--3),
\renewcommand{\labelenumi}{{(A\arabic{enumi})}}
\begin{enumerate}
\setcounter{enumi}{3}
\item \label{A4} $\Omega$ is bounded, \quad $r=2$, \quad $\rho(x)\equiv \rho_0>0$, \quad $\delta, \kappa\ge 0$ on $\Gamma_1$.
\end{enumerate}
This one is exactly the case  studied in \cite{beale} when $\sigma\equiv 0$ (a case we explicitly exclude by means of (A1)). While \cite{beale} focuses on unusual spectral properties arising as $\sigma\equiv 0$, in the present paper, due the regularizing effect of the second order elliptic operator on $\Gamma_1$, the problem exhibits more standard spectral properties. Hence we can concentrate on studying  the asymptotic  behavior of solutions.
The main difficulty one faces on this concern is constituted by several different types of trivial (to the authors) solutions that \eqref{1.1} may possess, whose number increases with the numbers of connected components of $\Gamma_1$ in which $\kappa$ and  $(\kappa,\delta)$ identically vanish.

For the sake of simplicity, we are going to present our results in detail only  when $\Gamma_1$ is connected, and offer only a brief overview of the situation in the disconnected case at the end of this introduction.
When $\Gamma_1$ is connected, the problem \eqref{1.1} possess only three types of vanishing velocity solutions, i.e., solutions for which $u(t,x)=u(t)$.
\begin{itemize}
\item The first one, arising independently on $\kappa$, is simply
$u(t,x)\equiv u_0\in\C$, $v(t,x)\equiv 0$. Since only $-\nabla u$ and $\rho_0u_t$ have a physical meaning, this fact is not surprising: adding a space-time constant to the solution has no physical effects and only corresponds to a different gauge.  The standard way to take into account this fact, used also in \cite{beale}, consists in disregarding  the average $\fint_\Omega u(t):=\frac 1{|\Omega|}\int_\Omega u(t)$ of $u(t)$ , but not that of $u_t(t)$. Since, as it will be evident in the sequel,  $\fint_\Omega u$ is not invariant under the flow,
one cannot formalize this fact by restricting to data $u_0$ with vanishing average. Indeed, completely disregarding the average of $u$ (as done in \cite{beale}) amounts to simultaneously changing the pressure. One can add to solutions only space-time constants. Instead we shall take the projection of $H^1(\Omega)$ onto
$\faktor{H^1(\Omega)}{\C}\simeq H^1_c(\Omega)$, where
\begin{equation}\label{5.1}
 H^1_c(\Omega):=\left\{u\in H^1(\Omega): {\textstyle\int}_\Omega u=0\right\},
\end{equation}
 in which  $\|\nabla (\cdot)\|_2$ is a norm, equivalent to the standard one. Hence we shall study, as far as $u$ is concerned, the behavior of the couple $(\nabla u, u_t)$,  disregarding the one of $\fint_\Omega u$.
Indeed, as it will be clear from the derivation of the physical model in \S~\ref{section7}, $\fint_\Omega u$ is chosen in a suitable (but physically meaningless) way to recover the fluid velocity and the incremental pressure from the velocity potential field $u$.
\item The second type of vanishing velocity solutions only arises when $\kappa\not\equiv 0$ and, to the best of our knowledge,  has never been noticed before,
even if it exists also when $\sigma\equiv 0$ (provided $1/\kappa\in H^1(\Gamma_1)$). Indeed, when $\kappa\not\equiv 0$, the elliptic equation
\begin{equation}\label{1.9}
-\DivGamma(\sigma \nabla_\Gamma v^*)+\kappa v^*+\rho_0=0\qquad \text{on $\Gamma_1$,}
\end{equation}
has a unique (trivially real valued) weak solution $v^*\in H^2(\Gamma_1)$ (see Lemma \ref{lemma5.4}  below), and for all $u_1\in \C$ problem \eqref{1.1} admits the solution $u(t,x)=u_1t$, $v(t,x)=u_1v^*(x)$. This type of solutions does not seem, at a first look, to have a physical interpretation. Indeed it is characterized by vanishing fluid velocity and arbitrarily large  incremental pressure and boundary deformation.
The discussion of their physical meaning, if any, is the main motivation of \S~\ref{section7}.
\item The third type of vanishing velocity solutions only arises when $\kappa\equiv 0$ and it is simply $u(t,x)\equiv 0$, $v(t,x)\equiv v_0\in \C$. These solutions, existing also when $\sigma\equiv 0$, were ignored in \cite{beale} since this case was not considered of interest, and allow for arbitrary constant boundary deformations. Physically they reflect the insensitivity on boundary constant deformations of the physical model.
\end{itemize}
The last two types of vanishing velocity solutions are exact solutions for data $(u_0,v_0,u_1,v_1)$ belonging to the one dimensional subspace $V_0$ of the phase space ${\cal H}$  (defined in \eqref{1.2})
 spanned by
\begin{equation}\label{1.10}
V^*:=\begin{cases}
(0,v^*,1,0)\qquad &\text{if $\kappa\not\equiv 0$,}\\
(0,1,0,0)  \qquad &\text{if $\kappa\equiv 0$,}
\end{cases}
\end{equation}
where $v^*$ is still defined as the solution of \eqref{1.9}.
Moreover one easily finds that ${\cal H}=V_0\oplus {\cal H}_1$, where
\begin{equation}\label{1.11}
{\cal H}_1:=\left\{(u,v,w,z)\in {\cal H}: \int_\Omega w=c^2\int_{\Gamma_1}v\right\},
\end{equation}
and that ${\cal H}_1$ is invariant under the flow induced by problem \eqref{1.1} on ${\cal H}$.

In the sequel, we shall give our main results for data in ${\cal H}_1$, since the behavior for data in ${\cal H}$ easily follows from them and, as we shall explain in the sequel, they are the most significant ones.

As it is clear by the energy identity \eqref{energyidentity} one can expect asymptotic stability, at least in the energy sense, when the system is damped, i.e., $\delta\not\equiv 0$. We shall call this type of stability {\em physical stability}.
This notion will be motivated in Corollary~\ref{corollary7.1} in \S~\ref{section7.3}.
By the contrary, when the system is undamped, i.e., $\delta\equiv 0$, the energy is conserved and one can expect to get pure oscillatory solutions.

Indeed, in the first case we get the following main result.
\begin{thm}[\bf Physical stability when $\boldsymbol{\delta\not\equiv 0}$]  \label{theorem1.3}
Let (A0--4) hold, $\Gamma_1$ be connected, $\delta\not\equiv 0$ and
$(u,v)$ denote the weak solution of \eqref{1.1} corresponding to data $(u_0,v_0,u_1,v_1)\in {\cal H}_1$. Then, as $t\to\infty$,
\begin{equation}\label{1.12}
\left\{
\begin{alignedat}{4}
&\nabla u(t)\to 0\quad&&\text{in $[L^2(\Omega)]^N$},\qquad  && u_t(t)\to 0\quad&&\text{in $L^2(\Omega)$}, \\
&u(t)-{\textstyle\fint}_\Omega u(t)\to 0\quad&&\text{in $H^1(\Omega)$},\qquad  && \tfrac 1t {\textstyle\fint}_\Omega u(t)\to 0\quad&&\text{in $\C$},\\
& v(t)\to 0\quad &&\text{in $H^1(\Gamma_1)$}, \qquad && v_t(t)\to 0\quad&&\text{in $L^2(\Gamma_1)$}.
\end{alignedat}\right.
\end{equation}
\end{thm}

\begin{rem} The limit behavior for solutions originating from data in $\cal{H} _1$ is what one can expect from the physical model.
It does not give, rigorously speaking, asymptotic stability in the mathematical sense, since the behavior of $\fint_\Omega u(t)$ as $t\to\infty$ is partially undetermined. To recover mathematical stability one has to consider the projection of the dynamical system on ${\cal H}_{01}={\cal H}_1\cap {\cal H}_0$, where
${\cal H}_0=H^1_c(\Omega)\times H^1(\Gamma_1)\times L^2(\Omega)\times L^2(\Gamma_1)$.
\end{rem}

The proof of Theorem~\ref{theorem1.3} is based on a far from trivial combination of the abstract strong stability result of \cite{ArendtBatty, LyubichVu} (see also \cite{EngelNagel}) with a recent rigidity result for partially overdetermined elliptic problems, based on the classical Unique Continuation Principle, given in \cite{FarinaValdinoci}.

In the undamped case we have the following result.
\begin{thm}[\bf Fourier decomposition when $\boldsymbol{\delta\equiv 0}$]  \label{theorem1.4}
Let (A0--4) hold, $\Gamma_1$ be connected and  $\delta\equiv 0$.
Then there is sequence $((u^n,v^n))_n$ of standing wave solutions of \eqref{1.1}, i.e., solutions of the form
\begin{footnote}{Here and in the sequel the imaginary unit will be denoted by $\mathfrak{i}$.}\end{footnote}
\begin{equation}\label{1.13}
u^n(t,x)=u^{0n}(x) \, e^{-\mathfrak{i}\lambda_nt},\qquad
v^n(t,x)=v^{0n} (x) \, e^{-\mathfrak{i}(\lambda_nt+\pi/2)},
\end{equation}
satisfying the following properties:
\renewcommand{\labelenumi}{{\roman{enumi})}}
\begin{enumerate}
\item $((u^{0n},v^{0n}))_n$ is a sequence of non-identically vanishing real-valued elements of $H^2(\Omega)\times H^2(\Gamma_1)$ such that
\begin{equation}\label{1.14}
\begin{cases} -c^2\Delta u^{0n}=\lambda_n^2u^{0n} \qquad &\text{in
$\Omega$,}\\
- \DivGamma (\sigma \nabla_\Gamma v^{0n})+\kappa v^{0n}+\rho_0\lambda_n u^{0n} =\mu \lambda_n^2v^{0n}\qquad
&\text{on $\Gamma_1$,}\\
\partial_\nu u^{0n}=0 &\text{on $\Gamma_0$,}\\
\partial_\nu u^{0n}=-\lambda_nv^{0n}\qquad
&\text{on
$\Gamma_1$;}
\end{cases}
\end{equation}
\item $(\lambda_n)_n$ is a real sequence such that $\lambda_{2n+2}=-\lambda_{2n+1}$ for all $n\in\N$, each $\lambda_{2n+1}$ repeating at most finitely many times,
$$0<\lambda_1\le \cdots\le \lambda_{2n-1}\le\lambda_{2n+1}\le \cdots,\qquad\qquad \lambda_{2n+1}\to\infty;$$
\item  for all data $U_0=(u_0,v_0,u_1,v_1)\in {\cal H}_1$  the corresponding weak solution  $(u,v)$ of \eqref{1.1}
admits the following Fourier decomposition:
\begin{equation}\label{1.15}
\left\{
\begin{alignedat}{2}
u =&\sum_{n=1}^\infty \alpha_n u^n+\alpha_0&&\qquad \text{in $C_b(\R; H^1(\Omega))\cap C^1_b(\R; L^2(\Omega))$,}\\
v =&\sum_{n=1}^\infty \alpha_n v^n&& \qquad         \text{in $C_b(\R; H^1(\Gamma_1))\cap C^1_b(\R; L^2(\Gamma_1))$,}\\
\end{alignedat}\right.
\end{equation}
where
\begin{equation}\label{1.19bis}
\begin{alignedat}2
&\alpha_n=\alpha_n(U_0)=&&\rho_0\int_\Omega \nabla u_0\nabla u^{0n}-\lambda_n \int_{\Gamma_1}\mu v_1v^{0n}+\mathfrak{i}\int_{\Gamma_1}\kappa v_0v^{0n}\\
&&&+\mathfrak{i}\int_{\Gamma_1}\sigma (\nabla_\Gamma v_0,\nabla_\Gamma v^{0n})_\Gamma+\frac{\mathfrak{i}\rho_0\lambda_n}{c^2}\int_{\Gamma_1}u_1u^{0n}\quad\text{for $n\in\N$,}\\
&\alpha_0=\alpha_0(U_0)=&&\fint_\Omega u_0-\sum_{n=1}^\infty \alpha_n\fint_\Omega u^{0n}.
\end{alignedat}
\end{equation}
\end{enumerate}

\end{thm}
\begin{rem}\label{remark1.2}
It was remarked in \cite[pp.~261--262]{morseingard} that, when there is no energy absorption at the boundary, reflecting acoustic plane waves
have a pure-imaginary acoustic impedance $z:=\rho_0u_t/\partial_\nu u=\rho_0u_t/v_t$ corresponding to a vanishing acoustic resistance $\Real  z$ and to a (frequency independent) phase difference of $\pi/2$.
The standing waves in \eqref{1.13}, of frequency $\lambda_n$ and period $\tau_n=2\pi/\lambda_n$,  enjoy the same property, that is they are at quadrature.
For this reason $v^{0n}\not=v^n(0,\cdot)$ , but $v^{0n}=v^n(-\tau_n/4,\cdot)$.
The complex nature of the Fourier coefficients $\alpha_n$ naturally arises from this fact. Indeed from the interaction between $u$ and $v$ we get waves of frequency $\lambda_n$ which are combination of two waves at quadrature.
It is worth observing that since $u$ has a physical meaning only up to a space-time constant,  the constant $\alpha_0$ in \eqref{1.15}$_1$ is given for the sake of completeness only. For practical purposes one can write \eqref{1.15}  as $(-\nabla u,\rho_0 u_t,v) =\sum_{n=1}^\infty \alpha_n (-\nabla u^n,\rho_0 u^n_t,v^n)$. See Corollary~\ref{corollary7.2} in \S~\ref{section7.3}.
Hence, also in the undamped case, solutions originating from data in $\cal{H}_1$, have a pure oscillatory nature, as expected.
\end{rem}

The proof of Theorem~\ref{theorem1.4} is based on standard spectral decomposition theory for skew-symmetric operators with compact resolvent and on an {\em ad-hoc} argument to get a basis of real-valued elements satisfying \eqref{1.14}, which gives the spatial profile of the standing waves. It is possible to give it explicitly when $\Omega$ is a ball.

To understand the important role of  vanishing velocity solutions solutions of the second and third types we now give our results on  solutions behavior for general data in $\cal H$. Once the projectors of ${\cal H}=V_0\oplus {\cal H}_1$ onto its addends are explicitly written, they follows by linearity from Theorems~\ref{theorem1.3}--\ref{theorem1.4}  and by the already given form of vanishing velocity solutions solutions. In \S~\ref{section6} we shall prove them together with Theorems~\ref{theorem1.3}--\ref{theorem1.4}.

We preliminarily  remark that, by \eqref{1.9}, multiplying by $v^*$ and integrating by parts one easily gets that $$|\Omega|-c^2\int_{\Gamma_1}v^*=|\Omega|+\tfrac{c^2}{\rho_0}\int_{\Gamma_1}\left(\kappa|v^*|^2+\sigma |\nabla_\Gamma v^*|_\Gamma^2\right)>0.$$
In the damped case we have the following result.
\begin{cor}[\bf Asymptotic behavior when $\boldsymbol{\delta\not\equiv 0}$]  \label{corollary1.1}
Let (A0--4) hold, $\Gamma_1$ be connected, $\delta\not\equiv 0$ and
$(u,v)$ denote the weak solution of \eqref{1.1} corresponding to data $(u_0,v_0,u_1,v_1)\in {\cal H}$. The following limits hold true as $t\to\infty$.
\renewcommand{\labelenumi}{{\Roman{enumi}.}}
\begin{enumerate}
\item If $\kappa\not\equiv 0$ then
$$\left\{\begin{alignedat}{4}
&\nabla u(t)\to 0\,\, && \text{in $[L^2(\Omega)]^N$}, \quad &\text{or} &\quad u(t)-{\textstyle\fint}_\Omega u(t)\to 0\quad &&\text{in $H^1(\Omega)$},\\
&\tfrac 1t {\textstyle\fint}_\Omega u(t)\to \tfrac{\int_\Omega u_1-c^2\int_{\Gamma_1}v_0}{|\Omega|-c^2\int_{\Gamma_1}v^*} \quad &&\text{in $\C$,}&& \quad  v_t(t)\to 0\,\, &&\text{in $L^2(\Gamma_1)$},\quad \\
& u_t(t)\to \tfrac{\int_\Omega u_1-c^2\int_{\Gamma_1}v_0}{|\Omega|-c^2\int_{\Gamma_1}v^*}\,\, &&\text{in $L^2(\Omega)$},  && \quad v(t)\to \tfrac{\int_\Omega u_1-c^2\int_{\Gamma_1}v_0}{|\Omega|-c^2\int_{\Gamma_1}v^*}\,v^* \,\, &&\text{in $H^1(\Gamma_1)$}.
\end{alignedat}\right.
$$
\item If $\kappa\equiv 0$ then
$$\left\{\begin{alignedat}{4}
&\nabla u(t)\to 0\,\, && \text{in $[L^2(\Omega)]^N$}, \quad &\text{or} &\quad u(t)-{\textstyle\fint}_\Omega u(t)\to 0\quad &&\text{in $H^1(\Omega)$},\\
&\tfrac 1t {\textstyle\fint} u(t)\to 0 &&\text{in $\C$,}&& \quad  u_t(t)\to 0\,\, &&\text{in $L^2(\Omega)$},\quad \\
& v(t)\to \tfrac{c^2\int_{\Gamma_1}v_0-\int_\Omega u_1}{c^2{\cal H}^{N-1}(\Gamma_1)}\quad &&\text{in $H^1(\Gamma_1)$}, &&\quad  v_t(t)\to 0\quad &&\text{in $L^2(\Gamma_1)$}.
\end{alignedat}\right.
$$
\end{enumerate}
\end{cor}
Corollary~\ref{corollary1.1} shows that the asymptotic behavior at $t\to\infty$ of solutions originating from data outside ${\cal H}_1$   is exactly the one of vanishing velocity solutions solutions of the second and third types.
It is then evident that the conclusions of Corollary~\ref{corollary1.1} for data outside $\cal{H}_1$ are different from those of Theorem~\ref{theorem1.3} for data inside $\cal{H}_1$. Indeed, in the physical sense,  while in the first case solutions stabilize to zero, in the second one they stabilize to a physical configuration characterized by different excess pressure $\rho_0u_t$ (when $\kappa\not\equiv 0$) and non-vanishing boundary deformation $v$. So, for example, after an acoustic perturbation inside a balloon one can get a larger balloon with a fall of pressure. This behavior does not correspond to physical (or even common) experience, so motivating a further  analysis of the physical model.

Also in the undamped case vanishing velocity solutions solutions of the second or third types play an important role for data outside $\cal{H}_1$, as an inexplicable addendum in the decomposition.
\begin{cor}\label{corollary1.2}
Let (A0--4) hold, $\Gamma_1$ be connected, $\delta\equiv 0$,
and $(u,v)$ denote the weak solution of \eqref{1.1}  corresponding to data $(u_0,v_0,u_1,v_1)\in {\cal H}$.
Then, keeping the notation in Theorem~\ref{theorem1.4}, the following  decompositions hold.
\renewcommand{\labelenumi}{{\Roman{enumi}.}}
\begin{enumerate}
\item If $\kappa\not\equiv 0$ then
$$\begin{alignedat}{2}
u(t) =&\sum_{n=1}^\infty \alpha_n u^n(t)+\alpha_0+\tfrac{\int_\Omega u_1-c^2\int_{\Gamma_1}v_0}{|\Omega|-c^2\int_{\Gamma_1}v^*}t&&\quad \text{in $H^1(\Omega)$ and uniformly in $t\in\R$,}\\
 u_t=&\sum_{n=1}^\infty \alpha_nu^n_t+\tfrac{\int_\Omega u_1-c^2\int_{\Gamma_1}v_0}{|\Omega|-c^2\int_{\Gamma_1}v^*}\,\,&&\quad\text{in $C_b(\R;L^2(\Omega))$}, \\
v=&\sum_{n=1}^\infty \alpha_n v^n+\tfrac{\int_\Omega u_1-c^2\int_{\Gamma_1}v_0}{|\Omega|-c^2\int_{\Gamma_1}v^*}\,v^*&& \quad \text{in $C_b(\R;H^1(\Gamma_1))\cap C^1_b(\R; L^2(\Gamma_1))$.}
\end{alignedat}
$$
\item If $\kappa\equiv 0$ then
$$\begin{alignedat}{2}
u =&\sum_{n=1}^\infty \alpha_n u^n+\alpha_0&&\qquad \text{in $C_b(\R; H^1(\Omega))\cap C^1_b(\R; L^2(\Omega))$,}\\
v =&\sum_{n=1}^\infty \alpha_n v^n+\tfrac{c^2\int_{\Gamma_1}v_0-\int_\Omega u_1}{c^2\cal{H}^{N-1}(\Gamma_1)}&& \qquad         \text{in $C_b(\R; H^1(\Gamma_1))\cap C^1_b(\R; L^2(\Gamma_1))$.}\\
\end{alignedat}
$$\end{enumerate}
\end{cor}

\section{Physical discussion of the model}\label{PhysicalDiscussion}
The results in \S~\ref{mainresultsII} motivate a detailed analysis of the derivation of problem \eqref{1.1} in the physical literature. Its aim is  to understand if either these vanishing velocity solutions   or the constraint
\begin{equation}\label{eq:pohpoh}
\int_\Omega u_t=c^2\int_{\Gamma_1}v
\end{equation}
on solutions (or on data), are physically meaningful.

In particular, in \S~\ref{section7.1}, we recall the mainstream derivation of the acoustic wave equation in the Eulerian  framework, the standard one in Fluid Mechanics, and the interpretation of boundary conditions in  \eqref{1.1}. We refer to two famous textbooks, chosen among others. The first one is \cite{morseingard}, a classical treatise on Theoretical Acoustics. The second one is  \cite[Chapter VIII]{Landau6}, a reference textbook for Theoretical Physics. Here we shall summarize the main outcomes of our analysis.

Problem \eqref{1.1} is a first order approximation of the nonlinear motion equations of the fluid, namely the Continuity  and the Euler Equations. It models small amplitude perturbations of the rest state.  Despite first order approximations are usually  considered inadequate in Theoretical Physics, such a first order approximation (see \cite[p.257]{morseingard})
\begin{center}\em{"represents surprisingly well many acoustic phenomena. Sounds of ordinary intensity  involve acoustic pressures of $10^{-4}$ atm or smaller".}\end{center}

In this derivation, the generalized Hooke's law is used at an infinitesimal level. It  says that the pressure perturbation corresponding to  a small volume strain is proportional to it, the proportionality constant being the so-called bulk modulus $B$. The acoustic wave equation $u_{tt}-c^2\Delta u=0$ is then derived (with $c^2=B/\rho_0$). The same Hooke's law, written at a macroscopic level (where it is experimentally verified), reads exactly as \eqref{eq:pohpoh}, thanks to boundary conditions (see \S~\ref{section7.1}). Consequently, although \eqref{eq:pohpoh} is not a consequence of problem \eqref{1.1} (it gets lost in its physical derivation), it should be added to problem \eqref{1.1} when $\Omega$ is bounded. It is an open problem to generalize such a condition to unbounded domains.

To further confirm this conclusion we  recall in \S~\ref{section7.2} a second  derivation of the acoustic wave equation in the Lagrangian framework, given in several textbooks.  We refer to \cite[Ch. 7]{Achterberg}, \cite[Ch. Three]{ElmoreHeald}, \cite[Ch. 47]{Feynman} and \cite{GoldsteinMechanics}. In this framework it makes sense to consider the displacement vector field, describing the displacement of fluid particles. Still using the Hooke's law at a infinitesimal level, one gets the acoustic wave equation $u_{tt}-c^2\Delta u=0$ (with $c^2=B/\rho_0$) as a first order perturbation. However, in this framework, by properly formulating the boundary conditions in terms of the displacement, one gets \eqref{eq:pohpoh} simply as a consequence of them and of the Divergence Theorem.

Moreover, as shown in \cite[Ch. 7]{Achterberg}, the Eulerian and Lagrangian approaches are equivalent as first order approximations of the problem. Hence, as far as this order of approximation is adequate, they are interchangeable.
\begin{footnote}{This conclusion appears mathematically paradoxical: the descriptions given by physical fields at a given point and at the displaced ones are equivalent! The explanation is given by the sentence quoted before: the displacement is very  small due to the different scales between acoustic pressure perturbations and the pressure of normal fluids at rest.}\end{footnote}

Hence we conclude that, to model physical phenomena, problem \eqref{1.1} should be completed when $\Omega$ is bounded  by the integral condition \eqref{eq:pohpoh}, and we would like to propose this small addition to the scientific community.
Consequently the vanishing velocity solutions  of the second and third type, which are ruled out by \eqref{eq:pohpoh}, are physically meaningless. Moreover Theorems~\ref{theorem1.1}--\ref{theorem1.2} have the following (expected) physical meanings: if the physical system is linearly damped it stabilizes to the rest state, while when it is undamped all solutions are superpositions of eigenstates. See Corollaries~\ref{corollary7.1} and \ref{corollary7.2}.

\section{The case of $\boldsymbol{\Gamma_1}$  disconnected and paper organization }\label{conclusion}
When $\Gamma_1$ is disconnected the solutions behavior can be  more complex.  Our partial  results in  this case will be given in detail in \S~\ref{subsection6.3}. Roughly speaking, two main cases occur, depending on the number $\mathfrak{n}_0$ of the connected components of $\Gamma_1$ where $\kappa$ vanish identically.
\renewcommand{\labelenumi}{{\roman{enumi})}}
\begin{enumerate}
\item If $\mathfrak{n}_0\le 1$ the conclusions of Theorems~\ref{theorem1.3}--\ref{theorem1.4} continue to hold, and also those of Corollaries~\ref{corollary1.1}--\ref{corollary1.2} in a slightly modified form. See Theorem~\ref{theorem5.3} below.
\item If $\mathfrak{n}_0\ge 2$  the conclusions of Theorems~\ref{theorem1.1}--\ref{theorem1.2} are no longer true. The precise behavior of solutions also depends on the number $\mathfrak{n}_{00}\le \mathfrak{n}_0$ of the connected components of $\Gamma_1$ where both $\kappa$ and $\delta$ identically  vanish.
    \renewcommand{\labelenumii}{{\roman{enumi}.\arabic{enumii})}}
    \begin{enumerate}
    \item When $\mathfrak{n}_{00}\le 1$ (so the system is damped) every solution of \eqref{1.1} and \eqref{eq:pohpoh} possess an asymptotic limit,
    with the component $v$ being linear combinations  of the characteristic functions  of connected components of $\Gamma_1$ where $\kappa$ vanish. The coefficients essentially reflect the dynamics of these boundary pieces.  The model is geometrically consistent only for small data.
    \item When $\mathfrak{n}_{00}\ge 2$ (so covering the undamped case) problem \eqref{1.1} and \eqref{eq:pohpoh} possess further trivial solutions, with every connected component $\Gamma'$ of $\Gamma_1$ where $\kappa$ and $\delta$ vanish linearly advancing or retreating. When $\Omega$ is a spherical shell (and $\Gamma_0=\emptyset$) they are easily calculated, and the model becomes geometrically inconsistent in finite time, also for initial data  of arbitrarily small amplitude.
         \end{enumerate}
\end{enumerate}
The general conclusion is that when $\Gamma_1$ is disconnected each boundary component (but one) has  to be retained by a spring-type term, also in a small subset, to avoid non-acoustic phenomena and  model's inconsistency.

The paper is organized as follows:
\begin{itemize}
\item[-] in \S~\ref{section2} we make precise our regularity assumption on $\Gamma$ and we recall some background material on $C^r$ manifolds;
\item[-] in \S~\ref{section3} we treat Sobolev Spaces on $\Gamma$ and boundedness properties of linear operators involved in \eqref{1.1};
\item[-] \S~\ref{section4} is devoted to  the proof of Theorem~\ref{theorem1.1};
\item[-] in \S~\ref{section5} we give regularity results and we prove Theorem~\ref{theorem1.2};
\item[-]  in \S~\ref{section6} we study the qualitative behavior of solutions in the bounded case, proving Theorems~\ref{theorem1.3}--\ref{theorem1.4} and Corollaries~\ref{corollary1.1}--\ref{corollary1.2};
\item[-]  \S~\ref{section7} contains our discussion of the physical model.
\end{itemize}

\chapter{Background and preliminaries} \label{section2}
\section{Notation}
We shall denote by $I$ the identity function on any set. Given any two sets $D$ and $E$,  $E^D=\{u: D\to E\}$. Moreover $\N=\{1,2,\ldots\}$ and $   \N_0:=\N\cup\{0\}$.
For any $l\in\N$, $x\in\R^l$ and $\eps>0$ we shall denote by $B_\eps(x)$ the ball of $\R^l$ centered ad $x$ of radius $\eps$.
Moreover $Q^l=(-1,1)^l$. For $x=(x_1,\ldots,x_l),y=(y_1,\ldots,y_l)\in\C^l$, $xy=\sum_{i=1}^lx_iy_i$ and  $\overline{x}$ is the conjugate of $x$.
Moreover $d(\cdot,C)$ denotes the distance from a nonempty closed set $C\subseteq\R^l$, with $d(\cdot,\emptyset)=\infty$.

For any open $\vartheta\subset \R^l$ and $\eps>0$
 \begin{equation}\label{3.11}
   \vartheta^{\eps}:=\{x\in\vartheta: B_{\eps}(x)\subseteq\vartheta\}=\{x\in\vartheta: d(x,\vartheta^c)\ge \eps\}.
 \end{equation}
Moreover we shall use the standard notation for (complex) Lebesgue and Sobolev spaces of real order on $\vartheta$, referring to \cite{adams}  for details. We shall also use the standard notation for the $\C^l$-valued version of them. For simplicity
$$\|\cdot\|_{\tau,\vartheta}=\|\cdot\|_{L^\tau(\vartheta)},\quad \|\cdot\|_{\tau}=\|\cdot\|_{\tau,\Omega},\quad
\|\cdot\|_{s,\tau,\vartheta}=\|\cdot\|_{W^{s,\tau}(\vartheta)},\quad \|\cdot\|_{s,\tau}=\|\cdot\|_{s,\tau,\Omega}.
$$

We shall use the standard notations $C^m(\vartheta)$ for $m\in\N\cup\{\infty\}$, and  the subscripts ``$c$" and ``$b$"  will respectively denote subspaces of compactly supported functions and of functions with bounded derivatives  up to order $m$. Moreover for any multiindex $\alpha=(\alpha_1,\ldots, \alpha_l)\in\N_0^l$ we shall use the standard notation $D^\alpha=\partial_1^{\alpha_1}\ldots \partial_l^{\alpha_l}$,  $|\alpha|=\alpha_1+\ldots\alpha_l$, $\alpha!=\alpha_1!+\ldots\alpha_l!$ and $v^\alpha=v_1^{\alpha_1}\ldots v_l^{\alpha_l}$ for any $v=(v_1,\ldots, v_l)\in \R^l$. Given two multiindex $\alpha$ and $\beta$ we shall
write  $\beta\le \alpha$ to mean $\beta_i\le \alpha_i$ for $i=1,\ldots,l$, and for any such couple
we  write  $\binom\alpha\beta=\alpha!/\beta!(\alpha-\beta)!$. By $J$ we shall denote the Jacobian matrix.

We shall use the standard notation $\ell^\tau$ for Lebesgue spaces on $\N$ with respect to the counting measure.
Given a Fréchet space $X$ we shall denote by $X'$ its dual, by $\langle\cdot,\cdot\rangle_X$ the duality product and, given another Fréchet space $Y$,
by ${\cal L}(X,Y)$ we shall indicate the space of continuous linear operators between them. When $X,Y$ are Banach spaces the standard norm on it will be denoted as $\|\cdot\|_{{\cal L}(X,Y)}$. Moreover, for any Banach space $X$ we shall use standard notations for Bochner--Lebesgue and Sobolev Spaces of $X$-valued functions on any open subset of $\R^l$.
\section{Assumptions}\label{section2.2}
We make precise our structural assumption (A0) mentioned  in \S~\ref{intro}:
\renewcommand{\labelenumi}{(A{\arabic{enumi})}}
\begin{enumerate}
\setcounter{enumi}{-1}
\item $\Omega\subset\R^N$, $N\ge
2$, is an open domain with boundary $\Gamma=\partial\Omega=\Gamma_0\cup\Gamma_1$, $\overline{\Gamma_0}\cap\overline{\Gamma_1}=\emptyset$, uniformly of class $C^r$, $r\in\N\cup\{\infty\}$, in the sense of \cite{Stein1970}.
More precisely, following \cite[pp.~423--424]{LeoniSobolev2}, we assume that there exist $\eps_0>0$, $M_i>0$ for $i\in\N_0$, $i\le r$, $N_0\in\N$ and a  countable locally finite (see \cite[p.~637]{LeoniSobolev2}) open cover $\{\Omega_n\}_n$ of $\Gamma$   such that
\renewcommand{\labelenumii}{({\roman{enumii})}}
\begin{enumerate}
\item if $y\in \Gamma$ then $B_{\eps_0}(y)\subseteq \Omega_n$ for some $n\in\N$;
\item no point of $\R^N$ is contained in more than $N_0$ of the $\Omega_n$'s;
\item for each $n\in\N$ there exists a rigid motion $T_n:\R^N\to\R^N$ (that is an affine transformation of the form $T_nx=\xi_n+O_nx$, $\xi_n\in\R^N$, $O_n$ a rotation) and $f_n\in C^r(\R^{N-1})$ with
    \begin{equation}\label{2.1}
    \|D^\alpha f_n\|_{\infty,\R^{N-1}}\le M_m\quad\text{for all $n\in\N$ and $|\alpha|\le m$, $m\in\N_0$, $m\le r$,}
    \end{equation}
     such that $\Omega_n\cap \Omega=\Omega_n\cap T_n(V_n)$, where
     \begin{equation}\label{2.1bis}
     V_n=\{(y',y_N)\in \R^{N-1}\times \R: y_N>f_n(y')\}.
     \end{equation}
\end{enumerate}
\end{enumerate}

\noindent{\bf Notice.} In the paper we shall introduce several positive constants, denoted by $c_1,c_2,\ldots$, which may depend (without further notice) on $\Omega$, $\Gamma_1$, $N$, $N_0$, $\eps_0$, $(M_i)_i$, which are all considered to be fixed. These constants will not depend neither on the cover $\{\Omega_n\}_n$ of $\Gamma$
or on the collections $(f_n)_n$ and $(T_n)_n$.
Some of these constants will also depend on other quantities, say $\Upsilon_1,\ldots,\Upsilon_s$, introduced in the paper. In this case this dependence will be recorded by writing $c_i=c_i(\Upsilon_1,\ldots,\Upsilon_s)$.

Trivially we can assume, without restriction, that the $M_i$'s are increasing with respect to $i$,  with $M_0\ge 1$ and that all $\Omega_n$'s have nonempty intersection with $\Gamma$.
We shall denote by ${\cal I}_\Gamma\subseteq\N$ the index set for the cover $\{\Omega_n\}_n$.

Clearly $\Omega$ is also of class $C^{r-1,1}$, and hence uniformly Lipschitz, in the sense of \cite[p.~423]{LeoniSobolev2}.

Before proceeding we point out the following  result, which also shows that  when $r=\infty$, thanks to  assumption (A0),  $\Gamma$ is an uniformly regular manifold in the sense of \cite{Amann2013}. Its elementary  proof is given  for the reader's convenience.
\begin{lem}\label{lemma2.1}
In assumption (A0) we can take, without restriction,
\begin{equation}\label{forma}
 \Omega_n=T_n(W_n),\quad\text{and}\quad W_n=\{(y',y_N)\in Q_0\times\R: |f_n(y')-y_N|<\tau_0\},
\end{equation}
for all $n\in {\cal I}_\Gamma$, where $Q_0=2\tau_0\, Q^{N-1}=(-2\tau_0,2\tau_0)^{N-1}$ and $\tau_0=2NM_1\eps_0$.
\end{lem}
\begin{proof}
 We denote $\widetilde{h}_n(y')=(y', f_n(y'))$ and $\overline{h}_n=T_n\cdot \widetilde{h}_n$ on $\R^{N-1}$. By \eqref{2.1} $f_n$ is Lipschitz continuous in $\R^{N-1}$, with Lipschitz constant $\sqrt{N-1}M_1$. Hence, being $T_n$ a rigid motion, and $M_1\ge 1$, we have, for any $y',y''\in \R^{N-1}$ and $n\in {\cal I}_\Gamma$,
\begin{equation}\label{A.1}
  |y'-y''|\le |\widetilde{h}_n(y')-\widetilde{h}_n(y')|=|\overline{h}_n(y')-\overline{h}_n(y')|\le  \sqrt{N} M_1 \,|y'-y''|.
\end{equation}
We introduce, for $m\in\Z^{N-1}$,  the open cubes $Q_m=\tau_0 (m+2Q^{N-1})$ and $Q_m'=\tau_0 (m+Q^{N-1})$, where $\tau_0=\eps_0/6NM_1$, so that $Q_m'\subset Q_m$,  $\{Q_m, m\in\Z^{N-1}\}$ and  $\{Q_m', m\in\Z^{N-1}\}$ are locally finite open covers of $\R^{N-1}$ and no point of it is contained in more than $4^{N-1}$ of the  $Q_m$'s. Moreover, we introduce the rigid motions $S_m$, defined by $S_m y=(\tau_0 m,0)+y$, and $T_{n,m}=T_n\cdot S_m$.

We also set, recalling the notation \eqref{3.11}, for $n\in {\cal I}_\Gamma$, $J_n=\{m\in\Z^{N-1}: \overline{h}_n(Q_m')\cap \Omega_n^{\eps_0}\not=\emptyset\}$ and $I^*_\Gamma=\{(n,m)\in {\cal I}_\Gamma\times\Z^{N-1}: m\in J_n\}$. Since assumption (A0)--(i) reads as $\Gamma\subset\bigcup\limits_{n\in {\cal I}_\Gamma} \Omega_n^{\eps_0}$, we have
\begin{equation}\label{A.2}
 \Gamma\subset \bigcup_{(n,m)\in {\cal I}_\Gamma^*} \overline{h}_n(Q_m').
\end{equation}
We now claim that
\begin{equation}\label{A.3}
  \overline{h}_n(Q_m)\subseteq \Omega_n^{\eps_0/2}\quad\text{for all $(n,m)\in {\cal I}_\Gamma^*$.}
\end{equation}
To prove our claim we take, for any $y'\in Q_m$, $y''=\tau_0 m +\frac 12 (y'-\tau_0 m)$, so $|y'-y''|<\sqrt{N-1}\tau$ and $y''\in Q_m'$. Since $\overline{h}_n(Q_m')\cap \Omega_n^{\eps_0}\not=\emptyset$, there is $y'''\in Q_m'$ such that $\overline{h}_n(y''')\in \Omega_n^{\eps_0}$. Since $\text{diam}\,Q_m'=2\sqrt{N-1}\,\tau_0$, we have $|y'''-y'|<3\sqrt{N-1}\,\tau_0$
and so, by \eqref{A.1}, $|\overline{h}_n(y''')-\overline{h}_n(y')|<3NM_1\tau_0$. Hence, as $d(\overline{h}_n(y'''),\Omega_n^c)\ge \eps_0$, it follows that $d(\overline{h}_n(y'),\Omega_n^c)\ge \eps_0-3NM_1\tau_0=\eps_0/2$, proving our claim.
We now set, for $(n,m)\in {\cal I}_\Gamma^*$,
\begin{equation}\label{A.4}
W_{n,m}=\{(y',y_N)\in Q_m\times\R: |f_n(y')-y_N|<\tau_0\},\quad \Omega_{n,m}=T_n(W_{n,m}),
\end{equation}
so, also setting $f_{n,m}(y')=f_n(y'+\tau_0 m)$, we have
\begin{equation}\label{A.4bis}
W_{n,m}=S_m(W_{n,m}'), \quad W_{n,m}'=\{(y',y_N)\in Q_0\times\R: |f_{n,m}(y')-y_N|<\tau_0\}.
\end{equation}
We also denote
\begin{equation}\label{A.4ter}
V_{n,m}=\{(y',y_N)\in \R^{N-1}\times \R: y_N>f_{n,m}(y')\}=S_m^{-1}V_n.
\end{equation}
By \eqref{A.3}, as $T_n$ is a rigid motion and $M_1\ge 1$, for any $y=(y',y_N)\in W_{n,m}$ we have
\begin{align*}
d(T_ny,\Omega_n^c)\ge & d(\overline{h}_n(y'),\Omega_n^c)-|T_ny-\overline{h}_n(y')|= d(\overline{h}_n(y'),\Omega_n^c)-|y-\widetilde{h}_n(y')|\\\ge & d(\overline{h}_n(y'),\Omega_n^c)-|y_N-f_n(y')|\ge \tfrac {\eps_0}2- \tau_0=\left(\tfrac 12 -\tfrac 1{6NM_1}\right)\eps_0>0.
\end{align*}
Consequently, also using  \eqref{A.4}--\eqref{A.4ter},
\begin{equation}\label{A.5}
  \Omega_{n,m}\subseteq \Omega_n,\quad\text{so}\quad \Omega_{n,m}\cap\Omega=\Omega_{n,m}\cap T_{n,m}(V_{n,m})\quad\text{for all $(n,m)\in {\cal I}_\Gamma^*$.}
\end{equation}
We now set $\eps_0'=\tau_0/2NM_1$ and we claim that
\begin{equation}\label{A.6}
 \overline{h}_n(Q_m')\subseteq \Omega_{n,m}^{\eps_0'}\quad\text{for all $(n,m)\in {\cal I}_\Gamma^*$,}
\end{equation}
which, by \eqref{A.4}, as $T_n$ is a rigid motion, is equivalent to claim that $\widetilde{h}_n(Q_m')\subseteq W_{n,m}^{\eps_0'}$. We then have to prove that for any $\tilde{y}'\in Q_m'$ we have $d(\widetilde{h}_n(\tilde{y}'), W_{n,m}^c)\ge \eps_0'$. We suppose by contradiction that $d(\widetilde{h}_n(\tilde{y}'), W_{n,m}^c)< \eps_0'$, so there is $y=(y',y_N)\in W_{n,m}^c$ such that
$|\widetilde{h}_n(\tilde{y}')-y|< \eps_0'$, which trivially yields
$|\tilde{y}'-y'|< \eps_0'$ and $|f_n(\tilde{y}')-y_N|< \eps_0'$.
Now, since $y\not\in W_{n,m}$, by \eqref{A.4} two cases can occur:
\renewcommand{\labelenumi}{{\roman{enumi})}}
\begin{enumerate}
\item  $y'\not\in Q_m$. In this case, denoting $y'=(y_1',\ldots, y_{N-1}')$, $\tilde{y}'=(\tilde{y}_1',\ldots, \tilde{y}_{N-1}')$ and $m=(m_1,\ldots, m_{N-1})$, there is  $i=1,\ldots, N-1$ such that $|y_i'-m_i\tau_0|\ge 2\tau_0$. The, as $|\tilde{y}_i'-m_i\tau_0|<\tau_0$
    and $|\tilde{y}_i'-y'_i|\le |\tilde{y}'-y'|< \eps_0'=\tau_0/2NM_1<\tau_0$, we get a contradiction;
\item $y'\in Q_m$ and $|f_n(y')-y_N|\ge \tau_0$.  In this case, using the Lipschitz continuity of $f_n$ and as $M_1\ge 1$ and $N\ge 2$, we have
    \begin{align*}
    |f_n(y')-y_N|\le &|f_n(\tilde{y}')-y_N|+\sqrt{N-1} M_1 \,|y'-\tilde{y}'|\\
    \le &(1+\sqrt{N-1} M_1)\eps_0'=\tfrac{(1+\sqrt{N-1} M_1)\tau_0}{2NM_1}<\tfrac{\tau_0}2,
    \end{align*}
    so getting a contradiction and proving our claim.
\end{enumerate}
 Now, since ${\cal I}_\Gamma^*$ is countable, we write it as ${\cal I}_\Gamma^*=\{(n_l,m_l), l\in {\cal I}_\Gamma'\}$ for some   ${\cal I}_\Gamma'\subseteq\N$ and we set $\Omega_l'=\Omega_{n_l,m_l}$ for $l\in {\cal I}_\Gamma'$.
 Hence, by \eqref{A.2}, \eqref{A.4} and \eqref{A.6}, $\{\Omega_n', n\in {\cal I}_\Gamma'\}$ is an open cover of $\Gamma$ and condition (i) holds with $\eps_0'$ instead of $\eps_0$. By \eqref{A.5} also condition (ii) holds provided we replace $N_0$ with $N_0'=4^{N-1}N_0$, and $\{\Omega_n', n\in {\cal I}_\Gamma'\}$ is locally finite since both $\{\Omega_n, n\in {\cal I}_\Gamma\}$ and $\{Q_m, m\in\Z^{N-1}\}$ are locally finite. Moreover, recalling the definition of the $f_{n,m}$'s, using \eqref{A.4}--\eqref{A.4ter} and setting
 $$f'_l=f_{n_l,m_l},\quad W_l'=W_{n_l,m_l},\quad T_l'=T_{n_l,m_l},\quad\text{and}\quad V_l=V_ {n_l,m_l}\quad\text{for $l\in {\cal I}_\Gamma'$}$$
 also (iii) hold provided $f_n$, $T_n$ and $V_n$ are respectively replaced by $f_n'$, $T_n'$ and $V_n'$.
By \eqref{A.4}--\eqref{A.4bis} one also has $\Omega_n'=T_n'(W_n')$, while \eqref{forma} follows from \eqref{A.4bis}, concluding the proof.
\end{proof}

In the sequel we shall take, without further mention, the cover $\{\Omega_n, n\in {\cal I}_\Gamma\}$ in the form prescribed in Lemma~\ref{lemma2.1} and we shall denote $h_n(y')=T_n(y', f_n(y'))$ for $y'\in \R^{N-1}$ and $U_n=\Omega_n\cap\Gamma$, so  $h_n$ is an homeomorphism from $Q_0$ onto $U_n$. We shall denote $\zeta_n=({h_n}_{|Q_0})^{-1}$.
Hence, $\Gamma$ is a regularly embedded $C^r$ submanifold of $\R^N$,  an atlas being given by ${\cal A}:=\{(U_n, \zeta_n), n\in {\cal I}_\Gamma\}$.
Moreover, since the $U_n$'s are connected, given any
relatively clopen $\Gamma'\subseteq\Gamma$
\footnote{ By definition of relatively openness and relative closedness, a relatively clopen subset of $\Gamma$ is any union of connected components of $\Gamma$.}
 any $U_n$ can intersect either $\Gamma'$ or $\Gamma\setminus\Gamma'$. Consequently ${\cal A}_{\Gamma'}:=\{(U_n, \zeta_n), n\in {\cal I}_{\Gamma'}\}$, where
 ${\cal I}_{\Gamma'}=\{n\in {\cal I}_\Gamma: \Omega_n\cap\Gamma'\not=\emptyset\}$,
 is an atlas for $\Gamma'$.

\section{Geometric preliminaries}  We now recall some preliminary notions of Differential Geometry  and Geometric Analysis, well-known when $r=\infty$ (see \cite{Boothby, hebey, jost, taylor}). Since we shall apply them on any relatively open subset $\Gamma'$ of $\Gamma$, which only possesses a $C^r$ structure, we are also going to show how these notions adapt to the case $r<\infty$, referring to \cite{sternberg} for basic knowledge on the $C^r$ case.
Recalling these notions we shall also remark some uniform estimates in the atlas $\cal A$ introduced above.
In the sequel of the paper $\Gamma'$, $\Gamma''$ and so on, will denote (sometimes without further notice) relatively open subsets of $\Gamma$.
At first we shall consider the case $r\ge1$, and then we shall restrict to the case $r\ge 2$.
\subsection{Preliminaries when $r\ge 1$}
For all $p,q\in\N_0$ we shall denote by $T^p_q(V)$ the standard space of $p$-times contravariant and $q$-times covariant, or of type $(p,q)$, tensors on a complex vector space $V$. We shall use for all tensor--like objects (without further notice) the standard contraction conventions $T^0_0(V)=\C$, $T^p=T^p_0$ and $T_q=T^0_q$.
We shall respectively denote by $T(\Gamma')$, $T^*(\Gamma')$, and $T^p_q(\Gamma')$,  the standardly fiber--wise complexified (see \cite{roman}) tangent, cotangent and tensor bundles of type $(p,q)$, which are just topological ones when $r=1$, and by
$T_y(\Gamma')$, $T_y^*(\Gamma')$, and $T^p_q(T_y(\Gamma'))$ their fibers at $y\in\Gamma'$. We shall denote by $\overline{\phantom{a}}$ the conjugation on them, and by $\Real $ and $\Ima $ the real and imaginary parts. Since $T(\Gamma')=T^1(\Gamma')$ and $T^*(\Gamma')=T_1(\Gamma')$ in the sequel we shall not separately treat these cases, when possible.

A tensor field on $\Gamma'$ will be $u:\Gamma'\to T^p_q(\Gamma')$ with $u(y)\in T^p_q(T_y(\Gamma'))$ for all $y\in\Gamma'$, and we shall denote by  ${\cal F}^p_q(\Gamma')$ the space of tensor fields of type $(p,q)$ on $\Gamma'$. Since $T^0_0(\Gamma')=\Gamma'\times\C$ we shall conventionally identify, as usual, ${\cal F}(\Gamma')=\C^{\Gamma'}$. By $\otimes$ we shall indicate the tensor product.

All these notation and conventions will be also used on any open subset $\vartheta$  of the trivial manifold $\R^{N-1}$.
In this case, however, tensor bundles are trivial, since by identifying, as usual (see \cite[p.~73]{sternberg}), $T(\vartheta)=\vartheta\times\C^{N-1}$, one has $T^p_q(\vartheta)=\vartheta\times T^p_q(\C^{N-1})$. Hence in this case we shall also identify  ${\cal F}^p_q(\vartheta)=[T^p_q(\C^{N-1})]^\vartheta$.

Dealing with tensors it will be also useful to denote, for $p\in\N$,
$\J_p=\{1,\ldots,N-1\}^p$, $\J_0=\emptyset$ and $(i)=(i_1,\ldots,i_p)$ for $(i)\in\J_p$, $p\in\N$, $(i)$ disappearing when $p=0$.
In the sequel we shall use the Einstein summation convention for repeated indices. The standard frame fields for $T(\vartheta)$ and $T^*(\vartheta)$ will be respectively denoted by $\{\partial_1,\ldots,\partial_{N-1}\}$ and $\{dy^1,\ldots,dy^{N-1}\}$, so setting $\partial_{(i)}= \partial_{i_1}\otimes\ldots\otimes \partial_{i_p}$,  and $dy^{(j)}=dy^{j_1}\otimes\ldots\otimes dy^{j_{N-1}}$, the standard frame field for $T^p_q(\vartheta)$  is given by
 $\{\partial_{(i)}\otimes dy^{(j)}, (i)\in\J_p, (j)\in J_q\}$ when $p+q>0$, while $\{1\}$ is the standard field for $\vartheta\times \C$.
Each tensor $u\in T^p_q(\C^{N-1})$ (respectively each tensor field $u\in {\cal F}^p_q(\vartheta)$) can be uniquely written as $u=u^{(i)}_{(j)}\,\partial_{(i)}\otimes y^{(j)}$. The $u^{(i)}_{(j)}$'s will be called the components of $u$ .

 Arbitrarily fixing and implicitly using  in the sequel a flattening bijection from $\{1,\ldots,N_1\}\to \J_p\times \J_q$, $N_1=(N-1)^{p+q}$, for any
 $u\in T^p_q(\C^{N-1})$ we can take $\big(u^{(i)}_{(j)}\big)_{(i),(j)}\in \C^{N_1}$ and the map $u\mapsto \big(u^{(i)}_{(j)}\big)_{(i),(j)}$ establishes a bijective isomorphism $T^p_q(\C^{N-1})\simeq\C^{N_1}$.
 We shall  simply denote  by $|\cdot|$  the norm on $T^p_q(\C^{N-1})$ and  the consequent bundle norm on $T^p_q(\vartheta)$ induced by the standard Hermitian norm of $\C^{N_1}$ trough this isomorphism.
 Moreover, also using the trivial bijective isomorphism  associating to a $\C^{N_1}$-valued function $u$ on $\vartheta$ the vector $(u_1,\ldots, u_{N_1})$ of its components, the isomorphism  $T^p_q(\C^{N-1})\simeq\C^{N_1}$ trivially induces  the chain of bijective isomorphisms
 \begin{equation}\label{N1}
 {\cal F}^p_q(\vartheta) \simeq [\C^{N_1}]^\theta\simeq  [{\cal F}(\vartheta)]^{N_1},\qquad N_1=(N-1)^{p+q},
\end{equation}
which allows to regard tensor fields on $\vartheta$ as vector valued functions or as  vectors of scalar fields. When using in the sequel the notation $\simeq$
between spaces of tensors and spaces of functions on $\vartheta$, and recalling \eqref{N1}, we shall always refer to isomorphisms induced by those in \eqref{N1}.

Given a chart $(U,\varphi)$ for $\Gamma'$  we shall use, when useful, the standard notations $\varphi_*$ and $\varphi^*$ for the push--forward and pull-back bundle isomorphisms, and for the corresponding bijections between tensor fields and Hermitian forms on $\R^{N-1}$ and on $\Gamma'$, recalling that $\varphi_*u=u\cdot \varphi^{-1}$ and $\varphi^*v=v\cdot \varphi$ for scalar fields.
However, as common in Geometrical Analysis, when writing formulas "in local coordinates" we shall identify tensor fields on $U$ with their counterparts on $\varphi(U)$. In particular spaces of functions and tensor fields on $U$ and their counterparts on $\varphi(U)$ will be equal in this sense.
\begin{footnote}{When dealing with different charts this identification will not be made.}\end{footnote}
As a first example we shall denote the frame fields on $T^p_q(\varphi(U))$ induced by the chart and the components of a tensor field $u$ with respect to it with the same notations used in $\R^{N-1}$.

Moreover, for any $p\in\N$ and $q\in\N_0$, we set the  trace (or contraction) \label{contractionoperator}operator with respect to the last $p$ couples of covariant and contravariant indices  ${\cal T}_p:{\cal F}_{p+q}^p(\Gamma')\to {\cal F}_q(\Gamma')$
defined in local coordinates by ${\cal T}_p u^{(i)}_{(j,i)}=u^{(i)}_{(j,i)}$ for any $(j)\in\J_q$. By tensor transformation laws this operator is well-defined even if it is given in local coordinates.

 As usual we shall denote by $C^m(\Gamma')$  the space of (complex-valued) functions of class $C^m$ on $\Gamma'$ for $m=0,\ldots,r$, and by
$C^{m,p}_q(\Gamma')$ the space of tensor fields which components in any chart are of class $C^m$, with no contraction convention for the order of derivatives, so  $C^{0,p}_q$ will be never contracted to $C^p_q$ to avoid ambiguity. The same notation will be used on open subsets of $\R^{N-1}$. Due to transformation laws for tensors this latter notion is well-defined on $\Gamma'$ only  for $m=0,\ldots,r-1$ when $r<\infty$. To unify the notation it is convenient to set the number \begin{equation}\label{2.2}
r_1=r_1(r,p,q)=\begin{cases}
\infty,\qquad &\text{if $r=\infty$},\\
r,\qquad      &\text{if $r<\infty, \quad p+q=0$},\\
r-1,\qquad      &\text{if $r<\infty, \quad p+q>0$},\\
\end{cases}
\end{equation}
so the notion is well-defined for $m=0,\ldots, r_1$.
Using the same notion on any open subset $\vartheta$ of $\R^{N_1}$, by \eqref{N1} one trivially has, for all $m\in\N_0$,
\begin{equation}\label{N1regolare}
  C^{m,p}_q(\vartheta)\simeq C^m(\vartheta; \C^{N_1})\simeq [C^m(\vartheta)]^{N_1}.
\end{equation}
Trivially $\Gamma$ (and hence $\Gamma'$) inherits from $\R^N$ a Riemannian metric, uniquely extended to an Hermitian one on $T(\Gamma)$, in the sequel denoted by
$(\cdot,\cdot)_\Gamma$, given in local coordinates by $(u,v)_\Gamma=g_{ij}u^i\overline{v^j}$, with the $g_{ij}$'s  of class $C^{r-1}$ as  $(\cdot,\cdot)_\Gamma$.
When dealing with charts $(U_n,\zeta_n)$  in the atlas $\cal A$ given in Lemma~\ref{lemma2.1} we shall denote $g_{ij}$ as $g^n_{ij}$. Since, by the definition of $h_n$, $\partial_i h_n=e_i+\partial_i f_n e_N$, where $\{e_1,\ldots,e_N\}$ is the standard basis of $\R^N$, and by \cite[p.~362]{Folland}  we have
$(g^n_{ij})=(Jh_n)^tJh_n$, we get
\begin{equation}\label{2.3}
g_{ij}^n=\partial_ih_n\cdot \partial_jh_n=\delta_{ij}+\partial_if_n\partial_jf_n,\qquad\text{for $i,j=1,\ldots,N-1$, $n\in {\cal I}_\Gamma$},
\end{equation}
where $\delta_{ij}$ stands for the Kronecker symbol. By \eqref{2.3} each $g_{ij}^n$ is the restriction of $g_{ij}^n\in C^{r-1}(\R^{N-1})$. By \eqref{2.1} for $m\in\N$, $m\le r-1$ there is $c_1=c_1(m)>0$ such that
\begin{equation}\label{2.4}
|D^\alpha g^n_{ij}|\le c_1\quad\text{on $\R^{N-1}$ for all $|\alpha|\le m$ and $n\in {\cal I}_\Gamma$}.
\end{equation}
We shall also denote $(g^{ij})=(g_{ij})^{-1}$, $g=\Det (g_{ij})$, $(g_n^{ij})=(g^n_{ij})^{-1}$, $g_n=\Det (g^n_{ij})$. By \eqref{2.3} trivially
\begin{equation}\label{2.5}
|\xi|^2\le g_{ij}^n\xi_i\overline{\xi_j}\le \left[1+(N-1)M_1^2\right]\,|\xi|^2\quad\text{in $\R^{N-1}$,  $\quad \forall\xi\in\C^{N-1}$, $\forall n\in {\cal I}_\Gamma$.}
\end{equation}
Now, by Cauchy--Binet formula (see for example \cite[\S 4.6]{BroidaWilliamson}) and \eqref{2.3}, we have $\Det (g^n_{ij})=\sum_{i=1}^N \Det [(Jh_n)_i]^2$, where $(Jh_n)_1,\ldots (Jh_n)_N$ are the $N$ square submatrices of $Jh_n$ of order $N-1$. An elementary linear algebra calculation then implies (as it is well-known) that $g_n=1+|\nabla f_n|^2$. Hence, by \eqref{2.1}, we have
\begin{equation}\label{2.6}
  1\le g_n\le 1+(N-1)M_1^2\qquad\text{in $\R^{N-1}$,   \quad $\forall n\in {\cal I}_\Gamma$.}
\end{equation}
Hence, by Leibniz rule, for $m\in\N$, $m\le r-1$ there is $c_2=c_2(m)>0$ such that
\begin{equation}\label{2.7}
|D^\alpha g_n|\le c_2\quad\text{in $\R^{N-1}$ for all $|\alpha|\le m$ and $n\in {\cal I}_\Gamma$}.
\end{equation}
By \eqref{2.6}--\eqref{2.7} and generalized Faa di Bruno formula (see \cite[Theorem~11.54]{LeoniSobolev2} and \cite{HardyComb})
for $m\in\N$, $m\le r-1$ there is $c_3=c_3(m)>0$ such that
\begin{equation}\label{2.8}
|D^\alpha (1/g_n)|\le c_3\qquad\text{in $\R^{N-1}$ for all $|\alpha|\le m$ and $n\in\N$},
\end{equation}
so, by the classical adjugate formula for the inverse matrix, \eqref{2.4} and Leibniz rule,
for $m\in\N$, $m\le r-1$ there is $c_4=c_4(m)>0$ such that
\begin{equation}\label{2.9}
|D^\alpha g_n^{ij}|\le c_4\qquad\text{in $\R^{N-1}$ for all $|\alpha|\le m$ and $n\in {\cal I}_\Gamma$}.
\end{equation}
Now the metric $(\cdot,\cdot)_\Gamma$ induces the conjugate-linear (fiber-wise defined) Riesz isomorphism $\flat:T(\Gamma')\to T^*(\Gamma')$, with its inverse
$\sharp$, known as the musical isomorphisms in the real case, given by
\begin{equation}\label{2.10}
  \langle \flat u,v\rangle_{T(\Gamma')}=(v,u)_\Gamma,\qquad\text{for all $u,v\in T(\Gamma)$,}
\end{equation}
where $\langle \cdot,\cdot\rangle_{T(\Gamma')}$ denotes the fiber-wise defined duality pairing.
One then defines the induced bundle metric on $T^*(\Gamma')$ by the formula $(\alpha,\beta)_\Gamma=\langle \alpha,\sharp\beta\rangle_{T(\Gamma')}$ for all $\alpha,\beta\in T^*(\Gamma')$. By \eqref{2.10} one has
\begin{equation}\label{2.11}
  (\alpha,\beta)_\Gamma=(\sharp \beta, \sharp\alpha)_\Gamma,\qquad \text{for all $\alpha,\beta\in T^*(\Gamma')$},
\end{equation}
so, in local coordinates,
\begin{equation}\label{2.12}
  \flat u=g_{ij}\overline{u^j}dx^i, \quad\text{and}\quad \sharp u=g^{ij}\overline{\alpha_j}\partial_i, \qquad \text{for all $u\in T(\Gamma')$, $\alpha\in T^*(\Gamma')$.}
\end{equation}
By \eqref{2.12}, for all $0\le m\le r-1$, the isomorphisms $\flat$ and $\sharp$ induce bijections, denoted with the same symbols, between $C^{m,1}(\Gamma')$ and $C^m_1(\Gamma')$, and  $\sharp=\flat^{-1}$.

More generally $(\cdot,\cdot)_{\Gamma}$ extends to a bundle metric (still denoted by the same symbol) on $T^p_q(\Gamma')$ for any $p,q\in\N_0$ (see \cite[p.~442]{Amann2013} for details), given in local coordinates by
\begin{equation}\label{2.13}
  (u,v)_\Gamma=g_{(i) (i')}\,g^{(j)(j')}\, u^{(i)}_{(j)}\,\overline{v}^{(i')}_{(j')}\qquad\text{for all $u,v\in T^p_q(\Gamma')$,}
\end{equation}
where $g_{(i) (i')}=g_{i_1 i_1'}\ldots g_{i_p i_p'}$ and $g^{(j) (j')}=g^{j_1 j_1'}\ldots g^{j_q j_q'}$ for $(i), (i')\in\J_p$ and $(j), (j')\in\J_q$. In the sequel we shall denote $|\cdot|_\Gamma^2= (\cdot,\cdot)_\Gamma$ on any $T^p_q(\Gamma')$.
We remark that, by \eqref{2.13}, for all $u\in T^{p_1}_{q_1}(\Gamma')$ and $v\in T^{p_2}_{q_2}(\Gamma')$ one has
 \begin{equation}\label{modulotensoriale}
   |u\otimes v|_{\Gamma}=|u|_\Gamma |v|_\Gamma.
 \end{equation}
We also remark that, pointwise diagonalizing the matrix $(g^n_{ij})$, by \eqref{2.5} one gets
$$\left[1+(N-1)M_1^2\right]^{-1}\,|\xi|^2\le g^{ij}_n\xi_i\overline{\xi_j}\le |\xi|^2\quad\text{in $Q_0$,  $\quad \forall\xi\in\C^{N-1}$, $\forall n\in {\cal I}_\Gamma$,}
$$
so by \eqref{2.13}, for all $p,q\in\N_0$, $n\in\N$ and $u\in T^p_q(U_n)$ we have
\begin{equation}\label{2.14}
\left[1+(N-1)M_1^2\right]^{-q}\,|u|^2\le |u|_\Gamma^2\le \left[1+(N-1)M_1^2\right]^p\,|u|^2\quad\text{in $Q_0$.}
\end{equation}

The natural volume element associated to $(\cdot,\cdot)_\Gamma$ on $\Gamma$ is given, in local coordinates, by $\sqrt gdy_1\wedge\ldots\wedge dy^{N-1}$ so, in the atlas $\cal A$, by $\omega_n=\sqrt{1+|\nabla f_n|^2}\,dy_1\wedge\ldots\wedge dy^{N-1}$.

Hence the Riemannian measure associated to $\omega_n$ on the Borel $\sigma$-algebra $\beta(\Gamma)$ of $\Gamma$ coincides with the restriction to $\beta(\Gamma)$ of the Hausdorff measure ${\cal H}^{N-1}$ by the Area Formula (see \cite[Theorem~9.27]{LeoniSobolev2}). By extending it to the $\sigma$-algebra of measurable sets in $\Gamma$ by the Carathéodory method (see \cite{royden})  one then gets that the notion of ${\cal H}^{N-1}$ measurability on $\Gamma$ and the one of Lebesgue measurability in local coordinates coincide.
\begin{footnote}{
By the continuity of transition maps one can check the asserted property using the atlas ${\cal A}$. Since all $U_n$'s are Borel sets the check is then reduced to subsets of a fixed $U_n$. If $A\subseteq Q_0$ is Lebesgue measurable then $h_n(A)$ is ${\cal H}^{N-1}$ measurable by the already quoted Area Formula. Conversely, is $A\subseteq U_n$ is ${\cal H}^{N-1}$ measurable then being ${\cal H}^{N-1}$  Borel regular, there is $B\in \beta(\Gamma)$ such that $A\subseteq B$ and ${\cal H}^{N-1}(A)={\cal H}^{N-1}(B)$. Since ${\cal H}^{N-1}(\overline{U_n})<\infty$ by the Area Formula and Lemma~\ref{lemma2.1}, we have
${\cal H}^{N-1}(B\setminus A)=0$. Moreover, since $\zeta_n$ is trivially Lipschitz (with Lipschitz constant $1$)  we then get ${\cal H}^{N-1}(\zeta_n(B\setminus A))=0$ (see
(\cite[Proposition~C44]{LeoniSobolev2}). Since ${\cal H}^{N-1}$ is nothing but the Lebesgue measure on $\R^{N-1}$, $\zeta_n(B)$ is a Borel set and $\zeta_n(A)=\zeta_n(B)\setminus \zeta_n(B\setminus A)$, we conclude that $\zeta_n(A)$ is Lebesgue measurable.
}
\end{footnote}

We say that a tensor field on $\Gamma'$ is ${\cal H}^{N-1}$ measurable provided its components on any chart are Lebesgue measurable, and the class of ${\cal H}^{N-1}$ measurable tensor fields of type $(p,q)$, $p,q\in\N_0$, modulo ${\cal H}^{N-1}$-a.e.\ equivalence, will be denoted by $\cal M^p_q(\Gamma')$.

We shall use the same notation for tensor fields on an open subset $\vartheta$ of $\R^{N-1}$. In this case, however, the  isomorphisms in \eqref{N1} trivially restrict to
\begin{equation}\label{N2}
  \cal M^p_q(\vartheta)\simeq {\cal M}(\vartheta, \C^{N_1})\simeq [\cal M(\vartheta)]^{N_1},
  \end{equation}
 where as usual  $\cal M(\vartheta;\C^l)$ denotes for any $l\in\N$ the standard space of measurable  $\C^l$-valued functions modulo a.e.\ equivalence on $\vartheta$ and $\cal M(\vartheta)=\cal M(\vartheta;\C)$.

Trivially  the isomorphisms $\flat$ and $\sharp$, between $C^m$ vector and covector fields,  extend to bijections, denoted with the same symbols,
\begin{equation}\label{musicalimisurabili}
\flat: \cal M^1(\Gamma')\to \cal M_1(\Gamma'),\qquad \sharp: \cal M_1(\Gamma')\to M^1(\Gamma'),\qquad \sharp=\flat^{-1},
\end{equation}
and \eqref{2.12} continues to hold.

The Riemannian gradient operator $\nabla_\Gamma$ is defined  by
\begin{equation}\label{2.15}
 \nabla_\Gamma u=\sharp d_\Gamma \overline{u}\qquad\text{for $u\in C^1(\Gamma')$,}
\end{equation}
where $d_\Gamma$ denotes the differential on $\Gamma$ (see \cite{sternberg}), so $d_\Gamma u$ is nothing but the standard differential $du$ in local coordinates.
By \eqref{2.12} we have, in local coordinates,
\begin{equation}\label{2.16}
\nabla_\Gamma u=g^{ij}\partial_ju\partial_i \qquad \text{for $u\in C^1(\Gamma')$,}
\end{equation}
and, by \eqref{2.11} and \eqref{2.15},
\begin{equation}\label{2.17}
(\nabla_\Gamma u,\nabla_\Gamma v)_\Gamma=(d_\Gamma u,d_\Gamma v)_\Gamma\qquad \text{for $u,v\in C^1(\Gamma')$.}
\end{equation}
\subsection{Preliminaries when $r\ge 2$}\label{subsection2.3.1}In this case we define, as in \cite[pp.~305--312]{Boothby}, the covariant derivative $D_vu$ of $u\in C^1(\Gamma')$ at $y\in\Gamma'$, with respect to $v\in T_y(\Gamma')$ (here $\Gamma\subset\R^N$ is used), and we get that, in local coordinates, $D_vu=v^i(\partial_iuu^k+\Gamma_{ij}^k)\partial_k$, where (after a long but straightforward calculation) the Christoffel symbols $\Gamma_{ij}^k$ are still given, for $i,j,k=1,\ldots,N-1$, by the classical formula
\begin{equation}\label{2.18}
\Gamma_{ij}^k=\tfrac 12g^{kl}(\partial_jg_{il}+\partial_ig_{jl}-\partial_lg_{ij}),
\end{equation}
and the functions $\Gamma_{ij}^k$ are of class $C^{r-2}$. Moreover, denoting them as $\Gamma_{ij,n}^k$ when they are evaluated on the atlas ${\cal A}$, by \eqref{2.4} and \eqref{2.9}
for $m\in\N$, $m\le r-2$ there is $c_5=c_5(m)>0$ such that
\begin{equation}\label{2.19}
|D^\alpha \Gamma_{ij,n}^k|\le c_5\qquad\text{in $Q_0$, for $i,j,k=1,\ldots,N-1$, $|\alpha|\le m$ and $n\in {\cal I}_\Gamma$}.
\end{equation}
We hence define, as in \cite[Theorem~3.12 p.319 and pp.~391--397]{Boothby}, the parallel translation of $u\in C^{1,p}_q(\Gamma')$, $p+q>0$, its covariant derivative $D_vu$ at $y\in\Gamma'$, with respect to $v\in T_y(\Gamma')$ and finally we set the covariant derivative of u
\begin{footnote}{Usually denoted by $\nabla u$ in the literature: however, the latter symbol indicates throughout this article the gradient in $\Omega$.}\end{footnote}
as a tensor field of type $(p,q+1)$ by defining, for all $\alpha_1,\ldots,\alpha_p\in T^*(\Gamma')$, $v^1,\ldots,v^{q+1}\in T(\Gamma')$,
\begin{equation}\label{2.20}
  D_\Gamma u[\alpha_1,\ldots,\alpha_p, v^1,\ldots,v^{q+1}]=D_{v^1}u[\alpha_1,\ldots,\alpha_p, v^2,\ldots,v^{q+1}].
  \end{equation}
In local coordinates then, for all $(i)\in\J_p$ and $(j)\in\J_q$ one has
\begin{equation}\label{2.21}
  \begin{aligned}
  (D_\Gamma u)_{(j)}^{(i)}=\quad &\partial_{j_1} u^{(i)}_{j_2,\ldots,j_{q+1}}-\sum_{k=2}^{q+1}\, \Gamma_{j_1 j_k}^l u^{(i)}_{j_2,\ldots,j_{k-1},l, j_{k+1},\ldots, j_{q+1}}\\
  +&\sum_{k=1}^p \,\Gamma_{j_1 l}^{i_k} \,u^{i_1,\ldots,i_{k-1}, l, i_{k+1},\ldots,\i_p}_{j_2,\ldots,j_{q+1}},
\end{aligned}
\end{equation}
the corresponding summation being suppressed when $p=0$ or $q=0$, so setting $D_\Gamma=d_\Gamma$ on scalar fields, the last formula holds true for any $p,q\in\N_0$ and $D_\Gamma$ maps $C^{m,p}_q(\Gamma')$ into $C^{m-1,p}_{q+1}(\Gamma')$, also when $r=1$, provided $1\le m\le r_1$, $r_1$ being given by \eqref{2.2}.

By \eqref{2.20} or \eqref{2.21} one trivially gets by induction  the generalized Leibniz formula
\begin{equation}\label{2.22}
 D_\Gamma^m(u\otimes v)=\sum_{i=0}^m \binom mi \,D_\Gamma^i u\,\otimes \,D_\Gamma^{m-i}v
 \end{equation}
for $u\in C^{m,p_1}_{q_1}(\Gamma')$, $v\in C^{m,p_2}_{q_2}(\Gamma')$, $p_1,p_2,q_1,q_2\in\N_0$, $1\le m\le r_1$.
By the same reason, for any $p\in\N$, $q\in\N_0$ and $1\le m\le r-1$, the covariant derivative of order $m$ commutes with the contraction operator ${\cal T}_p$ defined at page~\pageref{contractionoperator}, that is
\begin{equation}\label{contractioncommuta}
 D^m_\Gamma\cdot {\cal T}_p={\cal T}_p\cdot D_\Gamma^m\qquad\text{in $C^{m,p}_{p+q}(\Gamma')$.}
\end{equation}

We now introduce the Riemannian divergence operator $\DivGamma$ by setting, for $u\in C^{1,1}(\Gamma')$, $\DivGamma u$ as the unique continuous function on $\Gamma'$ such that $(\DivGamma u)\omega=d_{\rm ext}(\omega\lrcorner  u)$ on $\Gamma'$, where $d_{\rm ext}$ denotes the exterior derivative on forms and $\lrcorner$ the interior product (see \cite{taylor} for details). Since in local coordinates
$$\omega\lrcorner u=\sum_{i=1}^{N-1} \sqrt g u^i\,dy^1\wedge \ldots\wedge dy^{i-1}\wedge dy^{i+1}\wedge\ldots \wedge dy^{N-1},$$
we get, in the same coordinates, the classical formula
\begin{equation}\label{2.23}
\DivGamma u=g^{-1/2}\partial_i (g^{1/2}u^i)\qquad\text{for all $u\in C^{1,1}(\Gamma')$.}
\end{equation}
Hence, integrating by parts on coordinate neighborhoods and using a $C^2$  partition of the unity (see \cite[Theorem~4.1, p.~57]{sternberg}),
for any for $u\in C^1(\Gamma')$, $v\in C^{1,1}(\Gamma')$ such that $uv$ is compactly supported we get
\begin{equation}\label{2.24}
\int_{\Gamma'} (\nabla_\Gamma u,v)_\Gamma=-\int_{\Gamma'}u \DivGamma \overline{v}.
\end{equation}
We also introduce the Laplace--Beltrami operator $\Delta_\Gamma:=\DivGamma\!\!\!\nabla_\Gamma$. By \eqref{2.16} and \eqref{2.23} we have, in local coordinates,
\begin{equation}\label{LaplaceBeltrami}
  \Delta_\Gamma u=g^{-1/2}\partial_i(g^{1/2}g^{ij}\partial_j u),\qquad\text{for all $u\in C^2(\Gamma')$}.
\end{equation}

\section{Lebesgue spaces on $\Gamma$} For $1\le \tau\le \infty$ we shall denote by $L^\tau (\Gamma')$ the (complex) Lebesgue space on $\Gamma'$ with respect to ${\cal H}^{N-1}$. More generally, for any $p,q\in\N_0$ we shall denote $L^{\tau,p}_q(\Gamma')=\{u\in {\cal M}^p_q(\Gamma'): |u|_\Gamma\in L^\tau(\Gamma')\}$, which is a Banach space with respect to the trivially related norm. For simplicity we shall denote all $L^\tau$-type norms on $\Gamma'$ by $\|\cdot\|_{\tau,\Gamma'}$. Moreover we shall also use the notation  $\Gamma''\subset\subset\Gamma'$ and $L_{q,\loc}^{\tau,p}(\Gamma')$ in the standard sense.

 The same notation will be used on open subsets $\vartheta$ of $\R^{N-1}$. Trivially the isomorphism  in \eqref{N2} restrict to bijective isomorphisms
 \begin{equation}\label{N3}
  L^{\tau,p}_q(\vartheta)\simeq L^\tau(\vartheta;\C^{N_1})\simeq [L^\tau(\vartheta)]^{N_1}.
 \end{equation}
When endowing the second and third space in \eqref{N3} with their standard norms, the first isomorphism is isometric, while the second one is bounded with bounded inverse.  For simplicity we shall endow the three spaces in \eqref{N3}  with the norm $\|\cdot\|_{\tau, \vartheta}=\|\cdot\|_{L^\tau(\vartheta,\C^{N_1})}$.

 We remark that, by \eqref{modulotensoriale}, the standard H\"{o}lder inequality extends to tensor fields by extending the product to the tensor product.

 The space  $L^{\tau,p}_q(\Gamma')$ is easily characterized in local coordinates when $\Gamma'$ is contained in a coordinate neighborhood $(U,\varphi)$, provided \begin{equation}\label{2.26}
\Gamma'\subset\subset U\text{ or }(U,\varphi)\in{\cal A},
 \end{equation}
 where ${\cal A}$ is the atlas defined after Lemma~\ref{lemma2.1}.
 Indeed, in the first case, using compactness instead that \eqref{2.1}, one still gets estimates \eqref{2.4}--\eqref{2.6}, \eqref{2.7}--\eqref{2.9} and \eqref{2.14}, with different positive constants depending on  $(U,\varphi)$, so
 $L^{\tau,p}_q(\Gamma')=L^{\tau,p}_q(\varphi(\Gamma'))$ with equivalent norms.

 In the second one, setting $\vartheta_n=\zeta_n(\Gamma')$,   by \eqref{2.5}--\eqref{2.6} and \eqref{2.14} we still get $L^{\tau,p}_q(\Gamma')=L^{\tau,p}_q(\vartheta_n)$ and the uniform estimate
 \begin{equation}\label{2.29}
\left[1+(N-1)M_1^2\right]^{-\frac q2}\,\|u\|_{\tau,\vartheta_n} \le \|u\|_{\tau, \Gamma'}\le \left[1+(N-1)M_1^2\right]^{\frac p2+\frac 1\tau}\,\|u\|_{\tau,\vartheta_n},
\end{equation}
where we mean $1/\infty=0$ when $\tau=\infty$, which also  holds for all $u\in {\cal M}^p_q(\Gamma')$.

\chapter{Sobolev spaces on $\Gamma$ and operators} \label{section3}
This chapter is devoted to Sobolev spaces, of both functions and tensor fields, of real nonnegative order on relatively open subsets $\Gamma'$ of $\Gamma$ and to boundedness
properties of linear operators related to problem \eqref{1.1} between them. The need of such a study  has been already motivated in \S~\ref{intro}. We shall generalize the approach in \cite{Amann2013}, simplifying it by skipping Sobolev spaces of tensors of negative order.
\begin{footnote}{They could be treated as well as in quoted reference, but we do not need them here.}\end{footnote}
 One of the goals, and difficulties,  in our study is to prove that the approaches in \cite{Amann2013} and in \cite{LeoniSobolev2} are equivalent. Moreover we shall also deal with the spaces $W^{m,\infty}$, $m\in\N$, used in  assumptions (A1--3) and not considered in \cite{Amann2013}.
\section{Spaces of integer order} Essentially following \cite{hebey} we let, for $\Gamma'\subseteq \Gamma$, $p,q\in\N_0$, $1\le \tau<\infty$ and $m\in\N_0$, $m\le r_1$, $r_1$ being given by \eqref{2.2},
\begin{equation}\label{3.1}
C^{m,\tau,p}_q(\Gamma'):=\{u\in C^{m,p}_q(\Gamma'): D^i_\Gamma u\in L^{\tau,p}_{q+i}(\Gamma')\quad\text{for $i=0,\ldots,m$}\}
  \begin{footnote}{The spaces $C^{m,p}_q(\Gamma')$ and $L^{\tau,p}_{q+i}(\Gamma')$ were respectively defined in \eqref{N1regolare} and \eqref{N3}.}\end{footnote},
\end{equation}
endowed with  the norm
\begin{equation}\label{3.2}
 \|u\|_{m,\tau,\Gamma'}:=\left( \sum_{i=0}^m\|D_\Gamma ^i u\|_{\tau,\Gamma'}^\tau \right)^{1/\tau}\negquad,
\end{equation}
and $W^{m,\tau,p}_q(\Gamma')$ as the completion of $X:=\left (C^{m,\tau,p}_q(\Gamma'),\|\cdot\|_{m,\tau,\Gamma'}\right)$.
According to standard contraction conventions on tensor orders, when $p=q=0$, $W^{m,\tau,p}_q(\Gamma')$ is simply denoted as $W^{m,\tau}(\Gamma')$.
The space $W^{m,\tau,p}_q(\Gamma')$ is a Banach space and we shall naturally identify it with the subspace of $L^{\tau,p}_q (\Gamma')$  consisting of the $L^\tau$-limits of Cauchy sequences in $X$. This identification is justified by the following two remarks:
\begin{itemize}
\item[--] any Cauchy sequence in $X$ trivially converges in $L^{\tau,p}_q(\Gamma')$;
\item[--] any Cauchy sequence $(u_n)_n$ in $X$ converging to zero in the $L^\tau$-norm  also does in the $W^{m,\tau}$-norm, as we are going to show. Trivially  $D_\Gamma^iu_n$ converges in $L^{\tau,p}_{q+i} (\Gamma')$ to some $\xi_i$ for all $i=0,\ldots, m$, with $\xi_0=0$. We now suppose, by induction, that $\xi_0=\ldots \xi_i=0$ for $i\le m-1$. Then, using \eqref{2.21}, we easily get that, in any coordinate neighborhood $(U,\varphi)$ and for all $(i')\in\J_p$, $(j')\in\J_q$, $\psi\in C^\infty_c(\varphi(U))$, $l=1,\ldots, N-1$
$$\int_{\varphi(U)}\negqquad (\xi_{i+1})^{(i')}_{(j')}\psi=\lim_n \int_{\varphi(U)}\negqquad\partial_l (D^i_\Gamma u_n)^{(i')}_{(j')}\psi
=-\lim_n \int_{\varphi(U)}\negqquad (D^i_\Gamma u_n)^{(i')}_{(j')}\partial_l\psi=0,
$$
so $\xi_{i+1}=0$.
\end{itemize}
Now for all $u\in W^{m,\tau,p}_q(\Gamma')$, taking $(u_n)_n$ a Cauchy sequence in $X$ that converges to it, we can set $D_{\Gamma'} u$ as the limit in $W^{m-1, \tau,p}_{q+1}(\Gamma')$ of the sequence $D_\Gamma u_n$, and $D_{\Gamma'} u$ is trivially the restriction of $D_\Gamma u$.
In this way the operator $D_\Gamma$ extends by construction to
\begin{equation}\label{3.2bis}
D_\Gamma\in {\cal L}(W^{m,\tau,p}_q(\Gamma'); W^{m-1,\tau,p}_{q+1}(\Gamma')),\qquad\text{provided $1\le m\le r_1$,}
\end{equation}
and the norm in $W^{m,\tau,p}_q(\Gamma')$ is still given by \eqref{3.2}.
Moreover the definition \eqref{2.15} of the Riemannian gradient $\nabla_\Gamma$ extends to
\begin{equation}\label{3.2ter}
\nabla_\Gamma\in {\cal L}(W^{1,\tau}(\Gamma');L^{\tau,1}(\Gamma')),
\end{equation}
with the same localization property of $D_\Gamma$, and \eqref{2.16}--\eqref{2.17} continue to hold.
Boundedness of $\nabla_\Gamma$ between higher order spaces will be proved in Lemma~\ref{lemma3.4} below.

We shall denote $H^{m,p}_q(\Gamma')=W^{m,2,p}_q(\Gamma')$, and in particular $H^m(\Gamma')=W^{m,2}(\Gamma')$ as usual,
which is an Hilbert space with respect to the inner product
\begin{equation}\label{3.3}
  (u,v)_{m,\Gamma'}=\sum_{i=0}^m \int_{\Gamma'}(D_\Gamma^i u,D_\Gamma^i v)_\Gamma.
\end{equation}
By the extended version of \eqref{2.17}, when $p=q=0$ and $m=1$, \eqref{3.3} reads as
\begin{equation}\label{3.4}
  (u,v)_{1,\Gamma'}=\int_{\Gamma'}u\overline{v}+\int_{\Gamma'}(\nabla_\Gamma u,\nabla_\Gamma v)_\Gamma\qquad\text{for all $u,v\in H^1(\Gamma')$.}
\end{equation}

In the same way we also set, for $m\in\N$, $p,q\in\N_0$ and $1\le \tau<\infty$, the spaces $W^{m,\tau,p}_q(\vartheta)$ and $H^{m,p}_q(\vartheta)=W^{m,2,p}_q(\vartheta)$,  denoting by $\|\cdot\|_{m,\tau,\vartheta}$ its norm,  on
any open subset $\vartheta$ of the trivial manifold $\R^{N-1}$. Since the covariant derivative is nothing but the standard differential, using \eqref{N1regolare} and  Meyers--Serrin's Theorem, the isomorphisms in \eqref{N3} restrict to isomorphisms
\begin{equation}\label{N4}
 W^{m,\tau,p}_q(\vartheta)\simeq W^{m,\tau}(\vartheta;\C^{N_1})\simeq [W^{m,\tau}(\vartheta)]^{N_1},
\end{equation}
and the norms induced by them are all equivalent.

As for Lebesgue spaces, the spaces $W^{m,\tau,p}_q(\Gamma')$ are easily characterized when $\Gamma'$ is contained in a coordinate neighborhood $(U,\varphi)$ and \eqref{2.26} holds. Since we need uniform estimates, we  explicitly state this result. We give its  proof for the reader's convenience.
\begin{lem}\label{lemma3.1} Suppose that $\Gamma'$ is contained in a coordinate neighborhood $(U,\varphi)$, that \eqref{2.26} holds and let  $p,q\in\N_0$, $1\le \tau<\infty$, $m\in\N_0$, $m\le r_1$.
Then
$W^{m,\tau,p}_q(\Gamma') =W^{m,\tau,p}_q(\varphi(\Gamma'))$, with equivalence of norms.
Moreover there is $c_6=c_6(m, p,q,\tau)>0$ such that for all $n\in {\cal I}_\Gamma$, when $(U,\varphi)=(U_n,\zeta_n)\in {\cal A}$,
\begin{equation}\label{3.5}
  c_6^{-1} \|\cdot\|_{m,\tau,\vartheta_n}\le \|\cdot\|_{m,\tau, \Gamma'} \le   c_6 \|\cdot\|_{m,\tau,\vartheta_n}\quad\text{on $W^{m,\tau,p}_q(\Gamma')$,}
 \end{equation}
 where we are still denoting $\vartheta_n=\zeta_n(\Gamma')$.
\end{lem}
\begin{proof}
In the proof, for the sake of clearness,  we shall not identify $\varphi_*u$ with $u$ and $\varphi^*v$ with $v$.
Consequently the assertion $W^{m,\tau,p}_q(\Gamma') =W^{m,\tau,p}_q(\varphi(\Gamma'))$ in the statement reads as
$\varphi_*[W^{m,\tau,p}_q(\Gamma')]\subseteq W^{m,\tau,p}_q(\varphi(\Gamma'))$ and
$\varphi^*[W^{m,\tau,p}_q(\varphi(\Gamma'))]\subseteq W^{m,\tau,p}_q(\Gamma')$. Moreover, being $\varphi$ fixed, we shall simply denote
$\varphi_*u=\widetilde{u}$. Moreover, to make the proof more transparent, we shall make explicit the dependence of some constants on $m$ and $q$.

We preliminarily remark that the norm $\|\cdot\|_{m,\tau,\Gamma'}$ defined in \eqref{3.2} and the norm $\|\cdot\|_{m,\tau,\varphi(\Gamma')}$ are respectively equivalent to the norms defined by
\begin{equation}\label{B1}
\n u\n_{m,\tau,\Gamma'}^\tau=\sum_{i=0}^m{\textstyle\binom mi}\|D^i_\Gamma u\|_{\tau,\Gamma'}^\tau
,\quad  \n u\n_{m,\tau,\varphi(U)}^\tau=\sum_{i=0}^m{\textstyle\binom mi}\|d^i u\|_{\tau,\varphi(U)}^\tau,
\end{equation}
the equivalence constant depending only on $m$ and $\tau$.
Consequently \eqref{3.5} is equivalently rewritten as
\begin{equation}\label{3.5rewritten}
  c_7^{-1}(m,q) \n\cdot\n_{m,\tau,\vartheta_n}\le \n\cdot\n_{m,\tau, \Gamma'} \le   c_7(m,q) \n\cdot\n_{m,\tau,\vartheta_n}\quad\text{on $W^{m,\tau,p}_q(\Gamma')$,}
 \end{equation}
 for some $c_7(m,q)=c_7(m,\tau,p,q)>0$. Trivially, these equivalent norms verify the recursive property
\begin{equation}\label{B2}
\begin{aligned}
\n u\n_{m+1,\tau,\Gamma'}^\tau=\quad &\n u\n_{m,\tau,\Gamma'}^\tau+\|D_\Gamma u\|_{m,\tau,\Gamma'}^\tau,  \\
\n u\n_{m+1,\tau,\varphi(U)}^\tau=\quad &\n u\n_{m,\tau,\varphi(U)}^\tau+\|d u\|_{m,\tau,\varphi(U)}^\tau,
\end{aligned}
\end{equation}
for $m=0,\ldots,r_1-1$.

We first consider the case $(U,\varphi)=(U_n,\varphi_n)\in\cal{A}$. By \eqref{2.21} for all $u\in C^{1,p}_q(\Gamma')$, $(i)\in\J_p$ and $(j)\in\J_q$ we have
 $$ (\widetilde{D_\Gamma u}-d\widetilde{u})_{(j)}^{(i)}=-\sum_{k=2}^{q+1}\, \Gamma_{j_1 j_k}^l u^{(i)}_{j_2,\ldots,j_{k-1},l, j_{k+1},\ldots, j_{q+1}}
  +\sum_{k=1}^p \,\Gamma_{j_1 l}^{i_k} \,u^{i_1,\ldots,i_{k-1}, l, i_{k+1},\ldots,i_p}_{j_2,\ldots,j_{q+1}},
$$
with the same caveats made after \eqref{2.21}. Consequently, using \eqref{2.19} and Leibniz formula, for any $m=0,\ldots,r_1-1$ there is
$c_8=c_8(m,p,q)$ such that for all $u\in C^{m+1,p}_q(\Gamma')$. $0\le |\alpha|\le m$, $(i)\in\J_p$, $(j)\in \J_q$ and $n\in\cal{I}_{\Gamma}$
\begin{equation}\label{B3}
|D^\alpha (\widetilde{D_\Gamma u}-d\widetilde{u})_{(j)}^{(i)}|\le c_8\sum_{(i')\in\J_p}\sum_{(j')\in\J_q}\sum_{|\beta|\le m}|D^\beta\widetilde{u} _{(j')}^{(i')}|\qquad\text{on $\vartheta_n$.}
\end{equation}
Consequently there is $c_9=c_9(m,p,q,\tau)\ge 1$ such that
\begin{equation}\label{B4}
\n \widetilde{D_\Gamma u}-d\widetilde{u}\n_{m,\tau.\vartheta_n}^\tau\le c_9\,\n \widetilde{u}\n_{m,\tau,\vartheta_n}^\tau
\quad\text{for all $u\in C^{1,p}_q(\Gamma')$ and $n\in\cal{I}_{\Gamma}$,}
\end{equation}
where the left and hence the right hand side may be infinite.

We now claim that \eqref{3.5rewritten} hold for all $u\in C^{m,p}_q(\Gamma')$, the norms being finite or not. We shall prove our claim by induction on $m=0,\ldots,r_1$. The case $m=0$ immediately follows from \eqref{2.29}, so we suppose that \eqref{3.5rewritten} holds for $m\le r_1-1$.
Hence, by also using \eqref{B2} and \eqref{B4},  for all $u\in C^{m+1,p}_q(\Gamma')$  we have
\begin{equation}\label{B4bis}
\begin{aligned}
&\n u\n_{m+1,\tau,\Gamma'}^\tau
\le  c_7^\tau(m,q)\n \widetilde{u}\n_{m,\tau,\vartheta_n}^\tau+c_7^\tau(m,q+1)\n \widetilde{D_\Gamma u}\n_{m,\tau,\vartheta_n}^\tau\\
&\le  c_7^\tau(m,q)\n \widetilde{u}\n_{m,\tau,\vartheta_n}^\tau+2^{\tau-1}c_7^\tau(m,q+1)\left(\n d\widetilde{u}\n_{m,\tau,\vartheta_n}+c_9\n \widetilde{u}\n_{m,\tau,\vartheta_n}^\tau\right)\\
&\le c_7^\tau(m+1,q)\,\n \widetilde{u}\n_{m+1,\tau,\vartheta_n}^\tau,
\end{aligned}
\end{equation}
setting $c_7^\tau(m+1,q)=c_7^\tau(m,q)+2^{\tau-1}\,c_9 \,\,c_7^\tau(m,q+1)$.
The same arguments also give the reverse inequality, proving our claim.

Since $C^{m,\tau,p}_q$ is dense in $W^{m,\tau,p}_q(\Gamma')$ by construction and $W^{m,\tau,p}_q(\vartheta_n)$ is complete, by the just proved claim we then get that ${\zeta_n}_*[W^{m,\tau,p}_q(\Gamma')]\subseteq W^{m,\tau,p}_q(\vartheta_n)$. By \eqref{N4} and Meyers--Serrin Theorem, $C^{m,\tau,p}_q(\vartheta_n)\cap W^{m,\tau,p}_q(\vartheta_n)$ is dense in $W^{m,\tau,p}_q(\vartheta_n)$. Hence, as $W^{m,\tau,p}_q(\Gamma')$ is complete and ${\zeta_n}_*$ is a bijective isomorphism between $C^{m,\tau,p}_q(\Gamma')$ and $C^{m,\tau,p}_q(\vartheta_n)$, with inverse $\zeta_n^*$, by our claim we also get that
$${\zeta_n}^*[W^{m,\tau,p}_q(\vartheta_n)]\subseteq W^{m,\tau,p}_q(\Gamma').$$ Previous arguments also show that \eqref{3.5rewritten} extends to $W^{m,\tau,p}_q(\Gamma')$.

We now consider the case $\Gamma'\subseteq U$, in which the estimate \eqref{2.19} is replaced by its non-uniform version
\begin{equation}\label{tridente}
|D^\alpha \Gamma_{ij}^k|\le c_{10}\quad\text{in $\varphi(U)$, \quad for $i,j,k=1,\ldots,N-1$ and $|\alpha|\le m\le r-2$,}
\end{equation}
where $c_{10}=c_{10}(m,\varphi,U)$. To prove \eqref{tridente} one uses \eqref{2.18} and the non-uniform versions of \eqref{2.4} and \eqref{2.9} already pointed out in the case \eqref{2.26}. Hence we simply repeat the proof of the previous case allowing the constants $c_7$, $c_8$ and $c_9$ to depend also on $\varphi$ and $\Gamma'$.
\end{proof}

It is convenient to introduce, for $p,q\in\N_0$, $1\le \tau<\infty$ and $m\in\N_0$, $m\le r_1$,
the standardly defined space $W^{m,\tau,p}_{q,\loc}(\Gamma')$ which, since $\Gamma'$ is trivially $\sigma$-compact, is a Fréchet space.
Then the operator $D_\Gamma$ in \eqref{3.2bis} trivially extends to
$$D_\Gamma\in {\cal L}(W^{m,\tau,p}_{q,\loc}(\Gamma'); W^{m-1,\tau,p}_{q+1,\loc}(\Gamma')),$$
for $p,q\in\N_0$, $1\le \tau<\infty$, $m\in\N_0$, $m\le r_1$,  the definition being consistent with possible different values of $m$ and $\tau$.

Moreover \eqref{2.21} extends, on any coordinate neighborhood $(U,\varphi)$, to $W^{1,\tau,p}_{q,\loc}(U)$.
We then have the following easy result. Its proof is given for the reader's convenience.
\begin{lem}\label{lemma3.2} For $p,q\in\N_0$, $1\le \tau<\infty$ and  $m\in\N_0$, $m\le r_1$ we have
\begin{equation}\label{3.7}
W^{m,\tau,p}_q(\Gamma') =\{u\in W^{m,1,p}_{q,\loc}(\Gamma'): \, D_\Gamma^i u \in L^{\tau,p}_{q+i}(\Gamma')\quad \text{for $i=0,\ldots, m$}\}.
\end{equation}
Moreover $W^{m,\tau,p}_q(\Gamma')\cap C^{r_1,p}_q(\Gamma')$ is dense in $W^{m,\tau,p}_q(\Gamma')$.
\end{lem}
\begin{proof} Since $\|\cdot\|_{m,\tau,\Gamma'}$ is defined, finite or not, for any $u\in W^{m,1,p}_{q,\loc}(\Gamma')$,
after setting $Y^{m,\tau,p}_q(\Gamma')=\{u\in W^{m,1,p}_{q,\loc}(\Gamma'): \|u\|_{m,\tau,\Gamma'}<\infty\}$
we can rewrite \eqref{3.7} as $Y^{m,\tau,p}_q(\Gamma')\subseteq W^{m,\tau,p}_q(\Gamma')$, the reverse inclusion being trivial. When referring to
\eqref{3.7} in this proof we shall refer to this equivalent form.

We first claim that \eqref{3.7} holds when $\Gamma'$ is contained in a coordinate neighborhood $(U,\varphi)$. To prove this claim  we remark that, by repeating (once again) the arguments in the proof of Lemma~\ref{lemma3.1} there is $c_{11}=c_{11}(m,p,q,\tau,\varphi,U)>0$ such that
\begin{equation}\label{B5}
  c_{11}^{-1} \|\cdot\|_{m,\tau,\varphi(U)}\le \|\cdot\|_{m,\tau, \Gamma'} \le   c_{11} \|\cdot\|_{m,\tau,\varphi(U)}\quad\text{on $W^{m,1,p}_{q,\loc}(\Gamma')$.}
 \end{equation}
Now let $u\in Y^{m,\tau,p}_q(\Gamma')$. By Lemma~\ref{lemma3.1} then  $u\in Y^{m,1,p}_{q,\loc}(\varphi(U))$ and, by \eqref{B5}, $\|u\|_{m,\tau,\varphi(U)}<\infty$. By \eqref{N4} and the standard definition of $W^{m,\tau}(\varphi(U))$ we then get that $u\in W^{m,\tau,p}_q(\varphi(U))$.
By Lemma~\ref{lemma3.1} then $u\in W^{m,\tau,p}_q(\Gamma')$, proving our claim.

We now consider the general case, in which both parts of the statement follow provided $W^{m,\tau,p}_q(\Gamma')\cap C^{r_1,p}_q(\Gamma')$ is dense in
$Y^{m,\tau,p}_q(\Gamma')$ with respect to $\|\cdot\|_{m,\tau,\Gamma'}$.

To prove this fact we recall (see \cite[Lemma~4.2,~p.~56]{sternberg}) that, being $\Gamma'$ a $C^r$-manifold, it admits an atlas $\{(\widetilde{V_n},\widetilde{\zeta_n}), n\in\N\}$ such that $\{\widetilde{V_n}, n\in\N\}$
is locally finite, $\widetilde{\zeta_n}(\widetilde{V_n})=B_3$ (in $\R^{N-1}$) for all $n$ and $\{\widetilde{U_n},n\in\N\}$ is an open cover of $\Gamma'$, where $\widetilde{U_n}:=\widetilde{\zeta_n}^{-1}(B_1)$. Trivially, $\widetilde{U_n}\subset\subset \widetilde{V_n}$ for all $n$.
By \cite[Theorem~4.1,~p.~57]{sternberg} there is a $C^r$-partition of the unity $\{\widetilde{\varphi_n}, n\in\N\}$ subordinated to $\{\widetilde{U_n},n\in\N\}$.

We now fix $u\in Y^{m,\tau,p}_q(\Gamma')$ and $\eps\in (0,1/2)$, so   $\widetilde{\varphi_n} u\in W^{m,1,p}_{q,\loc}(\Gamma')$ and $\|\widetilde{\varphi_n} u\|_{m,\tau,\Gamma'}<\infty$ for all $n$. By our claim then $\widetilde{\varphi_n} u\in W^{m,\tau,p}_q(\Gamma')$ and, by Lemma~\ref{lemma3.1}, $\widetilde{\varphi_n} u\in W^{m,\tau,p}_q(B_1)$, with $\text{supp\,}\widetilde{\varphi_n} u\subset\subset B_1$ for all $n$. We can then trivially extend $\widetilde{\varphi_n} u$ to an element of
$W^{m,\tau,p}_q(\R^{N-1})$ and, by mollifying it, there is $v_n\in C^{r_1,p}_q(\R^{N-1})$, with $\text{supp\,}v_n\subset\subset \widetilde{U_n}$, such that
$\|v_n-\widetilde{\varphi_n} u\|_{m,\tau,\widetilde{U_n}}\le \eps^n/2$ for all $n$.
Each $v_n$ can be trivially extended to $v_n\in C^{r_1,p}_q(\Gamma')$. We then set $v=\sum_nv_n$, the sum being locally finite, so $v\in C^{r_1,p}_q(\Gamma')$.
We now set $\Gamma_n''=\widetilde{U_1}\cup\ldots\cup\widetilde{U_n}$, so
$(\Gamma_n'')_n$ is increasing and  $\Gamma'=\bigcup_n \Gamma''_n$.
Since $\sum_n\widetilde{\varphi_n}=1$ and $\eps<1/2$ we thus have, for all $n$,
\begin{equation}\label{B6}
  \|v-u\|_{m,\tau,\Gamma_n''}\le \sum_l\|v_i-\widetilde{\varphi_i} u\|_{m,\tau,\widetilde{U_i}}\le \eps.
\end{equation}
Hence $\|v\|_{m,\tau,\Gamma_n''}\le \|u\|_{m,\tau,\Gamma'}+\eps$ for all $n$, so by Beppo Levi's Theorem we get $v\in W^{m,\tau,p}_q(\Gamma')$ and, arguing as in \eqref{B6}, $\|v-u\|_{m,\tau,\Gamma'}\le \eps$, concluding the proof.
\end{proof}

Lemma~\ref{lemma3.2} shows that the $C^m$ regularity in \eqref{3.1} can be equivalently replaced by the maximal possible one, which is $C^{r_1}$. Hence, when $p=q=0$ and $r(=r_1)=\infty$, the Sobolev spaces defined above coincide with the well-known Sobolev spaces $H^{m,\tau}(\Gamma')$ widely used in Geometric Analysis, see \cite{hebey}.
Moreover \eqref{3.7} also suggests how to extend previous definition to $\tau=\infty$, being clearly unavailable in this case the procedure used before.

We then set, for $p,q\in\N_0$, $m\in\N_0$, $m\le r_1$,
\begin{equation}\label{3.8}
W^{m,\infty ,p}_q(\Gamma') =\{u\in W^{m,1,p}_{q,\loc}(\Gamma'): \, D_\Gamma^i u \in L^{\infty,p}_{q+i}(\Gamma')\quad \text{for $i=0,\ldots, m$}\},
\end{equation}
 endowed with the norm
\begin{equation}\label{3.9}
\|u\|_{m,\infty,\Gamma'} =\max_{i=0,\ldots,m}\|D_\Gamma^i u\|_{\infty,\Gamma'},
\end{equation}
which is trivially a Banach space.

By repeating previous procedure we can also set,  for $m\in\N$ and $p,q\in\N_0$, the spaces $W^{m,\infty,p}_q(\vartheta)$ on
any open subset $\vartheta$ of the trivial manifold $\R^{N-1}$. Since
$$W^{m,\infty}(\vartheta;\C^l)=\{u\in W^{m,1}_{\loc}(\vartheta;\C^l): \, D^\alpha u \in L^{\infty}(\vartheta;\C^l)\quad \text{for $0\le |\alpha|\le m$}\}$$
for any $l\in\N$, the isomorphism \eqref{N4} extends to  $1\le \tau\le \infty$.
Moreover we have the following result.
\begin{lem}\label{lemma3.3} The statement of Lemma~\ref{lemma3.1}  also holds when $\tau=\infty$.
\end{lem}
\begin{proof} We keep the notation convention in the proof of Lemma~\ref{lemma3.1}. Since \eqref{2.21} extends to $W^{1,1,p}_{q,\loc}(U)$, formula \eqref{B3} extends to $u\in W^{m+1,1,p}_{q,\loc}(U)$. We first consider the case $(U,\varphi)=(U_n,\zeta_n)\in\cal{A}$. By previous remark for $m=0,\ldots,r_1-1$ there is $c_{11}=c_{11}(m,p,q)\ge 1$ such that
\begin{equation}\label{B7}
\| \widetilde{D_\Gamma u}-d\widetilde{u}\|_{m,\infty,\vartheta_n}\le c_{11}\,\| \widetilde{u}\|_{m,\infty,\vartheta_n}
\quad\text{for all $u\in W^{m+1,1,p}_{q,\loc}(\Gamma')$ and $n\in\cal{I}_{\Gamma}$,}
\end{equation}
the norms being finite or not.

Also in this case we proceed by induction on $m=0,\ldots,r_1-1$, the case $m=0$ being nothing but \eqref{2.9}, so we suppose that \eqref{3.5} holds on $W^{m,\infty,p}_q(\Gamma')$ for $m\le r_1-1$. Then, by \eqref{B7}, for all $u\in W^{m+1,\infty,p}_q(\vartheta_n)$ we have
\begin{align*}
&\|u\|_{m+1,\infty,\Gamma'}=\max\{\|u\|_{m,\infty,\Gamma'}, \|D_\Gamma u\|_{m,\infty,\Gamma'}\}\\
&\le \max\{c_6(m,q)\|\widetilde{u}\|_{m,\infty,\vartheta_n}, c_6(m,q+1)\|\widetilde{D_\Gamma u}\|_{m,\infty,\vartheta_n'}\}\\
&\le  \max\{c_6(m,q)\|\widetilde{u}\|_{m,\infty,\vartheta_n}, c_6(m,q+1)\|d\widetilde{u}\|_{m,\infty,\vartheta_n'}+c_{11}\,c_6(m,q)\|\widetilde{u}\|_{m,\infty,\vartheta_n}\}\\
&\le c_6(m+1,q)\, \|\widetilde{u}\|_{m+1,\infty,\vartheta_n},
\end{align*}
where we set $c_6(m+1,q)=\max\{c_6(m,q), \,c_{11} \,\,c_6(m,q+1)\}$. The same arguments allow us to prove the reverse inequality.
In the case $\Gamma'\subseteq U$ we repeat the proof, by replacing \eqref{2.19} with \eqref{tridente} and allowing $c_{11}$ and $c_6$ to depend also on $U$ and $\varphi$.
\end{proof}

We have already recalled that  \eqref{2.21} can be extended to Sobolev spaces: accordingly, we then extend  \eqref{contractioncommuta} to Sobolev spaces as well, thus getting that for any $p\in\N$, $q\in\N_0$, $1\le \tau\le \infty$, $m\in\N$, $m\le r-1$,
\begin{equation}\label{contractioncommutaestesa}
D^m_\Gamma\cdot {\cal T}_p u={\cal T}_p\cdot D_\Gamma^m u\qquad\text{for all $u\in W^{m,\tau,p}_{p+q,\loc}(\Gamma')$.}
\end{equation}
By the same reason for any $p_1,p_2,q_1,q_2\in\N_0$, $1\le \tau_1,\tau_2\le \infty$ such that $1/\tau_1+ 1/\tau_2\le 1$, $m\in\N_0$, $m\le r_2:=\min\{r_1(r,p_1,q_1), r_1(r,p_2,q_2)\}$, and
$u\in W^{m,\tau_1,p_1}_{q_1,\loc}(\Gamma')$, $v\in W^{m,\tau_2,p_2}_{q_2,\loc}(\Gamma')$
 \begin{equation}\label{2.22estesa}
 u\otimes v \in W^{m,\tau_3,p_1+p_2}_{q_1+q_2,\loc}(\Gamma'),\quad\text{and the generalized Leibniz formula \eqref{2.22} holds.}
 \end{equation}
 Consequently, by Lemma~\ref{lemma3.2} and H\"{o}lder inequality,  one trivially derives the following multiplier property, which we state for  easy reference.
 \begin{lem} For all $p_1,p_2,q_1,q_2\in\N_0$, $1\le \tau\le \infty$,  $m\in\N_0$, $m\le r_2$,
$u\in W^{m,\infty,p_1}_{q_1}(\Gamma')$ and  $v\in W^{m,\tau,p_2}_{q_2}(\Gamma')$ we have
 $u\otimes v \in W^{m,\tau,p_1+p_2}_{q_1+q_2}(\Gamma')$. Moreover there is $c_{12}=c_{12}(m,\tau,p_1,p_2,q_1,q_2)>0$ such that
 \begin{equation}\label{3.10}
 \|u\otimes v\|_{m,\tau,\Gamma'}\le \,c_{12} \,\|u\|_{m,\infty,\Gamma'} \,\|v\|_{m,\tau, \Gamma'}
 \end{equation}
 for all $u\in W^{m,\infty,p_1}_{q_1}(\Gamma')$ and  $v\in W^{m,\tau,p_2}_{q_2}(\Gamma')$.
 \label{Multiplierlemma}\end{lem}

Our aim is now to show that, when $\Gamma'$ is relatively clopen in $\Gamma$, the spaces $W^{m,\tau,p}_q(\Gamma')$ essentially behave like the corresponding spaces in $\R^{N-1}$, and that they coincide with the spaces introduced in \cite{LeoniSobolev2}.  Hence we are going to introduce two suitable uniform localization system subordinated to the atlas ${\cal A}_{\Gamma'}$ defined in Lemma~\ref{lemma2.1},
 the first one being appropriate to use the argument in \cite{Amann2013}, and the second one coming from \cite{LeoniSobolev2}.

  We recall, following \cite[pp.~425--426]{LeoniSobolev2}, some auxiliary regular functions from  \cite{Stein1970} related to assumption (A0),  correcting an inessential misprint.
 Using the notation \eqref{3.11} and denoting by $*$ the standard convolution in $\R^N$, by $\chi_E$ the characteristic function of a set $E$ and by $(\rho_\eps)_{\eps>0}$ the standard mollifiers (see \cite[p.~687]{LeoniSobolev2} or \cite[p.~29]{adams}), we set
 \begin{equation}\label{3.12}
   \phi_n=\rho_{\eps_0/8}*\chi_{\Omega_n^{3\eps_0/8}},\qquad\text{for all $n\in {\cal I}_\Gamma$,}
 \end{equation}
 where $\eps_0$ is the parameter in assumption (A0). Hence, for all $n\in {\cal I}_\Gamma$, we have
 \begin{equation}\label{3.13}
 \begin{alignedat}{2}
 &\phi_n\in C^\infty(\R^N),\quad &&\text{supp}\,\phi_n\subseteq \Omega_n^{\eps_0/4}, \quad 0\le \phi_n\le 1\,\text{in $\R^N$},\\
 &\phi_n=1\,\,\text{in $\Omega_n^{\eps_0/2}$}\negquad,\quad &&D^\alpha \phi_n=D^\alpha\rho_{\eps_{\eps_0/8}}*\chi_{\Omega_n^{3\eps_0/8}}\quad\text{for every multiindex $\alpha$,}
\end{alignedat}
 \end{equation}
 so for all $m\in\N$ there is $c_{13}=c_{13}(m)>0$ such that
\begin{equation}\label{3.14}
|D^\alpha \phi_n|\le c_{13}\qquad\text{in $\R^N$ for all $|\alpha|\le m$ and $n\in {\cal I}_\Gamma$}.
\end{equation}
Introducing for all $\delta'>0$ the open set $\Gamma^{\delta'}=\{x\in\R^N: d(x,\Gamma)<\delta'\}$, by assumption (A0)--(i) we have $\Gamma^{\eps_0/2}\subseteq\bigcup_{n\in{\cal I}_\Gamma}\Omega_n^{\eps_0/2}$, so by \eqref{3.13}
\begin{equation}\label{3.15}
  \sum_{n\in{\cal I}_\Gamma}\phi_n\ge 1,\qquad \sum_{n\in{\cal I}_\Gamma} \phi_n^2\ge 1\qquad\text{in $\Gamma^{\eps_0/2}$,}
\end{equation}
the summations being finite at each point by assumption (A0)--(ii) and \eqref{3.13}. By the already recalled Faa di Bruno formula then for all $m\in\N$
there is $c_{14}=c_{14}(m)>0$ such that
\begin{equation}\label{3.16}
\left|D^\alpha \left( \textstyle\sum_{n\in{\cal I}_\Gamma}\phi_n\right)^{-1}\right| \le c_{14},\quad\text{and}\quad
\left|D^\alpha \left( \textstyle\sum_{n\in{\cal I}_\Gamma}\phi_n^2\right)^{-1/2}\right|\le c_{14}
\end{equation}
in $\Gamma^{\eps_0/2}$, for all $|\alpha|\le m$. We now set
\begin{equation}\label{3.17}
 \phi_n'=\phi_n(\textstyle\sum_{n\in{\cal I}_\Gamma}\phi_n^2)^{-1/2},  \quad\text{and}\quad \phi_n''=\phi_n(\textstyle\sum_{n\in{\cal I}_\Gamma}\phi_n)^{-1} \quad\text{for all $n\in {\cal I}_\Gamma$.}
\end{equation}
By \eqref{3.13}--\eqref{3.17} we then have, for all $n\in {\cal I}_\Gamma$,
\begin{equation}\label{3.18}
 \phi_n',\phi_n''\in C^\infty(\Gamma^{\eps_0/2}),\, \text{supp}\,\phi_n', \text{supp}\,\phi_n''\subseteq \Omega_n^{\eps_0/4}, \, 0\le \phi_n',\phi_n''\le 1\,\text{in $\Gamma^{\eps_0/2}$},
 \end{equation}
 and for all $m\in\N$ there is $c_{15}=c_{15}(m)>0$ such that
\begin{equation}\label{3.19}
|D^\alpha \phi_n'|,\quad |D^\alpha \phi_n''|,\quad
|D^\alpha (\phi_n'/\phi_n'')|,\quad
|D^\alpha (\phi_n''/\phi_n')|,\quad\le c_{15}\quad\text{in $\Gamma^{\eps_0/2}$ }
\end{equation}
for all $|\alpha|\le m$ and $n\in {\cal I}_\Gamma$. Moreover
\begin{equation}\label{3.20}
  \sum\limits_{n\in{\cal I}_\Gamma}\phi_n''=\sum\limits_{n\in{\cal I}_\Gamma}{\phi_n'}^2=1\quad\text{in $\Gamma^{\eps_0/2}$.}
\end{equation}
Finally we set $\psi_n={\phi_n'}_{|\Gamma}$ and $\psi'_n={\phi_n''}_{|\Gamma}$ so, being  $U_n$ bounded by assumption (A0)--(ii) and Lemma~\ref{lemma2.1}, using \eqref{3.13} we have
\begin{gather}\label{3.21}
  \psi_n,\psi_n'\in C^r_c(U_n)\quad \text{and}\quad
  0\le \psi_n, \psi'_n\le 1\quad \text{on $\Gamma$}\qquad\text{for all $n\in {\cal I}_\Gamma$},\\
  \label{3.21bis}
  \sum\limits_{n\in{\cal I}_\Gamma}\psi_n'=\sum\limits_{n\in{\cal I}_\Gamma}\psi_n^2=1\quad \text{on $\Gamma$,}
\end{gather}
  and, by \eqref{2.1}, \eqref{3.19} and Faa di Bruno formula for all $m\in\N_0$, $m\le r$, there is $c_{16}=c_{16}(m)>0$ such that
 \begin{equation}\label{3.22}
|D^\alpha (\psi_n\cdot h_n)|,\quad
|D^\alpha [(\psi_n/\psi_n')\cdot h_n]|,\quad
|D^\alpha [(\psi_n'/\psi_n)]\cdot h_n|\le c_{16}\quad\text{in $Q_0$}
\end{equation}
 for all $|\alpha|\le m$ and $n\in {\cal I}_\Gamma$. Consequently, using \eqref{2.19} and \eqref{2.21} (the  $\Gamma_{ij,n}^k$ appear only after first derivation), for all $m\in\N_0$, $m\le r$, there is $c_{17}=c_{17}(m)>0$ such that
 \begin{equation}\label{3.23}
|D_\Gamma ^l \psi_n|_\Gamma\le c_{17}\quad\text{on $\Gamma$ \quad  for all $l=0,\ldots,m$ and $n\in {\cal I}_\Gamma$.}
\end{equation}
Moreover, using \eqref{2.7}, \eqref{2.9}, \eqref{2.16} and \eqref{LaplaceBeltrami}, when $r\ge 2$, for all $m\in\N_0$,
$m\le r-1$ there is $c_{18}=c_{18}(m)>0$ such that
 \begin{equation}\label{3.23bis-}
|D_\Gamma ^l (\nabla_\Gamma \psi_n)|_\Gamma\le c_{18}\quad\text{on $\Gamma$ \quad  for all $l=0,\ldots,m$ and $n\in {\cal I}_\Gamma$,}
\end{equation}
and for all $m\in\N_0$, $m\le r-2$, there is $c_{19}=c_{19}(m)>0$ such that
 \begin{equation}\label{3.23bis}
|D_\Gamma ^l (\Delta_\Gamma \psi_n)|_\Gamma\le c_{19}\quad\text{on $\Gamma$ \quad  for all $l=0,\ldots,m$ and $n\in {\cal I}_\Gamma$.}
\end{equation}

We now introduce, for any  $\Gamma'=\overline{\Gamma'}$ and $p,q\in \N_0$, the
linear operators
\begin{equation}\label{3.24}
  \Phi_{q,\Gamma'}^p: {\cal M}^p_q(\Gamma')\to \prod\limits_{n\in {\cal I}_{\Gamma'}}{\cal M}^p_q(\R^{N-1}),
  \quad
  \Psi_{q,\Gamma'}^p: \prod\limits_{n\in {\cal I}_{\Gamma'}}{\cal M}^p_q(\R^{N-1}) \to {\cal M}^p_q(\Gamma'),
\end{equation}
defined by
\begin{equation}\label{3.24bis}
\Phi_{q,\Gamma'}^p u=\Big({\zeta_n}_*(\psi_n u)\Big)_n,\quad\text{and}\quad \Psi_{q,\Gamma'}^p \left((v_n)_n\right)=\sum\limits_{n\in{\cal I}_{\Gamma'}} \psi_n \zeta_n^*{v_n}_{|Q_0},
\end{equation}
 where ${\zeta_n}_*(\psi_n u)$ and $\psi_n \zeta_n^*{v_n}_{|Q_0}$ are trivially extended as zero, respectively on $\R^{N-1}$ and $\Gamma'$, and the summation is locally finite on $\Gamma$ by assumption (A0) and \eqref{3.21}. By \eqref{3.24bis} and \eqref{3.21}--\eqref{3.21bis} for all $u\in {\cal M}^p_q(\Gamma')$ we have
 \begin{align*}
 \Psi_{q,\Gamma'}^p \cdot \Phi_{q,\Gamma'}^p u=&\Psi_{q,\Gamma'}^p \Big({\zeta_n}_*(\psi_n u)\Big)_n\\=&\sum\limits_{n\in {\cal I}_{\Gamma'}} \psi_n \zeta_n^*\,{\zeta_n}_*(\psi_n u)=\sum\limits_{n\in {\cal I}_{\Gamma'}}\psi_n^2 u=u,
 \end{align*}
 that is  $\Psi_{q,\Gamma'}^p \cdot\Phi_{q,\Gamma'}^p=I$.
 By \eqref{3.21} and \eqref{3.24bis} the operators $\Psi_{q,\Gamma'}^p$ and $\Phi_{q,\Gamma'}^p$ preserve $C^m$ regularity for $m\in\N_0$, $m\le r_1$.

To characterize the spaces $W^{m,\tau,p}_q$ we recall, following \cite[p.~120]{triebel}, given a sequence $(X_n)_n$ of Banach spaces, with corresponding norms $\|\cdot\|_{X_n}$, and  $1\le \tau\le\infty$, the Banach space
\begin{equation}\label{elltau}
\ell^\tau (X_n)=\ell^\tau((X_n)_n)= \left\{x=(x_n)_n,\,\, x_n\in X_n\,\forall n\in\N,\,\,(\|x_n\|_{X_n})_n\in\ell^\tau\right\},
\end{equation}
endowed with the norm $\|x\|_{\ell^\tau (X_n)}=\|(\|x_n\|_{X_n})_n\|_{\ell^\tau}$.

Then we set, for all $s\ge 0$,
\begin{equation}\label{spaziconcreti}
  \ell^{s,\tau,p}_q=\ell^\tau (X_n),\qquad\text{where}\quad X_n=
  \begin{cases}
  W^{s,\tau,p}_q(\R^{N-1})\quad&\text{if $n\in {\cal I}_{\Gamma'}$,}\\
  \phantom{aaaa}\{0\}\quad&\text{otherwise,}
    \end{cases}
 \end{equation}
so that
\begin{footnote}{When ${\cal I}_{\Gamma'}$ has  finite cardinality $l$,  $\ell^{s,\tau,p}_q$ is isomorphic to $[W^{s,\tau,p}_q(\R^{N-1})]^l$. This trivial case occurs when $\Gamma_1$ is compact.}\end{footnote}
\begin{equation}\label{normaspaziconcreti}
  \|(u_n)_n\|_{\ell^{s,\tau,p}_q}=
  \begin{cases}
  \left(\sum_{n\in {\cal I}_{\Gamma'}}\|u_n\|_{s,\tau,\R^{N-1}}^\tau\right)^{1/\tau}\quad&\text{if $1\le \tau<\infty$,}\\
  \sup_{n\in {\cal I}_{\Gamma'}}\|u_n\|_{s,\infty,\R^{N-1}}\quad&\text{if $\tau=\infty$.}
  \end{cases}
  \end{equation}
Disregarding all indices $n\not\in {\cal I}_{\Gamma'}$,  we can identify $\ell^{s,\tau,p}_q$ with the corresponding linear subspace of $\prod_{n\in {\cal I}_{\Gamma'}}{\cal M}^p_q(\R^{N-1})$.
We then have the following technical result. It says, in the language of interpolation theory  (see \cite[p.22]{triebel}), that
$\Phi_{q,\Gamma'}^p$ is a retraction and $\Psi_{q,\Gamma'}^p $ a coretraction belonging to it.
\begin{prop}\label{proposition3.1}For any $\Gamma'=\overline{\Gamma'}$, $p,q\in\N_0$, $1\le \tau\le\infty$ and $m\in\N_0$, $m\le r_1$ the linear operators in \eqref{3.24} restrict to
\begin{equation}\label{3.25prima}
  \Phi_{q,\Gamma'}^p\in {\cal L}(W^{m,\tau,p}_q(\Gamma'),\ell^{m,\tau,p}_q), \qquad \Psi_{q,\Gamma'}^p \in {\cal L}(\ell^{m,\tau,p}_q,W^{m,\tau,p}_q(\Gamma')),
  \end{equation}
and $\Psi_{q,\Gamma'}^p \cdot\Phi_{q,\Gamma'}^p=I$.
Consequently $W^{m,\tau,p}_q(\Gamma')=(\Phi_{q,\Gamma'}^p)^{-1}(\ell^{m,\tau,p}_q)$ and  two equivalent norms on it are given by setting, for any $u\in W^{m,\tau,p}_q(\Gamma')$,
\begin{equation}\label{3.25}
  \|\Phi_{q,\Gamma'}^p u\|_{\ell^{m,\tau,p}_q}=\left\|\Big({\zeta_n}_*(\psi_n u)\Big)_n\right\|_{\ell^{m,\tau,p}_q}\!\!,\quad\text{and}\,\,
  \left\|\Big({\zeta_n}_*(\psi'_n u)\Big)_n\right\|_{\ell^{m,\tau,p}_q}\!\!.
\end{equation}
\end{prop}
\begin{proof}
We first prove \eqref{3.25prima}, starting from the case $\tau<\infty$. Let $u\in W^{m,\tau,p}_q(\Gamma')$.
By \eqref{3.34}, Lemma~\ref{lemma3.1}, \eqref{3.2} and \eqref{2.22estesa} we have
\begin{align*}
\|\Phi^p_{q,\Gamma'}u\|_{\ell^{m,\tau,p}_q}^\tau=&\sum_{n\in\cal{I}_{\Gamma'}}\|{\zeta_n}_*(\psi_nu)\|_{m\,\tau,\R^{N-1}}^\tau\le
c_6^\tau \sum_{n\in\cal{I}_{\Gamma'}}\|\psi_nu\|_{m\,\tau,U_n}^\tau\\
= &c_6^\tau \int_{\Gamma'} \sum_{n\in\cal{I}_{\Gamma'}}\sum_{i=0}^m |D_\Gamma^i(\psi_nu)|_\Gamma^\tau
= c_6^\tau \int_{\Gamma'} \sum_{n\in\cal{I}_{\Gamma'}}\sum_{i=0}^m \left|\sum_{j=0}^i D_\Gamma^j \psi_n \otimes D_\Gamma^{i-j}u)\right|_\Gamma^\tau\negquad,
\end{align*}
so, as $\text{supp\,}\psi_n\in U_n$, by \eqref{modulotensoriale} and Beppo Levi's Theorem
\begin{align*}
\|\Phi^p_{q,\Gamma'}u\|_{\ell^{m,\tau,p}_q}^\tau\le & m^{\tau-1}c_6^\tau \int_{\Gamma'} \sum_{n\in\cal{I}_{\Gamma'}}\sum_{i=0}^m \sum_{j=0}^i |D_\Gamma^j \psi_n|_\Gamma^\tau \,|D_\Gamma^{i-j}u|_\Gamma^\tau\\
= & m^{\tau-1}c_6^\tau \int_{\Gamma'} \sum_{i=0}^m \sum_{j=0}^i \left(\sum_{n\in\cal{I}_{\Gamma'}}|D_\Gamma^j \psi_n|_\Gamma^\tau \right) |D_\Gamma^{i-j}u|_\Gamma^\tau.
\end{align*}
Consequently, by \eqref{3.23}, assumption (A0)--(ii) and \eqref{3.2},
\begin{align*}
\|\Phi^p_{q,\Gamma'}u\|_{\ell^{m,\tau,p}_q}^\tau
\le & N_0 m^{\tau-1}(c_6\,c_{16})^\tau \int_{\Gamma'} \sum_{i=0}^m \sum_{j=0}^i |D_\Gamma^{i-j}u|_\Gamma^\tau\\
\le & N_0 (m c_6\,c_{16})^\tau \|u\|_{m,\tau,\Gamma'}^\tau,
\end{align*}
so proving that $\Phi^p_{q,\Gamma'}\in \cal{L}(W^{m,\tau,p}_q(\Gamma'),\ell^{m,\tau,p}_q)$.

Now let $v=(v_n)_n\in\ell^{m,\tau,p}_q$ and denote $u_n=\zeta_n^*{v_n}_{|Q_0}$, so $\Psi^p_{q,\Gamma'}v=\sum\limits_{n\in\cal{I}_{\Gamma'}}\psi_nu_n$,
where $\psi_nu_n$ is trivially extended to $\Gamma'$. Since $\{U_n,n\in \cal{I}_{\Gamma'}\}$ is locally finite and $\text{supp\,}\psi_n\subset U_n$ we have $D^j_\Gamma \Big(\sum_{n\in\cal{I}_{\Gamma'}}\psi_n u_n\Big)=\sum_{n\in\cal{I}_{\Gamma'}} D^j_\Gamma (\psi_n u_n)$ for $j=0,\ldots, m$. Hence, by \eqref{3.2},
assumption (A0)--(ii), \eqref{2.22estesa} and \eqref{modulotensoriale}
\begin{align*}
\|\|\Psi^p_{q,\Gamma'}v\|_{m,\tau,\Gamma'}^\tau=&\int_{\Gamma'}\sum_{i=0}^m\Big|D^i_\Gamma \Big(\sum_{n\in\cal{I}_{\Gamma'}}\psi_n u_n\Big)\Big|_\Gamma^\tau
\le N_0^{\tau-1}\int_{\Gamma'}\sum_{i=0}^m \sum_{n\in\cal{I}_{\Gamma'}} |D^i_\Gamma (\psi_n u_n)|_\Gamma^\tau\\
=&N_0^{\tau-1}\int_{\Gamma'}\sum_{i=0}^m \sum_{n\in\cal{I}_{\Gamma'}} \Big|\sum_{j=0}^i D^j_\Gamma \psi_n \otimes D_\Gamma^{i-j}u_n\Big|_\Gamma^\tau\\
\le &(m N_0)^{\tau-1}\int_{\Gamma'}\sum_{i=0}^m \sum_{n\in\cal{I}_{\Gamma'}} \sum_{j=0}^i |D^j_\Gamma \psi_n|_\Gamma^\tau\, |D_\Gamma^{i-j}u_n|_\Gamma^\tau.
\end{align*}
Consequently, by \eqref{3.23}, Beppo Levi's Theorem and Lemma~\ref{lemma3.1}, as $\text{supp\,}\psi_n\subset U_n$,
\begin{align*}
\|\|\Psi^p_{q,\Gamma'}v\|_{m,\tau,\Gamma'}^\tau\le & (m N_0)^{\tau-1}c_{16}^\tau\int_{\Gamma'}\sum_{n\in\cal{I}_{\Gamma'}}\sum_{i=0}^m  \sum_{j=0}^i |D_\Gamma^{i-j}u_n)|_\Gamma^\tau\\
\le & N_0^{\tau-1}(mc_{16})^\tau\negquad\sum_{n\in\cal{I}_{\Gamma'}}\|u_n\|_{m,\tau,U_n}^\tau
\le  N_0^{\tau-1}(m c_6 c_{16})^\tau\negquad\sum_{n\in\cal{I}_{\Gamma'}}\|v_n\|_{m,\tau,\R^{N-1}}^\tau\\
=& N_0^{\tau-1}(m c_6 c_{16})^\tau\|v\|_{\ell^{m,\tau,p}_q}^\tau,
\end{align*}
so concluding the proof of \eqref{3.25prima} when $\tau<\infty$.

We now consider the case $\tau=\infty$. Let $u\in W^{m,\infty,p}_q(\Gamma')$. By \eqref{normaspaziconcreti}, Lemma~\ref{lemma3.1}, \eqref{3.10} and \eqref{3.23} we have
\begin{align*}
\|\Phi^p_{q,\Gamma'}u \|_{\ell^{m,\infty,p}_q}=&\sup_{n\in\cal{I}_{\Gamma'}}\|{\zeta_n}_*(\psi_nu)\|_{m,\infty,\R^{N-1}}\le
c_6\sup_{n\in\cal{I}_{\Gamma'}}\|\psi_nu\|_{m\,\infty,U_n}\\
\le &c_6\, c_5 \,c_{16}\sup_{n\in\cal{I}_{\Gamma'}}\|u\|_{m,\infty,U_n}=c_6\, c_5 \,c_{16}\|u\|_{m,\infty,\Gamma'}.
\end{align*}
so  $\Phi^p_{q,\Gamma'}\in \cal{L}(W^{m,\infty,p}_q(\Gamma'),\ell^{m,\infty,p}_q)$.
Now let  $v=(v_n)_n\in\ell^{m,\infty,p}_q$ and  $u_n$ as in case $\tau<\infty$. Then, for all $i=0,\ldots,m$, by \eqref{2.22estesa} and \eqref{modulotensoriale} we have
\begin{align*}
\Bigg|D^i_\Gamma \Bigg(\sum_{n\in\cal{I}_{\Gamma'}}\psi_n u_n\Bigg)\Bigg|_\Gamma=&\Bigg|\sum_{n\in\cal{I}_{\Gamma'}} D^i_\Gamma (\psi_n u_n)\Bigg|_\Gamma
=\Bigg|\sum_{n\in\cal{I}_{\Gamma'}} \sum_{j=0}^i D^j_\Gamma \psi_n \otimes D^{i-j}_\Gamma u_n\Bigg|_\Gamma\\
\le & \sum_{n\in\cal{I}_{\Gamma'}} \sum_{j=0}^i |D^j_\Gamma \psi_n|_\Gamma |D^{i-j}_\Gamma u_n|_\Gamma\\
\le & \max_{j=0,\ldots,m} \sup_{n\in\cal{I}_{\Gamma'}} |D^j_\Gamma u_n|_\Gamma
\sum_{n\in\cal{I}_{\Gamma'}} \sum_{j=0}^i |D^j_\Gamma \psi_n|_\Gamma
\end{align*}
on $\Gamma'$. So, by \eqref{3.23}, assumption (A0)--(ii),  Lemma~\ref{lemma3.1} and \eqref{normaspaziconcreti}
$$\Bigg|D^i_\Gamma \Bigg(\sum_{n\in\cal{I}_{\Gamma'}}\psi_n u_n\Bigg)\Bigg|_\Gamma\le mN_0c_6 c_{16}\|v\|_{\ell^{m\infty,p}_q}\quad\text{on $\Gamma'$,}$$
so concluding the proof of \eqref{3.25prima}.

By it  we immediately get the inclusion $W^{m,\tau,p}_q(\Gamma')\subseteq (\Phi_{q,\Gamma'}^p)^{-1}(\ell^{m,\tau,p}_q)$. Since $\Psi_{q,\Gamma'}^p \cdot\Phi_{q,\Gamma'}^p=I$ also the reverse inclusion follows from \eqref{3.25prima}. Then trivially an equivalent norm for $W^{m,\tau,p}_q(\Gamma')$ is given by $\|\Phi_{q,\Gamma'}^p(\cdot)\|_{\ell^{m,\tau,p}_q}$. Now by Leibniz formula for all $m\in\N$ and $1\le \tau\le \infty$ there is $c_{20}=c_{20}(m,\tau)$ such that $\|uv\|_{m,\tau,\R^{N-1}}\le c_{20}\|u\|_{m,\infty,\R^{N-1}}\|v\|_{m,\tau,\R^{N-1}}$ for all $u\in W^{m,\infty}(\R^{N-1})$ and $v\in W^{m,\infty,p}_q(\R^{N-1})$. Moreover, by \eqref{3.21} and \eqref{3.22},
$$\|{\zeta_n}_*(\psi_n/\psi_n')\|_{m,\infty,\R^{N-1}},\quad \|{\zeta_n}_*(\psi_n'\psi_n)\|_{m,\infty,\R^{N-1}}\le c_{15}\quad\text{for all $n\in\cal{I}_{\Gamma'}$}.$$
Hence, setting $c_{21}=c_{21}(m,\tau)=c_{15}c_{20}$,
$$c_{21}^{-1} \|{\zeta_n}_*(\psi_n u)\|_{m,\tau,\R^{N-1}}\le \|{\zeta_n}_*(\psi'_n u)\|_{m,\tau,\R^{N-1}}\le
c_{21} \|{\zeta_n}_*(\psi_n u)\|_{m,\tau,\R^{N-1}}$$
for all $u\in W^{m,\tau,p}_q(\Gamma')$ and $n\in\cal{I}_{\Gamma'}$. Then also the second norm in \eqref{3.25} is an equivalent norm for
$W^{m,\tau,p}_q(\Gamma')$, concluding the proof.
\end{proof}
The following result points out a more useful  characterization of the space $W^{m,\infty,p}_q(\Gamma')$ and, as a first consequence of Proposition~\ref{proposition3.1}, a further multiplier property for the bundle metric $(\cdot,\cdot)_\Gamma$ given in \eqref{2.13}, which trivially fails without the uniformity assured by condition (A0).
\begin{lem}\label{lemma3.5}For any $\Gamma'=\overline{\Gamma'}$, $p,q\in\N_0$  and $m\in\N_0$, $m\le r_1$
$$W^{m,\infty,p}_q(\Gamma')=\{u\in {\cal M}^p_q(\Gamma'): \sup\nolimits_{n\in {\cal I}_{\Gamma'}} \| {\zeta_n}_* u\|_{m,\infty,Q_0}<\infty \},$$
and there is
$c_{22}=c_{22}(m,p,q)>0$ such that
$$ c_{22}^{-1}  \|\cdot\|_{m,\infty, \Gamma'}\le
  \sup\nolimits_{n\in{\cal I}_{\Gamma'}} \| {\zeta_n}_* \cdot\|_{m,\infty,Q_0}
   \le   c_{22} \|\cdot\|_{m,\infty, \Gamma'}.$$
Moreover, for all $1\le \tau\le\infty$,
$u\in W^{m,\infty,p}_q(\Gamma')$, $v\in W^{m,\tau,p}_q(\Gamma')$ we have $(u,v)_\Gamma\in W^{m,\tau}(\Gamma')$ and there is
$c_{23}=c_{23}(m,\tau,p,q)>0$ such that
 \begin{equation}\label{3.26bis}
 \|(u, v)_\Gamma\|_{m,\tau,\Gamma'}\le \,c_{23} \,\|u\|_{m,\infty,\Gamma'} \,\|v\|_{m,\tau, \Gamma'}
 \end{equation}
 for all $u\in W^{m,\infty,p}_q(\Gamma')$ and  $v\in W^{m,\tau,p}_q(\Gamma')$.
\end{lem}
\begin{proof}The first statement follows by Lemma~\ref{lemma3.3}, the isomorphism \eqref{N4} and the trivial identity
$\|u\|_{\infty,\Gamma'}=\sup\nolimits_{n\in {\cal I}_{\Gamma'}} \|u\|_{\infty,U_n}$.
To prove the second one we note that, denoting  $g^n_{(i) (i')}=g^n_{i_1 i_1'}\ldots g^n_{i_p i_p'}$ and $g_n^{(j) (j')}=g_n^{j_1 j_1'}\ldots g_n^{j_q j_q'}$ for $(i), (i')\in\J_p$, $(j), (j')\in\J_q$ and $n\in {\cal I}_{\Gamma'}$, by \eqref{2.13} we have
${\zeta_n}_*(u,v)_\Gamma=g^n_{(i) (i')}g_n^{(j) (j')}u^{(i)}_{(j)}\overline{v}^{(i')}_{(j')}$. Hence, by \eqref{2.4}, \eqref{2.9}, \eqref{N4} and Leibniz rule  $$
\begin{aligned}
\|{\zeta_n}_*[\psi_n (u,v)_\Gamma]\|_{m,\tau,\R^{N-1}}
=\,&\|{\zeta_n}_*[\psi_n (u,v)_\Gamma]\|_{m,\tau,Q_0}\\
=\,&\left\|\psi_n\cdot h_n\, g^n_{(i) (i')}\,g_n^{(j) (j')}u^{(i)}_{(j)}\overline{v}^{(i')}_{(j')}\right\|_{m,\tau,Q_0}\\
\le\, &c_{24}\,\|{\zeta_n}_*u\|_{m,\infty,Q_0} \,\|\, {\zeta_n}_*(\psi_n v)\|_{m,\tau,Q_0},\\
\end{aligned}
$$
where $c_{24}=c_{24}(m,p,q)>0$. By the first statement and Proposition~\ref{proposition3.1} we conclude the proof.
\end{proof}
\section{Spaces of real order} We now extend the definition of the spaces $W^{s,\tau,p}_q(\Gamma')$ to $0\le s\le r_1$. When $s\not\in\N_0$ they are usually called (see \cite{triebel}) Sobolev--Slobodeckii spaces. For simplicity we shall restrict to the case $1<\tau<\infty$.

For $s\in (0,r_1)\setminus\N$, $p,q\in\N_0$,  $1<\tau<\infty$ we take $m=[s]$, $\theta=s-m$, so $m<s<m+1$, $0<\theta<1$, and we set
\begin{equation}\label{3.27}
W^{s,\tau,p}_q(\Gamma')=\left( W^{m,\tau,p}_q(\Gamma'),W^{m+1,\tau,p}_q(\Gamma')\right)_{\theta,\tau},\quad H^{s,p}_q(\Gamma')=W^{s,2,p}_q(\Gamma')
\end{equation}
where $(\cdot,\cdot)_{\theta,\tau}$ is the real interpolator functor (se \cite[pp.~39--46]{bergh} or \cite[p.~24]{triebel}).

This definition shows the independence of $W^{s,\tau,p}_q(\Gamma')$ on any auxiliary object,  makes interpolation theory available in the sequel and it is classical in the compact case. We define the spaces $W^{s,\tau,p}_q(\R^{N-1})$ is the same way.

The following characterization in the atlas ${\cal A}_{\Gamma'}$, when $\Gamma'=\overline{\Gamma'}$,  is nothing but the extension of Proposition~\ref{proposition3.1}.
\begin{thm}\label{theorem3.1}For any $\Gamma'=\overline{\Gamma'}$, $p,q\in\N_0$, $1<\tau<\infty$ and $s\in\R$, $0\le s\le r_1$, the linear operators in \eqref{3.24} restrict to
\begin{equation}\label{stellarossa}
  \Phi_{q,\Gamma'}^p\in {\cal L}(W^{s,\tau,p}_q(\Gamma'),\ell^{s,\tau,p}_q), \qquad \Psi_{q,\Gamma'}^p \in {\cal L}(\ell^{s,\tau,p}_q,W^{s,\tau,p}_q(\Gamma')),
  \end{equation}
and $\Psi_{q,\Gamma'}^p \cdot\Phi_{q,\Gamma'}^p=I$.
Consequently $W^{s,\tau,p}_q(\Gamma')=(\Phi_{q,\Gamma'}^p)^{-1}(\ell^{s,\tau,p}_q)$ and  two equivalent norms on it are given by setting, for any $u\in W^{s,\tau,p}_q(\Gamma')$,
\begin{equation}\label{3.25-bis}
  \|\Phi_{q,\Gamma'}^p u\|_{\ell^{s,\tau,p}_q}=\left\|\Big({\zeta_n}_*(\psi_n u)\Big)_n\right\|_{\ell^{s,\tau,p}_q}\!\!,\quad\text{and}\,\,
  \left\|\Big({\zeta_n}_*(\psi'_n u)\Big)_n\right\|_{\ell^{s,\tau,p}_q}\!\!.
\end{equation}
\end{thm}
\begin{proof}
By well-known interpolation properties of Sobolev spaces in $\R^{N-1}$ (see \cite[Def.~6.22 and Theorem~6.45]{bergh} or \cite[Def. 1~p.~169 and Remark 2 p.~185] {triebel}) we have $W^{s,\tau}(\R^{N-1})=\left( W^{m,\tau}(\R^{N-1}),W^{m+1,\tau}(\R^{N-1})\right)_{\theta,\tau}$, with equivalence of norms.
Consequently, by using the characterization of real interpolation spaces as trace spaces, see \cite[Theorem~3.12.2 and Exercise~4 p.~76]{bergh},
we have
\begin{equation}\label{primainterpolazione}
 [W^{s,\tau}(\R^{N-1})]^{N_1}=\left( [W^{m,\tau}(\R^{N-1})]^{N_1},[W^{m+1,\tau}(\R^{N-1})]^{N_1}\right)_{\theta,\tau},
\end{equation}
with equivalence of norms.
Hence, by interpolation, the isomorphisms in \eqref{N3} restrict, for  all $\tilde{s}\ge 0$ and $1<\tau<\infty$, to bijective isomorphisms
\begin{equation}\label{N5}
 W^{\tilde{s},\tau,p}_q(\R^{N-1})\simeq W^{\tilde{s},\tau}(\R^{N-1};\C^{N_1})\simeq [W^{\tilde{s},\tau}(\R^{N-1})]^{N_1},
\end{equation}
with equivalence of norms. By combining \eqref{primainterpolazione}--\eqref{N5}, since all isomorphisms in \eqref{N5} are restrictions of the ones in \eqref{N3}, we then get
\begin{equation}\label{secondainterpolazione}
 W^{s,\tau,p}_q(\R^{N-1})=\left( W^{m,\tau,p}_q(\R^{N-1}),W^{m+1,\tau,p}_q(\R^{N-1})\right)_{\theta,\tau},
\end{equation}
with equivalence of norms. Since trivially $(\{0\},\{0\})_{\theta,\tau}=\{0\}$, by combining \eqref{spaziconcreti} with \eqref{secondainterpolazione} and
applying \cite[Theorem p.121]{triebel} we then get $\ell^{s,\tau,p}_q=(\ell^{m,\tau,p}_q,\ell^{m+1,\tau,p}_q)_{\theta,\tau}$ with equivalence of norms.
Consequently, by Proposition~\ref{proposition3.1} and \eqref{3.27}, by interpolation we get \eqref{stellarossa} and the identity $\Psi_{q,\Gamma'}^p \cdot\Phi_{q,\Gamma'}^p=I$. Then $W^{s,\tau,p}_q(\Gamma')=(\Phi_{q,\Gamma'}^p)^{-1}(\ell^{s,\tau,p}_q)$, and the first norm in \eqref{3.25-bis} is equivalent to the norm of $W^{s,\tau,p}_q(\Gamma')$. To show the  equivalence of the second norm in \eqref{3.25-bis} we remark that, by \eqref{3.21} and \eqref{3.22}, there is $c_{25}=c_{25}(m,\tau)>0$ such that
\begin{equation}\label{equivv}
\|(\psi_n/\psi_n') u\|_{m,\tau,\R^{N-1}},\quad \|(\psi'_n/\psi_n) u\|_{m,\tau,\R^{N-1}}\le c_{25} \|u\|_{m,\tau,\R^{N-1}}
\end{equation}
for all $n\in  I_{\Gamma}$ and  $u\in W^{m,\tau}(\R^{N-1})$, so using interpolation, \eqref{N5} and  \eqref{secondainterpolazione}, there is
$c_{26}=c_{26}(s, \tau)>0$ such that
\begin{equation}\label{equiww}
c_{26}^{-1}\|\psi_n u\|_{s,\tau,\R^{N-1}}\le \|\psi'_n u\|_{s,\tau,\R^{N-1}}\le c_{26} \|\psi_n  u\|_{s,\tau,\R^{N-1}}
\end{equation}
for all $u\in W^{s,\tau,p}_q(\R^{N-1})$ and $n\in  {\cal I}_\Gamma$, so by \eqref{normaspaziconcreti} we complete the proof.
\end{proof}

As a first consequence of Theorem~\ref{theorem3.1}  we get the following density result, which may fails when dealing with more general  non-compact smooth Riemannian manifolds, as shown in \cite[Chapter 3]{hebey}.
\begin{cor}\label{corollary3.1} Under the assumptions of Theorem~\ref{theorem3.1} the space of compactly supported elements of $C^{r_1,p}_q(\Gamma')$ is dense in $W^{s,\tau,p}_q(\Gamma')$. In particular $C^r_c(\Gamma')$ is dense in $W^{s,\tau}(\Gamma')$.
\end{cor}
\begin{proof} Since $C_c^\infty(\R^{N-1})$ is dense in $W^{s,\tau}(\R^{N-1})$ (see \cite[p.~52]{lionsmagenesIII}), by \eqref{N1regolare} and \eqref{N5} compactly supported elements of $C^{\infty,p}_q(\R^{N-1})$ are dense in $W^{s,\tau,p}_q(\R^{N-1})$.
Since for any $u\in W^{s,\tau,p}_q(\Gamma')$ we have  $\Phi_{q,\Gamma'}^p u\in\ell^{s,\tau,p}_q$ by Theorem~\ref{theorem3.1}, for any $n\in {\cal I}_{\Gamma'}$ there is a sequence $(v^n_l)_l$ of compactly supported elements of $C^{\infty,p}_q(\R^{N-1})$ such that $v^n_l\to {\zeta_n}_*(\psi_n u)$ in $W^{s,\tau,p}_q(\R^{N-1})$ as $l\to\infty$. Hence, for any $\eps>0$ and $n\in {\cal I}_{\Gamma'}$   there is $l_n\in\N$ such that $\|v^n_{l_n}-{\zeta_n}_*(\psi_n u)\|_{s,\tau,\R^{N-1}}<2^{-n}\eps$. We now set, for  $n\in\N\setminus {\cal I}_{\Gamma'}$,  $v^n_{l_n}=0$ and, for all $n\in\N$, $v_n=(v^1_{l_1},\ldots,v^n_{l_n},0,\ldots)\in \ell^{s,\tau,p}_q$.
Hence, by \eqref{normaspaziconcreti}, we have
$$
\begin{aligned}
\|v_n-\Phi_{q,\Gamma'}^p u\|_{\ell^{s,\tau,p}_q}^\tau = &
\sum_{i\in {\cal I}_{\Gamma'}, i\le n}\|v_{l_i}^i-{\zeta_i}_*(\psi_i u)\|_{s,\tau,\R^{N-1}}^\tau+
\sum_{i\in {\cal I}_{\Gamma'}, i> n}\|{\zeta_i}_*(\psi_i u)\|_{s,\tau,\R^{N-1}}^\tau\\
\le &\quad\sum_{i=1}^\infty 2^{-\tau i}\eps^\tau+
\sum_{i\in {\cal I}_{\Gamma'}, i> n}\|{\zeta_i}_*(\psi_i u)\|_{s,\tau,\R^{N-1}}^\tau\\
= & \quad\frac{\eps^\tau}{2^\tau-1}+
\sum_{i\in {{\cal I}_{\Gamma'}}, i> n}\|{\zeta_i}_*(\psi_i u)\|_{s,\tau,\R^{N-1}}^\tau,
\end{aligned}
$$
so, as $\Phi_{q,\Gamma'}^p u\in\ell^{s,\tau,p}_q$, we get $\varlimsup_n \|v_n-\Phi_{q,\Gamma'}^p u\|_{\ell^{s,\tau,p}_q}^\tau\le \frac{\eps^\tau}{2^\tau-1}$. Since $\eps>0$ is arbitrary we thus have
$v_n\to \Phi_{q,\Gamma'}^p u$ in $\ell^{s,\tau,p}_q$ as $n\to\infty$.
Hence, by Theorem~\ref{theorem3.1}, $\Psi_{q,\Gamma'}^p v_n\to u$ in $W^{s,\tau,p}_q(\Gamma')$.
Since $\Psi_{q,\Gamma'}^p v_n$ is a compactly supported element of $C^{r_1,p}_q(\Gamma')$, the proof is complete.
\end{proof}

We conclude this section by setting  Sobolev spaces (of functions) of negative real order on $\Gamma'$.
For $s\in\R$, $-r\le s< 0$ and $1<\tau<\infty$ we set $W^{s,\tau}(\Gamma')=[W^{-s,\tau'}(\Gamma')]'$, where $1/\tau'+1/\tau=1$, with the standard identification $[L^\tau(\Gamma')]'=L^{\tau'}(\Gamma')$. By \eqref{3.27} and the Duality Theorem for real interpolation spaces (see \cite[Theorem 3.7.1, p.~54]{bergh} or \cite[Theorem, p.~69]{triebel}) the interpolation property \eqref{3.7} continues to hold also for $s\in (-r,0)$.
 As usual we denote $H^s(\Gamma')=W^{s,2}(\Gamma')$.
\section{Operators}\label{subsection3.3} We now give the boundedness properties of  the linear operators which we need in the sequel.
\subsection{The trace operator} Since for $s\in (0,\infty)\setminus\N$ and $1<\tau<\infty$ the space $W^{s,\tau}(\R^{N-1})$ coincides with the Besov space $B^{s,\tau}(\R^{N-1})$ in \cite{LeoniSobolev2}, with equivalence of norms (see \cite[pp.~189--190]{triebel}), by Theorem~\ref{theorem3.1} the space $W^{s,\tau}(\Gamma)$ coincides with the Besov space $B^{s,\tau}(\Gamma)$ defined in \cite[pp.~615 and 625]{LeoniSobolev2}, with equivalence of norms.

Hence, by \cite[Theorem~18.40]{LeoniSobolev2}, the standard trace operator $u\mapsto u_{|\Gamma}$ from $C(\overline{\Omega})$ to $C(\Gamma)$, when restricted to $C(\overline{\Omega})\cap W^{1,\tau}(\Omega)$, $1<\tau<\infty$, has a unique surjective extension $\Tr\in {\cal L}\left(W^{1,\tau}(\Omega),W^{1-\frac 1\tau,\tau}(\Gamma)\right)$. Moreover $\Tr$ has a bounded right--inverse, i.e., $R_1\in {\cal L}\left(W^{1-\frac 1\tau,\tau}(\Gamma),W^{1,\tau}(\Omega)\right)$
such that $\Tr\cdot R_1=I$.
We shall denote, as usual, $\Tr u=u_{|\Gamma}$.
Moreover, when clear, we shall omit  the subscript $_{|\Gamma}$.

\subsection{Normal derivatives operators}\label{subsection3.3.2}
Denoting  for $m\in\N$, $2\le m\le r$, and $u\in W^{m,\tau}(\Omega)$, $i=1,\ldots,m-1$,
$\partial_\nu^i u_{|\Gamma} =\sum_{|\alpha|=i}\frac 1{\alpha!}D^\alpha u_{|\Gamma}\nu^\alpha$, by
\cite[Theorem~18.51, p.~626]{LeoniSobolev2},
we can set the operator
$$\Tr_mu=(u_{|\Gamma},\partial_\nu u_{|\Gamma}\ldots,\partial_\nu^{m-1} u_{|\Gamma}),\quad \Tr_m\in {\cal L}\left(W^{m,\tau}(\Omega),\prod_{i=0}^{m-1}W^{m-i-\frac 1\tau,\tau}(\Gamma)\right).
$$

Clearly $\partial_\nu u_{|\Gamma}$ is defined in this trace sense only when $r\ge 2$ and $u\in W^{2,\tau}(\Omega)$. When $r=1$  we then set
$\partial_\nu u_{|\Gamma}$  in  distributional sense as follows. For any $u\in W^{1,\tau}(\Omega)$ such that $\Delta u\in L^\tau(\Omega)$ in the sense of distributions and any  $h\in L^\tau(\Gamma)$ we say that  $\partial_\nu u_{|\Gamma}=h$ in distributional sense provided
\begin{equation}\label{3.40}
  \int_\Gamma hv=\int_\Omega \nabla u \nabla v+\int_\Omega \Delta u v\qquad\text{for all $v\in W^{1,\tau'}(\Omega)$.}
\end{equation}
By using the operator $R_1$ one easily gets that $\partial_\nu u_{|\Gamma}$ is unique when it exists. Moreover, integrating by parts using \cite[Theorem~18.1, p.~592]{LeoniSobolev2}, one easily recognizes that the so--defined distributional derivative extends the one defined in trace sense.
By $u_{|\Gamma_0}$ and $u_{|\Gamma_1}$ we shall denote the restrictions of $u_{|\Gamma}$ to $\Gamma_0$ and $\Gamma_1$

\subsection{The Riemannian gradient $\nabla_\Gamma$}
We extend  \eqref{3.2ter} to higher order spaces.
\begin{lem}\label{lemma3.4}For any $\Gamma'=\overline{\Gamma'}$,  $0\le s\le r-1$,  $1\le\tau\le\infty$ when $s\in\N_0$, $1<\tau<\infty$ when $s\in\R\setminus\N_0$, the isomorphisms $\flat$ and $\sharp$ defined in \eqref{musicalimisurabili} respectively restrict to
$\flat\in{\cal L}\left(W^{s,\tau,1}(\Gamma'),W^{s,\tau}_1(\Gamma')\right)$ and
$\sharp\in{\cal L}\left(W^{s,\tau}_1(\Gamma'),W^{s,\tau,1}(\Gamma')\right)$.
Consequently the operator $\nabla_\Gamma$ defined in \eqref{3.2ter}  restricts to
\begin{equation}\label{3.31}
\nabla_\Gamma\in{\cal L}\left(W^{s+1,\tau}(\Gamma'),W^{s,\tau,1}(\Gamma')\right).
\end{equation}
\end{lem}
\begin{proof}By interpolation we can take $s=m\in\N_0$, $m\le r-1$.
We first claim that $\sharp\in{\cal L}\left(W^{m,\tau}_1(\Gamma'),W^{m,\tau,1}(\Gamma')\right)$.
To prove our claim we note that, by \eqref{2.12}, for all $u\in W^{m,\tau}_1(\Gamma')$ and  $n\in {\cal I}_{\Gamma'}$ we have
${\zeta_n}_*[\psi_n(\sharp u)]= {\zeta_n}_*[\sharp(\psi_n u)]$. Moreover, by \eqref{2.9} and \eqref{2.12} there is $c_{27}=c_{27}(m)>0$ such that
$$\|(\zeta_n)_*[\sharp(\psi_n u)]\|_{m,\tau,\R^{N-1}}\le c_{27}\|(\zeta_n)_*(\psi_n u)\|_{m,\tau,\R^{N-1}}\quad\text{for all $n\in {\cal I}_{\Gamma'}$.}$$
Hence our claim follows by Proposition~\ref{proposition3.1}. The asserted boundedness of $\flat$ follows by the same arguments and, by combining \eqref{2.15} and \eqref{3.2bis} with the previous claim we get \eqref{3.31}.
\end{proof}
\subsection{The Riemannian divergence $\DivGammaPrimo$} The Riemannian divergence operator defined in \S\ref{subsection2.3.1} when $r\ge 2$ is trivially linear from $C^{m,1}_c(\Gamma')$ to $C^{m-1}_c(\Gamma')$ for $m\in\N$, $m\le r-1$. Denoting this restriction by $\DivGammaPrimo$ we can extend it as follows
on any relatively clopen subset $\Gamma'$ of $\Gamma$.
\begin{lem}\label{lemma3.5bis}
Let $r\ge 2$. For any relatively open $\Gamma'=\overline{\Gamma'}$,  $1\le s\le r-1$,  $1\le\tau\le\infty$ when $s\in\N_0$, $1<\tau<\infty$ when $s\in\R\setminus\N_0$, $\DivGammaPrimo$ uniquely extends by density to
\begin{equation}\label{3.32}
\DivGammaPrimo\in{\cal L}\left(W^{s,\tau,1}(\Gamma'),W^{s-1,\tau}(\Gamma')\right).
\end{equation}
\end{lem}
\begin{proof}By interpolation we can take $s=m\in\N$, $m\le r-1$, and by Corollary~\ref{corollary3.1} we only have to prove that
$$\DivGammaPrimo: (C^{m,1}_c(\Gamma'),\|\cdot\|_{m,\tau,\Gamma'})\to (C^{m-1}_c(\Gamma'),\|\cdot\|_{m-1,\tau,\Gamma'})$$
is bounded. Since $\DivGammaPrimo=\DivGammaPrimo\cdot \Psi_{0,\Gamma'}^1 \cdot\Phi_{0,\Gamma'}^1$ and $\Phi_{0,\Gamma'}^1(C^{m,1}_c(\Gamma'))\subset X^{m,\tau}$, where
$$X^{m,\tau}=\{u=(u_n)_n\in \ell^{m,\tau,1}_0: u_n\in C^{m,1}(\R^{N-1})\, \text{and supp}\, u_n\subset Q_0\,\,\text{for all $n\in  {\cal I}_{\Gamma'}$}\},$$
by Proposition~\ref{proposition3.1} we can just show that
$$\DivGammaPrimo\cdot \Psi_{0,\Gamma'}^1: (X^{m,\tau},\|\cdot\|_{\ell^{m,\tau,1}_0})\to (W^{m-1,\tau}(\Gamma'),\|\cdot\|_{m-1,\tau,\Gamma'})$$
is bounded. Now, for any $u=(u_n)_n\in X^{m,\tau}$, by \eqref{2.23} and \eqref{3.24bis} we have
\begin{equation}\label{3.33}
 \DivGammaPrimo\cdot \Psi_{0,\Gamma'}^1 u=\sum_{n\in {\cal I}_{\Gamma'}}\psi_n\DivGammaPrimo (\zeta_n^*u_n)+
 \sum_{n\in {\cal I}_{\Gamma'}}\partial_i\psi_n (\zeta_n^*u_n)^i.
 \end{equation}
 We shall now estimate both summations in \eqref{3.33}, starting from the first one. Now
 $\sum_{n\in {\cal I}_{\Gamma'}}\psi_n\DivGammaPrimo (\zeta_n^*u_n)=\Psi_{0,\Gamma'}^1 v$, where $v=(v_n)_n$ and $v_n=\DivGammaPrimo(\zeta_n^*u_n)$ when $n\in {\cal I}_{\Gamma'}$, $v_n=0$ otherwise. To  estimate $\|v\|_{\ell^{m,\tau,1}_0}$ we note that,
 by \eqref{2.7} and \eqref{2.23}, there is $c_{28}=c_{28}(m)>0$ such that
 $\|v_n\|_{m-1,\tau,\R^{N-1}}\le c_{28}\, \|u_n\|_{m,\tau,\R^{N-1}}$ for all $n\in {\cal I}_{\Gamma'}$, so
 $\|v\|_{\ell^{m,\tau,1}_0}\le c_{28}\,\|u\|_{\ell^{m,\tau,1}_0}$. Hence, by Proposition~\ref{proposition3.1}, there is $c_{29}=c_{29}(m,p,q,\tau)>0$ such that
 \begin{equation}\label{3.34}
 \Big\|\sum\nolimits_{n\in {\cal I}_{\Gamma'}}\psi_n\DivGammaPrimo(\zeta_n^*u_n)\Big\|_{m-1,\tau,\Gamma'}\le c_{29}\,\|u\|_{\ell^{m,\tau,1}_0}.
 \end{equation}
 To estimate the second summation in \eqref{3.33} we note that, by assumption (A0), Lemma~\ref{lemma3.1} and \eqref{3.22}, as
 $\text{supp}\, u_n\subset Q_0$ for all $n\in  {\cal I}_{\Gamma'}$,  there is $c_{30}=c_{30}(m,\tau)>0$ such that
$$\begin{aligned}
 \Big\|\sum_{n\in {\cal I}_{\Gamma'}}\partial_i\psi_n (\zeta_n^*u_n)^i\Big\|_{m-1,\tau,\Gamma'}^\tau
 =&\int_{\Gamma'}\sum_{i=0}^{m-1}
 \Big|D_\Gamma^i\Big(\sum_{n\in {\cal I}_{\Gamma'}}\partial_i\psi_n (\zeta_n^*u_n)^i\Big)\Big|_\Gamma^\tau\\
\le&N_0^{\tau-1}\int_{\Gamma'}\sum_{i=0}^{m-1}\sum_{n\in {\cal I}_{\Gamma'}}
 \left|D_\Gamma^i\left(\partial_i\psi_n (\zeta_n^*u_n)^i\right)\right|_\Gamma^\tau\\
 \le&c_6^\tau\,N_0^{\tau-1}\sum_{n\in {\cal I}_{\Gamma'}}
 \left\|\partial_i(\psi_n\cdot h_n)\, u_n^i\right\|_{m-1,\tau,\R^{N-1}}^\tau\\
\le&c_{30}\sum_{n\in {\cal I}_{\Gamma'}}
 \|u_n\|_{m-1,\tau,\R^{N-1}}^\tau\\
 =&c_{30}\,\|u\|_{\ell^{m,\tau,1}_0}^\tau.
\end{aligned}
$$
By combining the last estimate with \eqref{3.33} and \eqref{3.34} we complete the proof.
 \end{proof}
 We now remark that, by Corollary~\ref{corollary3.1}, the integration by parts formula \eqref{2.24} extends by density to  $u\in W^{1,\tau}(\Gamma')$, $w\in W^{1,\tau',1}(\Gamma')$, $1<\tau<\infty$.

 Hence, when $r\ge 2$ the operator $\DivGammaPrimo$ -- defined in Lemma~\ref{lemma3.5bis} on any relatively clopen $\Gamma'$ -- extends by density to
 $\DivGammaPrimo\in{\cal L}\left(L^{\tau,1}(\Gamma'),W^{-1,\tau}(\Gamma')\right)$ defined by
 \begin{equation}\label{3.36}
 \langle \DivGammaPrimo u,v\rangle_{W^{1,\tau'}(\Gamma')}=-\int_{\Gamma'}(u,\nabla_\Gamma\overline{v})_\Gamma\qquad\text{for all $v\in W^{1,\tau'}(\Gamma')$},
 \end{equation}
 and clearly \eqref{3.36} defines $\DivGammaPrimo$, in a weak sense, also when $r=1$.
 By interpolation
 \begin{equation}\label{3.37}
 \DivGammaPrimo\in{\cal L}\left(W^{s,\tau,1}(\Gamma'),W^{s-1,\tau}(\Gamma')\right)\quad\text{for $s\in\R$, $0\le s\le r-1$, $1<\tau<\infty$.}
\end{equation}
\subsection{Splitting and the operator $-\DivGamma(\sigma\nabla_\Gamma)$.}\label{subsection3.3.5} Mainly to simplify the notation we remark that, by Theorem~\ref{theorem3.1},
as $ {\cal I}_\Gamma={\cal I}_{\Gamma_0}\sqcup {\cal I}_{\Gamma_1}$,
by identifying each $u\in W^{s,\tau,p}_q(\Gamma_i)$, $i=0,1$, $p,q\in\N_0$, $1<\tau<\infty$, $s\in\R$, $0\le s\le r_1$ with its trivial extension to $\Gamma$, we can identify $W^{s,\tau,p}_q(\Gamma_i)$ with its isomorphic image in $W^{s,\tau,p}_q(\Gamma)$. We shall constantly make this identification during the paper. As a consequence we have the splitting
\begin{equation}\label{eq:splitting}
W^{s,\tau,p}_q(\Gamma)=W^{s,\tau,p}_q(\Gamma_0)\oplus W^{s,\tau,p}_q(\Gamma_1)\quad p,q\in\N_0, \,1<\tau<\infty, \,0\le s\le r_1,
\end{equation}
which extends, by duality, also when $p=q=0$ and $s\in[-r,0]$.

Consequently the operator $\DivGamma$ defined in \eqref{3.36} splits to $\DivGamma=(\DivGammaZero,\DivGammaUno)$, and $\DivGamma=\DivGammaUno$ on
$W^{s,\tau,1}(\Gamma_1)$  for $s\in\R$, $0\le s\le r-1$, $1<\tau<\infty$. In the sequel we shall only use the simpler notation $\DivGamma$. In particular \eqref{3.37}  yields
\begin{equation}\label{3.37bis}
 \DivGamma\in{\cal L}\left(W^{s,\tau,1}(\Gamma_1),W^{s-1,\tau}(\Gamma_1)\right)\quad\text{for $0\le s\le r-1$, $1<\tau<\infty$.}
\end{equation}
Hence, by assumption (A1), Lemma~\ref{Multiplierlemma} and interpolation we have
\begin{equation}\label{3.38}
 \DivGamma(\sigma\nabla_\Gamma)\in{\cal L}\left(W^{s+1,\tau}(\Gamma_1),W^{s-1,\tau}(\Gamma_1)\right)\quad\text{for  $0\le s\le r-1$, $1<\tau<\infty$.}
\end{equation}
Since by \eqref{2.24} we have
\begin{equation}\label{3.39}
\int_{\Gamma_1} -\DivGamma(\sigma\nabla_\Gamma u)v=\int_{\Gamma_1}\sigma (\nabla_\Gamma u,\nabla_\Gamma \overline{v})_\Gamma
=\int_{\Gamma_1} -\DivGamma(\sigma\nabla_\Gamma v)u
\end{equation}
for all $u\in W^{1,\tau}(\Gamma_1)$, $v\in W^{1,\tau'}(\Gamma_1)$, $1<\tau<\infty$, we can extend by transposition the operator $\DivGamma(\sigma\nabla_\Gamma)$ to get
\begin{equation}\label{3.41}
 \DivGamma(\sigma\nabla_\Gamma)\in{\cal L}\left(W^{s+1,\tau}(\Gamma_1),W^{s-1,\tau}(\Gamma_1)\right)\quad\text{for  $|s|\le r-1$, $1<\tau<\infty$.}
\end{equation}

\chapter{Well-posedness of problem \eqref{1.1}} \label{section4}
\section{Abstract analysis} We endow the phase space $\cal H$ defined in \eqref{1.2} with the inner product given, for $V_i=(u_i,v_i,w_i,z_i)$, $i=1,2$, by
\begin{equation}\label{4.2}
\begin{aligned}
(V_1,V_2)_{\cal H}= &\int_\Omega \nabla u_1\nabla \overline{u_2}+\int_\Omega u_1\overline{u_2}+\int_{\Gamma_1}\frac\sigma\rho(\nabla_\Gamma v_1,\nabla_\Gamma v_2)_\Gamma\\
+&\int_{\Gamma_1}\frac{v_1\overline{v_2}}\rho+\frac 1{c^2}\int_\Omega w_1\overline{w_2}+\int_{\Gamma_1}\frac\mu\rho z_1\overline{z_2}
\end{aligned}
\end{equation}
which is well-defined and is equivalent to the standard inner product of $\cal H$ by assumptions (A1--2). To study problem \eqref{1.1} we standardly reduce it to a first order problem by setting $w=u_t$ and $z=v_t$, so formally getting
\begin{equation}\label{4.3}
\begin{cases}
u_t-w=0 \qquad &\text{in
$\R\times\Omega$,}\\
v_t- z =0\qquad
&\text{on $\R\times \Gamma_1$,}\\
w_t-c^2\Delta u=0 \qquad &\text{in $\R\times\Omega$,}\\
\mu z_t- \DivGamma (\sigma \nabla_\Gamma v)+\delta z+\kappa v+\rho w =0\qquad
&\text{on $\R\times \Gamma_1$,}\\
z =\partial_\nu u\qquad
&\text{on
$\R\times \Gamma_1$,}\\
\partial_\nu u=0 &\text{on $\R\times \Gamma_0$,}\\
u(0,x)=u_0(x),\quad w(0,x)=u_1(x) &
 \text{in $\Omega$,}\\
v(0,x)=v_0(x),\quad z(0,x)=v_1(x) &
 \text{on $\Gamma_1$.}
\end{cases}
\end{equation}
We then define the unbounded operator $A:D(A)\subset {\cal H}\to {\cal H}$ by
\begin{multline}\label{4.4}
  D(A)=\{(u,v,w,z)\in [H^1(\Omega)\times H^1(\Gamma_1)]^2: \Delta u\in L^2(\Omega),\\  \partial_\nu u_{|\Gamma_0}=0, \quad
  \partial_\nu u_{|\Gamma_1}=z, \quad  \DivGamma(\sigma\nabla_\Gamma v)\in L^2(\Gamma_1)\},
\end{multline}
where $\Delta u$, $\partial_\nu u_{|\Gamma_0}$ and $\partial_\nu u_{|\Gamma_1}$ are taken in the sense made precise in \S\ref{subsection3.3.2} and
$\DivGamma(\sigma\nabla_\Gamma v)$ in the one made precise in \S\ref{subsection3.3.5}, and
\begin{equation}\label{4.5}
 A\begin{pmatrix}u\\v\\w\\z\end{pmatrix} =
\begin{pmatrix}-w\\-z\\-c^2\Delta u\\
\frac 1\mu\left[-\DivGamma(\sigma\nabla_\Gamma v)+\delta z+\kappa v+\rho w_{|\Gamma_1}\right]
\end{pmatrix}.
\end{equation}
So setting $U=(u,v,w,z)$, problem \eqref{4.3} can be formally written as
\begin{equation}\label{4.6}
  U'+AU=0,\qquad\text{in $\cal H$}, \qquad U(0)=U_0:=(u_0, v_0,u_1,v_1).
\end{equation}
We now recall, for the reader's convenience, the classical definitions of strong and generalized (or mild) solution of \eqref{4.6} (see \cite[pp.~4 and 105]{pazy} or \cite[Chapter II, pp.~145--150]{EngelNagel}) which we shall use in the sequel.
\begin{definition}\label{definition4.1}
We say that
\renewcommand{\labelenumi}{{\roman{enumi})}}
\begin{enumerate}
\item $U\in C^1(\R,{\cal H})$ is a {\em strong solution} of $U'+AU=0$ provided $U(t)\in D(A)$ and $U'(t)+AU(t)=0$ for all $t\in\R$;
\item $U\in C(\R,{\cal H})$ is a {\em generalized solution} of $U'+AU=0$ if it is the limit in $C(\R,{\cal H})$ of a sequence of strong solutions of it;
\item $U$ is a strong or generalized solution of \eqref{4.6} provided it is a solution of $U'+AU=0$ of the same type and $U(0)=U_0$.
\end{enumerate}
\end{definition}

Three important properties of the operator $A$ are given by the following lemmas.
\begin{lem}\label{lemma4Delio}Under assumptions (A0--3) for any $U=(u,v,w,z)\in D(A)$ we have
\begin{equation}\label{lemma4.1Delio1}
\begin{aligned}
(AU,U)_{\cal H}&=\,2\mathfrak{i}\, \Ima  \left\{\int_{\Omega}\nabla u\nabla\overline{w}+\int_{\Gamma_1} w\overline{z}+\int_{\Gamma_1}\frac{v\overline{z}}\rho+\int_{\Gamma_1}\frac \sigma\rho (\nabla_\Gamma v,\nabla_\Gamma z)_\Gamma\right\}\\
&-\int_\Omega w\overline{u}-\int_{\Gamma_1}\frac{\sigma \overline{z}}{\rho^2} (\nabla_\Gamma v,\nabla_\Gamma z)_\Gamma+\int_{\Gamma_1}\frac {\kappa-1}\rho v\overline{z}+\int_{\Gamma_1}\frac \delta\rho |z|^2,
\end{aligned}
\end{equation}
and consequently
\begin{equation}\label{lemma4.1Delio2}
\begin{split}
\Real (AU,U)_{\cal H}&=-\Real \int_\Omega w\overline{u}-\Real \int_{\Gamma_1}\frac{\sigma \overline{z}}{\rho^2} (\nabla_\Gamma v,\nabla_\Gamma \rho)_\Gamma\\
&\quad +\Real \int_{\Gamma_1}\frac {\kappa-1}\rho v\overline{z}+\int_{\Gamma_1}\frac \delta\rho |z|^2.
\end{split}
\end{equation}
\end{lem}
\begin{proof} By \eqref{4.2} and \eqref{4.5} for any $U=(u,v,w,z)\in D(A)$ we have
\begin{equation}\label{4.6.1}
\begin{aligned}
(AU,U)_{\cal H}=&-\int_\Omega\nabla w\nabla \overline{u}-\int_\Omega u\overline{w}-\int_\Omega\Delta u\overline w-\int_{\Gamma_1}\frac \sigma\rho (\nabla_\Gamma z,\nabla_\Gamma v)_\Gamma\\
&-\int_{\Gamma_1}\frac{z\overline{v}}\rho+\int_{\Gamma_1}\frac 1\rho[-\DivGamma (\sigma\nabla_\Gamma v)+\kappa v+\delta z+\rho w]\overline{z}.
\end{aligned}
\end{equation}
Since, by assumption (A2) and Lemma~\ref{lemma3.5}, $1/\rho\in W^{1,\infty}(\Gamma_1)$, by Lemma~\ref{Multiplierlemma} we have $z/\rho\in H^1(\Gamma_1)$.
Hence, by \eqref{3.40},  \eqref{3.36} and \eqref{4.6.1}, we get
$$
\begin{aligned}
(AU,U)_{\cal H}=&-\int_\Omega\nabla w\nabla \overline{u}-\int_\Omega u\overline{w}+\int_\Omega\nabla u\nabla \overline w-\int_{\Gamma_1}\frac \sigma\rho (\nabla_\Gamma z,\nabla_\Gamma v)_\Gamma-\int_{\Gamma_1}\frac{z\overline{v}}\rho\\
&-\int_{\Gamma_1}z\overline{w}+\int_{\Gamma_1} \sigma\left(\nabla_\Gamma v,\nabla_\Gamma\left(\tfrac 1 \rho z\right)\right)_\Gamma +\int_{\Gamma_1}\frac \kappa\rho v\overline{z}+
\int_{\Gamma_1}\frac \delta\rho |z|^2+ \int_{\Gamma_1}w\overline{z}.
\end{aligned}
$$
Using Leibniz formula $\nabla_\Gamma(z/\rho)=\nabla_\Gamma (1/\rho)z+\nabla_\Gamma z/\rho$, which follows from \eqref{2.15}, in the previous identity, we get
\begin{align*}
(AU,U)_{\cal H}=&\,2\mathfrak{i}\, \Ima  \left\{\int_{\Omega}\nabla u\nabla\overline{w}+\int_{\Gamma_1} w\overline{z}+\int_{\Gamma_1}\frac{v\overline{z}}\rho\right\}
-\int_\Omega w\overline{u}\\
&-\int_{\Gamma_1}\frac \sigma\rho (\nabla_\Gamma z,\nabla_\Gamma v)_\Gamma+\int_{\Gamma_1}\sigma \left(\nabla_\Gamma v,\nabla_\Gamma \left(\tfrac 1\rho\right)\right)_\Gamma \overline{z}\\
&+\int_{\Gamma_1}\frac \sigma\rho (\nabla_\Gamma v,\nabla_\Gamma z)_\Gamma+\int_{\Gamma_1}\frac {\kappa-1}\rho v\overline{z}+\int_{\Gamma_1}\frac \delta\rho |z|^2.
\end{align*}
Since, by \eqref{2.15} and the Chain Rule in $H^1(Q_0)$, we have $\nabla_\Gamma (1/\rho)=-\nabla_\Gamma \rho/\rho^2$,  we can rewrite it as \eqref{lemma4.1Delio1}. A simple calculation then gives \eqref{lemma4.1Delio2}.
\end{proof}
\begin{lem}\label{lemma4.1}Under assumptions (A0--3) there is $\Lambda_0=\Lambda_0(\sigma, \rho, \kappa, \delta, \mu, c)>0$ such that for any $\lambda\ge \Lambda_0$ the operators $A+\lambda I$ and $-A+\lambda I$are accretive on $\cal H$, i.e.,
$$\Real  \big(\,\,(\pm A+\lambda I)U,\,\,U\big)_{\cal H}\ge 0\qquad\text{for all $U\in D(A)$.}$$
\end{lem}
\begin{proof} We use H\"{o}lder, Young  and Cauchy--Schwarz inequalities in the identity \eqref{lemma4.1Delio2} given of Lemma~\ref{lemma4Delio}. For any $U=(u,v,w,z)\in D(A)$ we thus get the estimate
\begin{align*}
|\Real (AU,U)_{\cal H}|\le &\tfrac 12 \|u\|_2^2+\tfrac 12 \|w\|_2^2+\int_{\Gamma_1}\frac \sigma{\rho^2} |\nabla_\Gamma v|_\Gamma|\nabla_\Gamma \rho|_\Gamma |z|\\
+&\tfrac {\|\kappa-1\|_{\infty,\Gamma_1}}2  \left[\int_{\Gamma_1}\tfrac {|v|^2+|z|^2}\rho\right]+
{\scriptstyle\|\delta\|_{\infty,\Gamma_1}}\int_{\Gamma_1}\tfrac {|z|^2}\rho.
\end{align*}
Hence, using assumptions (A1--2) and Young inequality,
\begin{align*}
|\Real (AU,U)_{\cal H}|\le& \tfrac 12 \|u\|_2^2+\tfrac 12 \|w\|_2^2+{\scriptstyle\|\nabla_\Gamma \rho\|_{\infty,\Gamma_1}}\int_{\Gamma_1}\tfrac \sigma{\rho^2} |\nabla_\Gamma v|_\Gamma |z|\\
+&\tfrac 12 {\scriptstyle\|\kappa-1\|_{\infty,\Gamma_1}} \int_{\Gamma_1}\frac {|v|^2}\rho+\tfrac{\|\kappa-1\|_{\infty,\Gamma_1}+2\|\delta\|_{\infty,\Gamma_1}}{2\mu_0}
\int_{\Gamma_1}\tfrac \mu \rho|z|^2\\
\le &\tfrac 12 \|u\|_2^2+\tfrac 12 \|w\|_2^2+\tfrac{\|\nabla_\Gamma \rho\|_{\infty,\Gamma_1}}2\left[\int_{\Gamma_1}\tfrac \sigma{\rho^2} |\nabla_\Gamma v|_\Gamma^2 +\int_{\Gamma_1}\tfrac \sigma{\rho^2}|z|^2\right]\\
+&\tfrac 12 {\scriptstyle\|\kappa-1\|_{\infty,\Gamma_1}} \int_{\Gamma_1}\frac {|v|^2}\rho+\tfrac{\|\kappa-1\|_{\infty,\Gamma_1}+2\|\delta\|_{\infty,\Gamma_1}}{2\mu_0}
\int_{\Gamma_1}\tfrac \mu\rho|z|^2\\
\le &\tfrac 12 \|u\|_2^2+\tfrac 12 \|w\|_2^2+\tfrac{\|\nabla_\Gamma \rho\|_{\infty,\Gamma_1}}{2\rho_0} \int_{\Gamma_1}\tfrac \sigma\rho |\nabla_\Gamma v|_\Gamma^2 +\tfrac 12 {\scriptstyle\|\kappa-1\|_{\infty,\Gamma_1}} \int_{\Gamma_1}\tfrac {|v|^2}\rho\\
+&\tfrac{\|\nabla_\Gamma \rho\|_{\infty,\Gamma_1}\|\sigma\|_{\infty,\Gamma_1}} {2\rho_0\mu_0}\int_{\Gamma_1}\tfrac \mu\rho|z|^2
+\tfrac{\|\kappa-1\|_{\infty,\Gamma_1}+2\|\delta\|_{\infty,\Gamma_1}}{2\mu_0}
\int_{\Gamma_1}\frac {\mu }\rho|z|^2\\
\le &\Lambda_0\left[\|u\|_2^2+\tfrac 1{c^2}\|w\|_2^2+\int_{\Gamma_1}\tfrac \sigma\rho |\nabla_\Gamma v|_\Gamma^2+\int_{\Gamma_1}\tfrac {|v|^2}\rho
+\int_{\Gamma_1}\tfrac {\mu }\rho |z|^2\right],
\end{align*}
where
$$\Lambda_0=\tfrac 12 \max\left\{1, c^2, \tfrac{\|\nabla_\Gamma \rho\|_{\infty,\Gamma_1}}{\rho_0},{\scriptstyle \|\kappa-1\|_{\infty,\Gamma_1}},
\tfrac{\|\nabla_\Gamma \rho\|_{\infty,\Gamma_1}\|\sigma\|_{\infty,\Gamma_1}+\rho_0\|\kappa-1\|_{\infty,\Gamma_1}+2\rho_0\|\delta\|_{\infty,\Gamma_1}}{\rho_0\mu_0}
\right\},$$
and then, by \eqref{4.2}, we have $|\Real (AU,U)_{\cal H}|\le\Lambda_0\|U\|_{\cal H}^2$.
Consequently, for any $\lambda\ge \Lambda_0$, we have $\Real  \big(\,\,(\pm A+\lambda I)U,\,\,U\big)_{\cal H}\ge (\lambda-\Lambda_0)\|U\|_{\cal H}^2\ge 0$.
\end{proof}

\begin{lem}\label{lemma4.2}Under assumptions (A0--3) there is $\Lambda_1=\Lambda_1(\sigma, \rho, \kappa, \delta, \mu, c)\ge\Lambda_0$ such that for any $\lambda\in\R$ with  $|\lambda|\ge \Lambda_1$ the operator $A+\lambda I$ is surjective.
\end{lem}
\begin{proof}
By \eqref{1.2}, \eqref{4.4} and \eqref{4.5} we have to prove that for $\lambda\in\R$ with  $|\lambda|$ large enough and for all $h_1\in H^1(\Omega)$, $h_2\in H^1(\Gamma_1)$, $h_3\in L^2(\Omega)$ and $h_4\in L^2(\Gamma_1)$ the system
\begin{equation}\label{4.7}
\begin{cases}
-w+\lambda u=h_1 \qquad &\text{in
$\Omega$,}\\
- z+\lambda v =h_2\qquad
&\text{on
$\Gamma_1$,}\\
-c^2\Delta u+\lambda w=h_3 \qquad &\text{in
$\Omega$,}\\
- \DivGamma (\sigma \nabla_\Gamma v)+\delta z+\kappa v+\rho w+\lambda\mu z =\mu h_4\qquad
&\text{on
$\Gamma_1$,}\\
\partial_\nu u=z\qquad
&\text{on
$\Gamma_1$,}\\
\partial_\nu u=0 &\text{on $\Gamma_0$,}
\end{cases}
\end{equation}
has a solution $u\in H^1(\Omega)$, $v\in H^1(\Gamma_1)$, $w\in H^1(\Omega)$, $z\in H^1(\Gamma_1)$, where
$\Delta u$, $\partial_\nu u_{|\Gamma}$ and $\DivGamma(\sigma\nabla_\Gamma v)$ are taken in the sense made precise in \S\ref{subsection3.3.2} and \S\ref{subsection3.3.5}. Since $u,w,h_1\in H^1(\Omega)$ and $v,z,h_2\in H^1(\Gamma_1)$ we can use the first two equations in \eqref{4.7} to eliminate $w$ and $z$, so writing \eqref{4.7} as
$$
\begin{cases}
-c^2\Delta u+\lambda^2u=\lambda h_1+h_3 \qquad &\text{in
$\Omega$,}\\
- \DivGamma (\sigma \nabla_\Gamma v)+\lambda\rho u+(\lambda^2\mu+\lambda\delta +\kappa)v =\rho h_1+(\lambda\mu+\delta)h_2+\mu h_4\qquad
&\text{on $\Gamma_1$,}\\
\partial_\nu u-\lambda v=h_2\qquad
&\text{on
$\Gamma_1$,}\\
\partial_\nu u=0 &\text{on $\Gamma_0$.}
\end{cases}
$$
Hence the statement holds true provided there is $\Lambda_1=\Lambda_1(\sigma, \rho, \kappa, \delta, \mu, c)\ge \Lambda_0$ such that the two coupled  problems
\begin{equation}\label{4.9}
\begin{cases}
-c^2\Delta u+\lambda^2u=h_1 \qquad &\text{in
$\Omega$,}\\
\partial_\nu u-\lambda v=h_2\qquad
&\text{on
$\Gamma_1$,}\\
\partial_\nu u=0 &\text{on $\Gamma_0$.}
\end{cases}
\end{equation}
\begin{equation}\label{4.10}
- \DivGamma (\sigma \nabla_\Gamma v)+\lambda\rho u+(\lambda^2\mu+\lambda\delta +\kappa)v =h_3\qquad
\text{on $\Gamma_1$,}
\end{equation}
have a solution for all $h_1\in L^2(\Omega)$ and $h_2,h_3\in L^2(\Gamma_1)$, provided $|\lambda|\ge \Lambda_1$.

By \eqref{3.40} and \eqref{3.39}  problems \eqref{4.9} and \eqref{4.10} respectively mean
\begin{equation}\label{4.11}
c^2\int_\Omega \nabla u \nabla \overline{\varphi}+\lambda^2 \int_\Omega u \overline{\varphi}-\lambda c^2\int_{\Gamma_1} v\overline{\varphi}=
\int_\Omega h_1\overline{\varphi}+c^2\int_{\Gamma_1} h_2\overline{\varphi}
\end{equation}
for all $\varphi\in H^1(\Omega)$, and
\begin{equation}\label{4.12}
\int_{\Gamma_1} \sigma (\nabla_\Gamma v,\nabla_\Gamma \psi)_\Gamma+\lambda\int_{\Gamma_1}\rho u\overline{\psi}+\int_{\Gamma_1}(\lambda^2\mu+\lambda\delta+\kappa)v\overline{\psi}=\int_{\Gamma_1}h_3\overline{\psi}
\end{equation}
for all $\psi\in H^1(\Gamma_1)$. Since, by assumption (A2), $\rho\in W^{1,\infty}(\Gamma_1)$ and, as remarked in the proof of previous Lemma,  $1/\rho\in W^{1,\infty}(\Gamma_1)$, by Lemma~\ref{Multiplierlemma}  we can replace the test functions $\psi$ in  equation \eqref{4.12} with $\frac\psi\rho$,  $\psi$ being an arbitrary element of $H^1(\Gamma_1)$, so equivalently writing \eqref{4.12} as
\begin{equation}\label{4.13}
\int_{\Gamma_1} \sigma \left(\nabla_\Gamma v, \nabla_\Gamma \left(\tfrac\psi\rho\right)\right)_\Gamma+\lambda \int_{\Gamma_1} u\overline{\psi}+\int_{\Gamma_1}\tfrac{\lambda^2\mu+\lambda\delta+\kappa}\rho v\overline{\psi}=\int_{\Gamma_1}\frac{h_3\overline{\psi}}\rho.
\end{equation}
Since $\varphi$ and $\psi$ are independent, the system composed by the equations \eqref{4.11} and \eqref{4.13} is equivalent to the single equation obtained by summing \eqref{4.11} and   \eqref{4.13} multiplied by $c^2$, that is to
\begin{multline}\label{4.14}
c^2\int_\Omega \nabla u \nabla \overline{\varphi}+\lambda^2 \int_\Omega u\overline{\varphi}-\lambda c^2\int_{\Gamma_1} v\overline{\varphi}
+c^2\int_{\Gamma_1} \sigma \left(\nabla_\Gamma v,\nabla_\Gamma \left(\tfrac\psi\rho\right) \right)_\Gamma\\
+c^2\lambda \int_{\Gamma_1} u\overline{\psi}+c^2\int_{\Gamma_1}\tfrac{\lambda^2\mu+\lambda\delta+\kappa}\rho v\overline{\psi}
=\int_\Omega h_1\overline{\varphi}+c^2\int_{\Gamma_1} h_2\overline{\varphi}+c^2\int_{\Gamma_1}\frac{h_3\overline{\psi}}\rho
\end{multline}
for all $(\varphi,\psi)\in H^1(\Omega)\times H^1(\Gamma_1)$.
To convince themselves of the last assertion, the reader can simply check \eqref{4.14} with the couples $(\varphi,0)$ and $(0, \psi)$.
Hence, setting for any $\lambda\in\R$ the sesquilinear  form $a_\lambda$ on $H^1(\Omega)\times H^1(\Gamma_1)$, by
\begin{equation}\label{4.15}
\begin{aligned}
a_\lambda(U,\Phi)&=c^2\int_\Omega \nabla u \nabla \overline{\varphi} +c^2\int_{\Gamma_1} \sigma \left(\nabla_\Gamma v,\nabla_\Gamma \left(\tfrac\psi\rho\right) \right)_\Gamma+\lambda^2 \int_\Omega u\overline{\varphi}\\
&- c^2\lambda\int_{\Gamma_1}v\overline{\varphi}
+c^2\lambda \int_{\Gamma_1} u\overline{\psi}+c^2\int_{\Gamma_1}\tfrac{\lambda^2\mu+\lambda\delta+\kappa}\rho v\overline{\psi},
\end{aligned}
\end{equation}
and $\eta'\in [H^1(\Omega)\times H^1(\Gamma_1)]'$ defined by
\begin{equation}\label{4.16}
  \langle \eta',\Phi\rangle_{H^1(\Omega)\times H^1(\Gamma_1)}=\int_\Omega\overline{h_1}\varphi+c^2\int_{\Gamma_1} \overline{h_2}\varphi+c^2\int_{\Gamma_1}\frac{\overline{h_3}\psi}\rho,
\end{equation}
where $U=(u,v)$ , $\Phi=(\varphi,\psi)$, equation \eqref{4.14} can be written as
\begin{equation}\label{4.17}
 a_\lambda(U,\Phi)=\overline{\langle \eta',\Phi\rangle}_{H^1(\Omega)\times H^1(\Gamma_1)}\qquad\text{for all $\Phi\in H^1(\Omega)\times H^1(\Gamma_1)$.}
\end{equation}
By assumptions (A1--3) and the same arguments used in the proof of Lemma~\ref{lemma4.1} one easily gets that the form $a_\lambda$ is continuous.

We now claim that there is $\Lambda_1=\Lambda_1(\sigma, \rho, \kappa, \delta, \mu, c)\ge \Lambda_0$ such that $a_\lambda$ is also coercive provided $|\lambda|\ge \Lambda_1$. This claim then concludes the proof by using the complex version of the Lax--Milgram Theorem.

To prove our claim we note that, by \eqref{4.15}, using the formulas
$$\nabla_\Gamma(v/\rho)=\nabla_\Gamma (1/\rho)v+\nabla_\Gamma v/\rho, \qquad\text{and}\quad \nabla_\Gamma (1/\rho)=-\nabla_\Gamma \rho/\rho^2$$
recalled in the previous proof, for all $U=(u,v)\in H^1(\Omega)\times H^1(\Gamma_1)$ we have
\begin{align*}
\Real a_\lambda(U,U)=&c^2\|\nabla u\|_2^2+c^2\int_{\Gamma_1}\tfrac\sigma\rho |\nabla_\Gamma v|_\Gamma^2+\lambda^2\|u\|_2^2
+c^2\int_{\Gamma_1}\tfrac{\lambda^2\mu+\lambda\delta+\kappa}\rho |v|^2\\
&\quad -c^2\Real \int_{\Gamma_1}\tfrac\sigma{\rho^2}(\nabla_\Gamma v,\nabla_\Gamma \rho)_\Gamma \overline{v}.
\end{align*}
Consequently, by using assumptions (A1--3), H\"{o}lder and weighted Young inequality, for any $\eps>0$ we have
\begin{align*}
\Real a_\lambda(U,U)\ge
&c^2\|\nabla u\|_2^2+c^2\int_{\Gamma_1}\tfrac\sigma\rho |\nabla_\Gamma v|_\Gamma^2+\lambda^2\|u\|_2^2
-\tfrac{c^2}{\rho_0}\|\nabla_\Gamma\|_{\infty,\Gamma_1}\int_{\Gamma_1}\tfrac\sigma \rho|\nabla_\Gamma v|_\Gamma|v|\\
&\quad +c^2\Big(\mu_0\lambda^2-\|\delta\|_{\infty,\Gamma_1}|\lambda|-\|\kappa\|_{\infty,\Gamma_1}\Big)\int_{\Gamma_1}\tfrac 1\rho |v|^2\\
&\ge c^2\|\nabla u\|_2^2+c^2\Big(1-\tfrac{\|\nabla_\Gamma\rho\|_{\infty,\Gamma_1}}{2\rho_0}\eps\Big)\int_{\Gamma_1}\tfrac\sigma\rho |\nabla_\Gamma v|_\Gamma^2+\lambda^2\|u\|_2^2\\
&\quad +c^2\Big(\mu_0\lambda^2-\|\delta\|_{\infty,\Gamma_1}|\lambda|-\|\kappa\|_{\infty,\Gamma_1}-\tfrac{\|\nabla_\Gamma\rho\|_{\infty,\Gamma_1}\|\sigma\|_{\infty,\Gamma_1}}
{2\rho_0\eps}\Big)\int_{\Gamma_1}\tfrac 1\rho |v|^2\\
\ge &c^2\|\nabla u\|_2^2+\tfrac{c^2\sigma_0}{\rho_0}\Big(1-\tfrac{\|\nabla_\Gamma\rho\|_{\infty,\Gamma_1}}{2\rho_0}\eps\Big)
\|\nabla_\Gamma v\|_{2,\Gamma_1}^2+\lambda^2\|u\|_2^2\\
&\quad +\tfrac{c^2}{\rho_0}\Big(\mu_0\lambda^2-\|\delta\|_{\infty,\Gamma_1}|\lambda|-\|\kappa\|_{\infty,\Gamma_1}-\tfrac{\|\nabla_\Gamma\rho\|_{\infty,\Gamma_1}\|\sigma\|_{\infty,\Gamma_1}}
{2\rho_0\eps}\Big)\|v\|_{2,\Gamma_1}^2.
\end{align*}
Hence, fixing $\eps=\eps_1$ such that $\|\nabla_\Gamma\rho\|_{\infty,\Gamma_1}\eps_1\le \rho_0$ we have
\begin{align*}
\Real a_\lambda(U,U)\ge & c^2\|\nabla u\|_2^2+\tfrac{c^2\sigma_0}{2\rho_0}\|\nabla_\Gamma v\|_{2,\Gamma_1}^2+\lambda^2\|u\|_2^2\\
&\quad +\tfrac{c^2}{\rho_0}\Big(\mu_0\lambda^2-\|\delta\|_{\infty,\Gamma_1}|\lambda|-\|\kappa\|_{\infty,\Gamma_1}-\tfrac{\|\nabla_\Gamma\rho\|_{\infty,\Gamma_1}\|\sigma\|_{\infty,\Gamma_1}}
{2\rho_0\eps_1}\Big)\|v\|_{2,\Gamma_1}^2,
\end{align*}
so there is $\Lambda_1=\Lambda_1(\sigma, \rho, \kappa, \delta, \mu, c)\ge \Lambda_0$ such that when $|\lambda|\ge \Lambda_1$
we have $$\mu_0\lambda^2-\|\delta\|_{\infty,\Gamma_1}|\lambda|-\|\kappa\|_{\infty,\Gamma_1}-\tfrac{\|\nabla_\Gamma\rho\|_{\infty,\Gamma_1}\|\sigma\|_{\infty,\Gamma_1}}
{2\rho_0\eps_1}\ge \tfrac{\mu_0}2\lambda^2$$
and consequently
$$\Real a_\lambda(U,U)\ge  \min\left\{c^2, \tfrac{c^2\sigma_0}{2\rho_0}, \tfrac{c^2\mu_0\Lambda_1^2}{2\rho_0},\Lambda_1^2\right\}
\Big(\|\nabla u\|_2^2+\|\nabla_\Gamma v\|_{2,\Gamma_1}^2+\|u\|_2^2
+\|v\|_{2,\Gamma_1}^2\Big),$$
proving our claim and concluding the proof.
\end{proof}

The following result shows that problem \eqref{4.6} is well-posed and gives higher regularity of its solutions in an abstract sense.
\begin{thm}\label{Theorem4.1} Under assumptions (A0--3) the operators $-A-\Lambda_1I$ and  $A-\Lambda_1I$ are dissipative and densely defined, hence each of them generates a contraction semigroup on $\cal H$. Consequently
\renewcommand{\labelenumi}{{\roman{enumi})}}
\begin{enumerate}
\item the operator $-A$ generates a strongly continuous group $\{T(t), t\in\R\}$ on $\cal H$ and, for any $U_0\in\cal H$, problem \eqref{4.6} has a  unique generalized  solution $U$, given by $U(t)=T(t)U_0$, which is also strong provided $U_0\in D(A)$;
\item recursively defining $D(A^n)$, $n\in\N$, by $D(A^{n+1})=\{u\in D(A^n): Au\in D(A^n)\}$, then $D(A^n)$ is an Hilbert space when endowed with the inner product given, for all
    $V,W\in D(A^n)$, by
    \begin{equation}\label{innerDAN}
    (V,W)_{D(A^n)}=\sum_{i=0}^n (A^i V, A^i W)_{\cal H};
    \end{equation}
    \item for all $n\in\N$ the restriction of the operator $-A$ on $D(A^n)$, with domain $D(A^{n+1})$, generates a strongly continuous group on it.
    Consequently, denoting $D(A^0)=\cal H$,  the solution $U$ of \eqref{4.6} enjoys the further regularity $U\in \bigcap _{i=0}^n C^{n-i}(\R;D(A^i))$ provided $U_0\in D(A^n)$.
\end{enumerate}
\end{thm}
\begin{proof} All assertions follow by Lemmas \ref{lemma4.1}--\ref{lemma4.2} and standard semigroup theory. In particular, by Lemma \ref{lemma4.1}, using  \cite[Theorem 4.2,~p.~14]{pazy} or \cite[Chapter II, Proposition 3.23,~p.~88]{EngelNagel}, the operators $B_\pm:=\pm A-\Lambda_1I$ are dissipative. Moreover, by Lemma \ref{lemma4.2}, for any $\lambda>0$, the operators $\lambda I-B_\pm$ are surjective, so by Lumer and Philips Theorem (see \cite[Theorem 4.3~p.14 and Theorem 4.6~p.16]{pazy} or \cite[Chapter II, Corollary 3.20~p.86]{EngelNagel}), $D(A)$ is dense in $\cal H$ and $B_\pm$ generate two contraction semigroups on it.

To prove i) we then note that, by the standard rescaling argument (see \cite[p.12]{pazy} or \cite[Example 2.2~p.60]{EngelNagel}), the operators $\pm A$ generate two strongly continuous semigroups. So, by \cite[Chapter II, Generation Theorem for Groups, p.~79]{EngelNagel}, $-A$ generates a strongly continuous group.~The proof of i) is then completed by \cite[p.105]{pazy} or \cite[Chapter II, Proposition~6.4~and~Theorem~6.7,~pp.146--150]{EngelNagel}.

To prove ii) we note that, since by \cite[Chapter II, Theorem~1.4~p.~51]{EngelNagel} the operator $A$ is closed, one trivially gets by induction that $D(A^n)$ is complete with respect to the norm $\|\cdot\|_{D(A^n)}$ induced by \eqref{innerDAN}.

Moreover, since by Lemma~\ref{lemma4.1} the operator $A+2\Lambda_1 I$ is bijective and there is $c_{30}(n)=c_{30}(n,\sigma, \rho, \kappa, \delta, \mu, c)>0$ such that  $\|(A+2\Lambda_1 I)^n \cdot\|_{\cal H}\le c_{30} \|\cdot\|_{D(A^n)}$, an equivalent norm on $D(A^n)$ is given by  $\|(A+2\Lambda_1 I)^n \cdot \|_{\cal H}$, used in \cite[Chapter II, Definition~5.1~p.~124]{EngelNagel}. Hence, by \cite[Chapter II, Proposition~5.2~p.124]{EngelNagel}, we get that $\pm A$ generate two strongly continuous semigroup on $D(A^n)$ and then, by the already recalled Generation Theorem for Groups, $-A$ generates a group on it. To complete the proof of iii) we simply remark that, when $U_0\in D(A^n)$,  it follows that $U\in C(\R, D(A^i))$ and $U^{(n-i)}=(-A)^{n-i}U$ for $i=0,\ldots,n$.
\end{proof}

\section{Well-posedness for problem \eqref{1.1}} This section is devoted to prove Theorem~\ref{theorem1.1}. Hence we shall first make precise what we mean by a  strong, generalized and weak solution of it.

By a {\em strong} or {\em generalized solution} of \eqref{1.1} we simply mean a couple $(u,v)$ such that $u$ and $v$ are the first two components of a solution $U=(u,v,w,z)$ of \eqref{4.6} of the same type. Trivially strong solutions are also generalized ones.
By \eqref{4.5} and \eqref{4.6} one immediately gets that any strong solution $(u,v)$ belongs to
$C^1(\R;H^1(\Omega)\times H^1(\Gamma_1))\cap C^2(\R;L^2(\Omega)\times L^2(\Gamma_1))$,
and that the corresponding strong solution $U$ of \eqref{4.6} is nothing but $(u,v,u_t,v_t)$.

Consequently, given any generalized solution $(u,v)$, denoting by $U=(u,v,w,z)$ any corresponding generalized solution of \eqref{4.6}, taking a sequence $(u^n,v^n,w^n,z^n)=(u^n,v^n,u^n_t,v^n_t))$ of strong solutions converging to it
in the topology of $C(\R;{\cal H})$,   for any $\varphi\in C^\infty_c(\R)$ we have $\int_{-\infty}^\infty u^n_t\varphi=-\int_{-\infty}^\infty u^n\varphi'$ in $L^2(\Omega)$, so passing to the limit as $n\to\infty$ we get $\int_{-\infty}^\infty w\varphi=-\int_{-\infty}^\infty u\varphi'$ in $L^2(\Omega)$ and consequently $w=u_t$. By the same arguments we prove that $z=v_t$, so $U=(u,v,u_t,v_t)$ and
$(u,v)\in C(\R;H^1(\Omega)\times H^1(\Gamma_1))\cap C^1(\R;L^2(\Omega)\times L^2(\Gamma_1))$.

To define weak solution of \eqref{1.1} we first make precise the meaning of weak solutions of the evolution boundary value problem
\begin{equation}\label{1.1bis}
\begin{cases} u_{tt}-c^2\Delta u=0 \qquad &\text{in
$\R\times\Omega$,}\\
\mu v_{tt}- \DivGamma (\sigma \nabla_\Gamma v)+\delta v_t+\kappa v+\rho u_t =0\qquad
&\text{on
$\R\times \Gamma_1$,}\\
v_t =\partial_\nu u\qquad
&\text{on
$\R\times \Gamma_1$,}\\
\partial_\nu u=0 &\text{on $\R\times \Gamma_0$.}
\end{cases}
\end{equation}
\begin{definition}
\label{Definition4.2bis}
Given
$$u\in L^1_\loc(\R;H^1(\Omega))\cap W^{1,1}_\loc(\R;L^2(\Omega)),\,v\in L^1_\loc(\R;H^1(\Gamma_1))\cap W^{1,1}_\loc(\R;L^2(\Gamma_1)),$$
we say that $(u,v)$ is a {\em weak solution} of \eqref{1.1bis} provided
the distributional identities
\begin{equation}\label{4.B1}
 \int_{-\infty}^\infty\left[-\int_\Omega u_t\varphi_t+c^2\int_\Omega\nabla u\nabla \varphi-c^2\int_{\Gamma_1}v_t\varphi\right]=0,
\end{equation}
\begin{equation}\label{4.B2}
 \int_{-\infty}^\infty\int_{\Gamma_1}\left[-\mu v_t\psi_t+\sigma(\nabla_\Gamma v,\nabla_\Gamma \overline{\psi})_\Gamma+\delta v_t\psi
 +\kappa v\psi-\rho u\psi_t\right]=0,
 \end{equation}
hold for all $\varphi\in C^\infty_c(\R\times\R^N)$ and $\psi\in C^r_c(\R\times\Gamma_1)$.
\end{definition}
Trivially any weak solution  of \eqref{1.1bis} possesses  a (unique) representative  $(u,v)\in C(\R;L^2(\Omega)\times L^2(\Gamma_1))$.
In the sequel we shall always consider this representative, so $u$ and $v$ possess a pointwise meaning.
It is a bit more involved to give this meaning also to $u_t$ and $v_t$ and it   requires an additional light regularity requirement.
\begin{lem}\label{Lemma4B} Let assumptions (A0--3) hold and $(u,v)$ be  a weak solution of \eqref{1.1bis}. Then
\renewcommand{\labelenumi}{{\roman{enumi})}}
\begin{enumerate}
\item $u_t$ has a unique representative in $C(\R;[H^1(\Omega)]')$, which satisfies the alternative distributional identity
\begin{equation}\label{4.19}
 \int_s^t\left[-\int_\Omega u_t\varphi_t+c^2\int_\Omega\nabla u\nabla \varphi-c^2\int_{\Gamma_1}v_t\varphi\right]+
 \left[\langle u_t,\varphi\rangle_{H^1(\Omega)}\right]_s^t=0,
\end{equation}
for all $s,t\in\R$ and  $\varphi\in C(\R;H^1(\Omega))\cap C^1(\R;L^2(\Omega))$;
\item if $u_{|\Gamma_1}\in C(\R; H^{-1}(\Gamma_1))$ then also $v_t$ possesses  a unique representative in $C(\R;H^{-1}(\Gamma_1))$, which satisfies the alternative distributional identity
    \begin{equation}\label{4.20}
 \begin{aligned}
 \int_s^t\int_{\Gamma_1}\left[-\mu v_t\psi_t+\sigma(\nabla_\Gamma v,\nabla_\Gamma \overline{\psi})_\Gamma+\delta v_t\psi
 +\kappa v\psi-\rho u\psi_t\right]\\
 +\left[\langle \mu v_t+\rho u_{|\Gamma_1},\psi\rangle_{H^1(\Gamma_1)}\right]_s^t=0,
\end{aligned}
\end{equation}
for all $s,t\in\R$ and  $\psi\in C(\R;H^1(\Gamma_1))\cap C^1(\R;L^2(\Gamma_1))$.
\end{enumerate}
\end{lem}
\begin{rem} The requirement $u_{|\Gamma_1}\in C(\R; H^{-1}(\Gamma_1))$ may appear difficult to verify. On the other hand  it is automatic when
$u\in C(\R; H^1(\Omega))$. Trivially weak solutions of \eqref{1.1bis} are characterized by the alternative distributional identities \eqref{4.19} and \eqref{4.20}  among functions verifying the regularity prescribed by Definition~\ref{Definition4.2bis} and such that  $u_{|\Gamma_1}\in C(\R; H^{-1}(\Gamma_1))$.
\end{rem}
\begin{proof}[Proof of Lemma~\ref{Lemma4B}] Let $(u,v)$ be a weak solution of \eqref{1.1bis}. Taking test functions $\varphi$ for \eqref{4.B1}
in the separate form $\varphi(t,x)=\varphi_1(t)\varphi_0(x)$, with $\varphi_1\in C^\infty_c(\R)$ and $\varphi_0\in C^\infty_c(\R^N)$, we get
$$\int_{-\infty}^\infty\left[\left(-\int_\Omega u_t\varphi_0\right)\varphi_1'+c^2\left(\int_\Omega\nabla u\nabla \varphi_0-\int_{\Gamma_1}v_t\varphi_0\right)\varphi_1\right]=0.$$
Hence, for any fixed $\varphi_0\in C^\infty_c(\R^N)$ the function $t\mapsto \int_\Omega u_t(t)\varphi_0$ belongs to $W^{1,1}_\loc(\R)$
and $\left(\int_\Omega u_t\varphi_0\right)'=c^2\left(-\int_\Omega \nabla u\nabla\varphi_0+\int_{\Gamma_1}v_t\varphi_0\right)$ a.e.\ in $\R$. Consequently
\begin{equation}\label{4.B3}
\int_\Omega u_t(t)\varphi_0-\int_\Omega u_t(s)\varphi_0=\int_s^t c^2\left(-\int_\Omega \nabla u\nabla\varphi_0+\int_{\Gamma_1}v_t\varphi_0\right)
\end{equation}
for almost all $s,t\in\R$. Since restrictions of functions in $C^\infty_c(\R^N)$ are dense in $H^1(\Omega)$ (see \cite[Theorem~11.35 p.~330]{LeoniSobolev2}), taking for all $\varphi_0\in H^1(\Omega)$ a sequence $(\varphi_{0n})_n$ in $C^\infty_c(\R^N)$ with ${\varphi_{0n}}_{|\Omega}\to \varphi_0$ in $H^1(\Omega)$ we get that \eqref{4.B3} holds true for a.a. $s,t\in\R$ and all $\varphi_0\in H^1(\Omega)$.

Consequently, denoting by $\cal{B}\in\cal{L}(H^1(\Omega)\times L^2(\Gamma_1); [H^1(\Omega)]')$ the operator defined by
$\langle \cal{B}(\mathfrak{u},\mathfrak{v}),\varphi_0\rangle_{H^1(\Omega)}=-\int_\Omega \nabla \mathfrak{u}\nabla\varphi_0+\int_{\Gamma_1}\mathfrak{v}\varphi_0$, we have that $(u_t)'=c^2\cal{B}(u,v_t)\in L^1_\loc (\R;[H^1(\Omega)]')$ weakly.
Hence $u_t\in W^{1,1}_\loc (\R;[H^1(\Omega)]')$, so it has a (unique) representative in $C(\R;[H^1(\Omega)]')$.

By a standard density argument in the space $W^{1,1}(a,b; H^1(\Omega))$, $a,b\in\R$,  we then get that for any $\varphi\in C^1(\R; H^1(\Omega))\subset W^{1,1}_\loc(\R; H^1(\Omega))$ we have $\langle u_t,\varphi\rangle_{H^1(\Omega)}\in W^{1,1}_\loc (\R)$ and
$\langle u_t,\varphi\rangle_{H^1(\Omega)}'=\langle u_t',\varphi\rangle_{H^1(\Omega)}+\langle u_t,\varphi_t\rangle_{H^1(\Omega)}$ a.e.\ in $\R$,
from which \eqref{4.19} follows for test functions $\varphi\in C^1(\R; H^1(\Omega))$.
By standard time regularization then  \eqref{4.19} holds for test functions $\varphi\in C(\R;H^1(\Omega))\cap C^1(\R;L^2(\Omega))$, completing the proof of i).

To prove ii) we remark that, taking test functions $\psi$ for \eqref{4.B2}
in the separate form $\psi(t,x)=\psi_1(t)\psi_0(x)$, with $\psi_1\in C^\infty_c(\R)$ and $\psi_0\in C^r_c(\Gamma_1)$, and using the same arguments as before we get
\begin{multline}\label{4.B4}
 \int_{\Gamma_1} [\mu v_t(t)+\rho u_{|\Gamma_1}(t)]\psi_0-\int_{\Gamma_1} [\mu v_t(s)+\rho u_{|\Gamma_1}(s)]\psi_0\\
 =-\int_s^t \left[\int_{\Gamma_1} \sigma(\nabla_\Gamma v,\nabla_\Gamma \overline{\psi_0})_\Gamma+\int_{\Gamma_1}\delta v_t\psi_0+\int_{\Gamma_1}\kappa v\psi_0\right]
\end{multline}
 for all $\psi_0\in C^r_c(\Gamma_1)$ and a.a. $s,t\in\R$. Using the density of $C^r_c(\Gamma_1)$ in $H^1(\Gamma_1)$ (see Corollary~\ref{corollary3.1}) we  get that \eqref{4.B4} holds for all $\psi_0\in H^1(\Gamma_1)$. Consequently $\mu v_t+\rho u_{|\Gamma_1}\in W^{1,1}_\loc(\R;H^{-1}(\Gamma_1))$ and
 $(\mu v_t+\rho u_{|\Gamma_1})'=\DivGamma(\sigma\nabla_\Gamma v)-\delta v_t-kv$ weakly in this space. Hence $v_t$ has a (unique) representative such that
 $\mu v_t+\rho u_{|\Gamma_1}\in C(\R;H^{-1}(\Gamma_1)$.
 Then, as $u_{|\Gamma_1}\in C(\R;H^{-1}(\Gamma_1)$ and $1/\mu$, $\rho\in L^\infty(\Gamma_1)$, by Lemma~\ref{Multiplierlemma}
we have $v_t\in C(\R;H^{-1}(\Gamma_1)$. Using the same density argument used before for the space $W^{1,1}(a,b; H^{-1}(\Gamma_1))$
we then get
$$\langle \mu v_t+\rho u_{|\Gamma_1},\psi\rangle_{H^1(\Gamma_1)}'=\langle (\mu v_t+\rho u_{|\Gamma_1})',\psi\rangle_{H^1(\Gamma_1)}+\langle \mu v_t+\rho u_{|\Gamma_1},\psi_t\rangle_{H^1(\Gamma_1)}$$
a.e.\ in $\R$, from which \eqref{4.20} follows for test functions $\psi\in C^1(\R; H^1(\Gamma_1))$.
By standard time regularization then  \eqref{4.20} holds for test functions $\psi\in C(\R;H^1(\Omega))\cap C^1(\R;L^2(\Omega))$.
\end{proof}
By Lemma~\ref{Lemma4B} we can define weak solutions of problem \eqref{1.1}. In the proof of Theorem~\ref{theorem1.1} we shall only need data in $\cal{H}$, but to prove uniqueness it is convenient to consider a more general case.
\begin{definition}
\label{Definition4.2}For any $U_0\in L^2(\Omega)\times L^2(\Gamma_1)\times[H^1(\Omega)]'\times H^{-1}(\Gamma_1)$ we say that
$(u,v)$ is a {\em weak solution} of \eqref{1.1} provided it is a weak solution of \eqref{1.1bis}, $u_{|\Gamma_1}\in C(\R;H^{-1}(\Gamma_1))$ and its representative $(u,v)\in C(\R;L^2(\Omega)\times L^2(\Gamma_1))\cap C^1(\R;[H^1(\Omega)]'\times H^{-1}(\Gamma_1))$
satisfies $(u(0), v(0), u_t(0),v_t(0))=U_0$.
\end{definition}

Trivially, by \eqref{3.40}, \eqref{3.36}, \eqref{4.5} and \eqref{4.6}, strong solutions of \eqref{1.1} are also weak. More generally, since Definition~\ref{Definition4.2} is stable with respect to convergence in
$C(\R,H^1(\Omega)\times H^1(\Gamma_1))\cap C^1(\R,L^2(\Omega)\times L^2(\Gamma_1)),$
also generalized solutions are weak solutions. The following result allows us to prove that
the classes of weak and generalized solutions of problem \eqref{1.1} coincide when data are taken in $\cal{H}$.
\begin{lem}\label{lemma4.3} Under assumptions (A0--3), weak solutions of \eqref{1.1} are unique.
\end{lem}
\begin{proof} By linearity, proving uniqueness reduces to prove that $U_0=0$ implies $(u,v)=0$ in $\R$. Moreover, since when $U_0=0$ the couple $(\hat{u}, \hat{v})$ given by $\hat{u}(t)=-u(-t)$, $\hat{v}(t)=v(-t)$ is still a weak solution of \eqref{1.1} with vanishing data provided $\delta$ is replaced by $-\delta$, we shall just prove that $(u,v)=0$ in $[0,\infty)$.

To prove it, adapting the argument in \cite[Proof of Theorem~4,~p.~406]{Evans},
we fix $t>0$ and  test functions $\varphi$ and $\psi$, depending on $t$, given by
\begin{equation}\label{4.20bis}
\varphi(s)=
\begin{cases}
\int_s^t \overline{u}(\tau)\,d\tau\quad &\text{if $s\le t$},\\
\overline{u}(t)(t-s) &\text{if $s\ge t$},
\end{cases}
\qquad
\psi(s)=
\begin{cases}
\frac {c^2}\rho\int_s^t \overline{v}(\tau)\,d\tau\quad &\text{if $s\le t$},\\
\frac {c^2}\rho \overline{v}(t)(t-s) &\text{if $s\ge t$},
\end{cases}
\end{equation}
so that, since $1/\rho\in W^{1,\infty}(\Gamma_1)$ and using Lemma~\ref{Multiplierlemma}, $\varphi\in C(\R;H^1(\Omega))\cap C^1(\R;L^2(\Omega))$, $\psi\in C(\R;H^1(\Gamma_1))\cap C^1(\R;L^2(\Gamma_1))$,
$$\varphi_t=-\overline{u},\quad \quad \psi_t=-\tfrac {c^2}\rho\overline{v}\quad\text{in $[0,t]$,}\qquad  \varphi(t)=0,\quad\text{and}\quad \psi(t)=0.$$
We then apply Lemma~\ref{Lemma4B}. As $u_0=u_1=0$ and $v_1=0$, when $s=0$  the sum of the distribution identities \eqref{4.19}--\eqref{4.20} reads as
\begin{multline}\label{4.21}
  \int_0^t\left[\int_\Omega u_t\overline{u}-c^2\int_\Omega\nabla\varphi\nabla\overline{\varphi}_t-c^2\int_{\Gamma_1} v_t\varphi
 +c^2\int_{\Gamma_1}\tfrac\mu \rho v_t\overline{v}\right.\\
 \left.-\frac 1{c^2}\int_{\Gamma_1}\sigma (\nabla_\Gamma (\rho\overline{\psi}_t),\nabla_\Gamma\overline{\psi})_\Gamma+\int_{\Gamma_1}\delta  v_t\psi-\frac 1{c^2}\int_{\Gamma_1}\kappa \rho \psi\overline{\psi}_t+ c^2\int_{\Gamma_1}u\overline{v}\right]=0.
\end{multline}
Since $v_0=\psi(t)=0$, integrating by parts in time we have
$$-c^2\int_0^t\int_{\Gamma_1} v_t\varphi=-c^2\int_0^t\int_{\Gamma_1} \overline{u}v,\quad \text{and}\quad
\int_0^t\int_{\Gamma_1} \delta v_t\psi=c^2\int_0^t\int_{\Gamma_1}\frac\delta\rho|v|^2.$$
Moreover, by Leibniz formula $\nabla_\Gamma (\rho\overline{\psi}_t)=\overline{\psi}_t \nabla_\Gamma\rho +\rho\nabla_\Gamma\overline{\psi}_t=
-\frac{c^2}\rho v \nabla_\Gamma \rho+\rho\nabla_\Gamma\overline{\psi}_t$.
Hence we can rewrite \eqref{4.21} as
\begin{multline}\label{4.21bis}
  \int_0^t\left[\int_\Omega u_t\overline{u}-c^2\int_\Omega\nabla\varphi\nabla\overline{\varphi}_t
 +c^2\int_{\Gamma_1}\frac\mu \rho v_t\overline{v}-\frac 1{c^2}\int_{\Gamma_1}\sigma \rho(\nabla_\Gamma \overline{\psi}_t,\nabla_\Gamma\overline{\psi})_\Gamma\right.\\
 \left.+\int_{\Gamma_1}\frac \sigma \rho(\nabla_\Gamma \rho,\nabla_\Gamma\overline{\psi})_\Gamma v+c^2\int_{\Gamma_1}\frac \delta  \rho |v|^2-\frac 1{c^2}\int_{\Gamma_1}\kappa \rho \psi\overline{\psi}_t+ c^2\int_{\Gamma_1}(u\overline{v}-\overline{u}v)\right]=0.
\end{multline}
Taking the real part in \eqref{4.21bis} we get
\begin{multline*}
\frac 12\int_0^t\frac d{dt}\left[\|u\|_2^2-c^2\|\nabla\varphi\|_2^2+c^2\int_{\Gamma_1}\frac\mu\rho |v|^2- \frac 1{c^2}\int_{\Gamma_1}\sigma \rho|\nabla_\Gamma \psi|_\Gamma^2
-\frac 1{c^2}\int_{\Gamma_1}\kappa \rho|\psi|^2\right]\\ =
-\int_0^t\int_{\Gamma_1}\frac \sigma \rho\Real [(\nabla_\Gamma \rho, \nabla_\Gamma \overline{\psi})_\Gamma v]
-c^2\int_0^t \int_{\Gamma_1}\frac\delta \rho|v|^2.
\end{multline*}
Consequently, since $u_0=\varphi(t)=0$ and $v_0=\psi(t)=0$, we have
\begin{multline}\label{4.22}
\frac 12 \|u(t)\|_2^2+\frac{c^2}2\|\nabla\varphi(0)\|_2^2+\frac {c^2}2 \int_{\Gamma_1} \frac\mu\rho |v(t)|^2
+\frac 1{2c^2}\int_{\Gamma_1}\sigma\rho|\nabla_\Gamma\psi(0)|_\Gamma^2\\
=-\frac 1{2c^2}\int_{\Gamma_1}\rho\kappa|\psi(0)|^2- \int_0^t \int_{\Gamma_1}\frac \sigma \rho \Real [(\nabla_\Gamma \rho, \nabla_\Gamma \overline{\psi})_\Gamma v]-c^2 \int_0^t\int_{\Gamma_1}\frac \delta \rho|v|^2.
\end{multline}
Since test functions $\varphi$, $\psi$ in \eqref{4.22} change when $t$ changes, to let $t$ vary we simply skip the second term in its left--hand side (so getting an inequality) and we express $\psi$ in terms of the unvarying function $\Upsilon\in C(\R;H^1(\Gamma_1))\cap C^1(\R;L^2(\Gamma_1))$ defined by
$\Upsilon(t)=\frac {c^2}\rho\int_0^tv(\tau)\,d\tau$, so that $\overline{\psi}(\tau)=\Upsilon(t)-\Upsilon(\tau)$ for $\tau\in[0,t]$.

Proceeding in this way, since $\mu\ge \mu_0>0$, we get
\begin{multline}\label{4.23}
\frac 12\|u(t)\|_2^2+\frac {c^2\mu_0}2\int_{\Gamma_1} \frac {|v(t)|^2}\rho+\frac 1{2c^2}\int_{\Gamma_1} \sigma\rho |\nabla_\Gamma \Upsilon(t)|^2_\Gamma\le -\frac 1{2c^2} \int_{\Gamma_1}\rho\kappa|\Upsilon(t)|^2\\
-\int_0^t \int_{\Gamma_1}\frac\sigma\rho \Real [(\nabla_\Gamma \rho, \nabla_\Gamma \Upsilon(t)-\nabla_\Gamma \Upsilon(\tau))_\Gamma v(\tau)]\,d\tau-c^2\int_0^t\int_{\Gamma_1}\frac\delta \rho|v|^2.
\end{multline}
We now estimate, using assumptions (A1--3), the terms in the right--hand side of \eqref{4.23}. By H\"{o}lder inequality in time
\begin{equation}\label{4.26BIS}
\begin{aligned}
-\frac 1{2c^2} \int_{\Gamma_1}\rho\kappa|\Upsilon(t)|^2\le &\frac{\|\kappa\|_{\infty,\Gamma_1}}{2c^2}\left\|\int_0^t \frac {c^2}\rho v(\tau)\,d\tau\right\|
_{L^2(\Gamma_1;\rho d \cal{H}^{N-1})}^2\\
\le&\frac{\|\kappa\|_{\infty,\Gamma_1}}{2c^2}\left(\int_0^t \left\|\frac {c^2}\rho v(\tau)\right\|_{L^2(\Gamma_1;\rho d \cal{H}^{N-1})} d\tau\right)^2 \\
\le &c^2\|\kappa\|_{\infty,\Gamma_1}t \int_0^t \|v(\tau)/\rho \|_{L^2(\Gamma_1;\rho d \cal{H}^{N-1})}^2 d\tau\\
= & c^2\|\kappa\|_{\infty,\Gamma_1}t \int_0^t \int_{\Gamma_1}\frac {|v|^2}\rho.
\end{aligned}
\end{equation}
Trivially
\begin{equation*}\label{4.27-ex}
-c^2\int_0^t\int_{\Gamma_1}\frac\delta \rho|v|^2\le c^2\|\delta\|_{\infty,\Gamma_1}\int_0^t\int_{\Gamma_1}\frac{|v|^2}\rho.
\end{equation*}
Moreover
\begin{equation}\begin{aligned}\label{4.28}
  -&\int_0^t \int_{\Gamma_1}\frac\sigma\rho \Real [(\nabla_\Gamma \rho, \nabla_\Gamma \Upsilon(t)-\nabla_\Gamma \Upsilon(\tau))_\Gamma v(\tau)]\,d\tau
  \\\le &\int_0^t \int_{\Gamma_1}\frac\sigma\rho |\nabla_\Gamma \rho|_\Gamma |\nabla_\Gamma \Upsilon(t)-\nabla_\Gamma \Upsilon(\tau)|_\Gamma |v(\tau)|\,d\tau
  \\ \le
  &\|\nabla_\Gamma\rho\|_{\infty,\Gamma_1}\int_0^t \left[\int_{\Gamma_1}\frac\sigma\rho |v(\tau)|^2+\int_{\Gamma_1}\frac\sigma{2\rho}
  |\nabla_\Gamma \Upsilon(t)|_\Gamma^2+\int_{\Gamma_1}\frac \sigma{2\rho}|\nabla_\Gamma\Upsilon(\tau)|_\Gamma ^2\,\right]d\tau\\
  \le &\|\nabla_\Gamma\rho\|_{\infty,\Gamma_1}\|\sigma\|_{\infty,\Gamma_1}\int_0^t \int_{\Gamma_1}\frac{|v|^2}\rho +\frac{\|\nabla_\Gamma\rho\|_{\infty,\Gamma_1}}{2\rho_0^2}t\,\,\int_{\Gamma_1}\sigma\rho
  |\nabla_\Gamma \Upsilon(t)|_\Gamma^2\\
  +&\frac{\|\nabla_\Gamma\rho\|_{\infty,\Gamma_1}}{2\rho_0^2}\int_0^t\int_{\Gamma_1}\sigma\rho|\nabla_\Gamma\Upsilon|_\Gamma ^2.
\end{aligned}
\end{equation}
Plugging the estimates \eqref{4.26BIS}--\eqref{4.28} in \eqref{4.23} we get that for all $t\ge 0$
\begin{multline}\label{4.29BIS}
\frac 12\|u(t)\|_2^2+\frac {c^2\mu_0}2\int_{\Gamma_1} \frac {|v(t)|^2}\rho+\frac 1{2c^2}\int_{\Gamma_1} \sigma\rho |\nabla_\Gamma \Upsilon(t)|^2_\Gamma
\\
\le\Big[c^2\|\kappa\|_{\infty,\Gamma_1}t+c^2\|\delta\|_{\infty,\Gamma_1}+\|\nabla_\Gamma\rho\|_{\infty,\Gamma_1}\|\sigma\|_{\infty,\Gamma_1}\Big]\int_0^t \int_{\Gamma_1}\frac{|v|^2}\rho
\\+\frac{\|\nabla_\Gamma\rho\|_{\infty,\Gamma_1}}{2\rho_0^2}t\,\,\int_{\Gamma_1}\sigma\rho
  |\nabla_\Gamma \Upsilon(t)|_\Gamma^2+\frac{\|\nabla_\Gamma\rho\|_{\infty,\Gamma_1}}{2\rho_0^2}\int_0^t\int_{\Gamma_1}\sigma\rho|\nabla_\Gamma\Upsilon|_\Gamma ^2.
\end{multline}
When $\nabla_\Gamma\rho=0$ the estimate \eqref{4.29BIS} immediately yields,  by applying Grönwall inequality, that $(u,v)=0$ in $[0,\infty)$.
When $\nabla_\Gamma\rho\not=0$, fixing $t_1=\rho_0^2/2c^2\|\nabla_\Gamma \rho\|_{\infty,\Gamma_1}$,  \eqref{4.29BIS} implies
\begin{multline*}
\frac 12\|u(t)\|_2^2+\frac {c^2}{2\mu_0}\int_{\Gamma_1} \frac {|v(t)|^2}\rho+\frac 1{4c^2}\int_{\Gamma_1} \sigma\rho |\nabla_\Gamma \Upsilon(t)|^2_\Gamma
\le \frac{\|\nabla_\Gamma\rho\|_{\infty,\Gamma_1}}{2\rho_0^2}\int_0^t\int_{\Gamma_1}\sigma\rho|\nabla_\Gamma\Upsilon|_\Gamma ^2
\\ +\Big[c^2\|\kappa\|_{\infty,\Gamma_1}t+c^2\|\delta\|_{\infty,\Gamma_1}+\|\nabla_\Gamma\rho\|_{\infty,\Gamma_1}\|\sigma\|_{\infty,\Gamma_1}\Big]\int_0^t \int_{\Gamma_1}\frac{|v|^2}\rho
\end{multline*}
for $t\in [0,t_1]$. By applying Grönwall inequality again we then get  $(u,v)=0$ in $[0,t_1]$. The regularity of $(u,v)$ also gives
$(u,v,u_t,v_t)=0$ in $[0,t_1]$.

Now, by Definitions~\ref{Definition4.2bis} and \ref{Definition4.2}, for any $t^*\in\R$, the couple $(u_{t^*},v_{t^*})$ defined by $u_{t^*}(t)=u(t^*+t)$ and
$v_{t^*}(t)=v(t^*+t)$ for all $t\in\R$ is still a weak solution of problem \eqref{1.1} corresponding to initial data $(u(t^*), v(t^*),u_t(t^*),v_t(t^*))$.
Hence previous conclusion allows to conclude by induction on $n\in\N$ that $(u,v)=0$ in $[0,nt_1]$ for all $n\in\N$, concluding the proof.
 \end{proof}
We can now give the
\begin{proof}[Proof of Theorem~\ref{theorem1.1}] By Theorem~\ref{Theorem4.1}--i), problem \eqref{1.1} has a unique generalized solution, which is also a weak one as explained immediately before Lemma~\ref{lemma4.3}. By this Lemma this solution is also unique among weak solutions. Continuous dependence on data is an immediate consequence of the strong continuity of the group asserted in Theorem~\ref{Theorem4.1}--i). By \eqref{4.4}, \eqref{1.2ter} implies \eqref{1.2quater} and also implies that \eqref{1.1}$_1$ holds true in $L^2(\R\times\Omega)$ and consequently a.e.\ in $\R\times\Omega$. By the same type of argument also \eqref{1.1}$_2$--\eqref{1.1}$_3$
hold true a.e.\ in $\R\times\Gamma_1$ and \eqref{1.1}$_4$ does in $\R\times\Gamma_0$, with respect to the product measure on $\R\times\Gamma$, where
one uses Corollary~\ref{corollary3.1} to establish in the standard way the isomorphism $L^2(\R\times\Gamma_i)\simeq L^2(\R;L^2(\Gamma_i))$ for $i=0,1$.
Finally the energy identity \eqref{energyidentity} holds for strong solutions by a straightforward calculation and then for generalized (and so weak) solutions
by the density of $D(A)$ asserted in Theorem~\ref{Theorem4.1}.
\end{proof}

\chapter{Regularity when $r\ge2$.} \label{section5}
The spatial regularity of strong solutions of \eqref{1.1} asserted in Theorem~\ref{theorem1.1} can be improved provided $\Gamma$ is more regular, i.e., $r\ge 2$. Actually, when $r=1$, more regularity on $\Gamma$ would be meaningless.
To prove our assertion we need two preliminary regularity results. They are well-known when $\Gamma$ is compact.

The first one concerns the operator $\DivGamma(\sigma\nabla_\Gamma)$ defined in \eqref{3.41}, and hence in particular the Laplace--Beltrami operator
$\Delta_\Gamma$, to which it reduces when $\sigma\equiv 1$.
\begin{thm}\label{lemma4.5} If (A0--3) hold for any $s\in[-r+1,r-1]$ the operator $B_s:=-\DivGamma(\sigma\nabla_\Gamma)+I$ is an algebraic and topological isomorphism between $H^{s+1}(\Gamma_1)$ and $H^{s-1}(\Gamma_1)$.
\end{thm}
\begin{proof} When $s=0$, by \eqref{3.36} for any $u$ and $v$ in $H^1(\Gamma_1)$ we have
$$\langle B_0 u, \overline{v}\rangle_{H^1(\Gamma_1)} =\int_{\Gamma_1}\sigma(\nabla_\Gamma u,\nabla_\Gamma v)_\Gamma+\int_{\Gamma_1}u\overline{v}.$$
By assumption (A1) then $\langle B_0 u, \overline{v}\rangle_{H^1(\Gamma_1)}$ is an inner product in $H^1(\Gamma_1)$ inducing a norm equivalent to the one induced by \eqref{3.4}. Hence the statement follows by the Riesz--Fréchet Theorem.

When $r=1$ the proof is complete, so we suppose $r\ge 2$. We prove the statement for $s=m\in\N_0$, $0\le m\le r-1$, by induction on $m$, the basis case $m=0$ being true.  Hence we suppose that $m\ge 1$ and  $B_{m-1}^{-1}\in{\cal L}(H^{m-2}(\Gamma_1), H^m(\Gamma_1))$.
For any $h\in H^{m-1}(\Gamma_1)$ we then have $u=B_{m-1}^{-1}h\in H^m(\Gamma_1)$, $u$ being the unique solution of the equation $-\DivGamma(\sigma\nabla_\Gamma u)+u=h$ in $H^{-1}(\Gamma_1)$. By the splitting in \S\ref{subsection3.3.5} and \eqref{3.36} we one then easily gets that
$$-\Delta_\Gamma u=\tilde{h}:=\frac {h-u-(\nabla_\Gamma u,\nabla_\Gamma \sigma)_\Gamma}\sigma.$$
By Lemma~\ref{lemma3.4} and assumption (A1) we have $\nabla_\Gamma u,\nabla_\Gamma\sigma\in H^{m-1}(\Gamma_1)$, so  by Lemma~\ref{lemma3.5} one gets $(\nabla_\Gamma u,\nabla_\Gamma \sigma)_\Gamma\in H^{m-1}(\Gamma_1)$. Since by assumption (A1) and Lemma~\ref{lemma3.5} one also gets $1/\sigma \in W^{m,\infty}(\Gamma_1)$, by Lemma~\ref{Multiplierlemma} we have $\tilde{h}\in H^{m-1}(\Gamma_1)$. Hence the proof reduces to the case $\sigma\equiv1$, that is we can take $\tilde{h}=h$ in the sequel.

We now use the functions $\psi_n$ introduced in \eqref{3.21}. By \eqref{3.36} for all $n\in{\cal I}_{\Gamma_1}$ we have
\begin{equation}\label{4.24}
\begin{aligned}
-\Delta_\Gamma (\psi_nu)=&-u \Delta_\Gamma \psi_n -2(\nabla_\Gamma\psi_n,\nabla_\Gamma u)_\Gamma-\psi_n\Delta_\Gamma u\\
=&-u \Delta_\Gamma \psi_n -2(\nabla_\Gamma\psi_n,\nabla_\Gamma u)_\Gamma+\psi_n h:=\omega_n
\end{aligned}
\end{equation}
in $H^{-1}(\Gamma_1)$. Trivially
\begin{multline}\label{4.25}
\sum_{n\in {\cal I}_{\Gamma_1}}\|\omega_n\|_{m-1,2,\Gamma_1}^2 \le   3\left(\sum_{n\in {\cal I}_{\Gamma_1}}\|u \Delta_\Gamma \psi_n\|_{m-1,2,\Gamma_1}^2\right.\\
+ \left.\sum_{n\in {\cal I}_{\Gamma_1}}\|(\nabla_\Gamma\psi_n,\nabla_\Gamma u)_\Gamma\|_{m-1,2,\Gamma_1}^2+\sum_{n\in {\cal I}_{\Gamma_1}}\|\psi_n h\|_{m-1,2,\Gamma_1}^2\right).
\end{multline}
We now estimate the three summations in the right hand side of \eqref{4.25}, starting from the first one.
By Beppo Levi's Theorem,  \eqref{2.22estesa}, assumption (A0)--(ii)  and   \eqref{modulotensoriale} we get
\begin{align*}
\sum_{n\in {\cal I}_{\Gamma_1}}\|u \Delta_\Gamma \psi_n\|_{m-1,2,\Gamma_1}^2=&\sum_{n\in {\cal I}_{\Gamma_1}}\int_{\Gamma_1}\sum_{i=0}^{m-1}|D_\Gamma^i
(u \Delta_\Gamma \psi_n)|_\Gamma^2\\
=&\int_{\Gamma_1}\sum_{i=0}^{m-1}\sum_{n\in {\cal I}_{\Gamma_1}}\left|\sum_{j=0}^i \binom ij D_\Gamma^j
(\Delta_\Gamma \psi_n)\otimes D_\Gamma^{i-j}u\right|_\Gamma^2\\
\le &N_0
\int_{\Gamma_1}\sum_{i=0}^{m-1}\sum_{j=0}^i\binom ij\sum_{n\in {\cal I}_{\Gamma_1}} |D^j_\Gamma(\Delta_\Gamma \psi_n)|_\Gamma^2 |D^{i-j}_\Gamma u|_\Gamma^2.
\end{align*}
Since $m\le r-1$, by assumption (A0)--(ii) and \eqref{3.23bis} we have
$$\sum_{n\in {\cal I}_{\Gamma_1}} |D^j_\Gamma(\Delta_\Gamma \psi_n)|_\Gamma^2\le N_0 c_{19}^2,$$
so there is $c_{31}=c_{31}(m)>0$ such that
\begin{equation}\label{4.26}\sum_{n\in {\cal I}_{\Gamma_1}}\|u \Delta_\Gamma \psi_n\|_{m-1,2,\Gamma_1}^2\le
c_{31}\int_{\Gamma_1}\sum_{i=0}^{m-1}|D^i_\Gamma u|_\Gamma^2=c_{31}(m)\|u\|^2_{m-1,2,\Gamma_1}.
\end{equation}
To estimate the second one we note that, by Beppo Levi's Theorem,
\begin{equation}\label{4.26bis}
\sum_{n\in {\cal I}_{\Gamma_1}}\|(\nabla_\Gamma\psi_n,\nabla_\Gamma u)_\Gamma\|_{m-1,2,\Gamma_1}^2=
\int_{\Gamma_1}\sum_{i=0}^{m-1}\sum_{n\in {\cal I}_{\Gamma_1}}| D^i_\Gamma(\nabla_\Gamma\psi_n,\nabla_\Gamma u)_\Gamma|_\Gamma^2.
\end{equation}
To estimate the right--hand side  of the integrand we remark that, denoting  $\frak{g}(\cdot,\cdot):=(\cdot,\overline{\cdot})_\Gamma$,  clearly $\frak{g}\in C^{r-1}_2(\Gamma)$ and we can write $(\nabla_\Gamma\psi_n,\nabla_\Gamma u)_\Gamma={\cal T}_2(\frak{g}\otimes\nabla_\Gamma\psi_n\otimes\nabla_\Gamma \overline{u})$. Hence, by \eqref{contractioncommutaestesa},
\begin{equation}\label{4.26ter}
D^i_\Gamma (\nabla_\Gamma\psi_n,\nabla_\Gamma u)_\Gamma={\cal T}_2[D^i_\Gamma(\frak{g}\otimes\nabla_\Gamma\psi_n\otimes\nabla_\Gamma \overline{u})]\qquad
\text{for $i=0,\ldots,m-1$.}
\end{equation}
Moreover, since for any $v\in T_i(\Gamma)$ we have $|{\cal T}_2 v|\le (N-1)|v|$, by \eqref{2.14} we also have $|{\cal T}_2v|_\Gamma^2\le (N-1)[1+(N-1)M_1^2]^i|v|_\Gamma^2$. Hence, by \eqref{4.26ter},
\begin{equation}\label{4.26quater}
| D^i_\Gamma(\nabla_\Gamma\psi_n,\nabla_\Gamma u)_\Gamma|_\Gamma^2\le   c_{32}|D^i_\Gamma(\frak{g}\otimes\nabla_\Gamma\psi_n\otimes\nabla_\Gamma \overline{u})|_\Gamma^2.
\end{equation}
where $c_{32}=c_{32}(m)=(N-1)[1+(N-1)M_1^2]^{m-1}$.
By \eqref{2.18} and \eqref{2.21} one gets, with  a trivial calculation, that $D_\Gamma{\frak g}=0$, which is well-known in the smooth case. Moreover, by \eqref{2.13}, one has $|\frak{g}|_\Gamma=1$, so by \eqref{modulotensoriale} and \eqref{4.26quater} we have
\begin{equation}\label{4.26penta}
| D^i_\Gamma(\nabla_\Gamma\psi_n,\nabla_\Gamma u)_\Gamma|_\Gamma^2\le  c_{32}|D^i_\Gamma(\nabla_\Gamma\psi_n\otimes\nabla_\Gamma \overline{u})|_\Gamma^2.
\end{equation}
Then combining \eqref{4.26bis} and \eqref{4.26quater} and using \eqref{2.22estesa}, assumptions (A0)--(ii) and \eqref{modulotensoriale} we get as before
\begin{multline*}
\sum_{n\in {\cal I}_{\Gamma_1}}\|(\nabla_\Gamma\psi_n,\nabla_\Gamma u)_\Gamma\|_{m-1,2,\Gamma_1}^2\le c_{32}\int_{\Gamma_1}\sum_{i=0}^{m-1}\sum_{n\in {\cal I}_{\Gamma_1}}|D_\Gamma^i (\nabla_\Gamma\psi_n,\nabla_\Gamma \overline{u})_\Gamma|_\Gamma^2\\
\le c_{32}N_0
\int_{\Gamma_1}\sum_{i=0}^{m-1}\sum_{j=0}^i\binom ij\sum_{n\in {\cal I}_{\Gamma_1}} |D^j_\Gamma(\nabla_\Gamma \psi_n)|_\Gamma^2 |D^{i-j}_\Gamma (\nabla_\Gamma u)|_\Gamma^2.
\end{multline*}
Since $m\le r-1$, by assumption (A0)--(ii) and \eqref{3.23bis-} we have
$$\sum_{n\in {\cal I}_{\Gamma_1}} |D^j_\Gamma(\nabla_\Gamma \psi_n)|_\Gamma^2\le N_0 c_{18}^2,$$
so there is $c_{33}=c_{33}(m)>0$ such that
\begin{equation}\label{4.26seis}\sum_{n\in {\cal I}_{\Gamma_1}}\|(\nabla_\Gamma\psi_n,\nabla_\Gamma u)_\Gamma\|_{m-1,2,\Gamma_1}^2\le
c_{33}\int_{\Gamma_1}\sum_{i=0}^{m-1}|D^i_\Gamma (\nabla_\Gamma u)|_\Gamma^2.
\end{equation}
By combining \eqref{4.26seis} and Lemma~\ref{lemma3.4} hence there is $c_{34}=c_{34}(m)>0$ such that
\begin{equation}\label{4.27}\sum_{n\in {\cal I}_{\Gamma_1}}\|(\nabla_\Gamma\psi_n,\nabla_\Gamma u)_\Gamma\|_{m-1,2,\Gamma_1}^2\le
c_{34}\|u\|_{m,2,\Gamma_1}^2.
\end{equation}
By combining \eqref{4.25}, \eqref{4.26}, \eqref{4.27} and using Lemma~\ref{lemma3.1} and Proposition~\ref{proposition3.1} then
there is $c_{35}=c_{35}(m)>0$ such that
$$\sum_{n\in {\cal I}_{\Gamma_1}}\|\omega_n\|_{m-1,2,\Gamma_1}^2\le
c_{35}\left(\|u\|_{m,2,\Gamma_1}^2+\|h\|_{m-1,2,\Gamma_1}^2\right).
$$
and then, by Lemma~\ref{lemma3.1}, as $\text{supp\,} \omega_n\cdot h_n\subset Q_0$,  we have
\begin{equation}\label{4.27bis}
\sum_{n\in {\cal I}_{\Gamma_1}}\|\omega_n\cdot h_n\|_{m-1,2,\R^{N-1}}^2\le c_6^2\,c_{35}\left(\|u\|_{m,2,\Gamma_1}^2+\|h\|_{m-1,2,\Gamma_1}^2\right).
\end{equation}

Now we note that, by \eqref{LaplaceBeltrami}, \eqref{4.24} and Lemma~\ref{lemma3.1}, $v_n:=(\psi_nu)\cdot h_n\in H^m(Q_0)$ is a weak solution of the elliptic equation in the form of the divergence
\begin{equation}\label{4.28bis}
-\partial_i(g_n^{1/2}g_n^{ij}\partial_j v_n)=g_n^{1/2}\omega_n\cdot h_n\qquad\text{in $Q_0$,}
\end{equation}
 and that, by \eqref{3.21}, $\text{supp\,}v_n\subset\subset Q_0$.
Hence, since $g_n^{1/2}g_n^{ij}\in C^m(\R^{N-1})$, the operator  $\partial_i(g_n^{1/2}g_n^{ij}\partial_j )$ is uniformly elliptic in $\R^{N-1}$ by \eqref{2.5}--\eqref{2.6} and  $g_n^{1/2}\omega_n\cdot h_n\in H^{m-1}(Q_0)$, by standard elliptic regularity theory (see \cite[Theorem~1, p.~327]{Evans}) we have $v_n\in H^{m+1}(Q_0)$.
However, since the estimates in \cite{Evans} are not uniform with respect to $n$,  we also note that, since $m+1\ge 2$ and $\text{supp\,}v_n\subset\subset Q_0$, by \eqref{4.28bis} we have $-a^{ij}_n\partial_{ij}v_n-b_n^j\partial_jv_n=\omega_n\cdot h_n$ in $L^2(\R^{N-1})$, where $a_n^{ij}=g_n^{ij}$ and $b_n^j=g_n^{-1/2}\partial_i g_n^{1/2}g_n^{ij}+\partial_ig_n^{ij}$. We can thus apply \cite[Theorem~5, p.~28]{Krylov}, which gives in particular
(adding a $\widetilde{\phantom{.}}$ to the notation of the author) the estimate
\begin{equation}\label{4.29}
\|v_n\|_{m+1,2,\R^{N-1}}\le \widetilde{N_0}\|\omega_n\cdot h_n-\widetilde{\lambda_0}v_n\|_{m-1,2,\R^{N-1}},
\end{equation}
where $\widetilde{N_0}, \widetilde{\lambda_0}>0$ depend only on $\|a_n^{ij}\|_{C^{m-1}_b(\R^{N-1})}$, $\|b_n^j\|_{C^{m-1}_b(\R^{N-1})}$, $m$, the constant of uniform ellipticity of the operator $-a^{ij}_n\partial_{ij}$ and the Lipschitz constants of the functions  $a^{ij}_n$.
Hence, by \eqref{2.5}--\eqref{2.7} and \eqref{2.9}, for our problem $\widetilde{N_0}$ and
$\widetilde{\lambda_0}$ only depend on $m$, $M_n$ and $N$. Hence, by combining \eqref{4.27bis}, \eqref{4.29} and Proposition~\ref{proposition3.1} we get
$\sum_{n\in {\cal I}_{\Gamma_1}}\|v_n\|_{m+1,2,\R^{N-1}}^2<\infty$, which again by Proposition~\ref{proposition3.1} gives that $u\in H^{m+1}(\Gamma_1)$, concluding the induction argument.

When $s=-m$, $m\in\N_0$, $0\le m\le r-1$ the statement follows by transposition, while when $s\in [1-r,r-1]\setminus\Z$ it does by interpolation, concluding the proof.
\end{proof}

The second preliminary regularity result extends a classical  regularity estimate (see for example \cite{lionsmagenes1}) for elliptic problems with nonhomogeneous Neumann boundary conditions to the case of non-compact boundaries satisfying assumptions (A0).
In particular we shall deal with weak solutions of  the classical problem
\begin{equation}\label{4.30}
\begin{cases}
-\Delta u+u=f\qquad &\text{in
$\Omega$,}\\
\partial_\nu u=\gamma\qquad
&\text{on
$\Gamma$,}\\
\end{cases}
\end{equation}
where $f\in L^2(\Omega)$ and $\gamma\in L^2(\Gamma)$, that is $u\in H^1(\Omega)$ such that
\begin{equation}\label{4.31}
  \int_\Omega\nabla u\nabla \varphi+\int_\Omega u\varphi=\int_\Omega f\varphi+\int_\Gamma \gamma \varphi\qquad\text{for all $\varphi\in H^1(\Omega)$.}
\end{equation}
\begin{thm}\label{lemma4.5BISS} Let assumption (A0) holds and $f\in H^{s-2}(\Omega)$,  $\gamma\in H^{s-3/2}(\Gamma)$ with $s\in\R$, $2\le s\le r$.
Then the unique weak solution $u$ of \eqref{4.31} belongs to $H^s(\Omega)$. Moreover there is $c_{36}=c_{36}(s)>0$ such that
\begin{equation}\label{4.32}
\|u\|_{s,2}^2\le c_{36}\left(\|f\|_{s-2,2}^2+\|\gamma\|_{s-3/2,2,\Gamma}^2\right)
\end{equation}
for all $f\in H^{s-2}(\Omega)$ and  $\gamma\in H^{s-3/2}(\Gamma)$.
\end{thm}

Before giving the proof of Theorem~\ref{lemma4.5BISS} it is convenient to point out that classical estimates on the half--space $\R^N_+=\{(x',x_N)\in \R^N: x_N>0\}$ for a family of auxiliary problems depending on $n\in {\cal I}_{\Gamma}$ are uniform with respect to $n$. We then introduce the bijective maps $H_n\in C^r(\R^N;\R^N)$ and their inverses $K_n$ given by
\begin{equation}\label{4.37}
H_n(x',x_N)=T_n(x', x_N+f_n(x')),\quad K_n(T_n(y',y_N))=(y',y_N-f_n(y')),
\end{equation}
where $T_n$ and $f_n$ are given in Lemma \ref{lemma2.1}, and we set for all $x\in\R^N$
\begin{equation}\label{4.39}
{\cal A}_n(x)= (a_n^{ij}(x))_{ij}=[JK_n]^t [JK_n]\cdot H_n(x).
\end{equation}
We point out, for easy reference, that denoting $Q=Q_0\times (-\tau_0,\tau_0)$ and $Q_+=Q_0\times(0,\tau_0)$, by Lemma~\ref{lemma2.1} and direct computation, since $O_n$ is orthogonal,
\begin{gather}
\label{maps}H_n(Q)=\Omega_n,\qquad H_n(Q_+)=\Omega\cap\Omega_n,\qquad H_n(Q_0\times\{0\})=U_n,\\
\label{Jacobians}
JH_n=O_n\,\left(\begin{array}{cc}I & 0\\\nabla f_n^t & 1\\\end{array}\right), \qquad
JK_n=\left(\begin{array}{cc}I & 0\\-\nabla f_n^t & 1\\\end{array}\right)O_n^t,\\
\label{determinants}|\text{det}\,JH_n|=|\text{det}\,JK_n|=1\quad \text{in $\R^N$.}
\end{gather}
Moreover, by assumption (A0) and \eqref{Jacobians}, for all $m\in\N$, $m\le r$, there is $c_{37}=c_{37}(m)>0$ such that
\begin{equation}\label{Jacest}
|D^\alpha H_n|,\qquad  |D^\alpha K_n|\le c_{37}\qquad\text{for all $|\alpha|\le m$ and $m\in{\cal I}_\Gamma$,  in $\R^N$.}
\end{equation}
Consequently, by \eqref{4.39}, for all $m\in\N$, $m\le r-1$, there is $c_{38}=c_{38}(m)>0$ such that
\begin{equation}\label{4.40}
|D^\alpha a^{ij}_n|\le c_{38}\qquad\text{in $\R^N$ for all $i,j=1,\ldots,N$, $|\alpha|\le m-1$ and $n\in{\cal I}_\Gamma$.}
\end{equation}
We shall consider weak solutions of the elliptic  problem
\begin{equation}\label{4.30bis}
\begin{cases}
-\text{div} ({\cal A}_n\nabla u)=\mathfrak{f}\qquad &\text{in
$\R^N_+$,}\\
-{\cal A}_n\nabla u\cdot e_N =\mathfrak{g}\qquad
&\text{on
$\partial\R^N_+$,}\\
\end{cases}
\end{equation}
where $e_N=(0,\ldots,0,1)$, $\mathfrak{f}\in L^2(\R^N_+)$ and $\mathfrak{g}\in L^2(\partial\R^N_+)$, i.e., $u\in H^1(\R^N_+)$ such that
\begin{equation}\label{4.42}
\int_{\R^N_+}\nabla u^t {\cal A}_n \nabla\psi=\int_{\R^N_+}\mathfrak{f}\psi+\int_{\partial\R^N_+}\mathfrak{g}\psi\quad\text{for all $\psi\in  H^1(\R^N_+)$.}
\end{equation}
In the sequel we shall identify, when useful $\partial\R^N_+$ with $\R^{N-1}$ by the trivial chart.

\begin{lem}\label{lemma4.5bis}
Let assumption (A0) holds and $\mathfrak{f}\in H^{m-2}(\R^N_+)$,  $\mathfrak{g}\in H^{m-3/2}(\partial\R^N_+)$, $m\in\N$, $2\le m\le r$.
Then any weak solution $u$ of \eqref{4.30bis} belongs to $H^m(\R^N_+)$ and  there is $c_{39}=c_{39}(m)>0$ such that
\begin{equation}\label{4.43}
\|u\|_{m,2}^2\le c_{39}\left(\|\mathfrak{f}\|_{m-2,2,\R^N_+}^2+\|\mathfrak{g}\|_{m-3/2,2,\partial\R^N_+}^2+\|u\|_{1,2,\R^N_+}^2\right)
\end{equation}
for all $\mathfrak{f}\in H^{m-2}(\R^N_+)$, $\mathfrak{g}\in H^{m-3/2}(\partial\R^N_+)$ and $n\in {\cal I}_\Gamma$.
\end{lem}
\begin{proof} The statement trivially reduces to the real case. By \eqref{Jacobians} and  \eqref{2.1} we have
\begin{alignat*}2
|dH_n[\xi]|^2&=|\xi'|^2+|\xi_N+\nabla f_n\xi'|^2&&\le (1+2|\nabla f_n|^2)|\xi'|^2+2|\xi_N|^2\\
&\le 2(1+|\nabla f_n|^2)|\xi|^2&&\le 2NM_1|\xi|^2,
\end{alignat*}
in $\R^N$, for all $\xi=(\xi',\xi_N)\in\R^N$. Consequently
$$|JK_n\xi|^2=|dK_n[\xi]|^2 \ge (2NM_1)^{-1}|\xi|^2.$$
By \eqref{4.39} it follows that for all $\xi\in\R^N$ we have
\begin{equation}\label{4.41}
\xi^t{\cal A}_n\xi=\left|[JK_n\cdot H_n] \xi\right|^2\ge(2NM_1)^{-1}|\xi|^2\quad\text{in $\R^N$.}
\end{equation}
Since \eqref{4.41} trivially gives that $a_n^{NN}\ge (2NM_1)^{-1}$ in $\R^N$ for all $n\in{\cal I}_\Gamma$, by \eqref{4.40}
we also get that for all $m\in\N$, $m\le r$, there is $c_{40}=c_{40}(m)>0$ such that
\begin{equation}\label{4.41bis}
|D^\alpha (1/a^{NN}_n)|\le c_{40}\qquad\text{in $\R^N$, for all $|\alpha|\le m-1$ and $n\in{\cal I}_\Gamma$.}
\end{equation}
At first we claim that, for each $m$,  the statement reduces to the case $\mathfrak{g}=0$. To prove our claim we introduce the operator $\Tr_m$, already defined in \S\ref{subsection3.3.2} for $\Omega$, in the particular case when $\Omega=\R^N_+$ and $\tau=2$, which in the sequel will be denoted as $\Tr^+_m$.
Hence $\Tr_m^+\in {\cal L}\left(H^m(\R^N_+),\prod_{i=0}^{m-1}H^{m-i-1/2}(\partial\R^N_+)\right)$ and trivially it is given by $\Tr^+_mw=(w(\cdot,0), -\partial_Nw(\cdot,0),\ldots,(-1)^{m-1}\partial_N^{m-1}w(\cdot,0))$ for any $w\in H^m(\R^N_+)$.
It is well-known (see \cite[p.215]{adams} or \cite[Theorem 7.5, p.~38]{lionsmagenes1}) that
$\Tr_m^+$ is surjective and thus, when restricted to $[{\rm Ker}\Tr_m^*]^\bot$, is bijective. Hence, by the Closed Graph Theorem, $\Tr_m^+$ possesses a right inverse $R_m^+\in {\cal L}\left(\prod_{i=0}^{m-1}H^{m-i-1/2}(\partial\R^N_+),H^m(\R^N_+)\right)$, i.e., $R_m^+\cdot \Tr_m^+=I$.
Hence, setting $\cal{R}_{m,n}^+\mathfrak{g}= R_m^+(0,\mathfrak{g}/a_n^{NN},0,\ldots,0)$ for all $\mathfrak{g}\in H^{m-3/2}(\partial\R^N_+)$, since $H^{m-3/2}(\R^{N-1})=[H^{m-1}(\R^{N-1}),H^{m-2}(\R^{N-1})]_{1/2}$, by interpolation we have
$$\cal{R}_{m,n}^+\in {\cal L}\left(H^{m-3/2}(\partial\R^N_+),H^m(\R^N_+)\right)$$ and, by \eqref{4.41bis} and the exactness of the complex interpolator functor, there is $c_{41}=c_{41}(m)>0$ such that
\begin{equation}\label{Add1}
 \|\cal{R}_{m,n}^+\|_{{\cal L}\left(H^{m-3/2}(\partial\R^N_+),H^m(\R^N_+)\right)}\le c_{41}\quad\text{for all $n\in\cal{I}_\Gamma$.}
\end{equation}
Moreover, denoting $v=\cal{R}_{m,n}^+\mathfrak{g}$, we have $v(\cdot,0)=0$,  $-a_n^{NN}\partial_N v(\cdot,0)=\mathfrak{g}$ and
consequently
$$
-{\cal A}_n\nabla v\cdot e_N =-\sum_{j=1}^{N-1}a_n^{Nj}\partial_jv(\cdot,0)-a_n^{NN}\partial_Nv(\cdot,0)=\mathfrak{g}
\qquad\text{on $\partial\R^N_+$.}
$$
Hence, multiplying by any $\psi\in  H^1(\R^N_+)$ and integrating by parts we get
\begin{equation}\label{Add4}
\int_{\R^N_+}\nabla v^t {\cal A}_n \nabla\psi=\int_{\R^N_+}-{\rm div}({\cal A}_n \nabla v)\psi +\int_{\partial\R^N_+}\mathfrak{g}\psi\quad\text{for all $\psi\in  H^1(\R^N_+)$.}
\end{equation}
By \eqref{4.40} and \eqref{Add1} there is $c_{42}=c_{42}(m)>0$ such that
\begin{equation}\label{Add5}
\|{\rm div}({\cal A}_n \nabla v)\|_{m-2,2,\R^N_+}\le c_{42}\|\mathfrak{g}\|_{m-3/2,2,\partial\R^N_+}.
\end{equation}
Now, to prove our claim we suppose that it holds for a fixed $m$ when $\mathfrak{g}=0$. By \eqref{4.42} and \eqref{Add4} $u-v$ is a weak solution of problem \eqref{4.30bis} provided $\mathfrak{f}$ and  $\mathfrak{g}$
are respectively replaced by $\mathfrak{f}-{\rm div}({\cal A}_n \nabla v)$ and zero. Hence $u-v\in H^m(\R^N_+)$ and, by combining \eqref{4.43} for $u-v$, \eqref{Add1} and \eqref{Add5}, we complete the proof for the $m$ fixed before.

We now prove the statement by induction on $m$. By previous claim in the initial step $m=2$ and in the recursive one we can take $\mathfrak{g}=0$.

 When $m=2$, by using \eqref{4.41}--\eqref{4.41bis}, we simply repeat  \cite[Proof of Theorem~9.25, Case C$_2$, pp.~305--306]{brezis2}, based on difference quotients,  with two differences: we do not need to chose $\psi\in H^1_0(\Omega)$ and, in estimate (66),  we keep $\|w\|_{1,2,\R^N_+}$ since Poincaré inequality is not available. In the present case, by \eqref{4.41}, we can take the positive constant $\alpha$ appearing there at page 305  as  $(2NM_1)^{-1}$.
Moreover  all other constants denoted by $C$ in the quoted proof only depend on upper estimates of $|\nabla a^{ij}_n|$, so in our case, by \eqref{4.40}, they only depend on  $M_1$ and $N$, so completing the proof when $m=2$.

We now suppose, by induction, that the statement holds true for $m\ge 2$, $m+1\le r$. We then take
$u\in H^m(\R^N_+)$, $\mathfrak{f}\in H^{m-1}(\R^N_+)$ and $\mathfrak{g}=0$. We fix any $i=1,\ldots,N-1$.
Then, for any $\psi\in H^2(\R^N_+)$ we can take $\partial_i\psi$ as a test function in \eqref{4.42}. By standard calculations and integration by parts we get that $\partial_i u$ satisfies the distribution identity \eqref{4.42} for any $\psi\in H^2(\R^N_+)$ provided $\mathfrak{f}$ and $\mathfrak{g}$ are respectively replaced by
\begin{equation}\label{Add7}
  \mathfrak{f}_i={\rm div}(\partial_i {\cal A}_n\nabla u)+\partial_i \mathfrak{f},\quad
\text{and}\quad \mathfrak{g}_i=\partial_i{\cal A}_n\nabla u \cdot e_N.
\end{equation}
Then we get by density that $\partial_i u$ is a weak solution of \eqref{4.30bis} provided previous replacement is performed. By \eqref{4.40} and \eqref{Add7} there is $c_{43}=c_{43}(m)>0$ such that
\begin{equation}\label{Add8}
\|\mathfrak{f}_i\|_{m-2,2,\R^N_+}^2+\|\mathfrak{g}_i\|_{m-3/2,2,\partial\R^N_+}^2\le c_{43}\left(\|u\|_{m,2,\R^N_+}^2+\|\mathfrak{f}\|_{m-1,2,\R^N_+}^2\right).
\end{equation}
Hence, by the induction hypothesis, combining \eqref{4.43} and \eqref{Add8},
there is $c_{44}=c_{44}(m)>0$ such that for all $i=1,\ldots,N_1$,
\begin{equation}\label{Add9}
\|\partial_i u\|_{m,2,\R^N_+}^2\le c_{44}\left(\|u\|_{1,2,\R^N_+}^2+\|\mathfrak{f}\|_{m-1,2,\R^N_+}^2\right).
\end{equation}
Since
\begin{equation}\label{Add2}
\|u\|_{m+1,2,\R^N_+}^2=\|u\|_{1,2,\R^N_+}^2+2\sum_{i=1}^{N-1}\|\partial_i u\|_{m,2,\R^N_+}^2+\|\partial^2_{NN}u\|_{m-1,2,\R^N_+}^2,
\end{equation}
to complete the proof we only have to estimate the last term in the right--hand side of \eqref{Add2}. Since $u\in H^m(\R^N_+)$ and $m\ge 2$, integrating by parts in \eqref{4.42} we get, as $\mathfrak{g}=0$, that $-{\rm div}(\cal{A}_n\nabla u)=\mathfrak{f}$ a.e.\ in $\R^N_+$, and consequently
\begin{equation}\label{Add10}
  \partial^2_{NN}u=-\frac 1{a_n^{NN}} {\sum_{i,j=1}^N}_{\phantom{}_{\phantom{}_{\negqquad \negqquad \negquad (i,j)\not=(N,N)}}} a_n^{ij}\partial^2_{ij} u
  -\sum_{i,j=1}^N \partial_ia_n^{ij}\partial_j u -\mathfrak{f}\quad\text{a.e.\ in $\R^N_+$}.
\end{equation}
By combining \eqref{Add9}--\eqref{Add10} with \eqref{4.40} and using the induction hypothesis again we then get the existence of $c_{45}=c_{45}(m)>0$ such that
$$\|u\|^2_{m+1,2,\R^N_+}\le c_{45}\left(\|u\|_{1,2,\R^N_+}^2+\|\mathfrak{f}\|_{m-1,2,\R^N_+}^2\right)$$
completing the proof.
\end{proof}

\begin{proof}[Proof of Theorem~\ref{lemma4.5BISS}] The proof is organized in several steps.

\noindent{\bf Step 1.} {\em{Preliminaries.}}
We shall use the sequence $(\phi_n)_{n\in{\cal I}_\Gamma}$ introduced in \eqref{3.12}, supplemented  with
\begin{equation}\label{4.33}
\phi_0=\rho_{\eps_0/8}*\chi_{_{\Omega_0}},\qquad\text{where $\Omega_0:=\{x\in\Omega: d(x,\Gamma)>\eps_0/4\}$.}
\end{equation}
We shall denote ${\cal I}_\Gamma^0={\cal I}_\Gamma\sqcup\{0\}$.
Denoting $\Omega_{\eps_0/2}=\{x\in \R^N: d(x,\overline{\Omega})<\eps_0/2\}$ and recalling that
$\Gamma^{\eps_0/2}=\{x\in\R^N: d(x,\Gamma)<\eps_0/2\}$, we have
\begin{equation}\label{4.34}
\phi_0\in C^\infty(\R^N),\quad \text{supp}\,\phi_0\subset \Omega, \quad 0\le \phi_n\le 1\,\text{in $\R^N$},\quad
 \phi_0=1\,\text{in $\Omega\setminus \Gamma^{\eps_0/2}$}.
 \end{equation}
 Since by assumption (A0)--(i) we have $\Gamma^{\eps_0/2}\subseteq\bigcup_{n\in{\cal I}_\Gamma}\Omega_n^{\eps_0/2}$, by \eqref{3.15} and \eqref{4.34} we get $\sum_{n\in {\cal I}_\Gamma^0}\phi_n\ge 1$, so by combining \eqref{3.14} and \eqref{4.34} for all $m\in\N$ there is $c_{46}=c_{46}(m)>0$ such that
$$\left|D^\alpha \left( \textstyle\sum_{n\in{\cal I}_\Gamma^0}\phi_n\right)^{-1}\right|\le c_{46}\quad\text{in $\Omega_{\eps_0/2}$ for all $|\alpha|\le m$.}$$
Consequently, by setting $\phi_n^\natural=\phi_n/\sum_{n\in{\cal I}_\Gamma^0}\phi_n$ for $n\in {\cal I}_\Gamma^0$ in  $\Omega_{\eps_0/2}$, by \eqref{3.14}, \eqref{4.33} and \eqref{4.34} we have
\begin{equation}\label{4.35}
 \begin{gathered}
 \phi_n^\natural\in C^\infty(\R^N),\quad \text{supp}\,\phi_n^\natural\subseteq \Omega_n^{\eps_0/4}\,\text{and}\quad {\phi_n^\natural}_{|\Gamma}=\psi_n'\quad\text{when $n\in{\cal I}_\Gamma$},\\
 \text{supp}\,\phi_0^\natural\subset \Omega, \qquad \sum\nolimits_{n\in{\cal I}_\Gamma^0}\phi_n^\natural=1\quad \text{in $\Omega_{\eps_0/2}$,}
 \end{gathered}
 \end{equation}
and for all $m\in\N$ there is $c_{47}=c_{47}(m)>0$ such that
\begin{equation}\label{4.35bis}
|D^\alpha \phi_n^\natural|\le c_{47}\,\,\text{in $\Omega_{\eps_0/2}$  for all $|\alpha|\le m$ and $n\in {\cal I}_\Gamma^0$.}
\end{equation}
We also remark that by \eqref{2.1bis} and the identity $g_n=\sqrt{1+|\nabla f_n|^2}$ we have $\nu\cdot h_n=g_n^{-1/2}O_n(\nabla f_n,-1)$. Consequently, as $O_n$ is orthogonal, by \eqref{2.1} and \eqref{2.8} for $m\in\N$, $m\le r$, there is $c_{48}=c_{48}(m)>0$ such that
\begin{equation}\label{4.35ter}
|D^\alpha(\nu\cdot h_n)|\le c_{48}\quad\text{in $Q_0$ \,for all $|\alpha|\le m-1$ and $n\in \cal{I}_\Gamma$.}
\end{equation}
Denoting $\nu=(\nu_1,\ldots,\nu_N)$, starting from \eqref{4.35bis}--\eqref{4.35ter} and using the same arguments which leaded to \eqref{3.23}, we get for $m\in\N$, $m\le r$, the existence of
$c_{49}=c_{49}(m)>0$ such that
$$\big|D_\Gamma^j\big(\partial_i {\phi_n^\natural}_{|\Gamma}\big)\big|_{\Gamma},\, |D_\Gamma^j \nu_i|_\Gamma\le c_{49}\,\text{on $\Gamma$, for all $i=1,\ldots,N$, $j=0,\ldots,m-1$, $n\in \cal{I}_\Gamma$.}
$$
Consequently, since ${\partial_\nu\phi_n^\natural}_{|\Gamma}=\sum_{i=1}^N \partial_i {\phi_n^\natural}_{|\Gamma}\nu_i$, using Leibniz formula \eqref{2.22},  for $m\in\N$, $m\le r$, there is
$c_{50}=c_{50}(m)>0$ such that
\begin{equation}\label{4.35quater}
|D^i_\Gamma({\partial_\nu\phi_n^\natural}_{|\Gamma})|_\Gamma\le c_{50}\,\,\text{on $\Gamma$,\, for all $i=0,\ldots,m-1$ and $n\in \cal{I}_\Gamma$.}
\end{equation}
\noindent{\bf Step 2.} {\em{The main claim.}} We shall prove at first that when
$f\in H^{m-2}(\Omega)$ and $\gamma\in H^{m-3/2}(\Gamma)$ , with $m\in\N$, $2\le
m\le r$, then $u\in H^m(\Omega)$. Clearly this assertion can be extended to $m\in\N$, $1\le m\le r$ by respectively replacing $H^{m-2}(\Omega)$ and  $H^{m-3/2}(\Gamma)$ with $H^{\max\{m-2,0\}}(\Omega)$ and $H^{\max\{m-3/2,0\}}(\Gamma)$.
We shall prove this fact by induction on $m$, and clearly when $m=1$ there is nothing to prove. Hence we shall now suppose by induction that $m\ge 2$, $f\in H^{m-2}(\Omega)$,  $\gamma\in H^{m-3/2}(\Gamma)$ and
$u\in H^{m-1}(\Omega)$, and we claim that $u\in H^m(\Omega)$.
By \eqref{4.35} we have $u=\phi_0^\natural u+u_1$, where $u_1=\sum_{n\in {\cal I}_\Gamma}\phi_n^\natural u$. Hence our claim reduces to prove that $\phi_0^\natural u\in H^m(\Omega)$ and $u_1\in H^m(\Omega)$.
These claims require an inner and an outer estimate, respectively.

\noindent{\bf Step 3.} {\em{Inner estimate.}} By  \eqref{4.35bis}, given any $\psi\in H^1(\R^N)$, we can take the restriction to $\Omega$ of $\phi_0^\natural \psi$ as a test function $\varphi$  in \eqref{4.31}. Then, by \eqref{4.35}, a standard use of Leibniz formula and integration by parts gives
\begin{equation}\label{Add11}
\int_\Omega \nabla(\phi_0^\natural u)\nabla \psi+\int_\Omega \phi_0^\natural u\psi=\int_\Omega f_0^\natural \psi\qquad\text{for all $\psi\in H^1(\R^N)$,}
\end{equation}
where $f_0^\natural:=-\Delta \phi_0^\natural u-2\nabla \phi_0^\natural\nabla u+\phi_0^\natural f$.
By \eqref{4.35} and \eqref{4.35bis} the trivial extension to $\R^N$ of $\phi_0^\natural$ belongs to $H^{m-1}(\R^N)$, while those of $\Delta \phi_0^\natural u$, $\nabla \phi_0^\natural\nabla u$ and   $\phi_0^\natural f$  belong to $H^{m-2}(\R^N)$, and we can replace
$\Omega$ with $\R^N$ in  \eqref{Add11}. We can then directly apply regularity theory for the Dirichlet problem in $\R^N$, and in particular the trivial extension to the complex case of \cite[Theorem 9.26, Case A]{brezis2}. Consequently we get that $\phi_0^\natural u\in H^m(\Omega)$, concluding the inner estimate.

\noindent{\bf Step 4.} {\em{Uniform estimate of $\phi_n^\natural u$ with respect to $n\in{\cal I}_\Gamma$.}}
By \eqref{4.35bis}, for any $n\in{\cal I}_\Gamma$ and $\psi\in H^1(\Omega)$ we can take $\varphi=\phi_n^\natural\psi$ in \eqref{4.31}. Then by a standard use of Leibniz formula and integration by parts we get
\begin{equation}\label{Add12}
\int_\Omega \nabla(\phi_n^\natural u)\nabla \psi=\int_\Omega f_n^\natural \psi+\int_\Gamma \gamma_n\psi
\qquad\text{for all $\psi\in H^1(\Omega)$,}
\end{equation}
where
\begin{equation}\label{Add13}
f_n^\natural:=-\Delta \phi_n^\natural u-2\nabla \phi_n^\natural\nabla u-\phi_n^\natural u +\phi_n^\natural f,
\quad \text{and}\quad \gamma_n=\partial_\nu \phi_n^\natural u_{|\Gamma}+\psi'_n\gamma.
\end{equation}
Trivially, by \eqref{4.35bis}, $f_n^\natural\in H^{m-2}(\Omega)$ for all $n\in\cal{I}_\Gamma$. Moreover, by \eqref{4.35bis} and  Theorem~\ref{theorem3.1}  we have $\psi'_n\gamma\in H^{m-3/2}(\Gamma)$. Since $u\in H^{m-1}(\Omega)$, by the Trace Theorem $u_{|\Gamma}\in H^{m-3/2}(\Gamma)$. By \eqref{4.35bis} and Lemma~\ref{lemma3.5} we have $\partial_\nu \phi_n^\natural\in W^{m-1,\infty}(\Gamma_1)$, so by Lemma~\ref{Multiplierlemma} , we have
$\partial_\nu \phi_n^\natural u_{|\Gamma}\in H^{m-3/2}(\Gamma)$. Consequently $\gamma_n\in H^{m-3/2}(\Gamma)$. Next, by \eqref{4.35}, we can respectively replace $\Omega$ and $\Gamma$ with $\Omega\cap\Omega_n$ and $U_n$ in all integrals in \eqref{Add12}. Moreover we set for any $n\in\cal{I}_\Gamma$ the cut--off function $\eta_n=\rho_{\eps_0/32}*\chi_{\Omega_n^{\eps_0/8}}$, so that $\eta_n=1$ in $\Omega_n^{\eps_0/4}$ and $\text{supp}\,\eta_n\subset \Omega_n$. Hence, given any $\psi\in H^1(\Omega\cap\Omega_n)$ we have $\eta_n\psi=\psi$ in $\Omega_n^{\eps_0/4}$ and $\text{supp}\,\eta_n\psi\subset \Omega_n$. Hence, by the distributional characterization of the elements of $H^1$ (\cite[Proposition~9.3]{brezis2}), applying L.~Schwartz's {\em Principe du recollement des morceaux}
(\cite[Th\'{e}or\`{e}me~IV, p.~27]{schwartz}) to the open sets $\Omega\setminus \text{supp}\,\eta_n\psi$ and $\Omega\cap\Omega_n$, we get that the trivial extension of $\eta_n\psi$ belongs to $H^1(\Omega)$. Hence, putting it as a test function in \eqref{Add12} we get
\begin{equation}\label{Add14}
\int_{\Omega\cap\Omega_n} \nabla(\phi_n^\natural u)\nabla \psi=\int_{\Omega\cap\Omega_n} f_n^\natural \psi+\int_{U_n} \gamma_n\psi
\qquad\text{for all $\psi\in H^1(\Omega\cap\Omega_n)$.}
\end{equation}
Now we remark that by combining \eqref{maps} and \eqref{Jacest} with the classical change of variable results for Sobolev spaces (\cite[Theorem 3.35, p.~63]{adams} or \cite[Theorem~11.57]{LeoniSobolev2})
we get that for any $i\in\N_0$, $i\le r$, the operators $v\mapsto v\cdot H_n$ map $H^i(Q^+)$ onto $H^i(\Omega\cap\Omega_n)$, with inverse $w\mapsto w\cdot K_n$, and there is $c_{51}=c_{51}(i)>0$ such that
\begin{equation}\label{Add15}
c_{51}^{-1}\|v\|_{i,2,Q^+}\le \|v\cdot H_n\|_{i,2,\Omega\cap\Omega_n}\le c_{51}\|v\|_{i,2,Q^+}
\end{equation}
for all $v\in H^i(Q^+)$ and $n\in\cal{I}_\Gamma$.

Moreover, since (see \cite{adams})
$H^{j-1/2}(Q_0)=(H^{j-1}(Q_0),H^j(Q_0))_{1/2,2}$ and the real interpolator functor is exact (see \cite{bergh}), by \eqref{3.5} we get that for $j\le r$ the operators
$v\mapsto v\cdot h_n$ map $H^{j-1/2}(Q_0)$ onto $H^{j-1/2}(U_n)$, with inverse $w\mapsto w\cdot \xi_n$, and there is $c_{52}=c_{52}(j)>0$ such that
\begin{equation}\label{Add16}
c_{52}^{-1}\|v\|_{j-1/2,2,Q_0}\le \|v\cdot H_n\|_{j-1/2,2,U_n}\le c_{52}\|v\|_{j-1/2,2,Q_0}
\end{equation}
for all $v\in H^{j-1/2}(Q_0)$ and $n\in\cal{I}_\Gamma$.

Hence, changing variables in all integrals in \eqref{Add14}, using \eqref{determinants} and \eqref{Jacest} we get that $v_n:=(\phi_n^\natural u)\cdot H_n\in H^{m-1}(Q^+)$ and
\begin{equation}\label{Add17}
\int_{Q^+}\nabla v_n^t {\cal A}_n \nabla\psi=\int_{Q^+} \dot f_n\psi+\int_{Q_0}\dot\gamma_n\psi\quad\text{for all $\psi\in  H^1(Q^+)$,}
\end{equation}
where
\begin{equation}\label{Add18}
\dot f_n=f_n^\natural\cdot H_n\in H^{m-2}(Q^+),\quad \text{and}\quad \dot\gamma_n=g_n(\gamma_n\cdot h_n)\in H^{m-3/2}(Q_0).
\end{equation}
By \eqref{maps}, \eqref{4.35}, \eqref{Add13} and \eqref{Add18}, $\text{supp}\,v_n, \text{supp}\,\dot f_n\subset Q_0\times [0,\tau_0)$ and $\text{supp}\,\dot\gamma_n\subset Q_0$. Moreover, by using the universal extension operator built in \cite[Proof of Theorem 13.17, p.~424 and p.~560]{LeoniSobolev2},
$v_n$, $\dot f_n$ and $\dot\gamma_n$ respectively extend to elements of $H^{m-1}(\R^N_+)$, $H^{m-2}(\R^N_+)$
and $H^{m-3/2}(\R^{N-1})$. Hence,  using suitable cut-off functions, their trivial extensions (denoted with the same symbols) belong to the same spaces listed above, and $v_n$ is a weak solution on problem \eqref{4.42} with $\mathfrak{f}=\dot f_n$ and $\mathfrak{g}=\dot\gamma_n$. Hence, by applying Lemma~\ref{lemma4.5bis} we get that $v_n\in H^m(Q^+)$ and
$$
\|v_n\|_{m,2}^2\le c_{39}\left(\|\dot f_n\|_{m-2,2,Q^+}^2+\|\dot\gamma_n\|_{m-3/2,2,Q_0}^2+\|v_n\|_{1,2,Q^+}^2\right)\quad\text{for all $n\in \cal{I}_\Gamma$.}
$$
Consequently, by \eqref{2.7}, \eqref{4.35}, \eqref{Add13}, \eqref{Add15}--\eqref{Add16} and \eqref{Add18} there is $c_{53}=c_{53}(m)>0$ such that, for all $n\in \cal{I}_\Gamma$ we have
\begin{multline}\label{Add21}
\|\phi_n^\natural u\|_{m,2}^2\le c_{53}\Big(\|\Delta\phi_n^\natural u\|_{m-2,2}^2+\|\nabla\phi_n^\natural\nabla u\|_{m-2,2}^2+\|\phi_n^\natural u\|_{m-1,2}^2\\+\|\phi_n^\natural f\|_{m-1,2}^2+\|\partial_\nu \phi_n^\natural u_{|\Gamma}\|_{m-3/2,2,\Gamma}^2+\|\psi'_n\gamma\|_{m-3/2,2,\Gamma}^2\Big).
\end{multline}
\noindent{\bf Step 5.} {\em Boundary estimate.} Starting from \eqref{Add21} we are going to  prove that $u_1\in H^m(\Omega)$. We note that, being the summations locally finite by assumption (A0), we have
$$\|D^mu_1\|_2^2=\int_\Omega \Big|D^m\Big(\sum_{n\in\cal{I}_\Gamma}\phi_n^\natural u\Big)\Big|^2=
\int_\Omega \Big|\sum_{n\in\cal{I}_\Gamma} D^m\Big(\phi_n^\natural u\Big)\Big|^2
$$
and then, by assumption (A0)--(ii) and Beppo Levi's Theorem,
$$\|D^mu_1\|_2^2\le N_0\int_\Omega \sum_{n\in\cal{I}_\Gamma}|D^m(\phi_n^\natural u)|^2=
N_0\sum_{n\in\cal{I}_\Gamma}\|D^m(\phi_n^\natural u)\|_2^2\le N_0\sum_{n\in\cal{I}_\Gamma}\|\phi_n^\natural u\|_{m,2}^2.
$$
Consequently, by \eqref{Add21}, setting $c_{54}=c_{54}(m)=N_0 c_{53}>0$, we have
\begin{equation}\label{Add55}
\begin{aligned}
\|D^mu_1\|_2^2\le c_{54}\Bigg(&\sum_{n\in\cal{I}_\Gamma} \|\Delta\phi_n^\natural u\|_{m-2,2}^2+\sum_{n\in\cal{I}_\Gamma}\|\nabla\phi_n^\natural\nabla u\|_{m-2,2}^2\\+&\sum_{n\in\cal{I}_\Gamma}\|\phi_n^\natural u\|_{m-1,2}^2+\sum_{n\in\cal{I}_\Gamma}\|\phi_n^\natural f\|_{m-1,2}^2\\+&\sum_{n\in\cal{I}_\Gamma}\|\partial_\nu \phi_n^\natural u_{|\Gamma}\|_{m-3/2,2,\Gamma}^2+\sum_{n\in\cal{I}_\Gamma}\|\psi'_n\gamma\|_{m-3/2,2,\Gamma}^2\Bigg).
\end{aligned}
\end{equation}
We shall prove that all summations in the right-hand side of \eqref{Add55} are finite, so proving our claim. By Beppo Levi's Theorem and Leibniz formula we have
\begin{align*}
\sum_{n\in\cal{I}_\Gamma} \|\Delta\phi_n^\natural u\|_{m-2,2}^2=
&\sum_{|\alpha\le m-2}\int_\Omega \sum_{n\in\cal{I}_\Gamma}\Big| \sum_{\beta\le\alpha}{\binom\alpha\beta} D^\beta(\Delta \phi_n^\natural)D^{\alpha-\beta}u\Big|^2\\
\le & N^{m-2}\sum_{|\alpha\le m-2}\int_\Omega \sum_{\beta\le\alpha}\sum_{n\in\cal{I}_\Gamma}|D^\beta(\Delta \phi_n^\natural)|^2|D^{\alpha-\beta}u|^2.
\end{align*}
Consequently, by \eqref{4.35}, \eqref{4.35bis} and assumption (A0--ii),
\begin{equation}\label{Add22}
  \begin{aligned}
\sum_{n\in\cal{I}_\Gamma} \|\Delta\phi_n^\natural u\|_{m-2,2}^2\le
&c_{47}^2 N^{m-2} N_0 \sum_{|\alpha\le m-2}\sum_{\beta\le\alpha}
\int_\Omega |D^{\alpha-\beta}u|^2\\
\le & c_{47}^2 N^{2m}N_0\|u\|_{m-2,2}^2<\infty.
\end{aligned}
\end{equation}
By using exactly the same arguments we get that there is $c_{55}=c_{55}(m)>0$ such that
\begin{equation}\label{Add23}
\left\{
\begin{alignedat}2
&\sum_{n\in\cal{I}_\Gamma}\|\nabla\phi_n^\natural\nabla u\|_{m-2,2}^2\le &&c_{55} \|u\|_{m-1,2}^2<\infty,\\
&\sum_{n\in\cal{I}_\Gamma}\|\phi_n^\natural u\|_{m-1,2}^2\le &&c_{55}\|u\|_{m-1,2}^2<\infty,\\
&\sum_{n\in\cal{I}_\Gamma}\|\phi_n^\natural f\|_{m-1,2}^2\le&&c_{55}\|f\|_{m-2,2}^2<\infty.
\end{alignedat}\right.
\end{equation}
Moreover, by Theorem~\ref{theorem3.1}, there is $c_{56}=c_{56}(m)>0$ such that
\begin{equation}\label{Add24}
\sum_{n\in\cal{I}_\Gamma}\|\psi'_n\gamma\|_{m-3/2,2,\Gamma}^2\le c_{56}\|\gamma\|_{m-3/2,2,\Gamma}^2<\infty.
\end{equation}
To estimate the fifth summation in the right-hand side of \eqref{Add55} we first remark that, by using \eqref{2.22estesa}, \eqref{4.35}, \eqref{4.35bis}, assumption (A0--ii) and Beppo--Levi's Theorem, for any $j\in\N_0$, $j+1\le r$, and any $v\in H^j(\Gamma)$ we have
\begin{align*}
\sum_{n\in\cal{I}_\Gamma}\|\partial_\nu \phi_n^\natural v\|_{j,2,\Gamma}^2=
&\sum_{i=0}^j \int_\Gamma \sum_{n\in\cal{I}_\Gamma}|D_\Gamma^i(\partial_\nu \phi_n^\natural v)|_\Gamma^2\\
=&\sum_{i=0}^j \int_\Gamma \sum_{n\in\cal{I}_\Gamma} \Bigg|\sum_{l=0}^i \binom il D_\Gamma^l(\partial_\nu \phi_n^\natural)\otimes D_\Gamma^{i-l} v)\Bigg|_\Gamma^2\\
\le &j \sum_{i=0}^j \sum_{l=0}^i \binom il \int_\Gamma \sum_{n\in\cal{I}_\Gamma} |D_\Gamma^l(\partial_\nu \phi_n^\natural)|^2|D_\Gamma^{i-l} v|_\Gamma^2.
\end{align*}
Consequently, by \eqref{4.35quater} and assumption (A0-ii), there is $c_{57}=c_{57}(j)>0$ such that
\begin{equation}\label{Add25}
\sum_{n\in\cal{I}_\Gamma}\|\partial_\nu \phi_n^\natural v\|_{j,2,\Gamma}^2\le N_0c_{50}^2
j \sum_{i=0}^j \sum_{l=0}^i \binom il \int_\Gamma |D_\Gamma^{i-l} v|_\Gamma^2\le c_{57}\|v\|_{j,2,\Gamma}^2.
\end{equation}
Now we consider the space $\ell^2(X_n^s)$ for $s\in\R$, $0\le s\le r$, recalled in \eqref{elltau}, with
$X_n^s=H^s(\Gamma)$ when $n\in {\cal I}_{\Gamma}$ and $X_n^s=\{0\}$ otherwise. We set the operator $\Phi_0:L^2(\Gamma)\to \ell^2(X_n^0)$
by $(\Phi_0 v)_n=\partial_\nu\phi_n^\natural v$ for $n\in {\cal I}_{\Gamma}$. By  \eqref{Add25} we have $\Phi_0\in \cal{L}(H^{m-1}(\Gamma);\ell^2(X_n^{m-1}))\cap\cal{L}(H^{m-2}(\Gamma);\ell^2(X_n^{m-2}))$. Since,  by \eqref{3.27} and
\cite[Theorem p.121]{triebel} we have  $\ell^2(X_n^{m-3/2})=(\ell^2(X_n^{m-1}),\ell^2(X_n^{m-2}))_{1/2,2}$ with equivalence of norms, by interpolation we get
$\Phi_0\in \cal{L}(H^{m-3/2}(\Gamma);\ell^2(X_n^{m-3/2}))$. Consequently, by the Trace Theorem, there is $c_{58}=c_{58}(m)>0$ such that
$$\sum_{n\in\cal{I}_\Gamma}\|\partial_\nu \phi_n^\natural u\|_{m-3/2,2,\Gamma}^2\le c_{58}\|u\|_{m-1,2}^2.$$
The last estimate, together with \eqref{Add22}--\eqref{Add25}, concludes the proof of the main claim made in Step 2.

\noindent{\bf Step 6.} {\em Conclusion.} Considering for $m\in\N$, $2\le
m\le r$, the operator $\cal{S}_m u= (-\Delta u +u, \partial_\nu u)$, by the Trace Theorem $\cal{S}_m u\in \cal{L}(H^m(\Omega);H^{m-2}(\Omega)\times H^{m-3/2}(\Gamma))$. Since, by Trace and Riesz--Fischer's Theorems for all $f\in L^2(\Omega)$ and $\gamma\in L^2(\Gamma)$ problem \eqref{4.30} has exactly one weak solution, our main claim gives that  $\cal{S}_m$ is bijective, so by the Closed Graph Theorem, $\cal{S}_m^{-1}$ is bounded, so concluding the proof when  $s=m\in\N$. When $s\in\R\setminus\N$, $2\le s\le r$, we take $m=[s]+1$ and $\theta=s-m+1$, so $m-1<s<m$, $3\le m\le r$. Then, denoting by $[\cdot,\cdot]_\theta$ the  complex interpolator functor (\cite[pp.~86]{bergh} or \cite[p.~59]{triebel}), using the well-known interpolation properties of spaces $H^t(\R^l)$ for $t\ge 0$ and $l\in\N$ (see \cite[(11), p.~185]{triebel}), Theorem~\ref{theorem3.1} and the extension theorem \cite[Theorem~13.17]{LeoniSobolev2}, one easily gets that
$H^s(\Omega)=[H^{m-1}(\Omega), H^m(\Omega)]_\theta$ and $$H^{s-2}(\Omega)\times H^{s-3/2}(\Gamma) =[H^{m-3}(\Omega)\times H^{m-5/2}(\Gamma); H^{m-2}(\Omega)\times H^{m-3/2}(\Gamma)]_\theta.$$
The proof is then completed by interpolation.
\end{proof}

Theorems~\ref{lemma4.5} and \ref{lemma4.5BISS} allow to characterize the subspaces $D(A^n)$ in Theorem~\ref{Theorem4.1} when $r\ge 2$. In particular, denoting $D_{n-1}=\{(u_0,v_0,u_1,v_1)\in \cal{H}^n: \,\,\text{\eqref{1.6} holds}\}$ for $n\in\N$, $2\le n\le r$,
we get the following result
\begin{lem}\label{lemma4.6} If assumptions (A0--3) hold and $r\ge 2$ then $D(A^n)=D_n$ for all $1\le n<r$. Moreover the norms $\|\cdot\|_{D(A^n)}$ and $\|\cdot\|_{\cal{H}^{n+1}}$ are equivalent on it.
\end{lem}
\begin{proof} At  first we more explicitly  rewrite the sets $D_n$. By identifying (as in \eqref{eq:splitting}) $L^2(\Gamma)=L^2(\Gamma_0)\oplus L^2(\Gamma_1)$ and using \eqref{3.38} together with the Trace Theorem, we introduce $\cal{E}^n\in\cal{L}(\cal{H}^{n+1};L^2(\Gamma))$ for $1\le n<r$, setting for $U=(u,v,w,z)$,
\begin{equation}\label{L1}
\left\{
\begin{alignedat}2
&\cal{E}^1(U) &&=(\partial_\nu u_{|\Gamma_0},\partial_\nu u_{|\Gamma_1}-z)\\
&\cal{E}^2(U) &&=(\partial_\nu w_{|\Gamma_0},\mu \partial_\nu w_{|\Gamma_1}-\DivGamma(\sigma \nabla_\Gamma v)+\delta \partial_\nu u_{|\Gamma_1} +\kappa v+\rho w_{|\Gamma_1})\\
&\cal{E}^3(U) &&=(\partial_\nu \Delta u_{|\Gamma_0},c^2\mu \partial_\nu \Delta u_{|\Gamma_1}-\DivGamma(\sigma \nabla_\Gamma \partial_\nu u)\\
&&&\phantom{=(}+\delta \partial_\nu w_{|\Gamma_1} +\kappa\partial_\nu u_{|\Gamma_1}+c^2\rho \Delta u_{|\Gamma_1})\\
&\cal{E}^4(U) &&=(\partial_\nu \Delta w_{|\Gamma_0},c^2\mu \partial_\nu \Delta w_{|\Gamma_1}-\DivGamma(\sigma \nabla_\Gamma \partial_\nu w)\\
&&&\phantom{=(}+c^2\delta \partial_\nu \Delta u_{|\Gamma_1} +\kappa\partial_\nu w_{|\Gamma_1}+c^2\rho \Delta w_{|\Gamma_1})\\
&\cal{E}^{2i+3}(U) &&= \cal{E}^3(\Delta^iu,v,\Delta ^iw,z),\qquad\text{for $i\in\N$, $2i+3<r$},\\
&\cal{E}^{2i+4}(U) &&= \cal{E}^4(\Delta^iu,v,\Delta ^iw,z),\qquad\text{for $i\in\N$, $2i+4<r$}.\\
\end{alignedat}
\right.
\end{equation}
We also denote $\cal{E}^n=(\cal{E}^n_{\Gamma_0},\cal{E}^n_{\Gamma_1})$. Trivially
\begin{equation}\label{L2}
 D_n=\{U\in \cal{H}^{n+1}: \cal{E}^iU=0\quad\text{for $i=1,\ldots,n$}\},
\end{equation}
so $D_n$ is a closed subspace of $\cal{H}^{n+1}$. We shall endow $D_n$ with the norm inherited from $\cal{H}^{n+1}$.

We now point out that, since $1/\mu\in W^{r-1,\infty}(\Gamma_1)$, by assumption (A1), Lemmas~\ref{lemma3.5} and \ref{Multiplierlemma} and \eqref{3.38} we have
\begin{equation}\label{L3}
  A\in \cal{L}(D_n;\cal{H}^n)\qquad\text{for $1\le n<r$.}
\end{equation}
We also remark that, since $\partial_\nu u_{|\Gamma_1}=z$ in $D_1$, we have $\cal{E}^2=-(\cal{E}^1_{\Gamma_0},\mu\cal{E}^1_{\Gamma_1})\cdot A$ in $\cal{H}^3\cap D_1$, $\cal{E}^3=-(\frac 1{c^2}\cal{E}^2_{\Gamma_0},\cal{E}^2_{\Gamma_1})\cdot A$ in $\cal{H}^4\cap D_1$,
$\cal{E}^4=-\cal{E}^3\cdot A$ in $\cal{H}^5\cap D_1$ and $\cal{E}^{i+2}=\frac 1{c^2}\cal{E}^i\cdot A^2$ in $\cal{H}^{i+3}\cap D_1$
for $3\le i<r-2$. Consequently for any $1\le n<r$ and $U\in \cal{H}^{n+1}\cap D_1$ we have $\cal{E}^n U=0$ if and only if $\cal{E}^{n-1} AU=0$.
Hence \eqref{L2} can be rewritten as $D_n=\{U\in \cal{H}^{n+1}: \cal{E}^1A^iU=0\quad\text{for $i=1,\ldots,n-1$}\}$. The recursive formula
\begin{equation}\label{L4}
  D_{n+1}=\{U\in \cal{H}^{n+2}\cap D_n: AU\in D_n\}, \qquad\text{for  $1\le n<r-2$}
\end{equation}
then immediately follows.

Now we claim that $D(A^n)=D_n$ for any $1\le n<r$. We shall prove it by induction.

When $n=1$, by \eqref{4.4} and  \eqref{L2}--\eqref{L3} we have $D_1\subseteq D(A)$. To prove the reverse inclusion we take $U=(u,v,w,z)\in D(A)$. By \eqref{4.4}--\eqref{4.5} the quadruplet $(u,v,w,z)$ solves the system \eqref{4.7} with $\lambda=0$ and $(h_1,h_2,h_3,h_4)\in \cal{H}$. Consequently, as in the proof of Lemma~\ref{lemma4.2}, we get that $u$ is a weak solution of problem \eqref{4.30} with $f=u+c^{-2}h_3$ and $\gamma=h_2$ on $\Gamma_1$, $\gamma=0$ on $\Gamma_0$, and that $v$ satisfies the equation $B_0v=(1-\kappa)v-\delta h_2 -\rho {h_1}_{|\Gamma_1}+\mu h_4$, where $B_0$ is the operator in Theorem~\ref{lemma4.5}.
Hence, as $f\in L^2(\Omega)$, $\gamma\in H^1(\Gamma_1)$ and $B_0v\in L^2(\Gamma_1)$,  by applying Theorems~\ref{lemma4.5} and \ref{lemma4.5BISS} we get $u\in H^2(\Omega)$ and $v\in H^2(\Gamma_1)$, so $U\in \cal{H}^2$ and then, by the remarks made in \S~\ref{subsection3.3.2}, $U\in D_1$, proving our claim when $n=1$.

We now suppose by induction that $D(A^n)=D_n$ for $n<r-1$. By \eqref{L4} we immediately get that $D_{n+1}\subseteq D(A^{n+1})$.
To prove the reverse inclusion let $U=(u,v,w,z)\in D(A^{n+1})$, so by the induction hypothesis $U, AU\in D_n\subset \cal{H}^{n+1}$.
Hence, by  the same arguments used in the case $n=1$ we get that $u$ and $v$ are solutions of the same problems. But in this case, by using assumptions (A1--3) and Lemma~\ref{Multiplierlemma},  we get
$f\in H^n(\Omega)$, $\gamma\in H^{n+1}(\Gamma_1)$and $B_0v\in H^n(\Gamma_1)$, so a further application of Theorems~\ref{lemma4.5} and \ref{lemma4.5BISS} gives $u\in H^{n+2}(\Omega)$ and $v\in H^{n+2}(\Gamma_1)$, so $U\in \cal{H}^{n+2}$. Since $U, AU\in D_n$, by \eqref{L4} we then get $U\in D_{n+1}$, concluding the proof of our claim.

To prove the stated equivalence of norms we remark that \eqref{L3} yields the existence of $c_{59}=c_{59}(m, \sigma, \delta, \kappa, \mu,\rho,c)>0$ such that
$\|\cdot \|_{D(A^n)}\|\le c_{59}\|\cdot\|_{\cal{H}^{n+1}}$ on $D_n$. To get the reverse inequality we remark that, by Lemmas~\ref{lemma4.1} and \ref{lemma4.2}, the operator $A+\Lambda_1I\in\cal{L}(D_1,\cal{H})$ is bijective, so being $D_n=D(A^n)$ also the operator $(A+\Lambda_1I)^n\in\cal{L}(D_n,\cal{H})$ is bijective for $n<r-1$. Hence, being $D_n$ complete, by the Closed Graph Theorem there is $c_{60}=c_{60}(m, \sigma, \delta, \kappa, \mu,\rho,c)>0$ such that
$\|\cdot \|_{\cal{H}^{n+1}}\|\le c_{60}\|(A+\Lambda_1I)^n(\cdot)\|_{\cal{H}}$ on $D_n$, from which the required inequality follows.
\end{proof}
We can finally give the
\begin{proof}[Proof of Theorem~\ref{theorem1.2}] The first part of the statement follows by combining Theorem~\ref{Theorem4.1} and Lemma~\ref{lemma4.6}. The second part of it immediately follows by the first one since,  when $r=\infty$,  by Morrey's Theorem we have $C^\infty_{L^2}(\overline{\Omega})=\bigcap_{n\in\N_0} H^n(\Omega)$, as remarked in the Introduction, and by the same result in coordinate neighborhoods we also have $C^\infty_{L^2}(\Gamma_1)=\bigcap_{n\in\N_0} H^n(\Gamma_1)$.
\end{proof}

\chapter{Qualitative behavior of solutions when $\Omega$ is bounded} \label{section6}
In the sequel we shall study the qualitative behavior as $t\to\infty$ of solutions of problem \eqref{1.1} when also assumption (A4) holds true. When $\delta=0$ we shall give the behavior also when $t\to-\infty$, as expected by symmetry.
As anticipated in \S~\ref{intro} the analysis depends on the number of the connected components of $\Gamma_1$, but to keep it as unitary as possible we start with a general discussion.  In the sequel we shall consider different cases.
\section{General discussion}\label{section6.1} We have to consider several different types of trivial solutions of \eqref{1.1}, but for a preliminary analysis it is enough to look for spatially constant solutions, i.e., $u(t,x)=u(t)$ and $v(t,x)=v(t)$. In this way one immediately finds two types of vanishing velocity solutions, the second one taking two different forms depending on $\kappa$:
\renewcommand{\labelenumi}{{(s\arabic{enumi})}}
\renewcommand{\labelenumii}{{(\roman{enumii})}}
\begin{enumerate}
\item \label{s1} if $u_0\in\C$, $u_1\equiv 0$, $v_0\equiv v_1\equiv 0$ then $u(t,x)\equiv u_0$ and $v(t,x)\equiv 0$;
\item \label{s2} when $\kappa(x)\equiv \kappa_0\in [0,\infty)$ then two cases occur:
\begin{enumerate}
\item if $\kappa_0>0$ and $u_0\equiv 0$, $u_1\in\C$, $v_0\equiv -\rho_0u_1/\kappa_0$ and $v_1\equiv 0$ then $u(t,x)=u_1t$ and $v(t,x)\equiv v_0$;
\item if $\kappa_0=0$ and $u_0\equiv u_1\equiv 0$, $v_0\in\C$ and $v_1\equiv 0$ then $u(t,x)\equiv0$ and $v(t,x)\equiv v_0$.
\end{enumerate}
\end{enumerate}
By resuming the quick discussion in \S~\ref{intro}, solutions of type (s1) are expected. Also solutions of type (s2--ii) are mathematically expected.
Both of them  correspond to a one dimensional subspace of $\text{Ker } A$. Unfortunately solutions of type (s2--ii) are not consistent with the physical model. Indeed, when $\kappa\equiv 0$, the piece $\Gamma_1$ of the boundary is not retained in its correct position, and we shall study this case mainly for the sake of  completeness.

Solutions of type (s2--i) do not seem to have an easy physical interpretation, due to their unboundedness in time. Indeed, although the term $\kappa_0v$ acts on $\Gamma_1$ like a spring restoring force, like in an harmonic oscillator, and it should force $\Gamma_1$ to goes back to its rest position, these solutions allow arbitrary large boundary deformation $v_0$ and excess pressure $\rho_0u_t$. They motivate the more detailed discussion of the physical derivation of problem \eqref{1.1} in \S~\ref{section7} below.  {Mathematically, the issue boils down to the unboundedness of both the forward and backward semigroups,
 $\{T(t),\,t\ge 0\}$, $\{T(-t),\,t\ge 0\}$, which follows from unboundedness (in time) of solutions of type (s2--i).}

This remark motivates the first step in our analysis. Indeed, to drop physically meaningless constant functions in $\Omega$ one can introduce the space
\[
\faktor{H^1(\Omega)}{\C}:=\{[u],u\in H^1(\Omega)\},
\]
where $[\cdot]$ denotes the equivalence class with respect to the equivalence relation $\sim$ defined by $u\sim v$ if and only  if $u-v\in\C$, endowed with the norm   $\|\nabla (\cdot)\|_2$, as mentioned in \S~\ref{intro} and done in \cite{beale}.
Unfortunately, when one replaces $H^1(\Omega)$ with $\faktor{H^1(\Omega)}{\C}$ in \eqref{1.2}, as in the quoted paper, that is when one sets the alternative space
$${\cal H}_{\text{alt}}=\faktor{H^1(\Omega)}{\C}\times H^1(\Gamma_1)\times L^2(\Omega)\times L^2(\Gamma_1)$$
and the operator $A_{\text{alt}}: D(A_{\text{alt}})\subseteq {\cal H}_{\text{alt}}\to{\cal H}_{\text{alt}}$ formally defined as
 $A$ in \eqref{4.5}, one gets an ill--defined operator on any linear manifold $D(A_{\text{alt}})$ of ${\cal H}_{\text{alt}}$.
Indeed, for any $U=(u,v,w,z)\in {\cal H}_{\text{alt}}$ the first component of $A_{\text{alt}}U$ must be $-w$. To get $A_{\text{alt}}U\in {\cal H}_{\text{alt}}$ one then has to require that $w\in \faktor{H^1(\Omega)}{\C}$. But, as $D(A_{\text{alt}})\subseteq {\cal H}_{\text{alt}}$, necessarily $w\in L^2(\Omega)$, and $\faktor{H^1(\Omega)}{\C}\cap L^2(\Omega)=\emptyset$.
In \cite{beale} the author implicitly replaces $w$ with its equivalence class $[w]$, i.e., he studies a slightly different operator. In this way one  gets solutions of a problem different from \eqref{1.1}.

Also by simultaneously replacing $H^1(\Omega)$ and $L^2(\Omega)$ with $\faktor{H^1(\Omega)}{\C}$ and $\faktor{L^2(\Omega)}{\C}$, i.e., by setting the
further alternative space
$$\widetilde{{\cal H}_{\text{alt}}}=\faktor{H^1(\Omega)}{\C}\times H^1(\Gamma_1)\times \faktor{L^2(\Omega)}{\C}\times L^2(\Gamma_1)$$
 and the operator $\widetilde{A_{\text{alt}}}: D(\widetilde{A_{\text{alt}})}\subseteq \widetilde{{\cal H}_{\text{alt}}}\to\widetilde{{\cal H}_{\text{alt}}}$ formally defined as $A$ in \eqref{4.5}, one gets an ill--defined operator on any linear manifold $D(\widetilde{A_{\text{alt}})}$ of $\widetilde{{\cal H}_{\text{alt}}}$.
 Indeed, for any $U=(u,v,w,z)\in \widetilde{{\cal H}_{\text{alt}}}$ the third component of $\widetilde{A_{\text{alt}}U}$ must be $-c^2\Delta u$.
 To get $\widetilde{A_{\text{alt}}U}\in \widetilde{{\cal H}_{\text{alt}}}$ one then has to require that $\Delta u\in \faktor{L^2(\Omega)}{\C}$. But
 $\Delta 1=0$, so one cannot add a free additive constant to $\Delta u$, unless one replace $\Delta u$ with $[\Delta u]$ in \eqref{4.5}. Also in this way one
 gets solutions of a problem different from \eqref{1.1}.

To rigorously fix the problem here we shall identify $\faktor{H^1(\Omega)}{\C}$ with its isomorphic image  $H^1_c(\Omega)=\left\{u\in H^1(\Omega): {\textstyle\int}_\Omega u=0\right\}$, already defined in \eqref{5.1}, using a Poincar\'{e} type inequality (see \cite[Theorem~13.2.7, p.~423]{LeoniSobolev2}).
Unfortunately solutions of type (s2--i) show that the average $\fint_\Omega u$ is not invariant under the flow induced by \eqref{1.1}.
To bypass the problem we are going to introduce a suitable quotient group (see \cite[Chapter I, p.~43 and Chapter II, p.~61]{EngelNagel}).

We start by recalling the standard orthogonal splitting $H^1(\Omega)=H^1_c(\Omega)\oplus\C$, the corresponding orthogonal projectors $P_0:H^1(\Omega)\to H^1_c(\Omega)$ and $P_\C:H^1(\Omega)\to \C$   given by $P_\C u=\fint_\Omega u$ and $P_0u=u-\fint_\Omega u$.
Since the factor spaces in \eqref{1.2} are mutually orthogonal with respect to the inner product $(\cdot,\cdot)_{\cal H}$ defined in \eqref{4.2}
we get the following corresponding orthogonal splitting  of $\cal{H}$:
\begin{equation}\label{5.2}
\begin{alignedat}2
&\cal{H}=\C_1\oplus\cal{H}_0 , \quad &&\text{where $\C_1$ and $\cal{H}_0$ are given by}\\
&\C_1=\C1_{\cal{H}}, \quad 1_{\cal{H}}=(1,0,0,0), \quad &&\cal{H}_0=H^1_c(\Omega)\times H^1(\Gamma_1)\times L^2(\Omega)\times L^2(\Gamma_1).
\end{alignedat}
\end{equation}
The corresponding orthogonal projections  $\Pi_{\C_1}:\cal{H}\to \C_1$  and $\Pi_0:\cal{H}\to \cal{H}_0$ are trivially given by
$\Pi_1(u,v,w,z)= (P_\C u,0,0,0)$ and $\Pi_0(u,v,w,z)=(P_0 u,v,w,z)$.

In the sequel we shall equip $\cal{H}_0$ with the inner product obtained from $(\cdot,\cdot)_{\cal{H}}$ in \eqref{4.2} by simply dropping the second term, i.e., with $(\cdot,\cdot)_{\cal{H}_0}$ defined for $V_i=(u_i,v_i,w_i,z_i)$, $i=1,2$, by
\begin{multline}\label{5.2ter}
(V_1,V_2)_{\cal{H}_0}= \int_\Omega \nabla u_1\nabla \overline{u_2}+\int_{\Gamma_1}\frac\sigma{\rho_0}(\nabla_\Gamma v_1,\nabla_\Gamma v_2)_\Gamma\\
+\int_{\Gamma_1}\frac{v_1\overline{v_2}}{\rho_0}+\frac 1{c^2}\int_\Omega w_1\overline{w_2}+\int_{\Gamma_1}\frac\mu{\rho_0} z_1\overline{z_2},
\end{multline}
and we shall denote $\|\cdot\|_{\cal{H}_0}=(\cdot,\cdot)_{\cal{H}_0}^{1/2}$ the related norm, which is trivially equivalent to  $\|\cdot\|_{\cal{H}}$ on $\cal{H}_0$. Solutions of type (s2--i) show that $\cal{H}_0$ is not in general invariant under the flow induced by \eqref{1.1}, while those of type (s1) show that $\C_1$ is invariant under this flow.
Since $\faktor{H^1(\Omega)}{\C}$ is identified with $H^1_c(\Omega)$ we then get the corresponding identification of $\cal{H}_0$ with $\faktor{H^1(\Omega)}{\C}\times H^1(\Gamma_1)\times L^2(\Omega)\times L^2(\Gamma_1)$.

Following \cite{EngelNagel} we then define on $\cal{H}_0$ the quotient group $\{T_0(t),t\in\R\}$   by
\begin{equation}\label{5.2bis}
T_0(t)[\Pi_0 U]=\Pi_0[T(t)U],\qquad \text{for all $U\in\cal{H}, t\in\R$}.
\end{equation}
Its generator is $-A_0: D(A_0)\subset \cal{H}_0\to \cal{H}_0$ given, using \eqref{4.4}--\eqref{4.5}, by
\begin{align}\label{5.3}
  D(A_0)=&\Pi_0(D(A))=D(A)\cap \cal{H}_0,\quad\text{and}\\
\label{5.4}
 A_0\begin{pmatrix}u\\v\\w\\z\end{pmatrix} =
&\begin{pmatrix}-P_0w\\-z\\-c^2\Delta u\\
\frac 1\mu\left[-\DivGamma(\sigma\nabla_\Gamma v)+\delta z+\kappa v+\rho w_{|\Gamma_1}\right]
\end{pmatrix}.
\end{align}
 Trivially for any $(u_0,v_0,u_1,v_1)\in \cal{H}_0$ we have
 \begin{equation}\label{5.5}
T_0(t)[(u_0,v_0,u_1,v_1)]=(P_0u(t), v(t), u_t(t), v_t(t)),\qquad t\in\R,
 \end{equation}
 where $(u,v)$ is the solution of \eqref{1.1} corresponding to these initial data. We remark that the projections of vanishing velocity solutions  of type (s1) vanish, as expected, while for those of type (s2), when $\kappa(x)\equiv \kappa_0\in [0,\infty)$ we have:
 \renewcommand{\labelenumi}{{(\roman{enumi})}}
\begin{enumerate}
\item if $\kappa_0>0$ and $u_0\equiv 0$, $u_1\in\C$, $v_0\equiv -\rho_0u_1/\kappa_0$ and $v_1\equiv 0$ then \\$T_0(t)[(0,v_0,u_1,0)]=(0,v_0,u_1,0)$ for $t\in\R$, while
\item if $\kappa_0=0$ and $v_0\in\C$  then $T_0(t)[(0,v_0,0,0]=(0,v_0,0,0)$ for $t\in\R$.
\end{enumerate}
The case (i) clearly shows how trajectories of $\{T_0(t), t\in\R\}$ are not in general solutions of \eqref{1.1}. In both cases they make evident fixed points of the group.

In the sequel we shall then study the behavior of trajectories of $\{T_0(t), t\in\R\}$. The first step in this direction would naturally be understanding if they are bounded or not. Unfortunately the answer to this question depends on several facts, and it is essentially equivalent to get the behavior of solutions of \eqref{1.1}, with several cases occurring. In the sequel we shall then give some general properties of $A_0$.

With this aim we now introduce on $\cal{H}_0$ the pseudo--inner product $[\cdot,\cdot]_{\cal{H}_0}$ suggested by the form of the energy function in \eqref{energyidentity}, that is for all $V_i=(u_i,v_i,w_i,z_i)$, $i=1,2$,
\begin{equation}\label{5.6}
[V_1,V_2]_{\cal{H}_0}\!=\! \rho_0\!\!\int_\Omega\!\!\! \nabla u_1\nabla \overline{u_2}+\!\!\!\int_{\Gamma_1}\negquad\sigma(\nabla_\Gamma v_1,\nabla_\Gamma v_2)_\Gamma
+\int_{\Gamma_1}\negquad \kappa v_1\overline{v_2}+\tfrac {\rho_0}{c^2}\int_\Omega \negquad w_1\overline{w_2}+\!\!\!\int_{\Gamma_1} \negquad \mu z_1\overline{z_2}.
\end{equation}
The main difference between $[\cdot,\cdot]_{\cal{H}_0}$ and the inner product $(\cdot,\cdot)_{\cal{H}_0}$ in \eqref{5.2ter} is that $\rho_0(\cdot,\cdot)_{\cal{H}_0}$, trivially equivalent to $(\cdot,\cdot)_{\cal{H}_0}$, contains the term $\int_{\Gamma_1}v_1\overline{v_2}$ while in $[\cdot,\cdot]_{\cal{H}_0}$ this term is replaced by $\int_{\Gamma_1}\kappa v_1\overline{v_2}$.
Its associated pseudo--norm and null space are defined by
\begin{equation}\label{5.7}
  \n\cdot\n_{\cal{H}_0}=[\cdot,\cdot]_{\cal{H}_0}^{1/2},\qquad \cal{N}_0=\{U\in \cal{H}_0: \n U\n_{\cal{H}_0}=0\}.
\begin{footnote}{ By using the parallelogram identity for $[\cdot,\cdot]_{\cal{H}_0}$ one easily sees that $\cal{N}_0$ is a closed subspace}\end{footnote}
\end{equation}
We shall characterize $\cal{N}_0$ later on.
We shall also denote by $[\cdot,\cdot]_{\cal{H}}$ the extension of $[\cdot,\cdot]_{\cal{H}_0}$  to the whole of $\cal{H}$, again defined by \eqref{5.6}.

 Our first result on $A_0$ is the following one.
\begin{lem}\label{lemma5.1} Let assumptions (A0--4) hold. Then:
\renewcommand{\labelenumi}{{\roman{enumi})}}
\begin{enumerate}
\item $\{T_0(t), t\in\R\}$ is strongly continuous, so $A_0$ has nonempty resolvent set;
\item $A_0$ has compact resolvent;
\item for all $U=(u,v,w,z)\in D(A_0)$ we have $\Real \, [A_0U,U]_{\cal{H}_0}=\int_{\Gamma_1}\delta |z|^2\ge 0$;
\item if $\delta\equiv 0$  then $A_0U=\lambda U\Rightarrow\Real \,\lambda=0$ or $U\in \cal{N}_0$ for all $\lambda\in\C$, $U\in D(A_0)$;
\item if $\delta\not\equiv 0$  then $A_0U=\mathfrak{i}\lambda U\Rightarrow\lambda=0$ or $U=0$ for all $\lambda\in\R$, $U\in D(A_0)$.
\end{enumerate}
\end{lem}
We stress that resolvent compactness of $A_0$ is a major difference with respect to the classical case $\sigma=0$ previously studied in the literature, see~\cite[Theorem~3.12 and Corollary~3.13]{mugnolo}.

\begin{proof} Trivially i) follows by \cite[Chapter I, 5.13, p.43, Chapter II, 2.4 p.61 and Theorem~1.10 p.55]{EngelNagel}. By combining \eqref{5.3} with Lemma~\ref{lemma4.6}
we get, since $\Omega$ is bounded, that the canonical injection $D(A_0)\hookrightarrow \cal{H}_0$ is compact.
 Hence, by
\cite[Chapter II, Proposition~4.25 p.~117]{EngelNagel}, $A_0$ has compact resolvent, proving ii). Next iii) follows since, by \eqref{3.40}, \eqref{3.36}, \eqref{5.4} and \eqref{5.6},
\begin{align*}
[A_0U,U]_{\cal{H}_0}=&-\rho_0\int_\Omega\nabla w\nabla \overline{u}-\rho_0\int_\Omega\Delta u\overline w-\int_{\Gamma_1}\sigma(\nabla_\Gamma z,\nabla_\Gamma v)_\Gamma-\int_{\Gamma_1}\kappa \overline{v}z\\
+&\int_{\Gamma_1}[-\DivGamma (\sigma\nabla_\Gamma v)+\kappa v+\delta z+\rho w]\overline{z}\\
=&-\rho_0\int_\Omega[\nabla w\nabla \overline{u}-\nabla u\nabla \overline w]+\int_{\Gamma_1}\sigma[(\nabla_\Gamma v,\nabla_\Gamma z)_\Gamma-(\nabla_\Gamma z,\nabla_\Gamma v)_\Gamma]\\
+&\int_{\Gamma_1} \kappa [ v\overline{z}-z\overline{v}]+\int_{\Gamma_1}\rho_0[ w\overline{z}-z\overline{w}]+\int_{\Gamma_1}\delta |z|^2.
\end{align*}
As to iv), we note that by \eqref{5.7} it is an immediate consequence of iii).

To prove v)  we take $U=(u,v,w,z)\in D(A_0)$ and $\lambda\in\R\setminus\{0\}$ such that $A_0U=i\lambda U$.
Hence, by Lemma~\ref{lemma4.6}, \eqref{5.2}, \eqref{5.3}, \eqref{5.4} we have $u\in H^2(\Omega)$, $v\in H^2(\Gamma_1)$, $w\in H^1(\Omega)$, $z\in H^1(\Gamma_1)$, $\int_\Omega u=0$ and
\begin{equation}\label{5.9}
\begin{cases}
-w+\fint_\Omega w =\mathfrak{i}\lambda u \qquad &\text{in
$\Omega$,}\\
- z=\mathfrak{i}\lambda v\qquad
&\text{on
$\Gamma_1$,}\\
-c^2\Delta u=\mathfrak{i}\lambda w \qquad &\text{in
$\Omega$,}\\
- \DivGamma (\sigma \nabla_\Gamma v)+\kappa v+\delta z+\rho_0 w =\mathfrak{i}\lambda \mu z\qquad
&\text{on
$\Gamma_1$,}\\
\partial_\nu u=z\qquad
&\text{on
$\Gamma_1$,}\\
\partial_\nu u=0 &\text{on $\Gamma_0$.}
\end{cases}
\end{equation}
Since $\Real \, [A_0U,U]_{\cal{H}_0}=\Real \, (\mathfrak{i}\lambda)\,[U,U]_{\cal{H}_0}=0$, by iii) we get $\delta |z|^2=0$ in $L^1(\Gamma_1)$.
Since, by assumption (A4), $r=2$, by (A3) we have $\delta\in W^{1,\infty}(\Gamma_1)$, so by applying Morrey's Theorem in local coordinates we have $\delta\in C(\Gamma_1)$. Hence, as $\delta\ge 0$ and $\delta\not\equiv 0$, the set $\widetilde{\Gamma_1}=\{y\in \Gamma_1: \delta(y)>0\}$ is nonempty and relatively open on $\Gamma_1$. Trivially, $z=0$ and  consequently, by \eqref{5.9}$_2$,  $v=0$ on $\widetilde{\Gamma_1}$, $\cal{H}^{N-1}$-a.e..

By denoting $\alpha=\fint_\Omega w$, using the first two equations to eliminate $w$ and $z$,  by \eqref{5.9} we get that $u$ and $v$ satisfy the  coupled elliptic problems
\begin{equation}\label{5.901}
\begin{cases}
\begin{aligned}
-c^2\Delta u&=\mathfrak{i}\lambda \alpha+\lambda^2u \qquad &&\text{in
$\Omega$,}\\
\partial_\nu u&=-\mathfrak{i}\lambda v\qquad
&&\text{on
$\Gamma_1$,}\\
\partial_\nu u&=0 &&\text{on $\Gamma_0$}
\end{aligned}
\end{cases}
\end{equation}
and
\begin{equation}\label{5.902}
- \DivGamma (\sigma \nabla_\Gamma v)+(\kappa-\mathfrak{i}\lambda\delta-\lambda^2\mu)v =\mathfrak{i}\rho_0\lambda u-\rho_0\alpha\qquad
\text{on $\Gamma_1$.}
\end{equation}
By standard elliptic regularity (see \cite[Theorem~3, p.~334]{Evans})  we have $u\in C^\infty(\Omega)$. Moreover we claim that
$(u,v)\in \bigcap_{q\in [1,\infty)}\mathbb{W}^{2,q}$, where
$$\mathbb{W}^{2,q}:=W^{2,q}(\Omega)\times W^{2,q}(\Gamma_1).$$
To prove our claim we adapt the standard bootstrap procedure to the coupled problems \eqref{5.901}--\eqref{5.902}. Indeed, suppose that
$(u,v)\in \mathbb{W}^{2,p_0}$ for some $1<p_0<\infty$. By the Sobolev Embedding Theorem, since $\Gamma_1$ is compact and $v\in W^{2,p_0}(\Gamma_1)$ we get  that  $v\in L^q(\Gamma_1)$ for all
$q\in [1,q_1)$, where $q_1=p_0(N-1)/(N-1-2p_0)$ if $p_0<(N-1)/2$, $q_1=\infty$ otherwise. Moreover, since $u\in W^{2,p_0}(\Omega)$, by the Trace Theorem
$u_{|\Gamma_1}\in W^{2-1/p_0,p_0}(\Gamma_1)$.
Using again the Sobolev Embedding Theorem we thus have  $u_{|\Gamma_1}\in L^q(\Gamma_1)$ for all
$q\in [1,p_1)$, where $p_1=p_0(N-1)/(N-2p_0)$ if $p_0<N/2$, $p_1=\infty$ otherwise.
Since trivially $p_1\le q_1$, then $u_{|\Gamma_1},v\in L^q(\Gamma_1)$ for all
$q\in [1,p_1)$. Then, since $\sigma,\kappa,\delta,\mu\in W^{1,\infty}(\Gamma_1)$, by applying elliptic regularity to \eqref{5.902}
we get that $v\in W^{2,q}(\Gamma_1)$ for the same values of $q$. Since, by the Sobolev Embedding Theorem on $\Omega$ we also have $u\in L^q(\Omega)$ for the same values of $q$, applying elliptic regularity to \eqref{5.901} we then get $u\in W^{2,q}(\Omega)$ for them, so $(u,v)\in\mathbb{W}^{2,q}$.
Now the sequence defined by recurrence by $p_{n+1}=p_n(N-1)/(N-2p_n)$, $p_0>1$,  trivially goes to $\infty$, so proving our claim.

Hence, by Morrey's Theorem we get that $u\in C^\infty(\Omega)\cap C^1(\overline{\Omega})$. Moreover, by \eqref{5.901}--\eqref{5.902}, since
 $v=0$ on $\widetilde{\Gamma_1}$, $u$ satisfies the partially overdetermined elliptic problem
\begin{equation}\label{5.10}
\begin{cases}
\begin{split}
-c^2\Delta u&=\mathfrak{i}\lambda\alpha+\lambda^2u \qquad &&\text{in
$\Omega$,}\\
\partial_\nu u&=0\qquad
&&\text{on
$\widetilde{\Gamma_1}$,}\\
u&=-\mathfrak{i}\alpha /\lambda &&\text{on $\widetilde{\Gamma_1}$,}
\end{split}
\end{cases}
\end{equation}
which possesses the trivial solution $\tilde{u}=-\mathfrak{i}\alpha/\lambda$. Now $u$ is regular enough to apply the enhancement of the classical Unique Continuation Principle recently given in \cite[Theorem~1]{FarinaValdinoci}, from which we get that $u=\tilde{u}$ in $\Omega$. Since $\int_\Omega u=0$,  we then get $\alpha=0$ and $u=0$ in $\Omega$. By \eqref{5.9} then $U=0$, concluding the proof.
\end{proof}

The last result clearly shows the usefulness of understanding when $[\cdot,\cdot]_{\cal{H}_0}$ is an inner product on $\cal{H}_0$ or, this assertion failing, on a suitable closed subspace of it, and if its associated norm is equivalent to $\|\cdot\|_{\cal{H}_0}$.  To avoid checking this fact on various subspaces we point out the following abstract Poincaré type inequality. In the existing literature we have been able to find it only in \cite[Proposition~1]{Graeser}, which deals with the real case only; unfortunately the approach in the quoted paper does not seem to extend to the complex case.

The following result is probably folklore in the theory of Krein spaces, but the authors were unable to find  appropriate references.  We prove it for the reader's convenience.
\begin{lem}\label{lemma5.2}
Let $[\cdot,\cdot]$ be a continuous pseudo--inner product on a real or complex Hilbert space $(H, (\cdot,\cdot))$, and denote
$\cal{N}=\{u\in H: [u,u]=0\}$. Then
\renewcommand{\labelenumi}{{\roman{enumi})}}
\begin{enumerate}
\item $\cal{N}=\{u\in H: [u,v]=0\,\,\text{for all $v\in H$}\}$ and $\cal{N}$ is a closed subspace of $H$;
\item $[\cdot,\cdot]$ is positive definite on a subspace $H_1$ of $H$ if and only if $H_1\cap \cal{N}=\{0\}$;
\item if $\text{dim}\,\cal{N}<\infty$ then $[\cdot,\cdot]$ is coercive on any closed subspace $H_1$ of $H$ such that $H_1\cap \cal{N}=\{0\}$ if and only if
it is coercive on $\cal{N}^\bot$.
\end{enumerate}
\end{lem}
\begin{rem} Taking $H=\ell^2$, $[(x_n)_n,(y_n)_n]=\sum_nx_{2n}\overline{y_{2n}}$ for all $(x_n)_n,(y_n)_n\in \ell^2$, $H_1=\{(x_n)_n\in\ell^2: x_{2n+1}=nx_{2n}\,\,\text{for all $n\in\N$}\}$ one easily sees that iii) fails when $\text{dim}\,\cal{N}=\infty$.
\end{rem}

\begin{proof}[Proof of Lemma~\ref{lemma5.2}] By using the Cauchy--Schwarz inequality for $[\cdot,\cdot]$ one easily get that $\cal{N}=\{u\in H: [u,v]=0\,\,\text{for all $v\in H$}\}$ (see also \cite[Lemma~4.4]{bognar}), and thus i) follows. Moreover ii) is trivial, as
the fact that coercivity on all closed subspaces $H_1$ having trivial intersection with $\cal{N}$ implies coercivity on $\cal{N}^\bot$.
Hence we suppose that $\text{dim}\,\cal{N}<\infty$, that $[\cdot,\cdot]$ is coercive on $\cal{N}^\bot$ and we take a closed subspace $H_1$ of $H$ with $H_1\cap \cal{N}=\{0\}$. Let $T\in\cal{L}(H,H_1^\bot)$ be the orthogonal projection, so $\text{Ker}\,T=H_1$. Hence $\text{Ker}\,T\cap \cal{N}=\{0\}$ and consequently the restriction of $T$ to $\cal{N}$ is injective. Then $\text{dim}\, H_1^\bot\ge \text{dim}\,\cal{N}$.

Actually we claim that we can assume, without restriction, that $\text{dim}\, H_1^\bot=\text{dim}\,\cal{N}$.
Indeed, if $\text{dim}\, H_1^\bot>\text{dim}\,\cal{N}$, we replace $H_1$ with the maximal closed subspace $H_2$ of $H$ containing $H_1$ and having trivial intersection with $\cal{N}$, that  is with $H_2=H_1\oplus (H_1\oplus \cal{N})^\bot$. Indeed trivially $H_1\subseteq H_2$ and for any $u\in H_2\cap\cal{N}$ there are $v\in H_1$ and $w\in (H_1\oplus \cal{N})^\bot$ such that $u=v+w$, so $w=u-v\in (H_1\oplus \cal{N})\cap (H_1\oplus \cal{N})^\bot=\{0\}$, from which $w=u-v=0$ and consequently $u=v\in H_1\cap\cal{N}=\{0\}$. Hence $H_2\cap \cal{N}=\{0\}$. Moreover, by well-known properties of orthogonal complements we have
$H_2^\bot=[H_1\oplus (H_1\oplus \cal{N})^\bot]^\bot=H_1^\bot\cap (H_1\oplus \cal{N})=H_1^\bot\cap \cal{N}\subseteq\cal{N}$ and consequently $H=H_2\oplus H_2^\bot\subseteq H_2\oplus \cal{N}$, so $H=H_2\oplus \cal{N}$. Hence $H_2\simeq \faktor{H}{\cal{N}} \simeq \cal{N}^\bot$ and consequently  $\text{dim}\, H_2^\bot=\text{dim}\,\cal{N}$, proving our claim.

Since $\text{dim}\, H_1^\bot=\text{dim}\,\cal{N}$ there is $\Phi_0: H_1^\bot\to \cal{N}$ linear and bijective . We now set $S=\Phi_0\cdot T\in \cal{L}(H,\cal{N})$, so $\Ker S=\Ker T=H_1$ and, as $H_1\cap\cal{N}=\{0\}$,  $S_{|\cal{N}}\in \cal{L}(\cal{N})$ is injective. Since $\text{dim\,}\cal{N}<\infty$ we then also have that $S_{|\cal{N}}$ is bijective and $S_{|\cal{N}}^{-1}\in \cal{L}(\cal{N})$.
Now, respectively denoting by $\pi_{\cal N}$ and $\pi_{{\cal N}^\bot}$ the orthogonal projections of $H$ onto $\cal N$ and $\cal {N}^\bot$, for any $u\in H_1$ we have $S\pi_{\cal N} u+S\pi_{{\cal N}^\bot} u=Su=0$, and then $\pi_{\cal{N}}u=-S_{|\cal{N}}^{-1}S\pi_{\cal{N}^\bot}u$. Consequently, denoting $\|\cdot\|=(\cdot,\cdot)^{1/2}$, for any $u\in H_1$ we have
\begin{equation}\label{5.11}
\|u\|^2= \|\pi_\cal{N} u\|^2+\|\pi_{\cal{N}^\bot} u\|^2\le \left(1+\|S_{|\cal{N}}^{-1}\|^2_{\cal{L}(\cal{N})}\|S\|^2_{\cal{L}(H)}\right)\|\pi_{\cal{N}^\bot} u\|^2.
\end{equation}
Now, since by i) one trivially get that for any $u\in H$ one has $[u,u]=[\pi_{\cal{N}^\bot} u,\pi_{\cal{N}^\bot} u]$, by combining \eqref{5.11} with the assumed coercivity of $[\cdot,\cdot]$ on  $\cal{N}^\bot$ we complete the proof.
\end{proof}

We are now going to characterize $\cal{N}_0$ in \eqref{5.7}. We introduce the notation
\begin{gather}\label{5.12}
\cal{C}(\Gamma_1)=\{\Gamma'\subset \Gamma_1: \text{$\Gamma'$ is a connected component of $\Gamma_1$}\},\\
\cal{C}^0(\Gamma_1)=\{\Gamma'\in \cal{C}(\Gamma_1): \kappa_{|\Gamma'}\equiv 0\},\qquad  \cal{C}^+(\Gamma_1)=\{\Gamma'\in \cal{C}(\Gamma_1): \kappa_{|\Gamma'}\not\equiv 0\},\\
\label{5.13}\mathfrak{n}=\# \cal{C}(\Gamma_1), \quad \mathfrak{n}_0=\#\cal{C}^0(\Gamma_1), \quad \text{and}\quad \mathfrak{n}_+=\#\cal{C}^+(\Gamma_1),
\end{gather}
where $\#$ is the cardinality. Since $\Gamma_1$ is compact $\mathfrak{n}_0+\mathfrak{n}_+=\mathfrak{n}<\infty$.
Moreover, for any $\Gamma'\subseteq \Gamma_1$ (relatively open) we shall denote by $[\cdot,\cdot]_{\Gamma'}$ the pseudo--inner product on $H^1(\Gamma')$ defined by
\begin{equation}\label{5.14}
[v_1,v_2]_{\Gamma'}=\int_{\Gamma'} \sigma (\nabla_\Gamma v_1,\nabla_\Gamma v_2)_\Gamma + \int_{\Gamma'} \kappa v_1\overline{v_2},
\end{equation}
and respectively by $\n \cdot\n_{\Gamma'}$, $\cal{N}(\Gamma')$ its associated pseudo--norm and null space. In analogy with \eqref{5.1} we shall also set
\begin{equation}\label{5.15}
 H^1_c(\Gamma'):=\left\{v\in H^1(\Gamma'): {\textstyle\int}_{\Gamma'} v=0\right\}.
\end{equation}
We then have
\begin{lem}\label{lemma5.3} Let assumptions (A0--4) hold. Then
\renewcommand{\labelenumi}{{\roman{enumi})}}
\begin{enumerate}
\item if $\Gamma'\in \cal{C}^+(\Gamma_1)$ then $[\cdot,\cdot]_{\Gamma'}$ is an inner product on $H^1(\Gamma')$ and $\n\cdot\n_{\Gamma'}$ is equivalent to $\|\cdot\|_{1,2,\Gamma'}$ on it;
\item if $\Gamma'\in \cal{C}^0(\Gamma_1)$ then $\cal{N}(\Gamma')=\C \chi_{\Gamma'}$ but $[\cdot,\cdot]_{\Gamma'}$ is an inner product on $H_c^1(\Gamma')$ and $\n\cdot\n_{\Gamma'}$ is equivalent to $\|\cdot\|_{1,2,\Gamma'}$ on it;
\item if $\mathfrak{n}_0=0$ then $[\cdot,\cdot]_{\Gamma_1}$ is an inner product on $H^1(\Gamma_1)$ and $\n\cdot\n_{\Gamma_1}$ is equivalent to $\|\cdot\|_{1,2,\Gamma_1}$ on it;
\item if $\mathfrak{n}_0\ge 1$  then $\cal{N}(\Gamma_1)=\text{span\,}\{\chi_{\Gamma'}, \Gamma'\in \cal{C}^0(\Gamma_1)\}\not=\{0\}$, but $[\cdot,\cdot]_{\Gamma_1}$ is an inner product on $H_\kappa^1(\Gamma_1)=\{v\in H^1(\Gamma_1): {\textstyle\int}_{\Gamma'}v=0\quad\forall \,\Gamma'\in \cal{C}^0(\Gamma_1)\}$ and $\n\cdot\n_{\Gamma_1}$ is equivalent to $\|\cdot\|_{1,2,\Gamma_1}$ on it;
\item if $\mathfrak{n}_0=0$ then $[\cdot,\cdot]_{\cal{H}_0}$ is an inner product on $\cal{H}_0$ and $\n\cdot\n_{\cal{H}_0}$ is equivalent to $\|\cdot\|_{\cal{H}_0}$ on it;
\item if $\mathfrak{n}_0\ge 1$  then $\cal{N}_0=\{0\}\times \cal{N}(\Gamma_1)\times\{0\}\times\{0\}\not=\{0\}$,
but $[\cdot,\cdot]_{\cal{H}_0}$ is an inner product on any closed subspace $\cal{H'}$ of $\cal{H}_0$ such that $\cal{H}'\cap\cal{N}_0=\{0\}$ and $\n\cdot\n_{\cal{H}_0}$ is equivalent to $\|\cdot\|_{\cal{H}_0}$ on it;
\end{enumerate}
\end{lem}
\begin{proof} To prove i--ii) we first recall that $\{v\in H^1(\Gamma'): \nabla_\Gamma v=0\}=\C \chi_{\Gamma'}$ for any $\Gamma'\in \cal{C}(\Gamma_1)$. Indeed for any
$v\in H^1(\Gamma')$ we have $\nabla_\Gamma v=0$ if and only if $v$ is $\cal{H}^{N-1}$-a.e.\ constant on any coordinate neighborhood. Moreover, since $\Gamma'$ is connected, a standard topological argument  shows that all these constant must coincide.
Hence, by assumptions (A1), (A4) and by \eqref{5.14}, we get that $v\in \cal{N}(\Gamma')$ if and only if $v=c\chi_{\Gamma'}$, $c\in\C$ and $c{\textstyle\int}_{\Gamma'}\kappa=0$. So for any $\Gamma'\in \cal{C}(\Gamma_1)$ we have $\cal{N}(\Gamma')\not=\{0\}$ if and only if $\Gamma'\in \cal{C}^0(\Gamma_1)$, and in this case $\cal{N}(\Gamma')=\C \chi_{\Gamma'}$, so $\cal{N}(\Gamma')\cap H^1_c(\Gamma')=\{0\}$.
Hence, since $[\cdot,\cdot]_{\Gamma'}$ is trivially a continuous pseudo-inner product, to prove the asserted equivalencies in i--ii) reduce to prove coercivity of $[\cdot,\cdot]_{\Gamma'}$.

To this aim, we apply in both cases a well-known Poincaré-type inequality, i.e., \cite[Lemma~4.1.3, p.178] {ziemer}. We take, in the author's notation, $X_0=L^2(\Gamma')$, $X=H^1(\Gamma')$, $\|\cdot\|_1=\|\nabla_\Gamma(\cdot)\|_{2,\Gamma'}$,  $Y=\{v\in H^1(\Gamma'): \|v\|_1=0\}=\C \chi_{\Gamma'}$ and $L\in\cal{L}(X,Y)$ given by
$$
Lv=
\begin{cases}
\left(\int_{\Gamma_1}\sqrt \kappa\right)^{-1}\int_{\Gamma_1}\sqrt \kappa v\quad \chi_{\Gamma'},\quad &\text{if $\Gamma'\in \cal{C}^+(\Gamma_1)$,}\\
\left(\cal{H}^{N-1}(\Gamma')\right)^{-1}\int_{\Gamma_1}v\quad \chi_{\Gamma'},\quad &\text{if $\Gamma'\in \cal{C}^0(\Gamma_1)$.}
\end{cases}
$$
A trivial check shows that $L^2=L$, so by the quoted result there is a positive constant $c_{61}=c_{61}(\kappa, \Gamma')$ such that
\begin{equation}\label{5.21}
  \|v-Lv\|_{2,\Gamma'}^2\le c_{61}\|\nabla_\Gamma v\|_{2,\Gamma'}^2\qquad\text{for all $v\in H^1(\Gamma')$.}
\end{equation}
Hence, by \eqref{5.14}--\eqref{5.15} and assumption (A1), $[\cdot,\cdot]_{\Gamma'}$ is coercive on $H^1_c(\Gamma')$ when $\Gamma'\in\cal{C}^0(\Gamma_1)$.
When $\Gamma'\in\cal{C}^+(\Gamma_1)$, by \eqref{5.21} and Cauchy--Schwarz inequality  we get
\begin{align*}
\|v\|_{2,\Gamma'}^2\le & 2 c_{61}\left(\|\nabla_\Gamma v\|_{2,\Gamma'}^2+\|Lv\|_{2,\Gamma'}^2\right)\\
=& 2 c_{61}\left[\|\nabla_\Gamma v\|_{2,\Gamma'}^2+\left({\textstyle\int}_{\Gamma'}\sqrt \kappa\right)^{-1}\left|{\textstyle\int}_{\Gamma'}\sqrt \kappa v\right|^2\cal{H}^{N-1}(\Gamma')\right]\\
\le & 2 c_{61}\left[\|\nabla_\Gamma v\|_{2,\Gamma'}^2+\left({\textstyle\int}_{\Gamma'}\sqrt \kappa\right)^{-1} \left(\cal{H}^{N-1}(\Gamma')\right)^2
{\textstyle\int}_{\Gamma'}\kappa|v|^2\right],
\end{align*}
from which the coercivity of $[\cdot,\cdot]_{\Gamma'}$ on $H^1(\Gamma')$ trivially follows, so proving i--ii).

To get iii)--iv) we remark that, by identifying for any $\Gamma'\in\cal{C}(\Gamma_1)$ the space $H^1(\Gamma')$ with its isometric image in $H^1(\Gamma_1)$ obtained by the trivial extension, we have \mbox{$H^1(\Gamma_1)=\negquad\bigoplus\limits_{\Gamma'\in \cal{C}(\Gamma_1)}H^1(\Gamma')$} and $[\cdot,\cdot]_{\Gamma_1}=\negquad\sum\limits_{\Gamma'\in \cal{C}(\Gamma_1)}[\cdot,\cdot]_{\Gamma'}$, so  $N(\Gamma_1)=\negquad\bigoplus\limits_{\Gamma'\in \cal{C}(\Gamma_1)}N(\Gamma')$. Hence, by i)--ii), $N(\Gamma_1)=\negquad\negquad\bigoplus\limits_{\Gamma'\in \cal{C}^0(\Gamma_1)}\negquad\C\chi_{\Gamma'}=\text{span\,}\{\chi_{\Gamma'},\,\, \Gamma'\in\cal{C}^0(\Gamma_1)\}$.
Moreover trivially
\mbox{$H^1_\kappa(\Gamma_1)=\Big[\bigoplus_{\Gamma'\in \cal{C}^+(\Gamma_1)} H^1(\Gamma')\Big]\oplus \Big[\bigoplus_{\Gamma'\in \cal{C}^0(\Gamma_1)}H^1_c(\Gamma')\Big]$} and the asserted equivalencies follow by i--ii) since $\cal{C}(\Gamma_1)$ is finite, proving iii)--iv).

To get v--vi) we remark that \eqref{5.6} trivially yields that $\cal{N}_0=\{0\}\times \cal{N}(\Gamma_1)\times\{0\}\times\{0\}\not=\{0\}$, so by iii) one gets v). When $\mathfrak{n}_0\ge 1$ we remark that one immediately get that $N(\Gamma_1)^\bot=H^1_\kappa(\Gamma_1)$ and consequently, by \eqref{5.6} , that
$\cal{N}_0^\bot=H^1_c(\Omega)\times H^1_\kappa(\Gamma_1)\times L^2(\Omega)\times L^2(\Gamma_1)$. Hence, by iv), $[\cdot,\cdot]_{\cal{H}_0}$ is coercive on $\cal{N}_0^\bot$. Since $\text{dim\,}\cal{N}_0=\mathfrak{n}_0$, by applying Lemma~\ref{lemma5.2}--iii) we complete the proof.
\end{proof}
Before starting our discussion on the behavior of trajectories of $\{T_0(t),t\in\R\}$ in various cases we are going to generalize vanishing velocity solutions of type (s2) to non-constant functions $\kappa$. By looking for weak solutions of the form $u(t,x)=u(t)$, $v(t,x)=v_0(x)$, corresponding to initial data $u(0)\equiv 0$, $u_t(0)\in\C$, $v(0)\equiv v_0$, $v_t(0)\equiv 0$, we find of interest to consider weak solutions $v_0\in H^1(\Gamma_1)$ of the elliptic equation
\begin{equation}\label{5.22}
-\DivGamma(\sigma \nabla_\Gamma v_0)+\kappa v_0+\rho_0u_1=0,\qquad \text{on $\Gamma_1$}
\end{equation}
with parameter $u_1\in\C$. Clearly \eqref{1.9} is nothing but \eqref{5.22} with parameter $1$. Such weak solutions are defined as $v_0\in H^1(\Gamma_1)$ such that
\begin{equation}\label{5.23}
  \int_{\Gamma_1}\sigma (\nabla_\Gamma v_0,\nabla_\Gamma \psi)_\Gamma +\int_{\Gamma_1}\kappa v_0\overline{\psi}+\rho_0u_1\int_{\Gamma_1}\overline{\psi}=0\quad\text{for all $\psi\in H^1(\Gamma_1)$,}
\end{equation}
and we shall also say that the couple $(u_1,v_0)\in\C\times H^1(\Gamma_1)$ solves \eqref{5.22}.
We remark that, by Theorem~\ref{lemma4.5}, automatically $v_0\in H^2(\Gamma_1)$.
We  have the following easy result
\begin{lem}\label{lemma5.4}Let assumptions (A0--4) hold. Then
\renewcommand{\labelenumi}{{\roman{enumi})}}
\begin{enumerate}
\item for any $u_1\in\C$ and $v_0\in H^1(\Gamma_0)$ the couple  $u(t)=u_1t$, $v(t)=v_0$ is a weak solution of \eqref{1.1} if and only if $v_0$ is a weak solution of \eqref{5.22} with parameter $u_1$;
\item if $\mathfrak{n}_0=0$  then for all $u_1\in \C$ the equation \eqref{5.22} has a unique weak solution $v_0=u_1v^*$, where $v^*$ is the solution of \eqref{1.9};
\item if $\mathfrak{n}_0\ge 1$  then equation \eqref{5.22} has solutions only when $u_1=0$, and in this case these solutions are nothing but the elements of $\cal{N}_0$;
\item in both cases $T_0(t)[(0,u_1,v_0,0)]=(0,u_1,v_0,0)$ for all $t\in\R$. Moreover
\begin{equation}\label{5.24}
\Ker A_0=\{(0,u_1,v_0,0)\in\cal{H}_0:\,\text{$(u_1,v_0)$ solves \eqref{5.22}}
\}.
\end{equation}
\end{enumerate}
\end{lem}
\begin{proof} By testing $u(t)=u_1t$, $v(t)=v_0$ with Definition~\ref{Definition4.2} one  gets i). When $\mathfrak{n}_0=0$, by Lemma~\ref{lemma5.3}--v) we can apply the Riesz--Fischer Theorem with inner product $[\cdot,\cdot]_{\cal{H}_0}$ and get ii).
When $\mathfrak{n}_0\ge 1$, taking any $\Gamma'\in \cal{C}^0(\Gamma_1)$ and $\psi=\chi_{\Gamma'}$ as a test function in \eqref{5.23} we get that \eqref{5.22} may admit solution only when $u_1=0$. Moreover, in this case, by Lemma~\ref{lemma5.3}--i), $v_0$ solves \eqref{5.23} is and only if $v_0\in \cal{N}(\Gamma_1)$. Then, by Lemma~\ref{lemma5.3}--vi) we complete the proof of iii). By \eqref{5.5} we have $T_0(t)[(0,u_1,v_0,0)]=(0,u_1,v_0,0)$ for all $t\in\R$. Finally \eqref{5.24}  follows by \eqref{5.4}.
\end{proof}
By combining the last two results we are now able to generalize vanishing velocity solutions of type (s2) to $\kappa\in W^{1,\infty}(\Gamma_1)$, $k\ge 0$, also non-constant,  as follows:
\renewcommand{\labelenumii}{{(\roman{enumii})}}
\begin{enumerate}
\item[(s2)']
\begin{enumerate}
\item
if $\mathfrak{n}_0=0$ and $u_0\equiv 0$, $u_1\in \C$, $v_0=u_1v^*$, $v_1\equiv 0$ then $u(t,x)=u_1t$, $v(t,x)=v_0(x)$. Moreover
$T_0(t)[(0,v_0,u_1,0)]=(0,v_0,u_1,0)$ for all $t\in\R$;
\item if $\mathfrak{n}_0\ge 1$ and $u_0=u_1\equiv 0$,  $v_0=\sum_{\Gamma'\in \cal{C}^0(\Gamma_1)}\alpha_{\Gamma'}\chi_{\Gamma'}$, $\alpha_{\Gamma'}\in \C$, $v_1\equiv 0$ then $u(t,x)\equiv 0$, $v(t,x)=v_0(x)$. Moreover
$T_0(t)[(0,v_0,0,0)]=(0,v_0,0,0)$ for all $t\in\R$.
\end{enumerate}
\end{enumerate}
Let us now briefly explain the general strategy for determining the behavior of trajectories of $\{T_0(t), t\in\R\}$ and  $\{T(t), t\in\R\}$ suggested by previous results. Since these behaviors are trivial on $\Ker A_0$ we shall look for a suitable complement $\cal{H}^*$ of it which is invariant under the flow, so the two groups split. By Lemmas~\ref{lemma5.1} and \ref{lemma5.3}--\ref{lemma5.4} the group restriction on any such $\cal{H}^*$ possesses good properties when  $\cal{H}^*$ is equipped with $[\cdot,\cdot]_{\cal{H}_0}$, which is an inner product on it with equivalent norm. The key point for the success of this strategy is clearly to find out $\text{dim Ker\,}A_0=\max\{\mathfrak{n}_0,1\}$ invariants.

For the sake of clearness, we shall first analyze, in the next section, the case  of $\Gamma_1$ connected, giving the proof of Theorems~\ref{theorem1.3} and \ref{theorem1.4}. In \S~\ref{subsection6.3} we shall then generalize the analysis to the general case.

\section{The simplest case:  $\Gamma_1$ connected}\label{subsection6.2} In this case the two cases arising from Lemma~\ref{lemma5.3} are $\mathfrak{n}_0=0$ and $\mathfrak{n}_0=1$, respectively  reducing to $\kappa\not\equiv 0$ and $\kappa\equiv 0$. We shall briefly analyze them separately, and later show how the two cases can be conveniently unified.
\subsection{The case $\kappa\not\equiv 0$} In this case, by Lemma~\ref{lemma5.3}, $[\cdot,\cdot]_{\cal{H}_0}$ is an inner product and $\n\cdot\n_{\cal{H}_0}$ is an equivalent norm on $\cal{H}_0$. We equip $\cal{H}_0$ with them. By Lemma~\ref{lemma5.4}, setting
\begin{equation}\label{5.25}
  V^*=(0,v^*,1,0),\qquad\text{and}\quad V_0=\C V^*,
\end{equation}
$V_0$ is entirely constituted by fixed points of the group.~ The trivial candidate for a complement of $V_0$ is clearly the orthogonal complements $V_0^{\bot_{\cal{H}_0}}$ with respect to the inner product $[\cdot,\cdot]_{\cal{H}_0}$. A straightforward calculation shows that
$$V_0^{\bot_{\cal{H}_0}}=\left\{(u,v,w,z)\in \cal{H}_0: {\textstyle\int}_\Omega w-c^2{\textstyle\int}_{\Gamma_1}v=0\right\}.$$
Now, even in the well-studied case $\sigma\equiv 0$, this subspace is invariant under the flow $\{T_0(t), t\in\R\}$, even when $\Gamma_1$ is disconnected and regardless of $\kappa$:  this fact was never remarked before, to the best our knowledge.

 The following result, which was not given in our general discussion since it originates from previous arguments, will  be used in all cases.
\begin{lem}\label{lemma5.5}Let assumptions (A0--4) hold. Then for any weak solution $(u,v)$ of \eqref{1.1} we have
\begin{equation}\label{5.27}
  \int_\Omega u_t(t)-c^2\int_{\Gamma_1} v(t)=\int_\Omega u_1-c^2\int_{\Gamma_1} v_0\qquad\text{for all $t\in\R$.}
\end{equation}
Consequently, setting
\begin{gather}\label{5.28}
  L_1\in \cal{H}'\quad\text{defined by}\quad L_1(u,v,w,z)=\int_\Omega w-c^2\int_{\Gamma_1}v,\\
\cal{H}_1=\Ker L_1,\qquad \cal{H}_{01}=\cal{H}_0\cap\cal{H}_1,
\end{gather}
then $\cal{H}_1$ and $\cal{H}_{01}$ are respectively invariant for $\{T_0(t), \,t\in\R\}$ and  $\{T(t),\, t\in\R\}$.
\end{lem}
\begin{proof} By \eqref{4.5}--\eqref{4.6}  for data $(u_0,v_0,u_1,v_1)\in D(A)$ we have
$$\frac d{dt}\left(\int_\Omega u_t(t)-c^2\int_{\Gamma_1} v(t)\right)=\int_\Omega u_{tt}-c^2\int_{\Gamma_1} v_t=
c^2\left(\int_\Omega \Delta u-\int_{\Gamma_1} \partial_\nu u\right)=0$$
for all $t\in\R$ by the Divergence Theorem in $H^2(\Omega)$, so \eqref{5.27} holds. Since, as shown in the proof of Theorem~\ref{theorem1.1}, weak solutions of \eqref{1.1} are also generalized solutions,  then \eqref{5.27}
extends by density to data in $\cal{H}$, thus proving the invariance of $\cal{H}_1$ with respect to $\{T(t),\, t\in\R\}$. By \eqref{5.5} we then get the invariance of $\cal{H}_{01}$ with respect to $\{T_0(t),\, t\in\R\}$.
\end{proof}
Resuming the discussion, recalling the splitting \eqref{5.2} and pointing out the trivial one $\cal{H}_1=\C_1\oplus \cal{H}_{01}$, we then have the $[\cdot,\cdot]_{\cal{H}}$-orthogonal splittings
\begin{equation}\label{5.30}
  \cal{H}_0=\cal{H}_{01}\oplus V_0,\qquad \cal{H}=\cal{H}_1\oplus V_0.
\end{equation}
It is useful to calculate the projections associated to them. By \eqref{5.6}, \eqref{5.23}, \eqref{5.25} and \eqref{5.28}
we have, for $U=(u,v,w,z)\in \cal{H}_0$, with $[\cdot,\cdot]_{\Gamma_1}$ defined in \eqref{5.14},
\begin{align}\label{5.31}
0<\n V^*\n_{\cal{H}_0}^2=&\tfrac {\rho_0}{c^2}{\textstyle\int}_\Omega 1+[v^*,v^*]_{\Gamma_1}=\tfrac {\rho_0}{c^2}\left({\textstyle\int}_\Omega 1-c^2{\textstyle\int}_{\Gamma_1}v^*\right)=\tfrac {\rho_0}{c^2}L_1V^*,\\
[U,V^*]_{\cal{H}_0}=&\tfrac {\rho_0}{c^2}{\textstyle\int}_\Omega w+[v,v^*]_{\Gamma_1}=\tfrac {\rho_0}{c^2}\left({\textstyle\int}_\Omega w-c^2{\textstyle\int}_{\Gamma_1}v\right)=\tfrac {\rho_0}{c^2}L_1U,
\end{align}

so the  projections $\Pi_{V_0}:\cal{H}_0\to V_0$  and $\Pi_{\cal{H}_{01}}:\cal{H}_0\to \cal{H}_{01}$ are  given by
\begin{equation}\label{5.33}
\Pi_{V_0}U=\tfrac{L_1U}{L_1V^*}V^*,\qquad\text{and}\quad \Pi_{\cal{H}_{01}}U=U-\tfrac{L_1U}{L_1V^*}V^*.
\end{equation}
Since, by \eqref{5.31}, $L_1V^*\not=0$, these  projections trivially extend, with the same laws, to the ones associated to the second splitting
in \eqref{5.30}, that is to $\Pi_{V_0}:\cal{H}\to V_0$  and $\Pi_{\cal{H}_{1}}:\cal{H}\to \cal{H}_1$, given by \eqref{5.33} provided $\cal{H}_{01}$ is replaced by $\cal{H}_1$. Trivially $\Pi_{V_0}\cdot\Pi_0=\Pi_{V_0}$, where we recall that $\Pi_0(u,v,w,z)=(u-\fint_\Omega u,v,w,z)$, while $\Pi_{V_0}$ is given in \eqref{5.33}.

Since, recalling vanishing velocity solutions of type (s2--i)' and Lemma~\ref{lemma5.5}, both $\cal{H}_{01}$ and $V_0$ are invariant under the flow of $\{T_0(t),\,t\in\R\}$ and $T_0(t)_{|V_0}=I$, the group can then be written as
\begin{equation}\label{5.35}
  T_0(t)=T_{01}(t)\cdot \Pi_{\cal{H}_{01}}+\Pi_{V_0}\qquad\text{for all $t\in\R$,}
\end{equation}
where $\{T_{01}(t),\,t\in\R\}$ denotes the subspace group induced by $\{T_0(t),\,t\in\R\}$ on $\cal{H}_{01}$ (see \cite[Chapter I, 5.13 p.~43 and Chapter II, 2.4 p.61]{EngelNagel}),
having generator $-A_{01}$, where $A_{01}: D(A_{01})\subset\cal{H}_{01}\to \cal{H}_{01}$ is given by
\begin{equation}\label{5.36}
 D(A_{01})=D(A)\cap\cal{H}_{01},\qquad A_{01}={A_0}_{|\cal{H}_{01}}.
\end{equation}
Moreover, denoting by $\{T_1(t),\,t\in\R\}$ the subspace group induced by $\{T(t),\,t\in\R\}$ on $\cal{H}_1$ and using the vanishing velocity solutions of type (s2)' ,
for all $U_0\in\cal{H}$ and $t\in\R$ we have
\begin{equation}\label{5.37}
  T(t)[U_0]=T_1(t)\cdot \Pi_{\cal{H}_1}[U_0]+\tfrac{L_1U_0}{L_1 V^*} \left(1_{\cal{H}}\,t+V^*\right) \qquad\text{for all $t\in\R$,}
\end{equation}
where $1_{\cal{H}}=(1,0,0,0)$ as in \eqref{5.2}.

\subsection{The case $\kappa\equiv 0$} In this case $[\cdot,\cdot]_{\cal{H}_0}$ is not an inner product and, by Lemma~\ref{lemma5.3},
$\cal{N}_0=\C_2:=\C(0,1,0,0)$. Since trivially $\cal{H}_1\cap \C_2=\cal{H}_{01}\cap \C_2=\{0\}$, by the same Lemma,  $[\cdot,\cdot]_{\cal{H}_0}$ restricts to an inner product  and $\n\cdot\n_{\cal{H}_0}$ is an equivalent norm on $\cal{H}_{01}$. We equip $\cal{H}_{01}$ with it.

Setting in this case
\begin{equation}\label{5.38}
  V^*=(0,1,0,0),\qquad\text{and}\quad V_0=\C V^*=\C_2,
\end{equation}
since $L_1 V^*=-c^2\cal{H}^{N-1}(\Gamma_1)<0$, we still have the  splittings \eqref{5.30}, with projections still given by \eqref{5.33}. They are still $[\cdot,\cdot]_{\cal{H}}$-orthogonal, since $V^*\in\cal{N}_0$.
Moreover, also in this case, by vanishing velocity solutions of type (s2--ii)' and Lemma~\ref{lemma5.5}, both $\cal{H}_{01}$ and $V_0$ are invariant for $\{T_0(t),\,t\in\R\}$ and $T_0(t)_{|V_0}=I$, so \eqref{5.35}--\eqref{5.36} still hold. By the different forms of vanishing velocity solutions of type (s2)' in this case the group splitting \eqref{5.37} is replaced, in this case, by the  simpler one
\begin{equation}\label{5.39}
 T(t)=T_1(t)\cdot \Pi_{\cal{H}_1}+\Pi_{V_0}.
\end{equation}
\subsection{Unified treatment of the cases $\kappa\not\equiv0$ and $\kappa\equiv0$} By the analysis made in previous two sections, and in particular by the group splittings  \eqref{5.35}, \eqref{5.37} and \eqref{5.39}, it is then clear than we can treat the two cases in an unified way, with $V^*$ given by \eqref{1.10} and unifying  \eqref{5.37} and \eqref{5.39} to
\begin{equation}\label{5.3739}
 T(t)[U_0]=T_1(t)\cdot \Pi_{\cal{H}_1}[U_0]+\tfrac{L_1U_0}{L_1 V^*} \left(s_\kappa 1_{\cal{H}}\,t+V^*\right) \qquad\text{for all $t\in\R$,}
\end{equation}
where $s_\kappa:=\text{sgn\,}\|\kappa\|_{\infty,\Gamma_1}$, that is, $s_\kappa=0$ if $\kappa=0$  a.e.\ while $s_\kappa=1$ otherwise.
 Moreover the analysis reduces to study the behavior of trajectories of the group $\{T_{01}(t),\,t\in\R\}$.

In the sequel we shall endow it with the restriction of  $[\cdot,\cdot]_{\cal{H}_0}$, which is an inner product on it with associated norm equivalent to $\|\cdot\|_{\cal{H}_0}$ and  $\|\cdot\|_{\cal{H}}$.
The first main step in our analysis is  the following result.
\begin{thm}[\bf Character  of $\mathbf{\{T_{01}(t),\,t\in\R\}}$ when $\Gamma_1$ is connected] \label{lemma5.6}
When \\ assumptions (A0--4) hold and $\Gamma
_1$ is connected, then $\{T_{01}(t),\,t\in\R\}$ is a strongly continuous group with compact resolvent and $\{T_{01}(t),\,t\ge 0\}$ is a contraction semigroup.
Moreover
\renewcommand{\labelenumi}{{\roman{enumi})}}
\begin{enumerate}
\item if  $\delta\not\equiv 0$, then $\{T_{01}(t),\,t\ge 0\}$ is strongly stable;
\item if $\delta\equiv 0$, then $\{T_{01}(t),\,t\in\R\}$ is unitary.
\end{enumerate}
\end{thm}
\begin{proof}
The fact that $\{T_{01}(t),\,t\in\R\}$ is a strongly continuous group is trivial. By \cite[Chapter II, Theorem~1.10--(i)~p.~55]{EngelNagel} its resolvent is compact being the restriction of the resolvent of $A_0$, which is compact by Lemma~\ref{lemma5.1}--ii). Moreover, by Lemma~\ref{lemma5.1}--iii) the operator $A_{01}$ is accretive. Since $-A_{01}$ generates a strongly continuous semigroup, by \cite[Chapter II, Theorem~1.10--(i)~p.~55 and Proposition~3.14~p.~82]{EngelNagel} its resolvent set contains $(0,\infty)$ and then by \cite[Chapter II, Corollary~3.20~p.~86]{EngelNagel} this semigroup is contractive.

Let now consider the alternative given in i)--ii). When $\delta\not\equiv 0$, using the previous statement, the semigroup $\{T_{01}(t),\,t\ge 0\}$ is bounded and, by \cite[Chapter IV, Corollary~1.19~p.~248]{EngelNagel}, the spectrum of $A_{01}$ reduces to its point spectrum. Moreover, by Lemma~\ref{lemma5.1}--(v), $A_0$ and {\em a fortiori} $A_{01}$ have no purely imaginary eigenvalues. By \eqref{5.24} and \eqref{5.30} we also have $\text{Ker }A_{01}=\{0\}$. The conclusion of i) then follows by applying a consequence of the Arendt, Batty, Lyubich and V\~{u} Theorem (see \cite[Chapter V, Corollary~2.22 p.~327]{EngelNagel}). When $\delta\equiv 0$ by Lemma~\ref{lemma5.1}--iii) also the operator  $A_{01}$ is accretive. Since $A_{01}$ generates the semigroup $\{T_{01}(-t), t\ge 0\}$, using the same arguments as before also $\{T_{01}(-t), t\ge 0\}$ is contractive. Since  $T_{01}^{-1}(t)=T_{01}(-t)$ for all $t\in\R$ the proof is complete.
\end{proof}
Theorem~\ref{lemma5.6} shows that the behavior of trajectories of the groups $\{T_{01}(t),\,t\in\R\}$, $\{T_1(t),\,t\in\R\}$, $\{T_0(t),\,t\in\R\}$ and $\{T(t),\,t\in\R\}$ depend (as expected) on the presence (or on the lack thereof) of the damping term.

In particular, in the first case, simply by combining Theorem~\ref{lemma5.6}--i) with our preliminary analysis we get the following result, which includes  Theorem~\ref{theorem1.3} and Corollary~\ref{corollary1.1}.
\begin{thm}[\bf Stability vs. instability when $\mathbf{\Gamma_1}$ is connected and $\mathbf{\delta\not\equiv 0}$] \label{theorem5.1}

When assumptions (A0--4) hold, $\Gamma_1$ is connected and $\delta\not\equiv 0$ the following conclusions hold.
\renewcommand{\labelenumi}{{\Roman{enumi}.}}
\begin{enumerate}
\item For any $U_0\in\cal{H}_{01}$ we have $T_{01}(t)[U_0]\to 0$ in $\cal{H}_0$ as $t\to\infty$. Consequently for any $U_0\in\cal{H}_1$ we have $\Pi_{0}T_1(t)[U_0]\to 0$ in $\cal{H}$ as $t\to\infty$ and the conclusion of Theorem~\ref{theorem1.3} holds.
\item For any $U_0\in\cal{H}_0$ we have $T_0(t)[U_0]\to \Pi_{V_0}U_0$ in $\cal{H}_0$ as $t\to\infty$. Consequently for any $U_0\in\cal{H}$ we have $\Pi_{0}T(t)[U_0]\to \Pi_{V_0}U_0$ in $\cal{H}$ as $t\to\infty$ and the conclusions of Corollary~\ref{corollary1.1} hold.
\end{enumerate}
\end{thm}
\begin{proof} The first statement in part I is simply a repetition of Theorem~\ref{lemma5.6}--i). By \eqref{5.2} we then get that for any $U_0\in\cal{H}_1$ we have $\Pi_{0}T_1(t)[U_0]\to 0$ in $\cal{H}$ as $t\to\infty$.  Hence all limits relations in \eqref{1.12}, but for the limit of $\frac 1t \fint_\Omega u(t)$, hold. Since $u_t(t)\to 0$ in $L^2(\Omega)$ as $t\to\infty$ we get that $\frac d{dt} \fint_\Omega u(t)\to 0$ as $t\to\infty$. An elementary application of the Mean Value Theorem then shows that  $\frac 1t \fint_\Omega u(t)\to 0$ as $t\to\infty$, concluding the proof of part I.

The first statement in part II follows by \eqref{5.35} and the first statement in part I.
For any $U_0\in\cal{H}$, using \eqref{5.2bis} and the first statement  we have $\Pi_{0}T(t)[U_0]=T_0(t)[\Pi_0 U_0]\to \Pi_{V_0}\cdot\Pi_0U_0$ in $\cal{H}$ as $t\to\infty$. By  the already remarked identity $\Pi_{V_0}\cdot\Pi_0=\Pi_{V_0}$ we then get $\Pi_{0}T(t)[U_0]\to \Pi_{V_0}U_0$ in $\cal{H}$ as $t\to\infty$. To get the conclusions of Corollary~\ref{corollary1.1} let $U_0=(u_0,v_0,u_1,v_1)\in\cal{H}$. By the group splittings \eqref{5.37} and \eqref{5.39} the weak solution $(u,v)$ of \eqref{1.1} corresponding to initial data $U_0$ is the sum of two solutions $(u^1,v^1)$ and $(u^2,v^2)$ of it, the first one corresponding to initial data $\Pi_{\cal{H}_1}U_0$ and the second one being the  vanishing velocity solution corresponding to initial data $\Pi_{V_0}U_0$. Hence, by part I, $(u^1,v^1)$ satisfies \eqref{1.12}, while
$$(u^2(t),v^2(t))=
\begin{cases}
\tfrac{\int_\Omega u_1-c^2\int_{\Gamma_1}v_0}{|\Omega|-c^2\int_{\Gamma_1}v^*}(t,v^*)\quad&\text{if $\kappa\not\equiv 0$,}\\
\tfrac{c^2\int_{\Gamma_1}v_0-\int_\Omega u_1}{c^2\cal{H}^{N-1}(\Gamma_1)}(0,1)\quad&\text{if $\kappa\equiv 0$.}
\end{cases}$$
Summing $(u^1,v^1)$ and $(u^2,v^2)$ we get Corollary~\ref{corollary1.1}.
\end{proof}
In the undamped case we have the following result, which trivially implies Theorem~\ref{theorem1.4} and Corollary~\ref{corollary1.2} in \S~\ref{intro}.
\begin{thm}[\bf Fourier expansions when $\mathbf{\Gamma_1}$ is connected and $\mathbf{\delta\equiv 0}$]\label{theorem5.2} When  (A0--4) hold, $\Gamma_1$ is connected and $\delta\equiv 0$ there is an Hilbert basis $\{W^{0n},n\in\N\}$ of $\cal{H}_{01}$ and a real sequence $(\lambda_n)_n$ satisfying Theorem~\ref{theorem1.4}--ii) such that $A_{01}W^{0n}=-\mathfrak{i}\lambda_nW^{0n}$ for all $n\in\N$.

Moreover $\{W^{0n},n\in\N\}$ can be chosen in such a way that
\begin{equation}\label{5.43}
  W^{0n}=\Pi_0V^{0n},\qquad V^{0n}=(u^{0n}, -\mathfrak{i} v^{0n}, -\mathfrak{i}\lambda_n u^{0n}, -\lambda_n v^{0n})\quad\text{for all $n\in\N$,}
\end{equation}
where $((u^{0n},v^{0n}))_n$ satisfies Theorem~\ref{theorem1.4}--i).

Finally, denoting $U_0=(u_0,v_0,u_1,v_1)$,  the following expansions hold:
\begin{alignat}2
\label{5.44} T_{01}(t)[U_0]=&\sum_{n=1}^\infty \alpha_n W^{0n} e^{-\mathfrak{i}\lambda_n t}\,&&\text{in $C_b(\R;\cal{H})$ for all $U_0\in\cal{H}_{01}$;}\\
\label{5.45} T_0(t)[U_0]=&\sum_{n=1}^\infty \alpha_n W^{0n} e^{-\mathfrak{i}\lambda_n t}+\frac{L_1U_0}{L_1 V^*}V^*\,&&\text{in $C_b(\R;\cal{H})$ for all $U_0\in\cal{H}_0$;}\\
\label{5.46} T_1(t)[U_0]=&\sum_{n=1}^\infty \alpha_n V^{0n} e^{-\mathfrak{i}\lambda_n t}+\alpha_0 1_{\cal{H}}\,&&\text{in $C_b(\R;\cal{H})$ for all $U_0\in\cal{H}_1$;}\\
\label{5.47} T(t)[U_0]=&\sum_{n=1}^\infty \alpha_n V^{0n} e^{-\mathfrak{i}\lambda_n t}+\alpha_0 1_{\cal{H}}+\frac{L_1U_0}{L_1 V^*}&&\left(s_\kappa 1_{\cal{H}}\,t+V^*\right)
\end{alignat}
for all $U_0\in\cal{H}_{01}$, uniformly in $t\in\R$, where $s_\kappa=\text{sgn\,}\|\kappa\|_{\infty,\Gamma_1}$ and $1_{\cal{H}}$, $V^*$, $L_1$, $\alpha_n$ are respectively given by \eqref{5.2}, \eqref{1.10}, \eqref{5.28} and
\eqref{1.19bis}.
\end{thm}
\begin{proof}
By Theorem~\ref{lemma5.6} and Stone Theorem (see \cite[Chapter II, Theorem~3.24, p.~89]{EngelNagel}) the operator $A_{01}$ is skew--adjoint, i.e., $A_{01}^*=-A_{01}$ and it has compact resolvent. Hence
by applying the standard Spectral Decomposition Theorem for self-adjoint operators with compact resolvent  (see \cite[Theorem~6, p.~38]{dautraylionsvol3}) to the operator $\mathfrak{i}A_{01}$,
 its spectrum reduces to its point spectrum $\Lambda=\{\lambda'_n,n\in \N\}$, each $\lambda_n'$ having finite multiplicity, with $|\lambda_n'|\to\infty$ as $n\to\infty$, $\Lambda\subset \mathfrak{i}\R$, and $\cal{H}_{01}=\bigoplus_{n=1}^\infty \Ker (A_{01}-\lambda_n'I)$, the sum being orthogonal.
By Lemma~\ref{lemma5.4}--iv) we have $\Ker A_{01}=\Ker A_0\cap\cal{H}_{01}=\{0\}$, so $0\not\in\Lambda$. Moreover, since by \eqref{5.4}  $\overline{A_{01}U}=A_{01}\overline{U}$ for all $U\in\cal{H}_{01}$, $\Lambda$ is invariant with respect to conjugation and hence symmetric with respect to the origin. Hence $\Lambda=\{\mathfrak{i}\lambda_n'',n\in\N\}$, where $(\lambda_n'')_n$ is a real sequence such that $\lambda''_{2n+2}=-\lambda_{2n+1}''$ for all $n\in\N$,
$$0<\lambda_1''< \cdots\le \lambda_{2n-1}''<\lambda_{2n+1}''< \cdots,\qquad\qquad \lambda_{2n+1}''\to\infty$$
and, denoting $H_n=\Ker (A_{01}-\mathfrak{i}\lambda_n''I)$ and $\mu_n=\text{dim\,} H_n$, we have $1\le \mu_n<\infty$ and
$\cal{H}_{01}=\bigoplus_{n=1}^\infty H_n$, the sum being orthogonal.

Any choice of an orthonormal basis for each $H_n$ then leads, by repeating $\mu_n$ times each $\lambda_n''$ , to the existence of  a real sequence $(\lambda_n)_n$ and an Hilbert basis $\{W^{0n},n\in\N\}$ as asserted in the first part of the statement. We now claim that this choice can be made in such a way that \eqref{5.43} holds true, with $(u^{0n},v^{0n})_n$ as asserted.

To prove our claim we fix $n\in\N$ and we remark that,  by \eqref{5.4} (here $\delta\equiv 0$), Lemma~\ref{lemma4.6}, \eqref{5.3} and \eqref{5.36},
$\dot U=(\dot u, \dot v, \dot w, \dot z)\in H_n$ if and only if
$\dot u\in H^2(\Omega)$, $\dot v\in H^2(\Gamma_1)$, $\dot w\in H^1(\Omega)$, $\dot z\in H^1(\Gamma_1)$, $\int_\Omega \dot u=0$, $\int_\Omega \dot w=c^2\int_{\Gamma_1}\dot v$ and
\begin{equation}\label{5.50}
\begin{cases}
-\dot w+\fint_\Omega \dot w=\mathfrak{i}\lambda_n\dot u \qquad &\text{in
$\Omega$,}\\
- \dot z =\mathfrak{i}\lambda_n \dot v\qquad
&\text{on
$\Gamma_1$,}\\
-c^2\Delta \dot u=\mathfrak{i}\lambda_n \dot w \qquad &\text{in
$\Omega$,}\\
- \DivGamma (\sigma \nabla_\Gamma \dot v)+\kappa\dot v+\rho_0\dot w =\mathfrak{i}\lambda_n\mu \dot z\qquad
&\text{on
$\Gamma_1$,}\\
\partial_\nu \dot u=\dot z\qquad
&\text{on
$\Gamma_1$,}\\
\partial_\nu \dot u=0 &\text{on $\Gamma_0$.}
\end{cases}
\end{equation}
Hence, setting in correspondence to $\dot U$, $U=(u,v,w,z)$ with $u\in H^2(\Omega)$, $v\in H^2(\Gamma_1)$, $w\in H^1(\Omega)$, $z\in H^1(\Gamma_1)$ given by
\begin{equation}\label{5.51}
  u=\dot u+\tfrac {\mathfrak{i}}{\lambda_n}{\textstyle\fint}_\Omega \dot w,\quad v=\mathfrak{i}\dot v,\quad w=\dot w,\quad z=\mathfrak{i}\dot z,
\end{equation}
we trivially get that
\begin{equation}\label{5.52}
  \dot u= u-\tfrac {\mathfrak{i}}{\lambda_n}{\textstyle\fint}_\Omega w,\quad \dot v=-\mathfrak{i} v,\quad \dot w= w,\quad \dot z=-\mathfrak{i} z,
\end{equation}
and, plugging it into \eqref{5.50}, that the last two components of $U$ are given by
\begin{equation}\label{5.53}
  w=-\mathfrak{i}\lambda_n u,\qquad z=-\mathfrak{i}\lambda_n v
\end{equation}
where the couple $(u,v)$ solves the problem \eqref{1.14}, that is
\begin{equation}\label{5.54}
\begin{cases} -c^2\Delta u=\lambda_n^2u \qquad &\text{in
$\Omega$,}\\
- \DivGamma (\sigma \nabla_\Gamma v)+\kappa v+\rho_0\lambda_n u =\mu \lambda_n^2v\qquad
&\text{on $\Gamma_1$,}\\
\partial_\nu u=0 &\text{on $\Gamma_0$,}\\
\partial_\nu u=-\lambda_nv\qquad
&\text{on
$\Gamma_1$.}
\end{cases}
\end{equation}
Conversely, setting
$$H_n'=\{(u,v,w,z)\in H^2(\Omega)\times H^2(\Gamma_1)\times H^1(\Omega)\times H^1(\Gamma_1): \text{\eqref{5.53}--\eqref{5.54} hold}\},$$
and for any  $U\in H_n'$ setting $\dot U$ given by \eqref{5.52}, trivially $\dot U$ satisfies \eqref{5.50}, $\int_\Omega \dot u=0$
and, integrating by parts and recalling that $\lambda_n\not=0$, $\int_\Omega \dot w=c^2\int_{\Gamma_1}\dot v$.

Hence the operator $Q_n\in\cal{L}(H_n,H_n')$ defined by $Q_n\dot U=U$, where $U$ is given by \eqref{5.51}, is a bijective isomorphism, with inverse $Q_n^{-1}$ given by $Q_n^{-1}U=\dot U$, where $\dot U$ is given by \eqref{5.52}.
Next we remark that, by \eqref{5.53}, the restricted projection $\pi_n\in\cal{L}(H'_n,H^2(\Omega)\times H^2(\Gamma_1))$ given by $\pi_n(u,v,w,z)=(u,v)$ is injective and, denoting by $H_n''$ its range, we have $H_n''=\{(u,v)\in H^2(\Omega)\times H^2(\Gamma_1): \text{\eqref{5.54} holds}\}$, its inverse $\pi_n^{-1}$
being given by $\pi_n^{-1}(u,v)=(u,v,-\mathfrak{i}\lambda_n u, -\mathfrak{i}\lambda_nv)$.
Consequently $S_n=\pi_n\cdot Q_n$ is a bijective isomorphism between $H_n$ and $H_n''$, with inverse $S_n^{-1}$ being given by
\begin{equation}\label{5.55}
S_n^{-1}(u,v)= (u-{\textstyle\fint}_\Omega u, -\mathfrak{i}v, -\mathfrak{i}\lambda_n u, -\mathfrak{i}\lambda_n v)=\Pi_0(u, -\mathfrak{i}v, -\mathfrak{i}\lambda_n u, -\mathfrak{i}\lambda_n v).
\end{equation}
Next we note that, by denoting by $H^2_\R(\Omega)$ and $H^2_\R(\Gamma_1)$ the real versions of the Sobolev spaces under consideration and by $V^\C$ the standard complexification of a real vector space $V$ (see \cite{roman}) we trivially have $[H^2_\R(\Omega)\times H^2_\R(\Gamma_1) ]^\C\simeq H^2(\Omega)\times H^2(\Gamma_1)$. Moreover $H_n''$ is closed with respect to conjugation  by \eqref{5.54}. Hence, by \cite[Proposition~7.4.1., p.~147]{roman} we have
$H_n''=H_{n,\R}''+\mathfrak{i} H_{n,\R}''$, where $H_{n,\R}''=\{(u,v)\in H_\R^2(\Omega)\times H_\R^2(\Gamma_1): \text{\eqref{5.54} holds}\}$ and
$\text{dim\,}_\R H_{n,\R}''=\text{dim\,} H_n''=\text{dim\,} H_n=\mu_n$.

Next, since $S_n$ is a bijective isomorphism, $(\cdot,\cdot)_{H_n''}:=[S_n^{-1}(\cdot),S_n^{-1}(\cdot)]_{\cal{H}_0}$
is an inner product on $H_n''$ and, using \eqref{5.6} and \eqref{5.55}, it is straightforward to check that its restriction to
$H_{n,\R}''$ is a real inner product on it. Hence, by the standard Gram--Schmidt orthonormalization process, we can construct an orthonormal basis $\{(u_{nj},v_{nj}), j=1,\ldots, \mu_n\}$ of $H_{n,\R}''$ with respect to it, which is also an orthonormal basis (with complex scalars) of $H_n''$ with respect to $(\cdot,\cdot)_{H_n''}$. Hence, by construction, $\{S_n^{-1}(u_{nj},v_{nj}), j=1,\ldots, \mu_n\}$ is an orthonormal basis of $H_n$ of the form \eqref{5.43}, with $(u_{nj},v_{nj})$ satisfying \eqref{5.54} and real-valued, so proving our claim.

To prove the expansions \eqref{5.44}--\eqref{5.47} we preliminarily remark that the coefficients $\alpha_n$ in \eqref{1.19bis} are nothing but $\alpha_n=[U_0,W^{0n}]_\cal{H}=[U_0,V^{0n}]_\cal{H}$ for all $U_0\in\cal{H}$ and $n\in\N$. Hence, since $-A_{01}$ generates $\{T_{01}(t), t\in\R\}$ and
$A_{01}W^{0n}=-\mathfrak{i}\lambda_nW^{0n}$ for all $n\in\N$, we immediately get \eqref{5.44}, the series being convergent in $C_b(\R;\cal{H})$ by Parseval's Identity.
By combining the group splitting \eqref{5.35}, also recalling \eqref{5.33},  with \eqref{5.44}, we then get that for all $U_0\in \cal{H}_0$ and $t\in\R$ we have
$$T_0(t)[U_0]= T_{01}(t)[\Pi_{\cal{H}_{01}}U_0]+\Pi_{V_0}U_0=\sum_{n=1}^\infty [\Pi_{\cal{H}_{01}}U_0,W^{0n}]_{\cal{H}_0}W^{0n}e^{-\mathfrak{i}\lambda_n t}+\tfrac{L_1U_0}{L_1 V^*}V^*,$$
the series being convergent in $C_b(\R;\cal{H})$. Since the first splitting in \eqref{5.30} is $[\cdot,\cdot]_{\cal{H}_0}$-orthogonal we have
$[\Pi_{\cal{H}_{01}}U_0,W^{0n}]_{\cal{H}_0}=[U_0,W^{0n}]_{\cal{H}_0}=[U_0,W^{0n}]_\cal{H}$ for all $n\in\N$, so proving \eqref{5.45}.

The proof of the expansions \eqref{5.46}--\eqref{5.47} is a bit different, since it essentially uses the previously constructed sequence $(u^{0n},v^{0n})_n$.
We preliminarily remark that for any $U_0\in\cal{H}_1$ the series appearing in \eqref{1.19bis}$_2$ is nothing but
$\sum_{n=1}^\infty [U_0,W^{0n}]_\cal{H}{\textstyle\fint}_\Omega u^{0n}$, and we claim that it converges in $\C$. Being $\Omega$ bounded this convergence is equivalent to the convergence in $H^1(\Omega)$, which will be used in the sequel.

To prove our claim we remark that, as the second splitting  in \eqref{5.30} is $[\cdot,\cdot]_{\cal{H}_0}$-orthogonal, we have
$[U_0,W^{0n}]_\cal{H}=[\Pi_0 U_0,W^{0n}]_{\cal{H}_0}$ for all $n\in\N$, with $\Pi_0U_0\in \cal{H}_{01}$. Hence, by Parseval's Identity,
$\sum_{n=1}^\infty [U_0,W^{0n}]_\cal{H}^2<\infty$. Since $\lambda_n\ge \lambda_1$ for all $n\in\N$ then  $\sum_{n=1}^\infty \lambda_n^{-2}[U_0,W^{0n}]_\cal{H}^2<\infty$, so the series $\sum_{n=1}^\infty [U_0,W^{0n}]_\cal{H}W^{0n}/\lambda_n$ converges in $\cal{H}_{01}$.
In particular, by \eqref{5.43}, then $\sum_{n=1}^\infty [U_0,W^{0n}]_\cal{H}u^{0n}$ converges in $L^2(\Omega)$. Since $u\mapsto\fint_\Omega u$ is bounded in
$L^2(\Omega)$ our claim is proved.

Now, using \eqref{1.14}, it is trivial to check that $(u^n,v^n)$ defined in \eqref{1.13} is for any $n\in\N$ the strong  solution of problem \eqref{1.1} corresponding to the initial data $V^{0n}$ in \eqref{5.43}, that is
\begin{equation}\label{5.56}
  T_1(t)[V^{0n}]=V^{0n}e^{-\mathfrak{i}\lambda_n t}\qquad\text{for all  $t\in\R$ and $n\in\N$.}
\end{equation}
Now using the orthogonal splitting \eqref{5.30}$_1$ any $U_0\in\cal{H}_1$ can be written as $U_0=\fint_\Omega u_0 1_{\cal{H}}+\Pi_0U_0$, with $\Pi_0U_0\in\cal{H}_{01}$, and we have $[\Pi_0 U_0,W^{0n}]_{\cal{H}_{01}}=[U_0,W^{0n}]_\cal{H}$ for all $n\in\N$.
Then, also using \eqref{5.43} we have
\begin{equation}\label{5.57}
\begin{aligned}
  U_0=&{\textstyle\fint}_\Omega u_0 1_{\cal{H}}+\sum_{n=1}^\infty [U_0,W^{0n}]_\cal{H}W^{0n}\\
  =&\left({\textstyle\fint}_\Omega u_0- \sum_{n=1}^\infty{\textstyle\fint}_\Omega u^{0n}\right)1_{\cal{H}}+\sum_{n=1}^\infty [U_0,W^{0n}]_\cal{H} V^{0n}.
\end{aligned}
\end{equation}
Since trivially $T_1(t)[1_\cal{H}]=1_\cal{H}$, by combining \eqref{5.56} and \eqref{5.57} we get \eqref{5.46}, the series being convergent in $C_b(\R;\cal{H})$
by the analogous property proved for \eqref{5.44}.
Using \eqref{5.3739} and \eqref{5.46} we thus have, for all $U_0\in \cal{H}$ and $t\in\R$,
\begin{align*}
T(t)[U_0]=&\left({\textstyle\fint}_\Omega u_0- \sum_{n=1}^\infty{\textstyle\fint}_\Omega u^{0n}\right)1_{\cal{H}}+\sum_{n=1}^\infty [\Pi_{\cal{H}_1}U_0,W^{0n}]_\cal{H} V^{0n}e^{-\mathfrak{i}\lambda_n t}\\
+&\tfrac{L_1U_0}{L_1 V^*} \left(s_\kappa 1_{\cal{H}}\,t+V^*\right),
\end{align*}
the series being convergent in $C_b(\R;\cal{H})$. Since the splitting \eqref{5.30}$_2$ is $[\cdot,\cdot]_\cal{H}$-orthogonal
in the last formula we have $[\Pi_{\cal{H}_1}U_0,W^{0n}]_\cal{H}=[U_0,W^{0n}]_\cal{H}$, so proving \eqref{5.47} and concluding the proof.
\end{proof}
\section{The general case:  $\Gamma_1$ possibly disconnected} \label{subsection6.3}
When $\Gamma_1$ is disconnected i.e., when $\mathfrak{n}\ge 2$, the analysis is more involved and we get different results depending on
the number $\mathfrak{n}_0$ of connected components where $\kappa$ vanishes, already introduced in \eqref{5.13}, and also
on the number $\mathfrak{n}_{00}$ of them where  $(\kappa,\delta)$  does, i.e., on
\begin{equation}\label{5.58}
\mathfrak{n}_{00}=\#\cal{C}^{00}(\Gamma_1), \quad\text{where}\quad \cal{C}^{00}(\Gamma_1)=\{\Gamma'\in \cal{C}^0(\Gamma_1): \delta_{|\Gamma'}\equiv 0\}.
\end{equation}
Trivially $\mathfrak{n}_{00}\le \mathfrak{n}_0\le \mathfrak{n}$.  We shall separately consider the two main cases $\mathfrak{n}_0\le 1$ and
$\mathfrak{n}_0\ge 2$. Clearly the case $\mathfrak{n}=1$ considered in \S~\ref{subsection6.2} is included in the first one.
\subsection{The case $\mathfrak{n}_0\le 1$} In this case we can essentially repeat the analysis made in \S~\ref{subsection6.2},  the cases $\mathfrak{n}_0=0$ and $\mathfrak{n}_0=1$ respectively corresponding to the cases $\kappa\not\equiv 0$ and  $\kappa\equiv 0$, with few adjustments indicated in the sequel.

When $\mathfrak{n}_0=0$  no changes are needed with respect to the case $\kappa\not\equiv 0$, since we can still define $V^*$ as in \eqref{5.25} and get the orthogonal splittings \eqref{5.30},
define the operator $A_{01}$ and the group $\{T_{01}(t),t\in\R\}$ and get the group reduction \eqref{5.37}.

When $\mathfrak{n}_0=1$ let us denote by $\Gamma_1^1$ the unique element of $\cal{C}^0(\Gamma_1)$. By Lemma~\ref{lemma5.3} we have $\cal{N}_0=\C(0,\chi_{\Gamma_1^1}, 0,0)$, and also in this case $\cal{H}_1\cap \cal{N}_0=\{0\}$, so $[\cdot,\cdot]_{\cal{H}_0}$ restricts
to an inner product on $\cal{H}_{01}$ inducing an equivalent norm. Hence, defining in this case $V^*=(0,\chi_{\Gamma_1^1}, 0,0)$ (so modifying \eqref{5.38}$_1$), since $L_1V^*=-c^2\cal{H}^{N-1}(\Gamma_1^1)<0$, we still have the orthogonal splittings \eqref{5.30} with projections given by \eqref{5.33}. Hence we also get \eqref{5.39}.

Hence unifying the two cases we still get \eqref{5.3739}, and by repeating the proofs of Theorems~\ref{lemma5.6}, \ref{theorem5.1} and \ref{theorem5.2} we get
\begin{thm}[The case $\mathfrak{n}_0\le 1$]\label{theorem5.3} When (A0--4) hold and $\mathfrak{n}_0\le 1$ the conclusions of Theorems~\ref{lemma5.6}, \ref{theorem5.1}, \ref{theorem5.2}, and consequently also those of Theorems~\ref{theorem1.3}--\ref{theorem1.4} and Corollaries~\ref{corollary1.1}--\ref{corollary1.2}, continue to hold provided we respectively mean the cases $\kappa\not\equiv 0$ and $\kappa\equiv 0$ as
 $\mathfrak{n}_0=0$ and $\mathfrak{n}_0=1$, and we replace $\frac{c^2\int_{\Gamma_1}v_0-\int_\Omega u_1}{c^2\cal{H}^{N-1}(\Gamma_1)}$ with
 $\frac{c^2\int_{\Gamma_1}v_0-\int_\Omega u_1}{c^2\cal{H}^{N-1}(\Gamma_1^1)}\chi_{\Gamma_1^1}$ in Corollaries~\ref{corollary1.1}--\ref{corollary1.2}.
\end{thm}
\subsection{The case $\mathfrak{n}_0\ge 2$} The  aim of this section it to complete the  mathematical study of \eqref{1.1}. We are also going to see that in this case
the physical model may become geometrically inconsistent.

We shall denote by $\Gamma_1^1,\ldots,\Gamma_1^{\mathfrak{n}_0}$ the elements of $\cal{C}^0(\Gamma_1)$. Moreover, when $\mathfrak{n}_{00}\ge 1$ we shall suppose, without restriction, that $\cal{C}^{00}(\Gamma_1)=\{\Gamma_1^1,\ldots,\Gamma_1^{{\mathfrak{n}_{00}}}\}$.

When $\mathfrak{n}_0\ge 2$ we have $\cal{N}_{01}:=\cal{N}_0\cap\cal{H}_{01}\not=\{0\}$. More precisely, by Lemma~\ref{lemma5.3},
\begin{equation}\label{N01}
\cal{N}_{01}=\cal{N}_0\cap\cal{H}_1=\{0\}\times\cal{N}_1(\Gamma_1)\times\{0\}\times\{0\},
\end{equation}
 where
\begin{equation}\label{calN1}
\cal{N}_1(\Gamma_1)=\left\{\sum_{i=1}^{\mathfrak{n}_0}\beta_i\chi_{\Gamma_1^i}\quad \beta_1,\ldots,\beta_{\mathfrak{n}_0}\in\C, \quad
\sum_{j=1}^{\mathfrak{n}_0}\cal{H}^{N-1}(\Gamma_1^i)\,\beta_i=0 \right\}.
\end{equation}
 By Lemma~\ref{lemma5.4} $\cal{N}_{01}\subset \Ker  A_0\subset\Ker  A$, so $\cal{N}_{01}$ is a subspace of dimension $\mathfrak{n}_0-1$ consisting of vanishing velocity solutions of type (s2--ii)' and fixed points of $\{T_{01}(t),t\in\R\}$.
 Hence, when $\delta\not\equiv 0$, the stability question can only consists in proving that all solution converge to elements of $\cal{N}_{01}$.
 We shall discuss the inconsistency of these vanishing velocity solutions with the physical model in \S~\ref{section7.3}.

Recalling our general strategy, explained above, since  $\cal{N}_0$ has dimension $\mathfrak{n}_0$,  we are looking for $\mathfrak{n}_0-1$ invariants of the system, independent on $L_1$. They are given by the following result.
\begin{lem}\label{lemma5.7} Let (A0--4) hold and $\mathfrak{n}_0\ge 2$. Then for any weak solution $(u,v)$ of \eqref{1.1}  and $i=2,\ldots,\mathfrak{n}_0$
we have
\begin{multline}\label{5.59}
\int_{\Gamma_1}[\mu v_t(t)+\delta v(t)+\rho_0 u(t)]\left[\frac {\chi_{\Gamma_1^1}}{\cal{H}^{N-1}(\Gamma_1^1)}-\frac {\chi_{\Gamma_1^i}}{\cal{H}^{N-1}(\Gamma_1^i)}\right]\\
=\int_{\Gamma_1}[\mu v_1+\delta v_0+\rho_0 u_0]\left[\frac {\chi_{\Gamma_1^1}}{\cal{H}^{N-1}(\Gamma_1^1)}-\frac {\chi_{\Gamma_1^i}}{\cal{H}^{N-1}(\Gamma_1^i)}\right]
\end{multline}
for all $t\in\R$. Consequently, setting for $i=2,\ldots,\mathfrak{n}_0$  the functionals $L_i\in\cal{H}'$ by
\begin{equation}\label{5.60}
  L_i(u,v,w,z)=\int_{\Gamma_1}(\mu z+\delta v+\rho_0 u)\left[\frac {\chi_{\Gamma_1^1}}{\cal{H}^{N-1}(\Gamma_1^1)}-\frac {\chi_{\Gamma_1^i}}{\cal{H}^{N-1}(\Gamma_1^i)}\right]
\end{equation}
and the subspaces of $\cal{H}$
\begin{equation}\label{5.61}
\cal{H}_{\mathfrak{n}_0}=\bigcap_{i=1}^{\mathfrak{n}_0}\Ker L_i,\qquad \cal{H}_{0,\mathfrak{n}_0}=\cal{H}_{\mathfrak{n}_0}\cap\cal{H}_0,
\end{equation}
$\cal{H}_{\mathfrak{n}_0}$ and $\cal{H}_{0,\mathfrak{n}_0}$ are respectively invariant for $\{T(t),\,t\in\R\}$ and $\{T_0(t),\,t\in\R\}$.
\end{lem}
\begin{proof} Since the $\Gamma_1^i$'s  are compact, by \eqref{4.5}--\eqref{4.6} and \eqref{2.24}, for all $(u_0,v_0,u_1,v_1)\in D(A)$ and $i=1,\ldots,\mathfrak{n}_0$ we have
$$\frac d{dt}\int_{\Gamma_1^i}\mu z(t)+\delta v(t)+\rho_0 u(t)=\int_{\Gamma_1^i}\DivGamma (\sigma \nabla_\Gamma v(t))=0$$
\eqref{5.59} trivially follows. By density then \eqref{5.59} follows for all data in $\cal{H}$, and consequently $\cal{H}_{\mathfrak{n}_0}$ is invariant for $\{T(t),\,t\in\R\}$. To show the second asserted invariance we remark that, using the just proved invariance of
$\cal{H}_{\mathfrak{n}_0}$, for any $(u_0,v_0,u_1,v_1)\in\cal{H}_{\mathfrak{n}_0}$  and $i=2,\ldots,\mathfrak{n}_0$ we have
\begin{multline}L_i(\Pi_0 T(t))=L_i(T(t))-L_i({\textstyle\fint}_\Omega u(t),0,0,0)\\
=-\rho_0{\textstyle\fint}_\Omega u(t)\int_{\Gamma_1}\left[\frac {\chi_{\Gamma_1^1}}{\cal{H}^{N-1}(\Gamma_1^1)}-\frac {\chi_{\Gamma_1^i}}{\cal{H}^{N-1}(\Gamma_1^i)}\right]=0,
\end{multline}
which by \eqref{5.2bis} concludes the proof.
\end{proof}
The second step in our strategy  then consists in checking if $\cal{H}_{0,\mathfrak{n}_0}\cap\cal{N}_0=\{0\}$  and, in this case, if $\cal{H}_0=\cal{H}_{0,\mathfrak{n}_0}\oplus\cal{N}_0$. The following result answers the question.
\begin{lem}\label{lemma5.8}
Let (A0--4) hold and $\mathfrak{n}_0\ge 2$. Then $\cal{H}_{0,\mathfrak{n}_0}\cap\cal{N}_0=\{0\}$  if and only if $\mathfrak{n}_{00}\le 1$. Moreover, in this case, we have the $[\cdot,\cdot]_\cal{H}$-orthogonal splittings
\begin{equation}\label{5.62}
\cal{H}_0=\cal{H}_{0,\mathfrak{n}_0}\oplus\cal{N}_0\qquad \cal{H}=\cal{H}_{\mathfrak{n}_0}\oplus\cal{N}_0.
\end{equation}
Next, introducing for $i,j=1,\ldots,\mathfrak{n}_0$ the nonnegative number
\begin{equation}\label{varid}
d_i=\prod_{\substack{{k=1}\\{k\ne i}}}^{\mathfrak{n}_0}\fint_{\Gamma_1^k}\delta,
\qquad
d_{ij}=\prod_{\substack{{k=1}\\{k\ne i,j}}}^{\mathfrak{n}_0}\fint_{\Gamma_1^k}\delta,
\end{equation}
($d_{ij}=1$ if $\mathfrak{n}_0=2$ and $i\not=j$ and the functional
$\dot L_i, \beta_i\in\cal{H}'$ given by $\dot L_i(u,v,w,z)={\textstyle\fint}_{\Gamma_1^i}\mu z+\delta v+\rho_0 u$ and
\begin{equation}\label{5.63}
  \beta_iU=\frac{c^2\sum_{j=1}^{\mathfrak{n}_0}\cal{H}^{N-1}(\Gamma_1^j)\,d_{ij}\,(\dot L_iU-\dot L_jU)-d_iL_1U}{c^2\sum_{j=1}^{\mathfrak{n}_0}d_j\cal{H}^{N-1}(\Gamma_1^j)},
\end{equation}
the projection operators on $\cal{N}_0$ associated to \eqref{5.62}, both denoted by $\Pi_{\cal{N}_0}$, are given by $\Pi_{\cal{N}_0}U=\sum_{i=1}^{\mathfrak{n}_0}(\beta_iU)\,(0,\chi_{\Gamma_1^i},0,0)$.
Finally $\Pi_{\cal{N}_0}\cal{H}_1=\cal{N}_{01}$.
\end{lem}
\begin{proof} When $\mathfrak{n}_{00}\ge 2$, by Lemma~\ref{lemma5.3} and \eqref{5.61}, for any $a,b\in\C$ such that $\cal{H}^{N-1}(\Gamma_1^1)+b \cal{H}^{N-1}(\Gamma_1^2)=0$ we have $(0,a\chi_{\Gamma_1^1}+b\chi_{\Gamma_1^2},0,0)\in \cal{H}_{\mathfrak{n},0}\cap\cal{N}_0\subset \cal{H}_{\mathfrak{n}}\cap\cal{N}_0$, so making these spaces nontrivial.

When $\mathfrak{n}_{00}\le 1$, by Lemma~\ref{lemma5.3} and \eqref{5.61}, the proof of \eqref{5.62}$_2$ and of the existence of the projection operator $\Pi_{\cal{N}_0}:\cal{H}\to\cal{N}_0$  consists in proving that for any $U\in\cal{H}$ there is a unique $(\beta_1,\ldots,\beta_{\mathfrak{n}_0})\in\C^{\mathfrak{n}_0}$ such that
\begin{equation}\label{5.64}
 L_i\left(U-\sum_{j=1}^{\mathfrak{n}_0}\beta_jV_j\right)=0,\,\,\text{for $i=1,\ldots,\mathfrak{n}_0$, \quad where $V_j=(0,\chi_{\Gamma_1^j},0,0)$.}
\end{equation}
Trivially $L_i=\dot L_1-\dot L_i$ for $i=2,\ldots,\mathfrak{n}_0$ and, denoting by $\delta_{jn}$ the Kronecker symbol,
 $\dot L_i V_j=\delta_{ij}{\textstyle\fint}_{\Gamma_1^i}\delta$  and $L_1(V_j)=-c^2\cal{H}^{N-1}(\Gamma_1^j)$
for $i,j=1,\ldots,\mathfrak{n}_0$. Consequently, denoting for shortness
\begin{equation}\label{5.64bis}
a_i=\cal{H}^{N-1}(\Gamma_1^i),\,\, \delta_i={\textstyle\fint}_{\Gamma_1^i}\delta,\quad -c^2e_0=L_1U,\quad e_i=\dot L_iU\,\,\text{for $i=1,\ldots,\mathfrak{n}_0$},
\end{equation}
we rewrite \eqref{5.64} as the linear system
\begin{equation}\label{5.65}
\sum_{j=1}^{\mathfrak{n}_0}a_j\beta_j=e_0, \qquad \delta_1\beta_1-\delta_i\beta_i=e_1-e_i,\quad\text{for $i=2,\ldots,\mathfrak{n}_0$,}
\end{equation}
where $a_1,\ldots,a_{\mathfrak{n}_0}>0$ and $\delta_2,\ldots,\delta_{\mathfrak{n}_0}>0$ since $\mathfrak{n}_{00}\le 1$. We  solve the last
$\mathfrak{n}_0-1$ equations with respect to $\beta_i$, getting
\begin{equation}\label{5.66}
\beta_i=\frac{\delta_1\beta_1-(e_1-e_i)}{\delta_i},\quad\text{for $i=2,\ldots,\mathfrak{n}_0$.}
\end{equation}
Plugging \eqref{5.66} into the first equation in \eqref{5.65}, multiplying by $d_1$ and taking into account \eqref{varid} we have
$$d_1a_1\beta_1+\sum_{j=2}^{\mathfrak{n}_0}d_{1j}a_j(\delta_1\beta_1-e_1+e_j)=d_1e_0,$$
 that is, using \eqref{varid} again,
\begin{equation}\label{5.67}
\beta_1 \sum_{j=1}^{\mathfrak{n}_0} d_ja_j
=\sum_{j=2}^{\mathfrak{n}_0}d_{1j}a_j(e_1-e_j)+d_1e_0.
\end{equation}
Because $\mathfrak{n}_{00}\le 1$, by \eqref{varid} we have $d_1>0$ and hence $\sum_{j=1}^{\mathfrak{n}_0} d_ja_j>0$: we can then solve \eqref{5.67} with respect to $\beta_1$,  getting
$$\beta_1=\frac{d_1e_0+\sum_{j=2}^{\mathfrak{n}_0}d_{1j}a_j(e_1-e_j)}{\sum_{j=1}^{\mathfrak{n}_0} d_ja_j},$$
which, by \eqref{5.64bis}, is nothing but \eqref{5.63} when $i=1$. Plugging it into \eqref{5.66} and taking into account \eqref{varid} we then get, for $i=2,\ldots,\mathfrak{n}_0$,
\begin{align*}
\beta_i=\quad& \frac{\delta_1\ldots\delta_{\mathfrak{n}_0}e_0+\sum_{j=2}^{\mathfrak{n}_0}d_ja_j(e_1-e_j)-\sum_{j=1}^{\mathfrak{n}_0}d_ja_j(e_1-e_i)}
{\delta_i\sum_{j=1}^{\mathfrak{n}_0} d_ja_j}\\
=&\frac {\delta_1\ldots\delta_{\mathfrak{n}_0}e_0-d_1a_1(e_1-e_i)-\sum_{j=2}^{\mathfrak{n}_0}d_ja_j(e_l-e_i)}
{\delta_i\sum_{j=1}^{\mathfrak{n}_0} d_ja_j}\\
=&\frac {\delta_i e_0+\sum_{j=1}^{\mathfrak{n}_0}a_jd_{ij}(e_i-e_j)}
{\sum_{j=1}^{\mathfrak{n}_0} d_ja_j}.
\end{align*}
Again by \eqref{5.64bis}, this is precisely \eqref{5.63} when $i\ge 2$, so getting \eqref{5.62}$_1$ and the existence of the projection operator $\Pi_{\cal{N}_0}:\cal{H}\to\cal{N}_0$  with the  $\beta_j'$s in \eqref{5.63}. By it we immediately get that
$c^2\sum_{i=1}^{\mathfrak{n}_0}\cal{H}^{N-1}(\Gamma_1^i)\beta_i=-L_1$ in $\cal{H}'$, so $\Pi_{\cal{N}_0}\cal{H}_1=\cal{N}_{01}$.
The restriction to  $\cal{H}_0$ is trivial.
\end{proof}
We shall then distinguish between the two cases $\mathfrak{n}_{00}\le 1<2\le\mathfrak{n}_0$  and $\mathfrak{n}_{00}\ge 2$.

In the first one, in which  clearly $\delta\not\equiv 0$, we get the following result.
\begin{thm}[\bf Boundedness and stability when $\mathfrak{n}_{00}\le 1<2\le\mathfrak{n}_0$]\label{theorem5.4}When (A0--4) hold and $\mathfrak{n}_{00}\le 1<2\le\mathfrak{n}_0$ there is an inner product on $\cal{H}_0$, inducing an equivalent norm on it, with respect to which $\{T_0(t),\,t\ge 0\}$ is a contraction semigroup.

Moreover, for any $U_0\in\cal{H}_0$ we have
$T_0(t)[U_0]\to \Pi_{\cal{N}_0} U_0$ in $\cal{H}_0$ as $t\to\infty$.
Consequently, for any $U_0\in\cal{H}$ we have $\Pi_0T(t)[U_0]\to \Pi_{\cal{N}_0} U_0$ in $\cal{H}_0$ as $t\to\infty$.

Finally, for any $U_0\in\cal{H}_1$, denoting by $(u,v)$ the weak solution of \eqref{1.1} corresponding to the initial data $U_0$, we have
$$\left\{\begin{alignedat}{4}
&\nabla u(t)\to 0\,\, && \text{in $[L^2(\Omega)]^N$}, \quad &\text{or} &\quad u(t)-{\textstyle\fint}_\Omega u(t)\to 0\quad &&\text{in $H^1(\Omega)$},\\
&\tfrac 1t {\textstyle\fint} u(t)\to 0 &&\text{in $\C$}&& \quad  u_t(t)\to 0\,\, &&\text{in $L^2(\Omega)$},\quad \\
& v(t)\to \sum_{i=1}^{\mathfrak{n}_0}(\beta_iU_0)\chi_{\Gamma_1^i}\quad &&\text{in $H^1(\Gamma_1)$}, &&\quad  v_t(t)\to 0\quad &&\text{in $L^2(\Gamma_1)$},
\end{alignedat}\right.
$$
where the $\beta_iU_0'$s are given by \eqref{5.63}.
\end{thm}
The main difference between the conclusions of Theorem~\ref{theorem1.3} (or Theorem~\ref{theorem5.3}) and those of Theorem~\ref{theorem5.4} is given by
the triviality or untriviality of the limit set $\cal{N}_{01}$. We shall discuss in \S~\ref{section7.3} the consequences of this untriviality for the physical model.
\begin{proof}[Proof of Theorem~\ref{theorem5.4}]By Lemma~\ref{lemma5.7}  we can introduce the subspace semigroup $\{T_{0,\mathfrak{n}}(t),\,t\ge 0\}$
induced by $\{T_0(t),\,t\ge 0\}$ on $\cal{H}_{0,\mathfrak{n}_0}$. Its generator is $-A_{0,\mathfrak{n}_0}$, where
$A_{0,\mathfrak{n}_0}:D(A_{0,\mathfrak{n}_0})\subset \cal{H}_{0,\mathfrak{n}_0}\to \cal{H}_{0,\mathfrak{n}_0}$ is given by
\begin{equation}\label{5.71}
D(A_{0,\mathfrak{n}_0})=D(A)\cap  \cal{H}_{0,\mathfrak{n}_0},  \qquad A_{0,\mathfrak{n}_0}={A_0}_{|\cal{H}_{0,\mathfrak{n}_0}}.
\end{equation}
By Lemmas~\ref{lemma5.2}--\ref{lemma5.3} and \ref{lemma5.8} we can endow $\cal{H}_{0,\mathfrak{n}_0}$ with $[\cdot,\cdot]_{\cal{H}_0}$, getting an equivalent norm. Then, by repeating the arguments in the proof of Theorem~\ref{lemma5.6}, we get that $\{T_{0,\mathfrak{n}}(t),\,t\ge 0\}$ is a contraction semigroup and that it is strongly stable, i.e., $T_{0,\mathfrak{n}}(t)[U_0]\to 0$ as $t\to\infty$ for all $U_0\in \cal{H}_{0,\mathfrak{n}_0}$. By Lemma~\ref{lemma5.4}--iv) and the splitting \eqref{5.62} we then get, denoting by $\Pi_{\cal{H}_{0,\mathfrak{n}_0}}$ the projection operator on $\cal{H}_{0,\mathfrak{n}_0}$,  that
\begin{equation}\label{5.72}
T_0(t)=T_{0,\mathfrak{n}}(t)\cdot\Pi_{\cal{H}_{0,\mathfrak{n}_0}}+\Pi_{\cal{N}_0}\qquad\text{for all $t\in\R$}.
\end{equation}
Hence $[\Pi_{\cal{H}_{0,\mathfrak{n}_0}}U,\Pi_{\cal{H}_{0,\mathfrak{n}_0}}V]_{\cal{H}_0}+[\Pi_{\cal{N}_0}U,\Pi_{\cal{N}_0}V]_{\cal{H}}$
defines an inner product on $\cal{H}_0$ which induces an equivalent norm and with respect to which $\{T_0(t),\,t\ge 0\}$ is contractive.
From \eqref{5.72} and from the proved strong stability of $\{T_{0,\mathfrak{n}}(t),\,t\ge 0\}$ we then get that $T_0(t)[U_0]\to \Pi_{\cal{N}_0} U_0$ in $\cal{H}_0$ as $t\to\infty$ for any $U_0\in\cal{H}_0$.
Consequently, for any $U_0\in\cal{H}$ we have $\Pi_0T(t)[U_0]\to \Pi_{\cal{N}_0} U_0$ in $\cal{H}_0$ as $t\to\infty$. By particularizing this conclusion to data in $\cal{H}_1$ the proof is completed.
\end{proof}

We now briefly discuss the case $\mathfrak{n}_{00}\ge 2$, in which our general strategy fails, since we are not able to find out other invariants.
 Theorem~\ref{theorem5.4} is no longer valid in this case, since the semigroup $\{T_0(t),\,t\ge 0\}$ is unbounded.
This conclusion follows from the appearance, in this case, of a new type of trivial solutions. Indeed, by looking for solutions of the form $u(t,x)=u_0(x)$ and $v(t,x)=\sum_{i=1}^{\mathfrak{n}_0}v_i(t)\chi_{\Gamma_1^i}$, suggested by the proof of Lemma~\ref{lemma5.8}, we get the following solutions.
\begin{enumerate}
\item[(s3)] If $\mathfrak{n}_{00}\ge 2$ then, for any $i=2,\ldots, \mathfrak{n}_{00}$, we introduce the unique solution $u_{0i}\in H^2(\Omega)$ of vanishing average  of the Neumann problem
\begin{equation}\label{5.73}
\begin{cases}
\Delta u_{0i}=0 \qquad &\text{in
$\Omega$,}\\
\partial_\nu u_{0i}=\cal{H}^{N-1}(\Gamma_1^i)\qquad
&\text{on
$\Gamma_1^1$,}\\
\partial_\nu u_{0i}=-\cal{H}^{N-1}(\Gamma_1^1)\qquad
&\text{on
$\Gamma_1^i$,}\\
\partial_\nu u_{0i}=0 &\text{on $\Gamma\setminus\left(\Gamma_1^1\cup\Gamma_1^i\right)$,}
\end{cases}
\end{equation}
the compatibility condition being trivially satisfied, and the corresponding function
\begin{equation}\label{5.74}
v_{1i}=\cal{H}^{N-1}(\Gamma_1^i)\chi_{\Gamma_1^1}-\cal{H}^{N-1}(\Gamma_1^1)\chi_{\Gamma_1^i}\in H^2(\Gamma_1).
\end{equation}
Taking initial data $U_{0j}=(u_{0i},0,0,v_{1i})$ the solution of \eqref{1.1} is
$$u(t,x)=u_{0i}(x),\qquad v(t,x)=v_{1i}(x)\,t.$$
\end{enumerate}
Trivially $(u(t),v(t),u_t(t), v_t(t))\in\cal{H}_{01}$ for all $t\in\R$ and all solutions originating from linear combinations of the $\mathfrak{n}_{00}-1$
given (linearly independent) data is a trivial solution, so originating an $\mathfrak{n}_{00}-1$-dimensional space of trivial solutions of type (s3).

Then, as for any $i=2,\ldots,\mathfrak{n}_{00}$
$$\|T_0(t)[U_{0i}]\|_{\cal{H}_0}^2=\|\nabla u_{0i}\|_2^2+\frac 1{\rho_0}\int_{\Gamma_1}|v_{1i}|^2(t^2+\mu)\to\infty\quad\text{as $t\to \infty$,}
$$
the semigroup $\{T_0(t),\,t\ge 0\}$ is unbounded.

A further analysis of problem \eqref{1.1} would still be possible by considering further quotient groups where averages of $v$ on the $\Gamma_1^i$   are partially neglected. Since this approach does not seem to have a  physical interpretation, we shall not proceed in this direction.

To better understand  the behavior of the trivial solutions of type (s3) we evaluate them in the simplest possible case.
\begin{example}\label{esempio}When assumptions (A1) and (A4) hold, $\delta\equiv \kappa\equiv 0$ and $\Omega$ is the spherical shell in $\R^3$ of radii $0<r<R$, that is $\Omega=B_R\setminus \overline{B_r}$, $\Gamma_0=\emptyset$, $\Gamma_1=\partial\Omega$, $\Gamma_1^1=\partial B_r$ and $\Gamma_1^2=\partial B_R$, the trivial solution of type (s3) is
\begin{equation}\label{5.75}
  u(t,x)=\frac{6\pi R^2r^2(R^2-r^2)}{r^3-R^3}+\frac {4\pi R^2r^2}{|x|},\quad v(t,x)=
  \begin{cases}
  \phantom{-}4\pi R^2 t\quad&\text{on $\partial B_r,$}\\
  -4\pi r^2 t\quad&\text{on $\partial B_R$.}
  \end{cases}
\end{equation}
Indeed the form of $v$ in \eqref{5.75} trivially follows from its definition, while to find the form of the  $u=u_{01}$  we remark that, in this case,  \eqref{5.73} is simply written as
\begin{equation}\label{5.76}
\Delta u=0\quad\text{in $\Omega$}, \quad \partial_\nu u=4\pi R^2 \quad\text{on $\partial B_r$},\quad \partial_\nu u=-4\pi r^2 \quad\text{on $\partial B_R$}.
\end{equation}
Now radial harmonic functions are of the form $u=a+b/|x|$, $a,b\in\C$, and both the boundary conditions reduce to $b=4\pi R^2r^2$. Then one finds $a$ by requiring  that
$$0=\int_\Omega\left[a+\frac{4\pi R^2r^2}{|x|}\right]=\frac 43 \pi(R^3-r^3)a+8\pi R^2r^2(R^2-r^2).$$
\end{example}

\chapter{The physical model} \label{section7}
The aim of this chapter is to show that the constraint
\begin{equation}\label{OurConstraint}
\int_\Omega u_t-c^2\int_{\Gamma_1}v=0
 \end{equation}
 appearing in Theorems~\ref{theorem1.3}--\ref{theorem1.4} naturally originates in the physical derivation of problem \eqref{1.1} and should be added as a constitutive part of it.

The discussion of this assertion goes beyond the scope of Mathematical Analysis, whose methods characterized the previous part of the paper. It naturally belongs to the realm of Theoretical Physics. For this reason the inherent discussion, performed in the following \S~\ref{section7.1}--\ref{section7.2}, will be methodologically different from the rest of the paper, since several approximations will be made, making a lack of mathematical rigor unavoidable. Moreover all fields will be tacitly assumed to be as smooth as needed.

 In the final section \S~\ref{section7.3} we shall  give the  physical interpretation of the results in the paper which takes into account this remark. In this section mathematical rigor will be recovered.

Since the physical meaning of the velocity potential $u$ is not self--evident, we shall recall two physical derivations of the acoustic wave equation and the motivation of boundary conditions in both of them. The first one, taken from references in Theoretical Acoustics and Fluid Mechanics, is given in the framework of an Eulerian description of the fluid. The second one, taken from several references in Physics, originates from a  more or less explicit Lagrangian description of it.

\section{The Eulerian approach}\label{section7.1} Referring to \cite[Chapter~6]{morseingard} and \cite[Chapter~VIII]{Landau6}, the main fields describing a fluid are the velocity $\mathbf{v}_{\cal{E}}$ of fluid particles, the hydrostatic pressure $P_{\cal{E}}$ and the fluid density $\rho_{\cal{E}}$, all of them being functions of time $t\in\R$ and position $x\in\Omega$, where $\Omega$ denotes a region of $\R^3$ filled by the fluid, bounded by $\Gamma=\partial\Omega$. Since, in our case, a part of $\Gamma$ may move, $\Omega$ should  depend on time, but this dependence is usually neglected, in this approach, for  small perturbations of $\Gamma$. In other terms the problem should be posed in a time--varying region $\Omega_t$, but for small perturbations one approximates $\Omega_t$ with $\Omega$. This hidden dependence will allow to understand the meaning of our constraint in the Eulerian approach.

The fluid is supposed to be homogeneous, ideal, compressible and to have zero heat conductivity, so thermal effects can be neglected, see \cite[Chapter VI, p.~230 and p.~242]{morseingard}. The fluid is at constant temperature $T_0$. One considers viscosity effects to be negligible, so the only energy involved is the mechanical one related to the fluid compression, and the motion is adiabatic. Finally, in absence of sound, the fluid is supposed to be at rest with uniform density $\rho_0$ and pressure $P_0$,  neglecting gravitational effects.

The two main equations governing the fluid are then the Continuity and the Euler's equations
\begin{equation}\label{7.1}
  \partial_t\rho_{\cal{E}}+\text{div\,}(\rho_{\cal{E}}\mathbf{v}_{\cal{E}})=0,\qquad \rho_{\cal{E}}[\partial_t \mathbf{v}_{\cal{E}}+(\mathbf{v}_{\cal{E}}\cdot\nabla) \mathbf{v}_{\cal{E}}]=-\nabla P_{\cal{E}},
\end{equation}
(the gradient being taken component-wise), respectively derived from conservation of matter and Newton Second Law (see \cite[p.2, (1.2) and (2.1)]{Landau6}.

Acoustic phenomena are physically described as small amplitude perturbations of density and pressure, respectively  denoted as $\rho'_{\cal{E}}$ and $p_{\cal{E}}$ in the sequel, from their equilibrium values $\rho_0$ and $P_0$. That is  $P_{\cal{E}}=P_0+p_{\cal{E}}$ and $\rho_{\cal{E}}=\rho_0+\rho'_{\cal{E}}$, with $\rho'_{\cal{E}}\ll\rho_0$ and $p_{\cal{E}}\ll P_0$.

 The motion  being adiabatic, $p_{\cal{E}}$ and $\rho'_{\cal{E}}$ are linearly related each other through the adiabatic compressibility factor $\kappa'$, depending on $T_0$ and $\rho_0$, that is (see \cite[p.~251]{Landau6} and \cite[(6.1.2)--(6.1.3), p.~230]{morseingard}) we have
 \begin{equation}\label{7.2}
 \rho'_{\cal{E}}=\rho_0\kappa'  p_{\cal{E}}.
 \end{equation}
 It is useful to recall another form of relation \eqref{7.2}, used in other references (see for example \cite[Chapter~5, p.~136]{ElmoreHeald}), which involves the infinitesimal volume strain. Indeed, if an infinitesimal volume $\mathcal{V}$ deforms to $\mathcal{V}+\mathcal{V}'$, by conservation of matter we should have
  $\rho_0\mathcal{V}=(\rho_0+\rho'_{\cal{E}})(\mathcal{V}+\mathcal{V}')$. Hence the volume strain $\theta=\cal{V}'/\cal{V}\ll 1$ can be approximated, up to the first order, by $-\rho'_{\cal{E}}/\rho_0$, that is \eqref{7.2} can be equivalently written as
  \begin{equation}\label{7.2Delio}
  \theta=-\kappa' p_{\cal{E}}
  \end{equation}
   when first order approximations are considered.
    This relation,  known as generalized Hooke's Law (see \cite[p.~136]{ElmoreHeald}), is an experimental law. It is generally expressed by saying  that a small volume strain $\theta$ is linearly related through the bulk modulus $B=1/\kappa'$ to a small pressure perturbation $p_{\cal{E}}$, that is
  \begin{equation}\label{7.3}
 p_{\cal{E}}=-B\theta.
  \end{equation}
Coming back to  the derivation of acoustic wave equation in the Eulerian approach, one then approximates the nonlinear equations \eqref{7.1} by a perturbation analysis, up to the first order, as in \cite[(64.2)--(64.3), p.251]{Landau6}, so getting the linearized equations
\begin{equation}\label{7.4}
\partial_t\rho'_{\cal{E}}+\rho_0\,\text{div\,}\mathbf{v}_{\cal{E}}=0,\qquad \rho_0 \partial_t \mathbf{v}_{\cal{E}}=-\nabla p_{\cal{E}}.
\end{equation}
By combining \eqref{7.2} and \eqref{7.4} to eliminate $\rho'_{\cal{E}}$, and using $B$ instead of $\kappa'$, we then get (see \cite[64.3)--(64.4), p.251]{Landau6} or \cite[p.~244]{morseingard})
\begin{equation}\label{7.5}
\partial_t p_{\cal{E}}+B\,\text{div\,}\mathbf{v}_{\cal{E}}=0,\qquad \rho_0 \partial_t \mathbf{v}_{\cal{E}}=-\nabla p_{\cal{E}},
\end{equation}
with $\rho'_{\cal{E}}$ is given in function of $p_{\cal{E}}$  as $\rho'_{\cal{E}}=\rho_0 B^{-1}p_{\cal{E}}$.  Such a first order approximation is considered to be adequate to describe small amplitude acoustic perturbations, as recalled in \S~\ref{intro}.  See \cite[pp.~244 and p.~257]{morseingard}.

In our problem $\Gamma=\Gamma_0\cup \Gamma_1$, the two surfaces having a different nature. Indeed $\Gamma_0$ cannot move and the fluid cannot penetrate in it, this condition being naturally expressed as
\begin{equation}\label{7.6}
\mathbf{v}_{\cal{E}}\cdot\nu=0\qquad\text{on $\R\times\Gamma_0$.}
\end{equation}
The surface $\Gamma_1$, as described in \cite{beale2} following \cite[p.~266]{morseingard}, models a membrane, with surface mass density $\mu$, which can move in reaction to the excess pressure $ p_{\cal{E}}$. Its normal displacement inside $\Omega$ (a scalar quantity) is denoted by $v$, and we can suppose that $v\ll 1$ as well.

In \cite{beale2} and \cite[p.~266]{morseingard} one supposes that $\Gamma_1$ reacts pointwise like an harmonic oscillator, so a force $\kappa v\nu$ acts on it. This type of reaction is called a local reaction. In \cite[p.~266]{morseingard} other types of reactions, called extended or local reactions, are mentioned. In the present paper we suppose that $\Gamma_1$ reacts like a membrane, a local component of the reaction being not excluded, so the force here is corrected to
$[-\DivGamma(\sigma \nabla_\Gamma v)+\kappa v]\nu$.

A damping effect on the membrane may be present, corresponding to the force $\delta v_t\nu$, where $\delta$ is the surface resistivity. Hence, by applying Newton Second Law on the motion of $\Gamma_1$ we get
\begin{equation}\label{7.7}
\mu v_{tt}- \DivGamma (\sigma \nabla_\Gamma v)+\delta v_t+\kappa v+p_{\cal{E}}=0\qquad
\text{on
$\R\times \Gamma_1$.}
\end{equation}
The second boundary condition on $\Gamma_1$ in \cite{beale2} express the fact that fluid particles cannot penetrate the membrane, so
by a classical argument of Fluid Mechanics, see \cite[p.~7]{Lamb}, they have to move with it. So we get
\begin{equation}\label{7.8}
\mathbf{v}_{\cal{E}}\cdot\nu=-v_t\qquad\text{on $\R\times\Gamma_1$.}
\end{equation}
Next, as the fluid is supposed to be non-viscous, we have
\begin{equation}\label{7.9}
\text{curl\,} \mathbf{v}_{\cal{E}}=0\qquad\text{in $\R\times\Omega$.}
\end{equation}
Putting together for future reference  \eqref{7.5}--\eqref{7.9} we then get the governing system
\begin{equation}\label{7.10}
\begin{cases}
 \partial_t p_{\cal{E}}+B\text{div\,}(\mathbf{v}_{\cal{E}})=0\qquad &\text{in
$\R\times\Omega$,}\\
 \rho_0 \partial_t \mathbf{v}_{\cal{E}}=-\nabla p_{\cal{E}}\qquad &\text{in
$\R\times\Omega$,}\\
\text{curl\,} \mathbf{v}_{\cal{E}}=0\qquad &\text{in
$\R\times\Omega$,}\\
\mu v_{tt}- \DivGamma (\sigma \nabla_\Gamma v)+\delta v_t+\kappa v+p_{\cal{E}}=0\qquad
&\text{on
$\R\times \Gamma_1$,}\\
\mathbf{v}_{\cal{E}}\cdot\nu=-v_t&\text{on
$\R\times \Gamma_1$,}\\
\mathbf{v}_{\cal{E}}\cdot\nu=0 &\text{on $\R\times \Gamma_0$.}
\end{cases}
\end{equation}
Since $\text{curl\,} \mathbf{v}_{\cal{E}}=0$, when $\Omega=\R^3$ (as in quoted textbooks), or more generally when $\Omega$ is
simply connected, one can then introduce  the velocity potential  (or, as in \cite{Landau6}, its opposite) as a scalar field $u_{\mathcal{E}}$ such that
$-\nabla u_{\mathcal{E}}=\mathbf{v}_{\cal{E}}$ in $\Omega$, so determining it (as a purely mathematical object) up to a free function of time $\varphi=\varphi(t)$. By the second equation in \eqref{7.10} one has $\nabla (\rho_0 \partial_t u_{\cal{E}}-p_{\cal{E}})=0$. Hence, by appropriately choosing $\varphi$ one then determines $u_{\mathcal{E}}$, up to a space-time constant, satisfying
\begin{equation}\label{7.11}
-\nabla u_{\mathcal{E}}=\mathbf{v}_{\cal{E}},\qquad  p_{\cal{E}}=\rho_0 \partial_t u_{\cal{E}}\qquad\text{in $\R\times\Omega$.}
\end{equation}
Consequently, denoting by $c=\sqrt{B/\rho_0}$ the sound speed, one gets that the couple $(u_{\mathcal{E}},v)$ satisfies the real version of problem \eqref{1.1}.
One then routinely considers complex-valued solutions for better mathematical handling.

The so described derivation of problem \eqref{1.1} is  rigorous, in terms of Theoretical Physics, as far as a first order perturbation can be. It allows to understand the meaning of our integral condition when $\Omega$ is bounded. Indeed the generalized Hooke's Law \eqref{7.3}, already used at an infinitesimal level, was experimentally verified at macroscopic level and when the fluid is at rest. In the present dynamic situation it is reasonable to replace $p_{\cal{E}}$ with its average $\fint_{\Omega_t}p_{\cal{E}}$, which is equivalent up to the first order with $\fint_{\Omega}p_{\cal{E}}$. Moreover, when $\cal{V}=|\Omega|$, being $v\ll 1$, the volume perturbation $\cal{V}'$ is approximated, up to the first order, by $-\int_{\Gamma_1}v$, so $\theta=-\frac 1{|\Omega|}\int_{\Gamma_1}v$. Hence generalized Hooke's Law reads as
\begin{equation}\label{7.12}
\int_{\Omega}p_{\cal{E}} = B\int_{\Gamma_1}v,
\end{equation}
or, by \eqref{7.11} (being $c^2=B/\rho_0$) as our constraint $\int_\Omega u_t=c^2\int_{\Gamma_1}v$.

Unfortunately this averaged form of generalized Hooke's Law is not automatically satisfied by solutions of \eqref{7.10}, which only have
$\int_{\Omega}p_{\cal{E}}-B\int_{\Gamma_1}v$ as a time invariant, as it is easily seen by the Divergence Theorem.
This fact does not matter when one uses \eqref{7.10} just to find wave-types solutions, but in a more theoretical analysis  one should add the integral condition \eqref{7.12} to \eqref{7.10} as an additional equation, at least when $\Omega$ is bounded.
To give a further confirmation of this assertion we now recall a second  derivation of the acoustic wave equation and our boundary conditions.

\section{The Lagrangian approach and its relation to the Eulerian one}\label{section7.2} A Lagrangian description of a fluid, without clearly distinguishing it from the Eulerian one, is given in mathematically less sophisticated textbooks in Physics like  \cite{ElmoreHeald,Feynman,GoldsteinMechanics}.  A mathematically more sophisticated treatment of this type is given in \cite[Chapter~7]{Achterberg}, dealing with gases. We shall essentially follow it in the sequel.

In the Lagrangian approach the spatial variable $x$ denotes the position of a fluid particle at some reference time $t_0$. The actual position of the particle at time $t$ is then a vector field, denoted in the sequel as $\mathbf{r}_1=\mathbf{r}_1(t,x)$. We also denote
$\mathbf{r}(t,x)=\mathbf{r}_1(t,x)-x$ the displacement of the particle from the original position. Consequently $\Omega$ denotes the region filled by the fluid at time $t_0$, while $\Omega_t=\{\mathbf{r}_1(t,x),\,x\in\Omega\}$ is the subset of $\R^3$  occupied by the fluid particles at time $t$.

The fluid is then described by the same fields as before, which are in the sequel denoted  as
$\mathbf{v}_{\cal{L}}$, $P_{\cal{L}}$ and $\rho_{\cal{L}}$, to distinguishing them from the Eulerian ones. The relation between any Lagrangian field $\cal{F}_{\cal{L}}$
and  the corresponding Eulerian one  $\cal{F}_{\cal{E}}$ is clearly given by
\begin{equation}\label{7.13}
\cal{F}_{\cal{L}}(t,x)=\cal{F}_{\cal{E}}(t,\mathbf{r}_1(t,x)).
\end{equation}
We also recall the well-known notion of Lagrangian time derivative $\frac {d\cal{F}_{\cal{L}}}{dt}$ (also known as co-moving or material derivative) of any  Lagrangian field $\cal{F}_{\cal{L}}$, which gives the rate of change of the field when it is measured following a specific fluid particle, i.e., moving with it.
Since $\mathbf{v}_{\cal{L}}(t,x)=\partial_t\mathbf{r}_1(t,x)=\mathbf{r}_t(t,x)$, by \eqref{7.13} we have \begin{equation}\label{7.14}
 \frac {d\cal{F}_{\cal{L}}}{dt}(t,x)=(\partial_t \cal{F}_{\cal{E}}+\mathbf{v}_{\cal{E}}\cdot \nabla \cal{F}_{\cal{E}})(t,\mathbf{r}_1(t,x)).
\end{equation}
One then makes the same physical assumptions on the fluid motion made in the Eulerian approach and one considers an acoustic perturbation of it, i.e., a small amplitude perturbation of density  and pressure  from the equilibrium values $\rho_0$, $P_0$,  here denoted as
$\rho'_{\cal{L}}$ and $p_{\cal{L}}$, that is  $\rho_{\cal{L}}=\rho_0+\rho'_{\cal{L}}$ and $P_{\cal{L}}=P_0+p_{\cal{L}}$.
Since the unperturbed fluid is supposed to be at rest, the Lagrangian perturbation of $\mathbf{v}_{\cal{L}}$ coincides with $\mathbf{v}_{\cal{L}}$. Moreover, from previous relations one trivially gets
\begin{equation}\label{7.15}
\rho'_{\cal{L}}(t,x)=\rho'_{\cal{E}}(x,\mathbf{r}_1(t,x))\qquad\text{and}\quad  p_{\cal{L}}(t,x)=p_{\cal{E}}(x,\mathbf{r}_1(t,x)).
\end{equation}
In \cite[\S~7.4]{Achterberg} a detailed perturbation analysis  with respect to the displacement $\mathbf{r}$ is given. We shall not repeat it here for shortness. Its conclusion  is that, since the unperturbed velocity vanishes and the unperturbed density and pressure are constant, see \cite[Table~7.1, p.~139, (7.4.40), p.~138 and (7.5.1) p.~143]{Achterberg},
\begin{equation}\label{7.16}
\mathbf{v}_{\cal{L}}\simeq \mathbf{v}_{\cal{E}}, \quad  \rho'_{\cal{L}}\simeq \rho'_{\cal{E}}, \quad p_{\cal{L}}\simeq p_{\cal{E}}\quad\text{and}\quad
 \rho'_{\cal{L}}\simeq -\rho_0\,\,\text{div\,}\mathbf{r}
\end{equation}
up to the first order. Hence, by \eqref{7.16} we can rewrite the infinitesimal Hooke's Law \eqref{7.2}, up to the first order, as
\begin{equation}\label{7.17}
 p_{\cal{L}}=-B\,\, \text{div\,}\mathbf{r}.
\end{equation}
Moreover, by directly applying Newton Second Law, or by rewriting Euler equation \eqref{7.1}$_2$ taking into account \eqref{7.13}--\eqref{7.14}, the main motion equation  for the fluid is
\begin{equation}\label{7.18}
  \rho_\cal{L}(t,x)\frac d{dt}\mathbf{v}_\cal{L}(t,x)=-\nabla p_\cal{E}(t,\mathbf{r}_1(t,x)).
\end{equation}
Since $\mathbf{v}_{\cal{L}}=\mathbf{r}_t$ we have $\frac d{dt}\mathbf{v}_\cal{L}=\mathbf{r}_{tt}$.
Moreover, since
$$\nabla p_\cal{L}(t,x)=\nabla p_\cal{E}(t,\mathbf{r}_1(t,x))\,J\mathbf{r}_1(t,x)
=\nabla p_\cal{E}(t,\mathbf{r}_1(t,x))\,[I+J\mathbf{r}(t,x)],$$
we also have $\nabla p_\cal{L}(t,x)\simeq \nabla p_\cal{E}(t,\mathbf{r}_1(t,x))$ us to the first order. Consequently, since $\rho_{\cal{L}}=\rho_0+\rho'_{\cal{L}}$, performing the same approximation
made to obtain \eqref{7.4}, the first order linearized version of \eqref{7.18}
is
\begin{equation}\label{7.19}
  \rho_0\mathbf{r}_{tt}=-\nabla p_\cal{L}.
\end{equation}
The system composed by \eqref{7.17} and \eqref{7.19} in the Lagrangian approach takes the place of the system \eqref{7.4} in the Eulerian one. It coincides with \cite[(5.13) and (5.16)]{ElmoreHeald}, and combining the two equations one gets the vector-valued wave equation
\begin{footnote}
The name will be justified in the sequel since we shall have $\text{curl } \mathbf{r}=0$ so that $\nabla\text{div\,}\mathbf{r}=\Delta \mathbf{r}$.
 \end{footnote}
$\mathbf{r}_{tt}-c^2\nabla\,\text{div\,}\mathbf{r}=0$, with $c=\sqrt{B/\rho_0}$ as before, given in \cite[Chapter~11]{GoldsteinMechanics} and \cite[(7.5.11)]{Achterberg}.
 We discuss in the sequel the translation of the already discussed boundary conditions \eqref{7.6}--\eqref{7.8} from the Eulerian to the Lagrangian approach.

Since fluid particles cannot penetrate neither $\Gamma_0$ nor $\Gamma_1$, as proved in \cite[p.7]{Lamb}, if they are at $\Gamma_0$ (or $\Gamma_1$) at the reference time $t_0$  they have to move with $\Gamma_0$ (or $\Gamma_1$). Hence, as  $\mathbf{v}_{\cal{L}}(t,x)=\mathbf{r}_t(t,x)$, \eqref{7.6}--\eqref{7.8} translate into
\begin{gather}\label{7.20}
 \mathbf{r}_t\cdot \nu=0\quad\text{on $\R\times\Gamma_0$},\quad\text{and}\quad \mathbf{r}_t\cdot \nu=-v_t\quad\text{on $\R\times\Gamma_1$,} \\
\label{7.21}\mu v_{tt}- \DivGamma (\sigma \nabla_\Gamma v)+\delta v_t+\kappa v+p_{\cal{L}}=0\qquad
\text{on
$\R\times \Gamma_1$.}
\end{gather}
We are now going to integrate from the reference time $t_0$. One has to be cautious about performing such an operation. Indeed  the model just describes  the propagation, and not the origin, of an acoustic perturbation, all information about it being contained in the initial data.  Hence the motion equations do not model the physical behavior from $t_0$ to $t$. Otherwise, the solution would  {\em ever} vanish. But we can integrate \eqref{7.20}, which only  express the impenetrability of $\Gamma_0$ and $\Gamma_1$.

Since $\mathbf{r}(t_0,\cdot)=0$, we have $\Omega_{t_0}=\Omega$ and  then $v(t_0,\cdot)\equiv 0$,
we get the {\em equivalent} boundary conditions
\begin{equation}\label{7.22}
\mathbf{r}\cdot \nu=0\quad\text{on $\R\times\Gamma_0$},\quad\text{and}\quad \mathbf{r}\cdot \nu=-v\quad\text{on $\R\times\Gamma_1$,}
\end{equation}
involving $\mathbf{r}$ instead of its time derivative. Such a simplification is impossible in the Eulerian approach, in which the displacement $\mathbf{r}$ is meaningless (see for example \cite[Chapter VI, p.~235]{morseingard}). We shall show in the sequel that it allows to {\em prove} that in the Lagrangian approach  our constraint is a mathematical consequence of the model.

We now assume, as in the Eulerian approach, that $\Omega$ is simply connected. The non-viscosity condition \eqref{7.9} then translate into
$\text{curl\,} \mathbf{v}_{\cal{E}}(t,\mathbf{r}_1(t,x))=0$. Since, by \eqref{7.13},
$$J\mathbf{v}_\cal{L}(t,x)=J\mathbf{v}_\cal{E}(t,\mathbf{r}_1(t,x))\,J\mathbf{r}_1(t,x)=
J\mathbf{v}_\cal{E}(t,\mathbf{r}_1(t,x))\,[I+J\mathbf{r}(t,x)]
$$
up to the first order we have $J\mathbf{v}_\cal{L}(t,x)\simeq J\mathbf{v}_\cal{E}(t,\mathbf{r}_1(t,x))$ and consequently
$\text{curl\,}\mathbf{r}_t(t,x)=\text{curl\,} \mathbf{v}_{\cal{L}}(t,x)\simeq \text{curl\,}\mathbf{v}_\cal{E}(t,\mathbf{r}_1(t,x))$. Hence we can assume that
$\text{curl\,}\mathbf{r}_t=0$. Also in this case we  integrate from the reference time $t_0$, since any acoustic perturbation, not considered in the model, must be irrotational. Hence, as $\mathbf{r}(t_0,\cdot)=0$ we get
\begin{equation}\label{7.23}
  \text{curl\,}\mathbf{r}=0\qquad\text{in $\R\times\Omega$,}
\end{equation}
obtaining the same accuracy of \eqref{7.9} as a first order approximation.

Putting together for future reference  \eqref{7.17}, \eqref{7.19} and \eqref{7.21}--\eqref{7.23} we then get the Lagrangian version of the governing system \eqref{7.10}
\begin{equation}\label{7.24}
\begin{cases}
 p_\cal{L}+B\,\text{div\,}\mathbf{r}=0\qquad &\text{in
$\R\times\Omega$,}\\
 \rho_0\,\mathbf{r}_{tt}=-\nabla p_\cal{L}\qquad &\text{in
$\R\times\Omega$,}\\
\text{curl\,} \mathbf{r}=0\qquad &\text{in
$\R\times\Omega$,}\\
\mu v_{tt}- \DivGamma (\sigma \nabla_\Gamma v)+\delta v_t+\kappa v+p_{\cal{L}}=0\qquad
&\text{on
$\R\times \Gamma_1$,}\\
\mathbf{r}\cdot\nu=-v&\text{on
$\R\times \Gamma_1$,}\\
\mathbf{r}\cdot\nu=0 &\text{on $\R\times \Gamma_0$,}
\end{cases}
\end{equation}
involving the displacement vector field instead of the velocity field, and formally appearing as the integrated in time version of \eqref{7.10}.

As in the Eulerian approach we then introduce a velocity potential $u_\cal{L}$ such that
$-\nabla u_\cal{L}=\mathbf{v}_\cal{L}=\mathbf{r}_t$. So, by \eqref{7.19}, we get $\nabla(\rho_0\partial_t u_\cal{L}-p_\cal{L})=0$ and, by appropriately choosing $\varphi$,
\begin{equation}\label{7.25}
-\nabla u_{\mathcal{L}}=\mathbf{r}_t,\qquad  p_{\cal{L}}=\rho_0 \,\partial_t u_{\cal{L}}\qquad\text{in $\R\times\Omega$,}
\end{equation}
so $u_\cal{L}$ satisfies the real version of problem \eqref{1.1} with $c=\sqrt{B/\rho_0}$.

In this approach the constraint \eqref{OurConstraint} is automatically satisfied. Indeed, when $\Omega$ is bounded, by \eqref{7.17}, \eqref{7.21}, \eqref{7.25} and the Divergence Theorem we have
$$\int_\Omega\partial_t u_\cal{L}=\frac 1{\rho_0}\int_\Omega p_\cal{L}=-c^2\int_\Omega \text{div\,}\mathbf{r}=-c^2\int_\Gamma \mathbf{r}\cdot \nu= c^2\int_{\Gamma_1}v.$$
Also in this case this equality is the macroscopic version of generalized Hooke's Law, but it appears as a mathematical consequence of the infinitesimal version \eqref{7.17}, thank to the use of $\mathbf{r}$ and to the transformation of \eqref{7.20} into \eqref{7.21}.

\section{Conclusions}\label{section7.3}
The two physical derivations of problem \eqref{1.1}  presented in previous sections show that its strong solutions originating from data in $\cal{H}_1$
represents, to some extent, solutions of \eqref{7.10}--\eqref{7.12}  and \eqref{7.24}, which are equivalent as  first order approximations of the physical problem.  It would be natural to make precise the meaning of strong (weak) solutions of \eqref{7.10}--\eqref{7.12}  and \eqref{7.24} and to discuss the relationships between them and strong (weak) solutions of \eqref{1.1}. Moreover, the discussion on weak solutions would naturally entail a discussion on the associated groups. These discussions are too long for the present paper.

Hence, to physically interpret our result, we limit to remark that, in both the presented derivations, the physical state of the fluid  is fully described by the excess pressure $p$, the density $\rho'$ and the fluid velocity $\mathbf{v}$. We do not make a choice between the Eulerian and Lagrangian meaning of these fields, since as shown by \eqref{7.16} they are equivalent as  first order approximation. The reader can make its own choice between them. Moreover, as shown by \eqref{7.2} and \eqref{7.16}--\eqref{7.17}, the density $\rho'$ and the excess pressure $p$ are linearly related through the relation $c^2\rho'=p$, so the physical state of the fluid  is fully described by the excess pressure $p$ and the fluid velocity $\mathbf{v}$, with $\text{curl\,}\mathbf{v}=0$.

We then introduce the space
$\nabla(H^1(\Omega))=\{\nabla\phi, \phi\in H^1(\Omega)\}$ and we recall
(see  \cite[Proposition~2, p.~219]{dautraylionsvol3}) that
since $\Omega$ is simply connected  $\nabla(H^1(\Omega))=\{\mathbf{v}\in [L^2(\Omega)]^3: \text{curl\, }\mathbf{v}=0\}$, where $\text{curl\,}\mathbf{v}$ is taken in the sense of distributions, and that $\nabla(H^1(\Omega))$ is a closed subspace of $[L^2(\Omega)]^3$. We shall equip $\nabla(H^1(\Omega))$ with the norm inherited from $[L^2(\Omega)]^3$.

By the way we also have to consider the boundary deformation $v$ and its time derivative. We mathematically formalize these remarks in the following
\begin{definition}\label{definition7.1} Suppose that (A0--4) hold,  $N=3$ and $\Omega$ is simply connected. Then we set $B=\rho_0c^2$ and for any weak solution $(u,v)$ of problem \eqref{1.1} with $U_0\in\cal{H}_1$  we define the excess pressure $p$ and  the fluid velocity $\mathbf{v}$ associated to the solution, and their initial values, by
\begin{equation}\label{7.26}
  p=\rho_0 u_t,\quad \mathbf{v}=-\nabla u \qquad\text{and}\quad
  p_0=\rho_0 u_1,\quad \mathbf{v}_0=-\nabla u_0.
\end{equation}
We shall respectively  call $\cal{U}=(p,\mathbf{v},v, v_t)$ and $\cal{U}_0=(p_0,\mathbf{v}_0,v_0,v_1)$
the {\em  state of the physical system} and   {\em its initial state}, while
$$\cal{X}=\left\{(p_0,\mathbf{v}_0,v_0,v_1)\in L^2(\Omega)\times \nabla(H^1(\Omega))\times H^1(\Gamma_1)\times L^2(\Gamma_1):{\textstyle\int}_\Omega p_0=B{\textstyle\int}_{\Gamma_1}v_0\right\}$$
is its {\em phase space}, endowed with the standard product norm.
\end{definition}

Clearly $\cal{U}_0\in \cal{X}$, any element of $\cal{X}$ can be written as $(\rho_0 u_1,-\nabla u_0, v_0,v_1)$ for some $U_0\in\cal{H}_1$, $\cal{U}\in C(\R;\cal{X})$ and $\cal{U}$ only depends on  $\cal{U}_0$, since for any choice of $u_0$ such that $-\nabla u_0=\mathbf{v}_0$ we have the same solution $(u,v)$, up to a space-time constant.

To better clarify  the physical meaning of our results we now make explicit their consequences for the just defined physical state behavior, starting from the  case  when $\Gamma_1$ is connected. When the system is damped we have

\begin{cor}[\bf Physical state stability for $\delta\not\equiv 0$, $\Gamma_1$ connected]\label{corollary7.1}
Suppose that (A0--4) hold, $\Gamma_1$ is connected, $\Omega$ is simply connected and $\delta\not\equiv 0$. Then
for all $\cal{U}_0\in\cal{X}$ we have $\cal{U}(t)\to 0$ in $\cal{X}$ as $t\to\infty$.
\end{cor}
\begin{proof}Combine Definition~\ref{definition7.1} and Theorem~\ref{theorem1.3}.
\end{proof}
This expected result motivates the notion of physical stability in  Theorem~\ref{theorem1.3}.
In the undamped case we have the following result, the meaning of which was already commented in Remark~\ref{remark1.2}.
\begin{cor}[\bf Fourier decomposition for $\delta\equiv 0$, $\Gamma_1$ connected]\label{corollary7.2}
Suppose that (A0--4) hold, $\Gamma_1$ is connected, $\Omega$ is simply connected and $\delta\equiv 0$.

Then there is sequence $((p^n,\mathbf{v}^n,v^n,v^n_t))_n$ of standing wave states of the physical system, i.e., states of the form
\begin{alignat*}2\label{7.29b}
p^n(t,x)=&p^{0n}(x)e^{-\mathfrak{i}(\lambda_nt+\pi/2)},\qquad &&\mathbf{v}^n(t,x)=\mathbf{v}^{0n}(x) e^{-\mathfrak{i}\lambda_nt},\\
v^n(t,x)=&v^{0n} (x) \, e^{-\mathfrak{i}(\lambda_nt+\pi/2)},&&v^n_t(t,x)=-\lambda_n v^{0n} (x)  e^{-\mathfrak{i}\lambda_nt},
\end{alignat*}
where $(\lambda_n)_n$ is a real sequence satisfying Theorem~\ref{theorem1.4}--i) and $(p^{0n}, \mathbf{v}^{0n}, v^{0n})_n$ is a sequence of non-identically vanishing real-valued elements of $H^2(\Omega)\times [H^1(\Gamma_1)]^3\times H^2(\Gamma_1)$ satisfying the problems
\begin{equation}\label{7.30}
\begin{cases} B\text{div\,}\mathbf{v}^{0n}=\lambda_np^{0n} \qquad &\text{in
$\Omega$,}\\
- \DivGamma (\sigma \nabla_\Gamma v^{0n})+\kappa v^{0n}+p^{0n} =\mu \lambda_n^2v^{0n}\qquad
&\text{on $\Gamma_1$,}\\
\mathbf{v}^{0n}\cdot\nu=0 &\text{on $\Gamma_0$,}\\
\mathbf{v}^{0n}\cdot\nu=\lambda_nv^{0n}\qquad
&\text{on
$\Gamma_1$.}
\end{cases}
\end{equation}
Moreover for all $\cal{U}_0\in\cal{X}$ we have $\cal{U} =\sum\limits_{n=1}^\infty \alpha_n \cal{U}^n$ in $C_b(\R; \cal{X})$, where
\begin{equation}\label{7.32}
\begin{alignedat}2
&\alpha_n=\alpha_n(\cal{U}_0)=&&\rho_0\int_\Omega \mathbf{v}^0\mathbf{v}^{0n}-\lambda_n \int_{\Gamma_1}\mu v_1v^{0n}+\mathfrak{i}\int_{\Gamma_1}\kappa v_0v^{0n}\\
&&&+\mathfrak{i}\int_{\Gamma_1}\sigma (\nabla_\Gamma v_0,\nabla_\Gamma v^{0n})_\Gamma+\frac{\mathfrak{i}c^2}{\rho_0}\int_{\Gamma_1}p_0p^{0n}.
\end{alignedat}
\end{equation}
\end{cor}
\begin{proof}Write Theorem~\ref{theorem1.4}, set $p^{0n}=\rho_0\lambda_n u^{0n}$, $\mathbf{v}^{0n}=-\nabla u^{0n}$
and use Definition~\ref{definition7.1}.
\end{proof}
We now consider the more involved case when $\Gamma_1$ can be disconnected. By Theorem~\ref{theorem5.3} it is evident that in Corollaries~\ref{corollary7.1}--\ref{corollary7.2} we can replace the assumption that  $\Gamma_1$ connected, i.e., $\mathfrak{n}=1$, to the more general one $\mathfrak{n}_0\le 1$, without differences. Hence problem \eqref{1.1} appears to give a good approximation of the acoustic phenomena when a local reaction is present in all connected components of $\Gamma_1$ but one.

The situation is different when $\mathfrak{n}_0\ge 2$, which we assume in the remainder of this section.
\begin{cor}[\bf Physical state stability when $\mathfrak{n}_{00}\le 1<2\le\mathfrak{n}_0$]\label{corollary7.3}
Suppose that (A0--4) hold, $\Omega$ is simply connected and $\mathfrak{n}_{00}\le 1<2\le\mathfrak{n}_0$.
Set, for $i=1,\ldots,\mathfrak{n}_0$, $\widetilde{\beta}_i\in\cal{X}'$ by
\begin{equation}\label{7.34}
\widetilde{\beta}_i\cal{U}_0\!\!=\!\!\frac{\sum_{j=1}^{\mathfrak{n}_0}\cal{H}^{N-1}(\Gamma_1^j)\,d_{ij}\!\!\left[\fint_{\Gamma_1^i}(\mu v_1+\delta v_0+\rho_0u_0)\!-\!\fint_{\Gamma_1^j}(\mu v_1+\delta v_0+\rho_0u_0)\right]}{\sum_{j=1}^{\mathfrak{n}_0}d_j\cal{H}^{N-1}(\Gamma_1^j)}
\end{equation}
for all $\cal{U}_0=(p_0,\mathbf{v}_0,v_0,v_1)\in\cal{X}$, where $u_0$ is any element of $H^1(\Omega)$ such that $-\nabla u_0=\mathbf{v}_0$.
Then for all $\cal{U}_0\in\cal{X}$ we have
\begin{equation}\label{7.35}
\cal{U}(t)\to  \sum\limits_{j=1}^{\mathfrak{n}_0} (\widetilde{\beta}_j \cal{U}_0)\,\,\, (0,0, \chi_{\Gamma_1^j},0)\qquad\text{in $\cal{X}$, as $t\to\infty$}.
\end{equation}
\end{cor}
\begin{proof}We combine  Theorem~\ref{theorem5.4} and Definition~\ref{definition7.1} and remark that the $\beta_i$'s in \eqref{5.63} actually do not depend on ${\textstyle\fint}_\Omega u_0$ and that $L_1U=0$ when $U\in\cal{H}_1$.
\end{proof}

By Lemma~\ref{lemma5.8} the limit set in \eqref{7.34}--\eqref{7.35} is nothing but
$$\cal{X}_{01}:=\{0\}\times\{0\}\times\cal{N}_1(\Gamma_1)\times\{0\}=\{(\rho_0u_1,-\nabla u_0, v_0,v_1)\in \cal{X}: (u_0,v_0,u_1,v_1)\in \cal{N}_{01}\}$$
 where $\cal{N}_{01}$ and $\cal{N}_1(\Gamma_1)$ are the spaces defined in \eqref{N01}--\eqref{calN1}.

Trivially, when one takes $p=p_\cal{E}$ and $\mathbf{v}=\mathbf{v}_\cal{E}$, each $\cal{U}_0\in\cal{X}_{01}$ is a stationary solution of \eqref{7.10} satisfying \eqref{7.12}. Also when one interpret $p$ and $\mathbf{v}$ as Lagrangian fields, for each $\cal{U}_0=(0,0,v_0,0)\in\cal{X}_{01}$ (i.e., $v_0\in\cal{N}(\Gamma_1)$), the Neumann problem
\begin{equation}\label{7.36}
  \Delta \phi_0=0\quad \text{in $\Omega$},\quad \partial_\nu\phi_0=0\quad \text{on $\Gamma_0$},\quad \partial_\nu\phi_0=v_0\quad\text{on $\Gamma_1$,}
\end{equation}
has a (unique up to an additive constant) solution $\phi_0\in H^2(\Omega)$, so that $\mathbf{r}_0=-\nabla \phi_0$, $p_\cal{L}=0$, $v=v_0$ is a stationary solution of \eqref{7.24} having $\cal{U}_0$ as physical state.

One can easily prove, this assertion being unessential for the present discussion, that all suitably smooth stationary solutions of \eqref{7.10} satisfying \eqref{7.12} and of \eqref{7.24} can be obtained in this way.

In both versions of the physical model these stationary solutions are characterized by the vanishing of the acoustic perturbation $p$ and of the (correctly defined) energy and by a possibly  large constant advance or retreat of  $\Gamma_1^1,\ldots,\Gamma_1^{\mathfrak{n}_0}$, compensating each other so that the first order volume perturbation vanishes.

These solutions are mathematically trivial, since no control (but for the time derivative) is present on  ${\textstyle\fint}_{\Gamma_1^1}v,\ldots {\textstyle\fint}_{\Gamma_1^{\mathfrak{n}_0}}v$ in the equation \eqref{7.7} (or the equivalent \eqref{7.21}) which governs the boundary dynamic, when \eqref{7.7} is rewritten for each $\Gamma_1^i$. But they are not consistent with the acoustic model.

Actually the asymptotic limits in \eqref{7.35} seem to be mostly related to the inertial (damped) dynamics of the
boundary pieces $\Gamma_1^1,\ldots,\Gamma_1^{\mathfrak{n}_0}$, which are not kept in position by the spring term, when the influence of the acoustic perturbation is neglected. Indeed, let us uncouple $u$ (and so $(p,\mathbf{v})$) and $v$ by neglecting the term $\rho_0 u_t$ in \eqref{1.1}$_2$ and the boundary condition \eqref{1.1}$_3$. Then, denoting $v^i=v_{|\Gamma_1^i}$, $v_0^i={v_0}_{|\Gamma_1^i}$ and $v^i_1={v_1}_{|\Gamma_1^i}$, each $v^i$, $i=1,\ldots,\mathfrak{n}_0$, satisfies the Cauchy problem
\begin{equation}\label{Ex1}
\begin{cases}
\mu v_{tt}^i-\DivGamma(\sigma \nabla_\Gamma v^i)+\delta v_t^i=0\quad &\text{on $\R\times \Gamma_1^i$,}\\
v^i(0)=v_0^i, \quad v_t^i(0)=v_1^i.
\end{cases}
\end{equation}
Then, Problem \eqref{Ex1} is easily solved, as it becomes a linear ODE's Cauchy problem, when $v_0^i, v_1^i$ and $\delta,\mu>0$  are constant.
Indeed in this case we have $v^i(t,x)=v^i(t)=\left(v_0^i+\frac \mu\delta v_1^i\right)-\frac \mu\delta v_1^ie^{-\frac \delta \mu t}$, and then
$v^i(t)\to v^i(\infty)$ as $t\to\infty$, where
\begin{equation}\label{Ex2}
v^i(\infty)=\left(v_0^i+\frac \mu\delta v_1^i\right),\qquad\text{for $i=1,\ldots, \mathfrak{n}_0$.}
\end{equation}
Now taking  $v_0^i, v_1^i$ and $\delta,\mu$ as before and choosing an initial   state $\cal{U}_0=(p_0,\mathbf{v}_0, v_0,v_1)\in\cal{X}$
such that $\text{supp\,}\mathbf{v}_0\subset\subset\Omega$, so we can take ${u_0}_{|\Gamma}=0$,
 and $v_0$, $v_1$ satisfying ${v_0}_{|\Gamma_1^i}=v_0^i$, ${v_1}_{|\Gamma_1^i}=v^i_1$ for $i=1,\ldots, \mathfrak{n}_0$,
the coefficients $\widetilde{\beta}_i\cal{U}_0$ in \eqref{7.35} are given by
\begin{equation}\label{Ex3}
\begin{aligned}
\widetilde{\beta}_i\cal{U}_0\quad=\quad
&\frac{\sum_{j=1}^{\mathfrak{n}_0}\cal{H}^{N-1}(\Gamma_1^j)\left[\frac\mu\delta v_1^i+ v_0^i-\frac\mu\delta v_1^j- v_0^j\right]}{\sum_{j=1}^{\mathfrak{n}_0}\cal{H}^{N-1}(\Gamma_1^j)}\\=\quad
&v^i(\infty)-\sum_{j=1}^{\mathfrak{n}_0}\frac{\cal{H}^{N-1}(\Gamma_1^j)}{\sum_{k=1}^{\mathfrak{n}_0}\cal{H}^{N-1}(\Gamma_1^k)}v^j(\infty)
\end{aligned}
\end{equation}
so they only depend by the asymptotic limits of solutions of \eqref{Ex1}.

To show that the model may become geometrically inconsistent for large data  we briefly analyze the case $\Omega=B_R\setminus\overline{B_r}$, $0<r<R$,
$\Gamma_0=\emptyset$, $\Gamma_1^1=\partial B_r$, $\Gamma_1^2=\partial B_R$ already considered in Example~\ref{esempio}, in which the real part of $\cal{N}(\Gamma_1)$ is $\{v_{0a},\,\, a\in\R\}$, where $v_{0a}=a(R^2\chi_{\partial B_r}-r^2\chi_{\partial B_R})$. In the Lagrangian framework the solution $\phi_{0a}$ of \eqref{7.36} associated to $v_{0a}$ is $\phi_a(x)=aR^2r^2/|x|$, so inducing the central displacement field $\mathbf{r}_{a}(x)=aR^2r^2|x|^{-3}x$ and actual position field $\mathbf{r}_{1a}(x)=(1+aR^2r^2|x|^{-3})x$, $x\in\Omega$.
We regard the deformed configuration $\mathbf{r}_{1a}(\overline{\Omega})$ of $\overline{\Omega}$ as geometrically consistent
 only when $\mathbf{r}_{1a}(\partial\Omega)=\partial [\mathbf{r}_{1a}(\Omega)]$ and neither $\mathbf{r}_{1a}(\partial B_r)$ or $\mathbf{r}_{1a}(\partial B_R)$ collapses to a point. Due to the radial symmetry of  $\mathbf{r}_{1a}$, the monotonicity requirement  is equivalent  to the strict monotonicity of the continuous real function $\tau\mapsto |\tau+aR^2r^2\tau^{-2}|$  in $[r,R]$. Taking the square, deriving with respect to $\tau$ and making the change of variable $\tilde{\tau}=\tau^3$ our monotonicity request is equivalent to the lack of roots $\tilde{\tau}$ of the equation $\tilde{\tau}^2-aR^2r^2\tilde {\tau}-2a^2R^4r^4=0$ in $(r^3,R^3)$, i.e., to $a\in [-r/R^2, r/2R^2]$. The non--collapsing request simply reduces to $a\in (-r/R^2, r/2R^2]$.
Hence our model is geometrically consistent only for small values of $|a|$. Since, by \eqref{7.34}, the asymptotic limit in \eqref{7.35} is small for small data, it is geometrically consistent, although physically unrealistic, only for small data.

By the way one can avoid this phenomenon by restricting \eqref{1.1} to initial data $U_0\in \cal{H}_{\mathfrak{n}_0}$, that is to initial physical states $(p_0,\mathbf{v}_0, v_0,v_1)\in \cal{X}$ such that, choosing any $u_0\in H^1(\Omega)$ with $-\nabla u_0=\mathbf{v}_0$, one has
\begin{equation}\label{7.37}
  \fint_{\Gamma_1^i}\mu v_1+\delta v_0+\rho_0u_0=\fint_{\Gamma_1^1}\mu v_1+\delta v_0+\rho_0u_0,\qquad i=2,\ldots,\mathfrak{n}_0.
\end{equation}
These conditions trivially hold if the initial state $(p_0,\mathbf{v}_0, v_0,v_1)$ satisfies $\fint_\Omega p_0=0$,
$\text{supp\,}\mathbf{v}_0\subset\subset\Omega$, so we can take ${u_0}_{|\Gamma}=0$, and $v_0=v_1=0$, that is if the acoustic perturbation, occurred in the past and originating inside $\Omega$, did not perturb the boundary. But in general the conditions \eqref{7.37}, as the asymptotic limits in \eqref{7.35} before, seem to be mostly related with  the inertial (damped) dynamics of the  boundary pieces $\Gamma_1^1,\ldots,\Gamma_1^{\mathfrak{n}_0}$ which are not kept in position by the spring term, when the influence of the acoustic perturbation is neglected. Indeed, all the coefficients $\widetilde{\beta}_i\cal{U}_0$ in \eqref{Ex3} vanish if and only if
$v^1(\infty)=v^2(\infty)=\ldots=v^{\mathfrak{n}_0}(\infty)$.

The just described phenomenon appears  to be more evident when $\mathfrak{n}_{00}\ge 2$. Indeed, when evaluating the physical state associated through \eqref{7.26} to the trivial solutions of type (s3) described at page~\pageref{5.73}, one gets for $i=2,\ldots,\mathfrak{n}_{00}$, $p=0$, $\mathbf{v}(t)=\mathbf{v}_{0i}$, where $\mathbf{v}_{0i}=-\nabla u_{0i}\in [H^1(\Omega)]^3$ solves the problem
$$\begin{cases}
\begin{split}
\text{div\,}\mathbf{v}_{0i}&=0\quad&&\text{in $\Omega$,}\\
\text{curl\,}\mathbf{v}_{0i}&=0\quad&&\text{in $\Omega$,}\\
\mathbf{v}_{0i}\cdot\nu&=-\cal{H}^{N-1}(\Gamma_1^i)\quad&&\text{on $\Gamma_1^1$,}\\
\mathbf{v}_{0i}\cdot\nu&=\cal{H}^{N-1}(\Gamma_1^1)\quad&&\text{on $\Gamma_1^i$,}\\
\mathbf{v}_{0i}\cdot\nu&=0\quad&&\text{on $\Gamma_1\setminus(\Gamma_1^1\cup\Gamma_1^i)$,}\\
\end{split}
\end{cases}$$
and $v(t)=v_{1i}t$, the $v_{1i}$'s being given by \eqref{5.74}.

Even in the Eulerian framework  it is clear that these solutions are characterized by a vanishing acoustic perturbation and an unbounded domain dislocation.
So, even  when one multiplies them by a small real constant (so having small initial data) these solutions become geometrically inconsistent for large $|t|$.

This fact is more evident in the Lagrangian framework. Fixing any $i=2,\ldots,\mathfrak{n}_{00}$, for any $a\in\R$ the solution of \eqref{7.24} associated to the trivial solution $u_a(t,x)=au_{0i}(x)$, $v_a(t,x)=av_{1i}(x)$, where $u_{0i}$ and $v_{1i}$ are given in \eqref{5.73}--\eqref{5.74},
is $\mathbf{r}_{a}(t,x)=-a\nabla u_{0i}(x)t$ and $v_a(t,x)=a\left[\cal{H}^{N-1}(\Gamma_1^i)\chi_{\Gamma_1^1}-\cal{H}^{N-1}(\Gamma_1^1)\chi_{\Gamma_1^i}\right]t$.
In the particular case of the spherical shell in Example~\ref{esempio} we then have
$\mathbf{r}_a(t,x)=4\pi a R^2 r^2|x|^{-3}xt$. Hence, by repeating the arguments used above, for any $a\not=0$ the configuration is geometrically inconsistent for $|t|$ large enough.

It is not clear if conditions \eqref{7.37}, of doubtful physical meaning, can avoid this phenomenon. One can prove that when $\mathfrak{n}=2$, in which \eqref{7.37} reduce to a single equation, trivial solutions of type (s3) do not verify it.  Also in this simplest case the authors were not able to find a further invariant of the evolution problem which allows to isolate them, following the already explained general strategy.

\backmatter
\newcommand{\etalchar}[1]{$^{#1}$}
\def\cprime{$'$}
\providecommand{\bysame}{\leavevmode\hbox to3em{\hrulefill}\thinspace}
\providecommand{\MR}{\relax\ifhmode\unskip\space\fi MR }
\providecommand{\MRhref}[2]{%
  \href{http://www.ams.org/mathscinet-getitem?mr=#1}{#2}
}
\providecommand{\href}[2]{#2}

\printindex
\end{document}